\documentclass[12pt]{amsart}

\usepackage{graphicx}
\usepackage{subcaption}
\usepackage{amsfonts}
\usepackage{epsf}
\usepackage{amssymb}
\usepackage{amsmath}
\usepackage{amscd}
\usepackage{hyperref}
\usepackage{color}

\newtheorem{theorem}{Theorem}[section]
\newtheorem{proposition}[theorem]{Proposition}

\theoremstyle{remark}
\newtheorem{definition}[theorem]{Definition}
\newtheorem{example}[theorem]{Example}

\begin{document}
\title{THE NON-ORIENTABLE 4-GENUS FOR  KNOTS with 10 CROSSINGS}
\begin{abstract}
Given a knot in the 3-sphere, the non-orientable 4-genus or 4-dimensional crosscap number of a knot is the minimal first Betti number of non-orientable surfaces, smoothly and properly embedded in the 4-ball, with boundary the knot. In this paper, we calculate the non-orientable 4-genus of knots with crossing number 10.
\end{abstract}


\author{Nakisa Ghanbarian}
\email{nakisa@unr.edu}
\address{Department of Mathematics and Statistics, University of Nevada, Reno NV 89557}

\thanks{The author gratefully acknowledges partial support from the National Science Foundation, DMS-1906413.}
\maketitle
\section{Main Results}
Given a knot $K$, the {\em non-orientable 4-genus} or {\em 4-dimensional crosscap number of $K$}, denoted $\gamma_4(K)$ and defined by Murakami and Yasuhara \cite{MurakamiYasuhara2000}, is the minimum first Betti number of non-orientable surfaces smoothly and properly embedded in the 4-ball $D^4$ and bounded by $K$. Among low-crossing knots, this knot invariant is currently only known for knots with crossing number up to 9, see \cite{Jabuka and Kelly, Knotinfo}.
The main result of this paper is a complete calculation of $\gamma_4$ for all 165 knots with 10 crossings. 
\begin{theorem} \label{MainTheorem}
Of the 165 knots with crossing number 10 there are exactly 104 knots with $\gamma_4$ equal to 1, namely

$$10_1, \,10_3, \, 10_4, \,10_6, \, 10_7, \, 10_8, \,10_{11}, \,10_{12}, \, 10_{15}, \, 10_{16}, \, 10_{17}, \, 10_{20}, \, 10_{21}, \,10_{22}, \,10_{23}, \,10_{24}, \, 10_{27}, \,10_{29},$$
$$10_{30}, \,10_{31}, \,10_{35}, \,10_{38}, \,10_{39}, \,10_{40}, \, 10_{41}, \,10_{42}, \, 10_{43}, \,10_{44}, \, 10_{45}, \, 10_{48}, \, 10_{49}, \,10_{50}, \, 10_{51}, \,10_{52},  \, 10_{54},$$
$$10_{55}, \,10_{57},\,10_{59}, \,10_{62}, \,10_{64},\,10_{65}, \, 10_{66}, \,10_{67},\,10_{68}, \,10_{69}, \,10_{70}, \,10_{73}, \,10_{74}, \,10_{75}, \,10_{77}, \,10_{78} , \, 10_{80},$$
$$10_{82}, \,10_{83}, \,10_{87}, \,10_{89}, \, 10_{91}, \, 10_{93}, \,10_{94}, \,10_{97}, \,10_{99}, \,10_{101}, \,10_{102}, \,10_{103}, \,10_{105}, \,10_{106}, \, 10_{108}, \, 10_{110},$$
$$10_{111}, \,10_{116}, \,10_{117}, \,10_{118}, \, 10_{121}, \,10_{122}, \,10_{123}, \, 10_{124}, \, 10_{125}, \, 10_{126}, \, 10_{127}, \,10_{128}, \, 10_{129}, \, 10_{130}, \, 10_{131},$$
$$ 10_{133}, \, 10_{134}, \, 10_{137}, \,10_{139}, \,10_{140}, \,10_{142}, \,10_{143}, \, 10_{144}, \, 10_{145}, \, 10_{146}, \, 10_{147}, \, 10_{148}, \, 10_{150}, \, 10_{151}, \,10_{152},$$
$$10_{153}, \, 10_{154}, \ 10_{155}, \,10_{160}, \, 10_{161}, 10_{165}. $$
\\
There are 61 knots with $\gamma_4$ equal to 2, which are

$$10_2, \, 10_5, \, 10_9, \, 10_{10}, \, 10_{13}, \, 10_{14}, \, 10_{18}, \, 10_{19}, \, 10_{25}, \, 10_{26}, \, 10_{28}, \, 10_{32}, \, 10_{33}, \, 10_{34}, \, 10_{36}, \, 10_{37}, \, 10_{46},$$
$$10_{47}, \,10_{53}, \, 10_{56}, \, 10_{58}, \, 10_{60}, \, 10_{61}, \, 10_{63}, \, 10_{71}, \, 10_{72}, \, 10_{76}, \, 10_{79}, \, 10_{81}, \,10_{84}, \, 10_{85}, \, 10_{86}, \, 10_{88}, \, 10_{90},$$
$$10_{92}, \,10_{95}, \,10_{96}, \,10_{98}, \,10_{100}, \,10_{104}, \,10_{107}, \,10_{109}, \,10_{112}, \, 10_{113}, \,10_{114}, \,10_{115}, \, 10_{119}, \,10_{120}, \,10_{132},$$
$$10_{135}, \,10_{136}, \,10_{138}, \, 10_{141}, \, 10_{149}, \,10_{156}, \, 10_{157}, \, 10_{158}, \,10_{159}, \, 10_{162}, \, 10_{163}, \, 10_{164}.$$.
\\
\end{theorem}%
\section{Background}{\label{TAT}}
In this section we review needed background material for computing $\gamma_4$. We largely follow the outline from \cite{Jabuka and Kelly} and refer the intersted reader to consult said reference for more details.
\subsection{Upper Bound for $\gamma_4(K)$}\label{Theorem and technique} 
Proposition \ref{band move} below outlines how to find upper bounds on $\gamma_4(K)$. I relies on non-orientable band moves and we digress to define those first.
\begin{definition}(Definition 2.3. in \cite{Jabuka and Kelly})
A non-oriented band move on an oriented knot $K$ is the operation of attaching an oriented band $h=[0, 1]\times[0, 1]$ to $K$ along $[0, 1]\times\partial[0, 1]$ in such a way that the orientation of the knot agrees with that of $[0, 1]\times{0}$ and disagrees with that of $[0,1]\times\partial{1}$ (or vice versa), and then performing surgery on $h$, that is replacing the arcs $[0, 1]\times\partial[0, 1]\subseteq{K}$ by the arcs $\partial[0, 1]\times[0, 1]$.
\end{definition}
\begin{proposition}{\label{band move}} (Proposition 2.4. in \cite{Jabuka and Kelly})
If the knots $K$ and $K'$ are related by a non-oriented band move then $$\gamma_4(K)\leq\gamma_4(K')+1.$$
If a knot $K$ is related to a slice knot $K'$ by a non-oriented band move, then $\gamma_4(K)=1$.
\end{proposition}

\begin{figure}[h]
\centering
\begin{subfigure}[b]{0.35\textwidth}
\includegraphics[width=\textwidth]{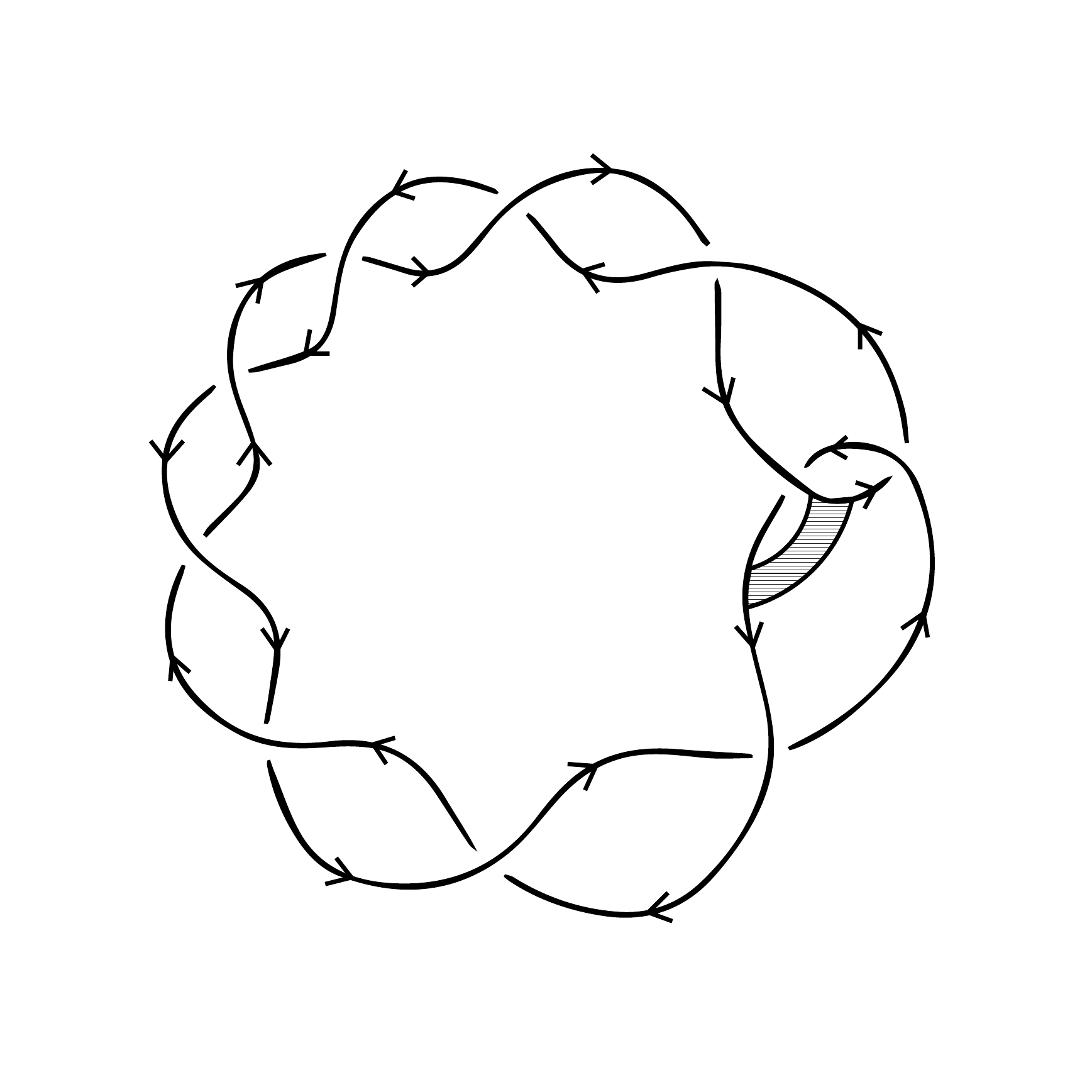}
\caption{the same orientation band move for $K=10_{1}$} \label{oa}
\end{subfigure}
\qquad \qquad \qquad \qquad
\begin{subfigure}[b]{0.35\textwidth}
\includegraphics[width=\textwidth]{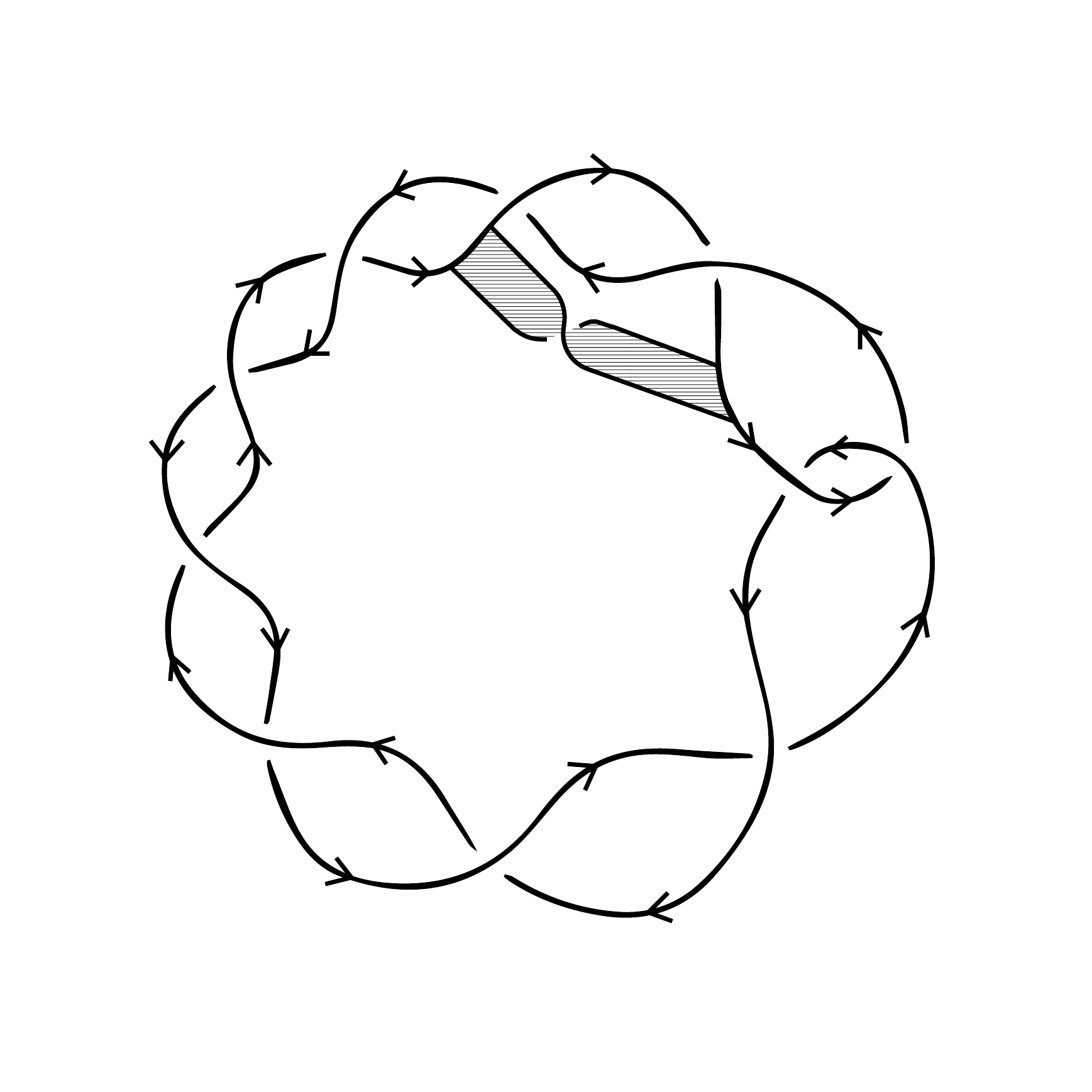}
\caption{opposite orientation band move for $K=10_{1}$}
\label{ob}
\end{subfigure}
\caption{Possible band moves for $K=10_{1}$}
\end{figure}
\begin{figure}[h]
\centering
\begin{subfigure}[b]{0.35\textwidth}
\includegraphics[width=\textwidth]{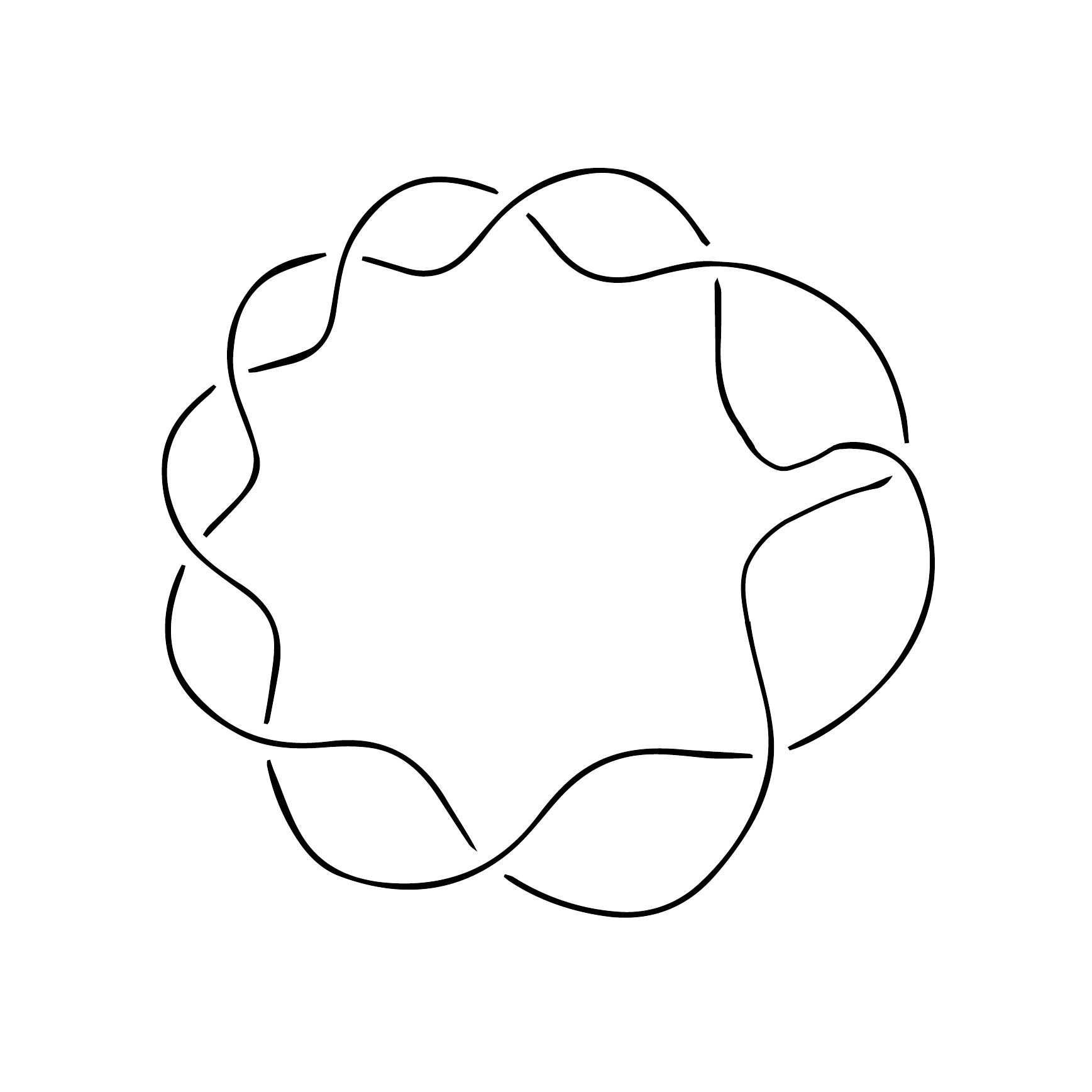}
\caption{$10_1\#h_1=9_1$} \label{fa}
\end{subfigure}
\qquad \qquad \qquad \qquad
\begin{subfigure}[b]{0.35\textwidth}
\includegraphics[width=\textwidth]{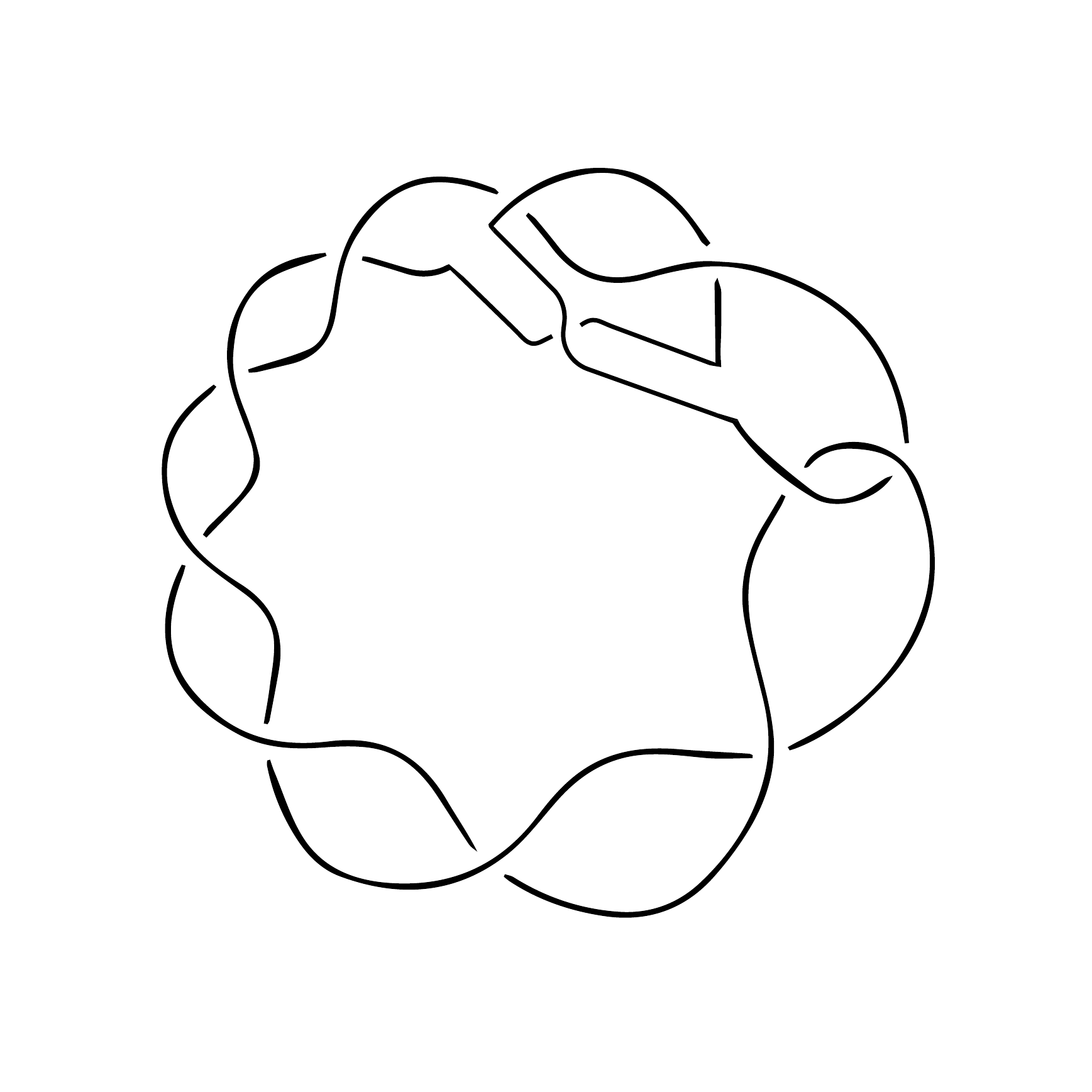}
\caption {$10_1\#h_2=6_1$}
\label{fb}
\end{subfigure}
\caption{Final results of adding band moves $h_1$ and $h_2$ to the knot $K=10_{1}$}
\end{figure}

\begin{example}
Figure \ref{oa} and \ref{ob} shows two different band moves of the knot $10_1$, leading to the knots $9_1$ and $6_1$ respectively, as in Figure \ref{fa} and \ref{fb}. The upper bound for $\gamma_4(10_1)$ resulting from these band moves are 2 and 1 respectively, according to Proposition \ref{band move}. Thus we can conclude that $\gamma_4(10_1)=1$. Thus, different band moves on the same knot $K$ can result in different upper bounds for $\gamma_4(K)$.
 \end{example}
\par
\begin{figure}[h]
\centering
\includegraphics[width=8cm]{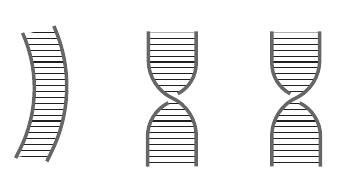}
\put(-220,-12){$\stackrel{0}{\longrightarrow}$}
\put(-130,-12){ $\stackrel{-1}{\longrightarrow}$}
\put(-45,-12){ $\stackrel{1}{\longrightarrow}$}
\caption{Different kinds of band moves.}
\label{bands}
\end{figure}

In this paper $0$ on top of an arrow represents a band move without twist, $1$ and $-1$ represent right-handed and left-handed twists respectively as shown in Figure \ref{bands}.
\par
As mentioned in Proposition \ref{band move} if we can find a non-oriented band move to a slice knot for a knot $K$, we conclude that $\gamma_4(K)=1$. In this paper, we will show that there are exactly 90 knots with 10 crossings with band moves to slice knots. Accordingly these 90 knots with all have $\gamma_4=1$. 
\subsection{Lower Bounds for $\gamma_4$}
Given a knot $K$, let $\sigma (K)$ denote its signature and let Arf$(K)$ denote the knot's Arf invariant. The following proposition is proved in \cite{GilmerLivingston}.
\begin{proposition}\label{lowerbound}
Given a knot $K$, if $\sigma(K)+4\cdot \text{Arf}(K)\equiv{4}\pmod 8$, then $\gamma_4(K)\geq{2}$.
\end{proposition}
We are able to use this proposition in conjunction with Proposition \ref{band move} to show that many 10-crossing knots have $\gamma_4$ equal to 2. Specifically, any knot $K$ which admits a band move to a knot $K'$ with $\gamma_4(K')=1$ and for which $\sigma(K)+4\cdot \text{Arf}(K)\equiv{4}\pmod 8$ will necessarily lead to $\gamma_4(K) = 2$. 

Among the 165 knots with 10 crossings there are exactly 43 knots which satisfy the congruence $\sigma(K)+4\cdot \text{Arf}(K)\equiv{4} \pmod 8$. We will show in the next section that all 43 knots in this group admit non-orientable band moves to knots with $\gamma_4$ equal to 1 and so all such knots have non-orientable 4-genus equal to 2. 

In the remaining cases of $\sigma(K)+4\cdot \text{Arf}(K)\equiv{0, \pm{2}}\pmod 8$, we need additional methods to bound  $\gamma_4$ from below. These extra methods are facilitated by the Goeritz Matrix. Our definition here follows that from \cite{Goeritz}, rather than Goertiz's original version. 

\begin{figure}[h]
\centering
\includegraphics[width=8cm]{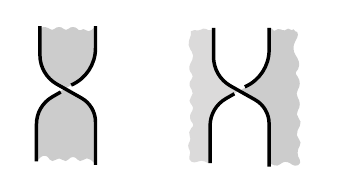}
\put(-204,-10){$\eta = -1.$}
\put(-88,-10){$\eta = +1.$}
\caption{Weights of a crossing in a checkerboard coloring of a knot projection.}
\label{crossing}
\end{figure}
Recall that a knot diagram admits two types of checkerboard colorings. Figure \ref{crossing} shows how to associated a weight of $\pm 1$ to a crossing in a checkerboard coloroing. The PreGoertiz matrix $PG$ asscociated to a chosen checkboard coloring of diagram of a knot $K$ is defined as follows. Let $A_1, A_2, ..., A_n$ be the white regions in said diagram, then $PG$ is the  $n\times n$ matrix 
$$PG = \begin{bmatrix}
g_{1,1}& g_{1,2} & \dots & g_{1,n} \\
\vdots & \vdots  & \ddots &  \vdots\\
g_{n,1} & g_{n,2} & \dots & g_{n,n} &
\end{bmatrix}$$
where for all $1\leq{i}\neq{j}\leq{n}$, $g_{ij} = -\sum _{p\in A_i \cap A_j} \eta (p)$ and $g_{ii}=-\sum_{j\neq{i}}g_{ij}$. Here $A_i\cap A_j$ denotes the set of double points in the knot projections that are incident to both regions $A_i$ and $A_j$, while $\eta(p)$ refers to the weight function $\eta$ from Figure \ref{crossing}. Note that $PG$ is symmetric. 

For any index $k\in \{1, \dots, n\}$, by removing the $k-th$ row and $k$-th column from $PG$, we obtain an $(n-1)\times (n-1)$ matrix $G$ called a {\em Goeritz matrix}:
$$G = \begin{bmatrix}
g_{1,1}& g_{1,2} &\dots & g_{1,k-1}&g_{1,k+1}&\dots & g_{1,n} \\
\vdots  & \vdots  & \vdots  & \vdots  & \vdots  & \vdots  & \vdots \\
g_{k-1,1}&g_{k-1,2}&\dots &g_{k-1,k-1}&g_{k-1,k-2}&\dots&g_{k-1,n}\\
g_{k+1,1}&g_{k+1,2}&\dots&g_{k+1,k-1}&g_{k+1,k+1}&\dots&g_{k+1,n}\\
\vdots & \vdots & \vdots & \vdots & \vdots & \vdots & \vdots  \\
g_{n,1} & g_{n,2}&...&g_{n,k-1}&g_{n,k+1}&...& g_{n,n} &
\end{bmatrix}.$$
It is well known that $\det(K)=|\det(G)|$. It is also well understood that if a knot is alternating, the Goeritz matrices associated to its two checkerboard colorings are definite, one positive definite and the other negative definite. Of the 165 knots with 10 crossings, there are 123 alternating knots whose Goertiz  matrices have said definiteness properties, something we shall rely on in Theorems \ref{group2,3} and \ref{last group}. We will explain our strategy for the remaining 42 non-alternating knots in Section \ref{NAK}.

Let $K$ be an alternating knot, and let $G_+$ and $G_-$ be Goertiz matrices associated to the two checkerboard colorings of some projection of $K$.  Theorems \ref{group2,3} and \ref{last group} rely on a branched covering construction that creates a pair 4-manifolds with boundary the 2-fold cover of $S^3$ branched along $K$, and with intersection forms given by $G_+$ and $G_-$. One can then boundary sum these 4-manifolds with the 2-fold cover of $D^4$ branched over a hypothetical non-oriented surface the knot bounds in $D^4$, obtaining a smooth, closed, definitie 4-manifold, to which  Donaldson's celebrated Diagonalization Theorem applies. The results of Theorems \ref{group2,3} and \ref{last group} are a direct consequences of this construction and the Diagonalization theorem.  We omit further details and refer the interested reader to \cite{Jabuka and Kelly}.
 
 \begin{theorem}{\label{group2,3}} (Theorem 2.10 in \cite{Jabuka and Kelly})(Case of $\sigma(K)+4\cdot \text{Arf}(K)\equiv\pm{2} \pmod 8$). Let $K$ be a knot with $\sigma(K)+4\cdot \text{Arf}(K)\equiv{2\epsilon} \pmod 8$ for a choice of $\epsilon\in\{\pm{1}\}$. Assume that $K$ admits a checkerboard coloring for which the associated Goeritz form $G$ is $-\epsilon$-definite. If no embedding exists of $G\oplus[l]$ into the $\epsilon$-definite diagonal form $(\mathbb{Z}^{rank(G)+1}, \epsilon Id)$ for any divisor $\ell$ of $\det K$ with $\det K/\ell$ a square, then $\gamma_4(K)\geq{2}$.
\end {theorem}
%
There are exactly 78 10-crossing knots which satisfy the congruence $\sigma(K)+4\cdot \text{Arf}(K)\equiv\pm{2} \pmod 8$. Precisely 36 knots with $\sigma(K)+4\cdot \text{Arf}(K)\equiv{2} \pmod 8$, and 42 knots with $\sigma(K)+4\cdot \text{Arf}(K)\equiv{-2} \pmod 8$. We will show in the next section that, among all 36 knots in the first congruence, we have only 5 knots with no embeddings and also among all 42 knots in the second congruence there are exactly 5 knots with no embeddings. The remaining knots yet again admite non-orientable band moves to slice knots and thus have $\gamma_4=1$. 
\begin{theorem}{\label{last group}} (Theorem 2.12 in \cite{Jabuka and Kelly}) (Case of $\sigma(K)+4\cdot \text{Arf}(K)\equiv{0} \pmod 8$). Let $K$ be a knot with $\sigma(K)+4\cdot \text{Arf}(K)\equiv{0} \pmod 8$ and assume that the Goeritz matrices $G_{\pm}$ of $K$, associated to the two possible checkerboard colorings of a knot projection $D$ of $K$, are positive and negative definite respectively (with the subscript $\pm$ indicating the definiteness type of $G_{\pm}$). 
If no embedding exists of $G_{+}\oplus[\ell]$ into the positive form $(\mathbb{Z}^{rank(G_{+})+1}, Id)$, and no embedding exists of $G_{-}\oplus[\ell]$ into the negative form $(\mathbb{Z}^{rank(G_{-})+1}, -Id)$, for any divisor $\ell$ of $\det K$ with $K/\ell$ a square, then $\gamma_4(K)\geq{2}$.
\end{theorem}
There are 44 knots with 10 crossings that satisfy the congruence $\sigma(K)+4\cdot \text{Arf}(K)\equiv{0} \pmod 8$, and in the next section we will prove that among these there are only 2 knots with no embeddings in this group. The remaining knots admit band moves to slice knots and thus have $\gamma_4=1$.  
\section{Computation of $\gamma_4$}
This section includes the computations of $\gamma_4(K)$ for all knots $K$ with 10 crossings. We study all 165 knots with 10 crossings in 6 different subsections. Section \ref{sliceknots} focuses on slice knots and knots related to slice knots with band moves. Section \ref{NAK} is about non-alternating knots. In Sections \ref{first group}, \ref{second group}, \ref{third group} and \ref{fourth group} we discuss knots $K$ based on the congruance $\sigma(K)+4\cdot \text{Arf}(K)\equiv{4}, {-2}, {2}$ or ${0}$ $\pmod 8$ respectively. Theorems and techniques we mentioned in Section \ref{TAT} and Proposition \ref{band move} lead into finding the exact value of $\gamma_4(K)$ for each knot. 
\subsection {Slice knots and band moves to slice knots} \label{sliceknots} Among all 10-crossing knots, we have 14 slice knots
$$10_3, \, 10_{22}, \, 10_{35}, \, 10_{42}, \, 10_{48}, \, 10_{75}, \, 10_{87}, \, 10_{99}, \, 10_{123}, \, 10_{129}, \, 10_{137}, \, 10_{140}, \, 10_{153}, \, 10_{155}$$
 $\gamma_4$ for each of these 14 knots is equal to 1.
 \par
\begin{figure}[h]

	\centering
	\begin{subfigure}[b]{0.3\textwidth}
		\includegraphics[width=\textwidth]{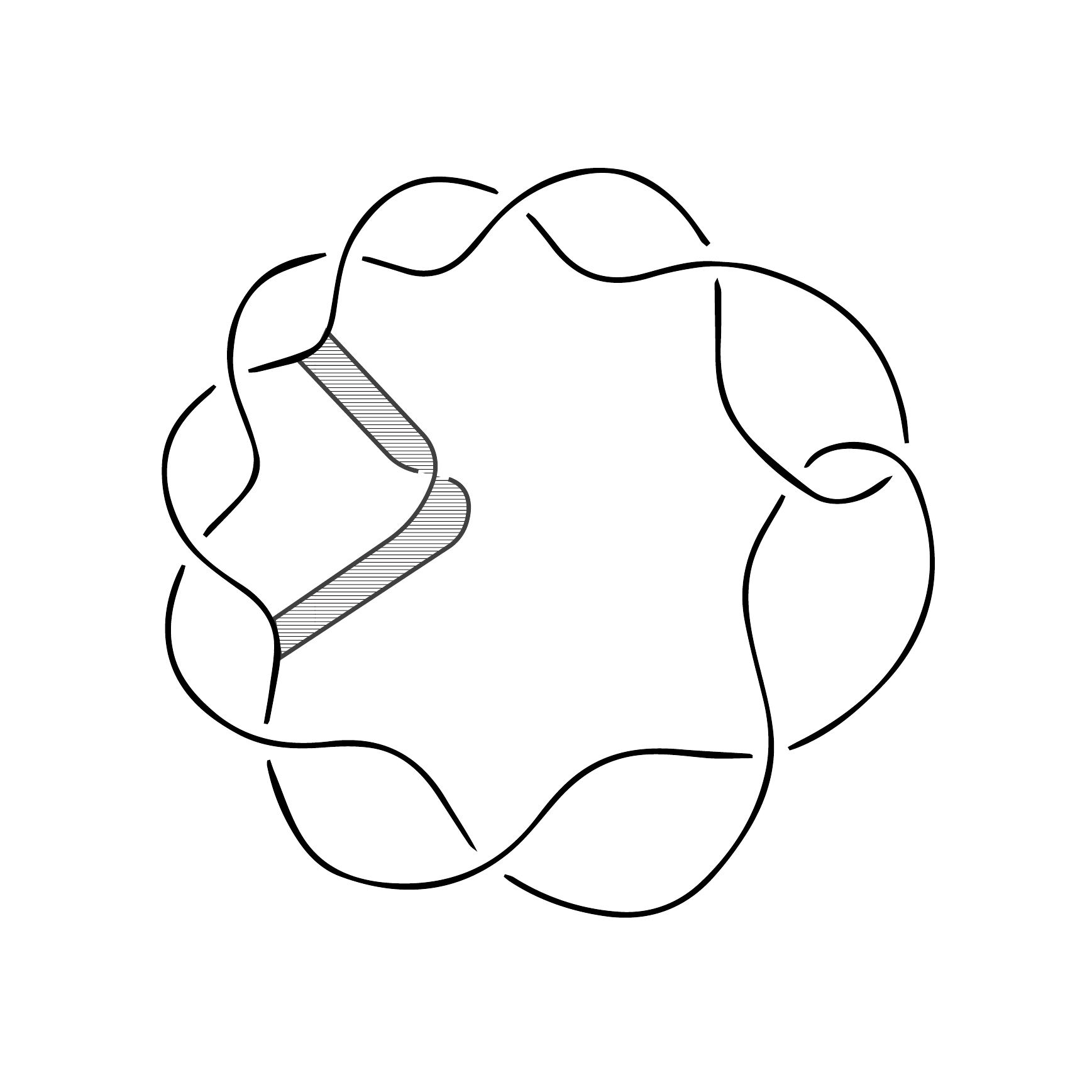}
		\caption{$10_{1}\stackrel{1}{\longrightarrow} 6_{1}$}
		\label{FigureFor10-1}
	\end{subfigure}
	~
	\begin{subfigure}[b]{0.3\textwidth}
		\includegraphics[width=\textwidth]{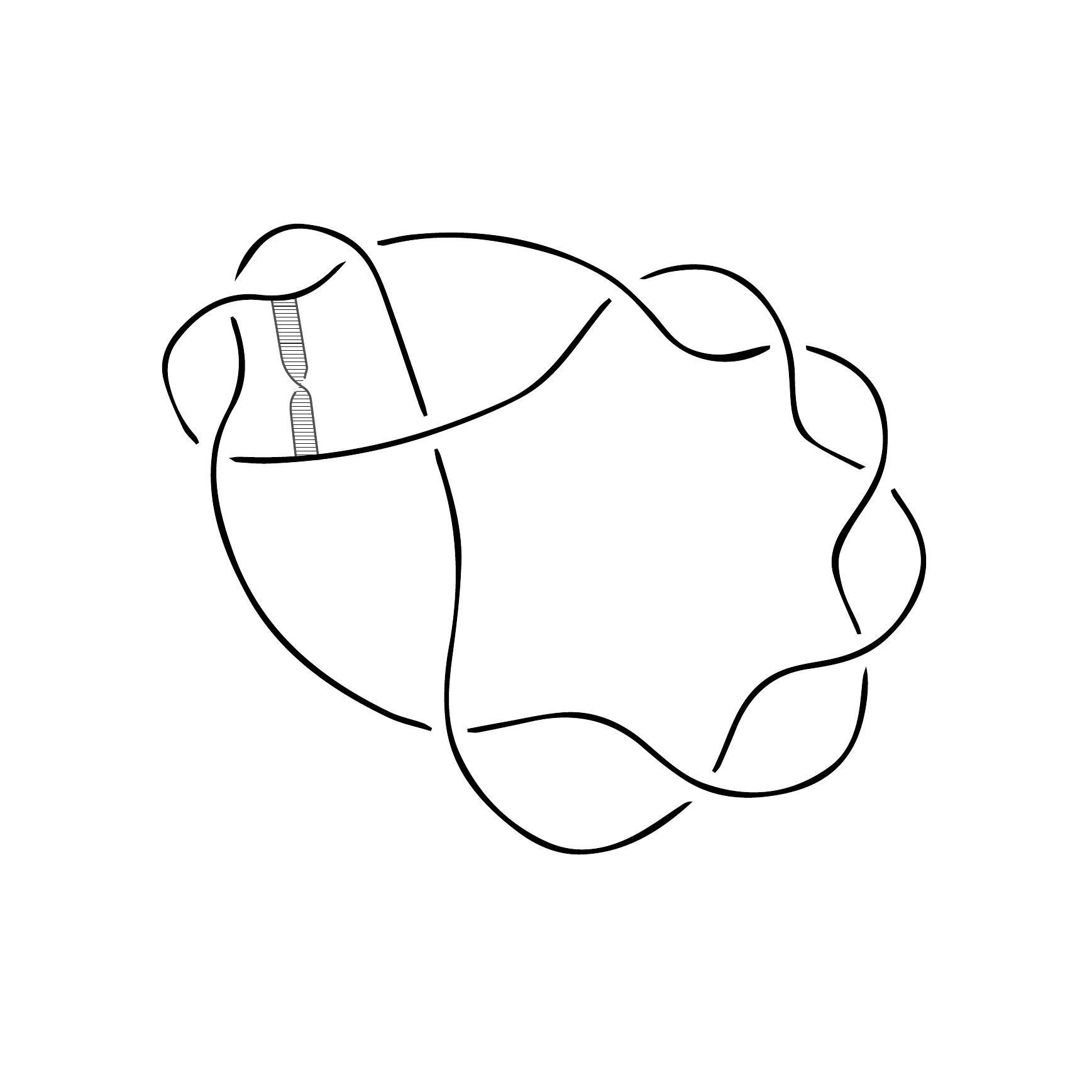}
		\caption{$10_{4}\stackrel{-1}{\longrightarrow} 0_{1}$}
		\label{FigureFor10-4}
	\end{subfigure}
	~
	\begin{subfigure}[b]{0.3\textwidth}
		\includegraphics[width=\textwidth]{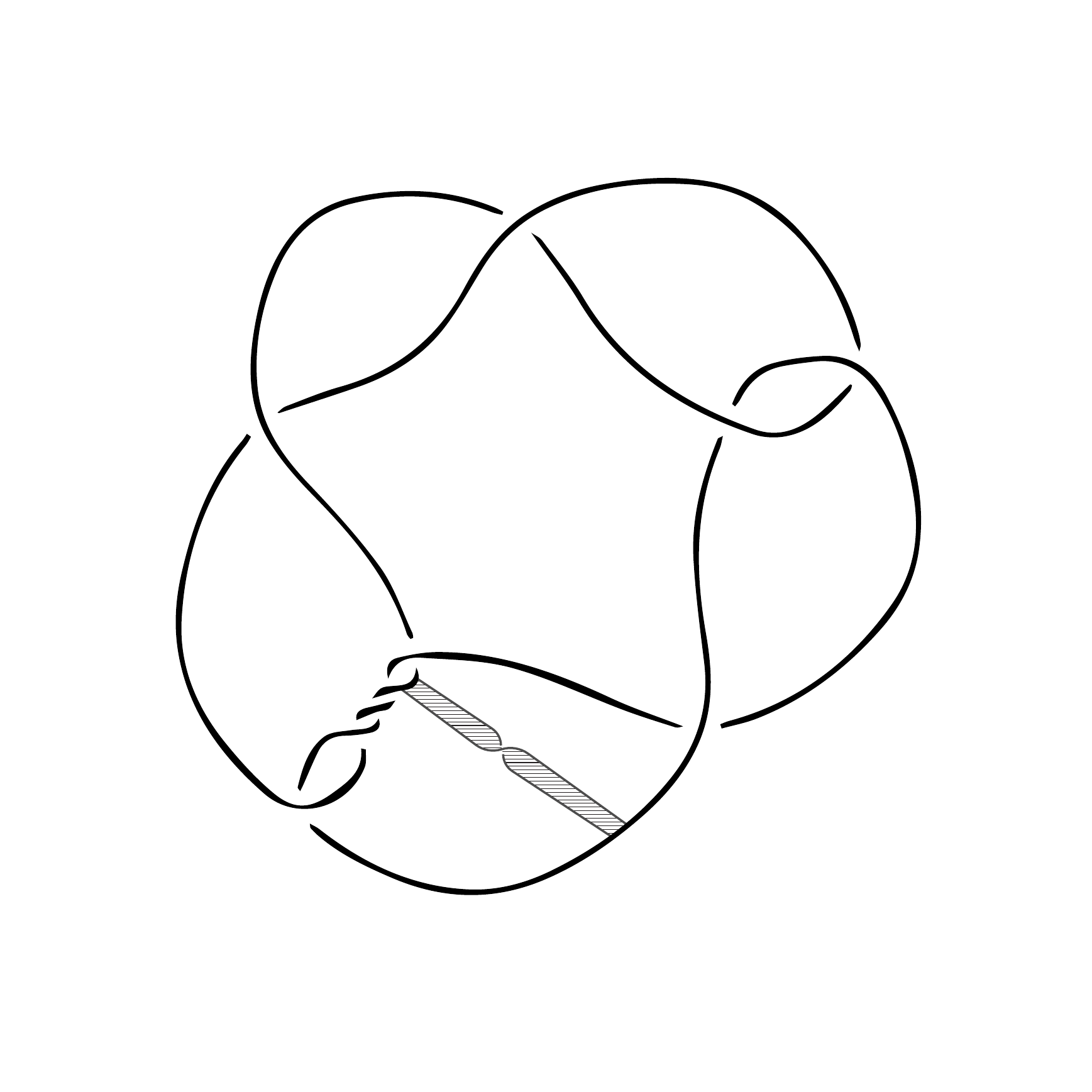}
		\caption{$10_{6}\stackrel{-1}{\longrightarrow} 0_{1}$}
		\label{FigureFor10-6}
	\end{subfigure}
	\vskip3mm
	\begin{subfigure}[b]{0.3\textwidth}
		\includegraphics[width=\textwidth]{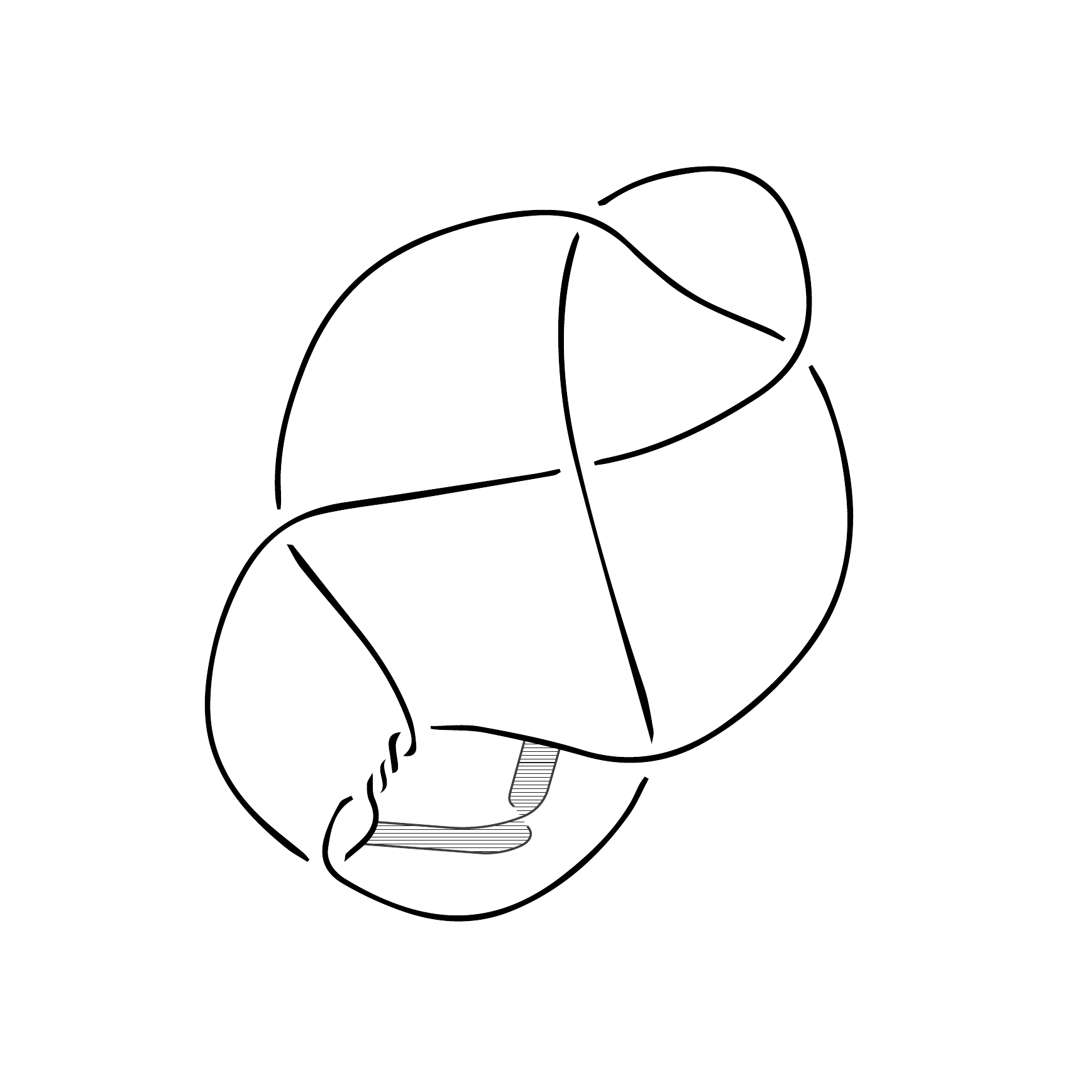}
		\caption{$10_{7}\stackrel{1}{\longrightarrow} 0_{1}$}
		\label{FigureFor10-7}
	\end{subfigure}
	~
	\begin{subfigure}[b]{0.3\textwidth}
		\includegraphics[width=\textwidth]{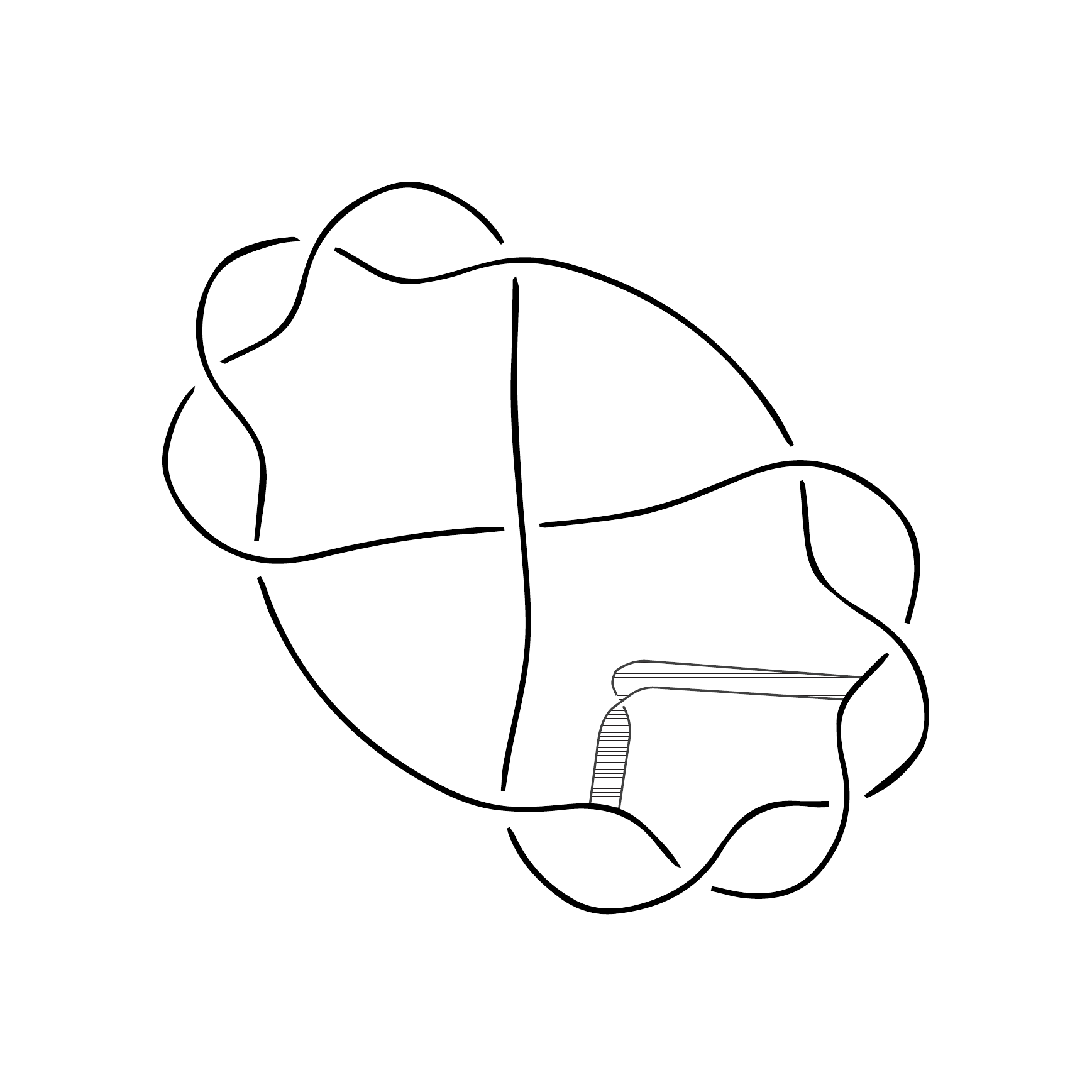}
		\caption{$10_{8}\stackrel{1}{\longrightarrow} 6_{1}$}
		\label{FigureFor10-8}
	\end{subfigure}%
	    \begin{subfigure}[b]{0.3\textwidth}
		\includegraphics[width=\textwidth]{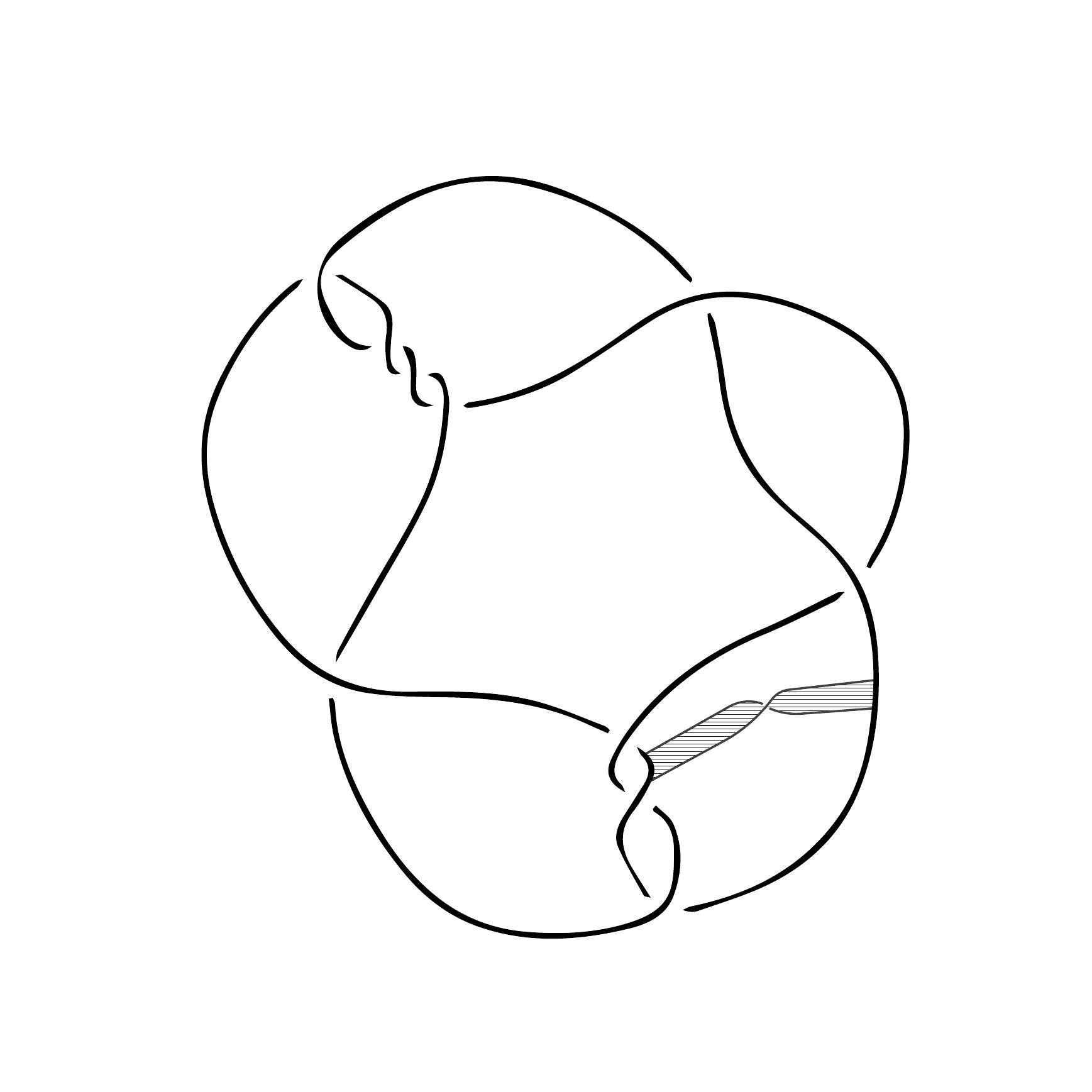}
		\caption{$10_{11}\stackrel{-1}{\longrightarrow} 6_{1}$}
		\label{FigureFor10-11}
	\end{subfigure}
	~
	\vskip3mm
	~
	\begin{subfigure}[b]{0.3\textwidth}
		\includegraphics[width=\textwidth]{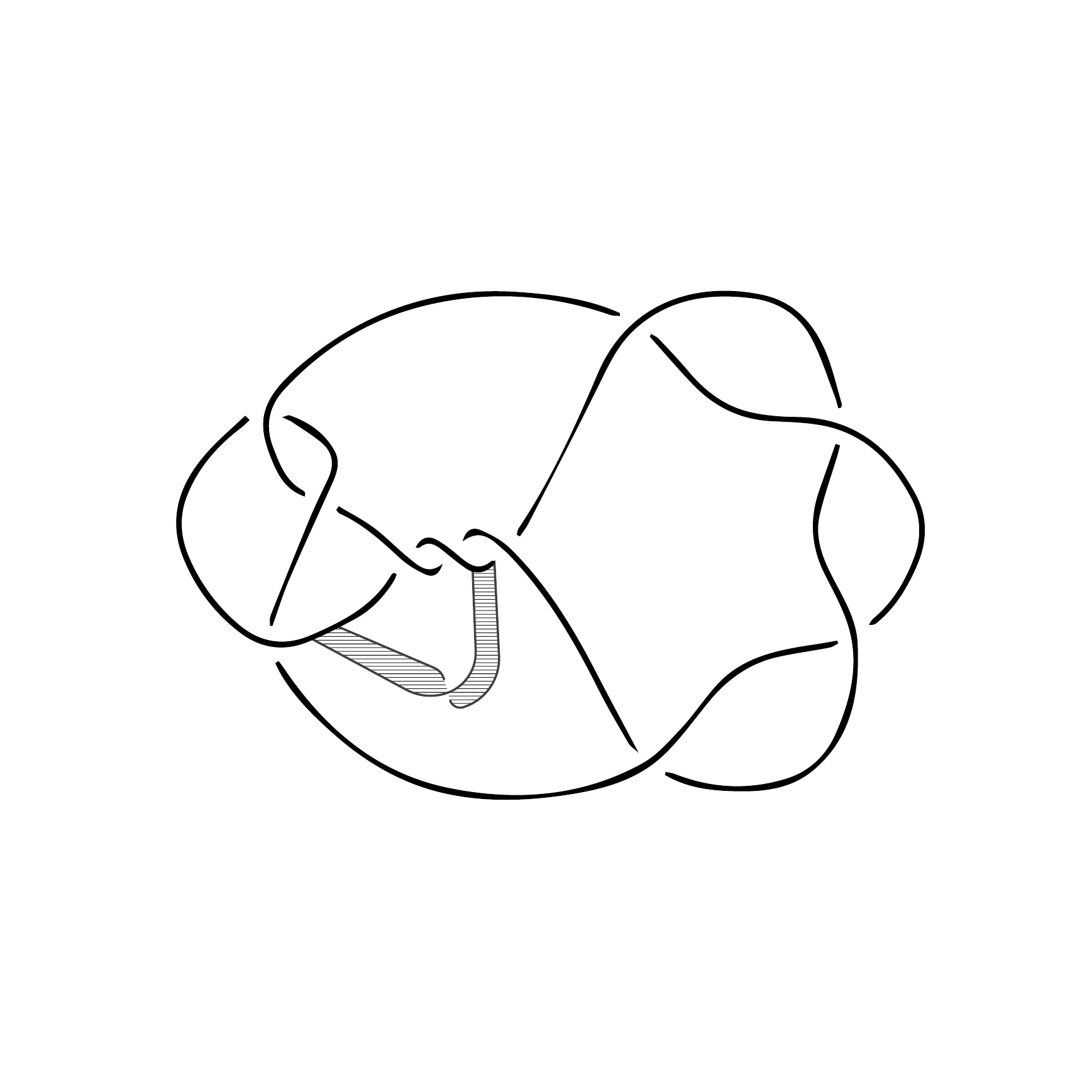}
		\caption{$10_{12}\stackrel{-1}{\longrightarrow} 0_{1}$}
		\label{FigureFor10-12}
	\end{subfigure}
	~
	\begin{subfigure}[b]{0.3\textwidth}
		\includegraphics[width=\textwidth]{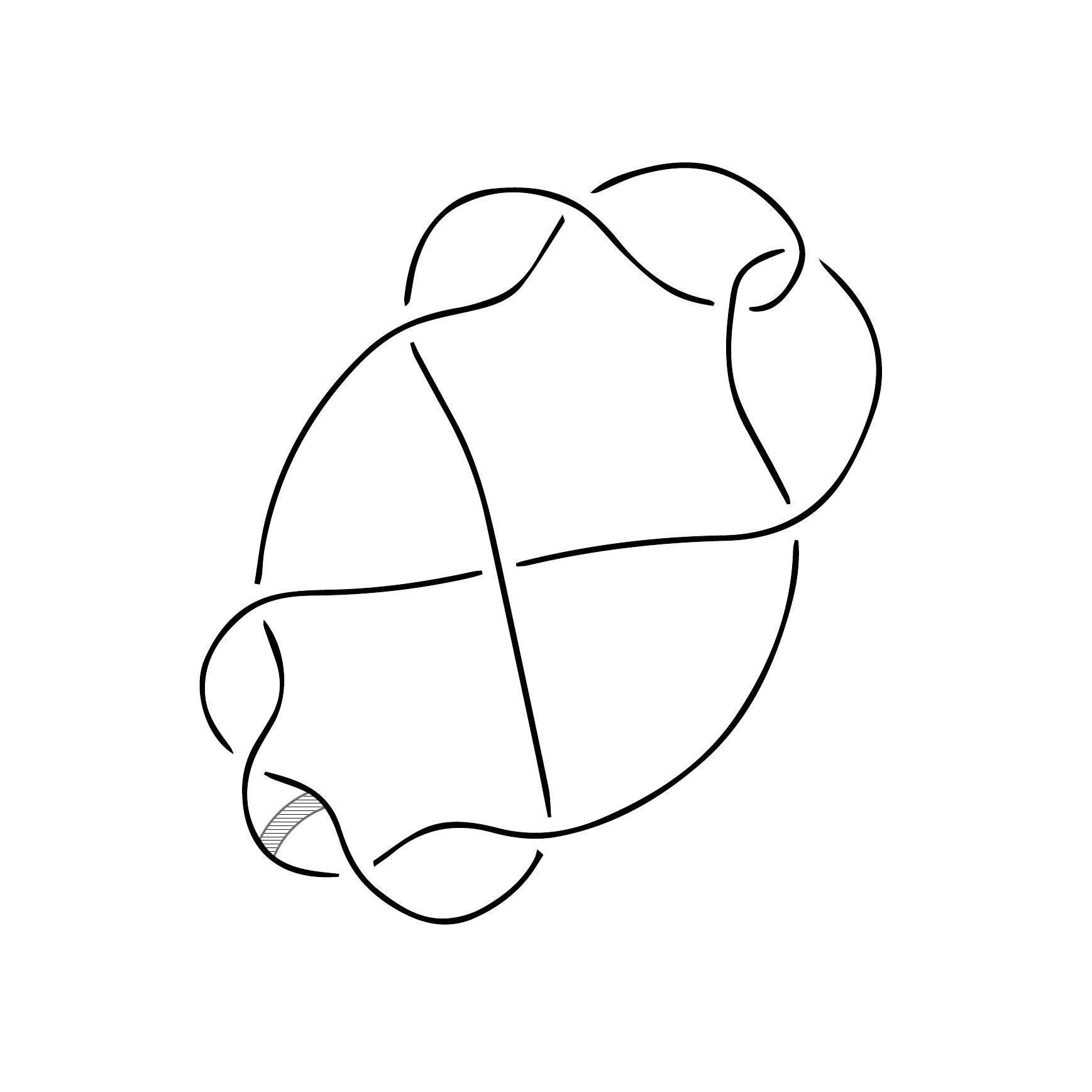}
		\caption{$10_{15}\stackrel{0}{\longrightarrow} 6_{1}$}
		\label{FigureFor10-15}
	\end{subfigure}
	~
	\begin{subfigure}[b]{0.3\textwidth}
		\includegraphics[width=\textwidth]{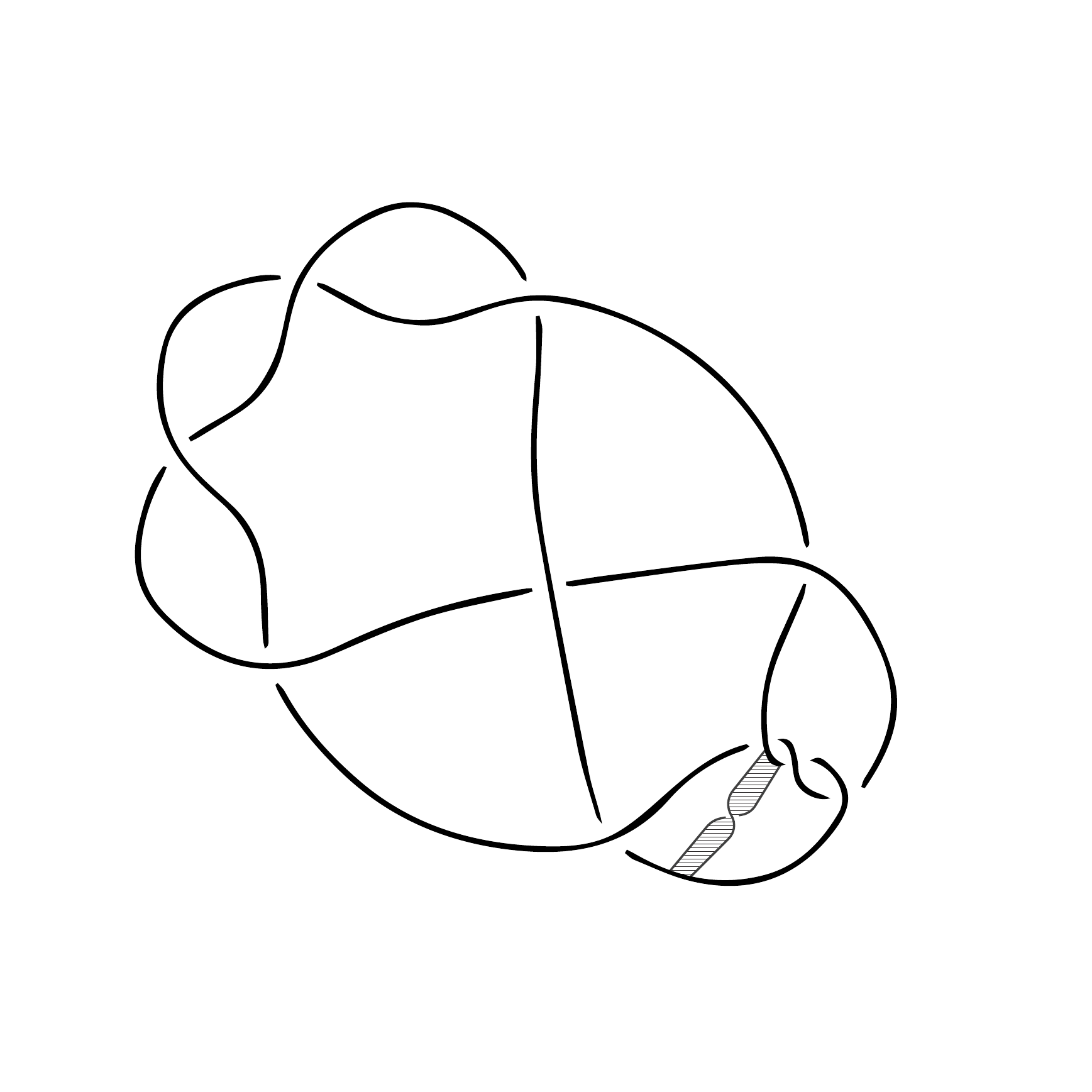}
		\caption{$10_{16}\stackrel{-1}{\longrightarrow} 6_{1}$}
		\label{FigureFor10-16}
	\end{subfigure}
	~
	\vskip3mm
	\caption{Non-oriented band moves from the knots $10_{1}$, $10_{4}$, $10_{6}$, $10_{7}$ , $10_{8}$, $10_{11}$, $10_{12}$, $10_{15}$, $10_{16}$ to slice knots.}
	\label{slice1}
\end{figure}
\begin{figure}[h]
	\centering
		\begin{subfigure}[b]{0.3\textwidth}
		\includegraphics[width=\textwidth]{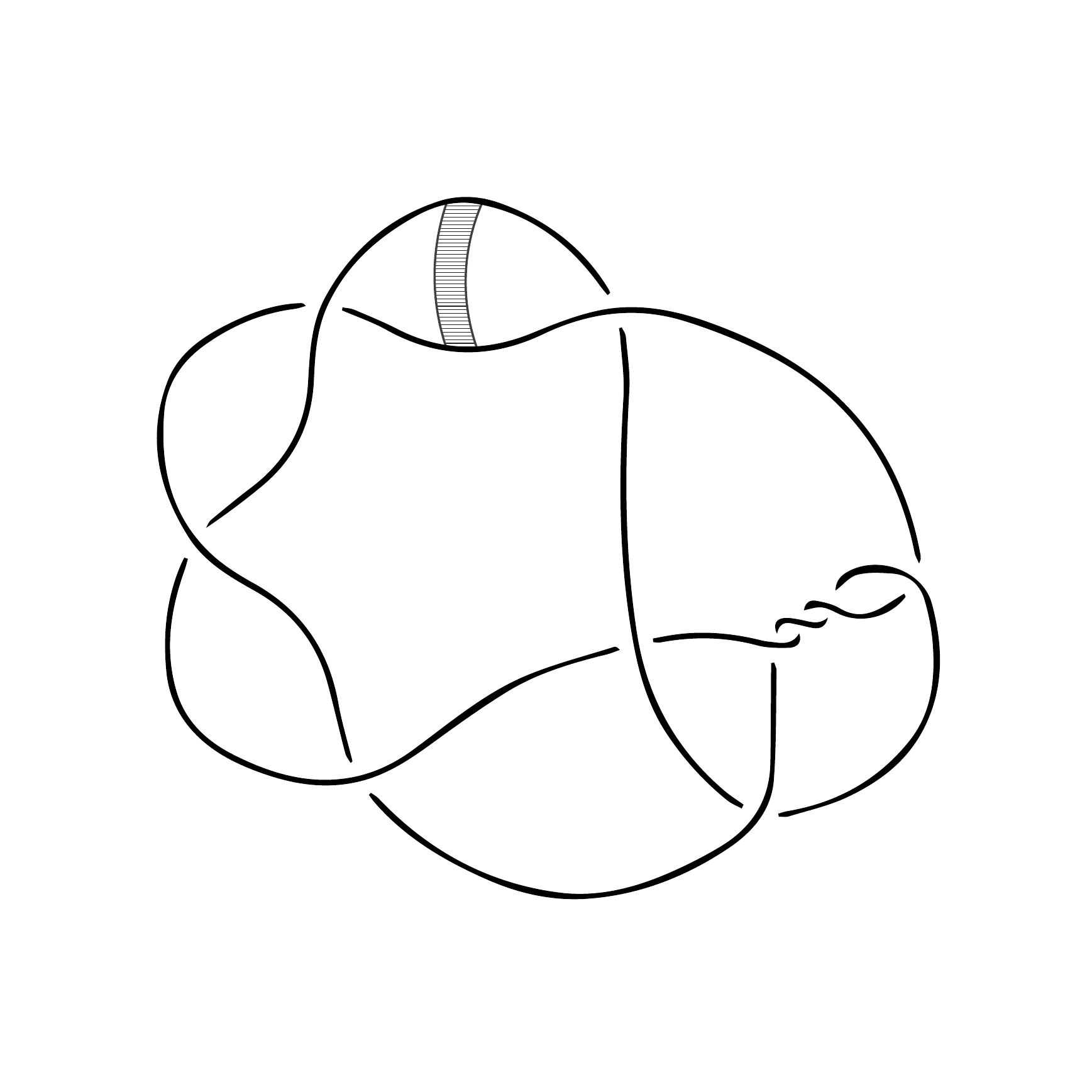}
		\caption{$10_{17}\stackrel{0}{\longrightarrow} 6_{1}$}
		\label{FigureFor10-17}
	\end{subfigure}
	~
	\begin{subfigure}[b]{0.3\textwidth}
		\includegraphics[width=\textwidth]{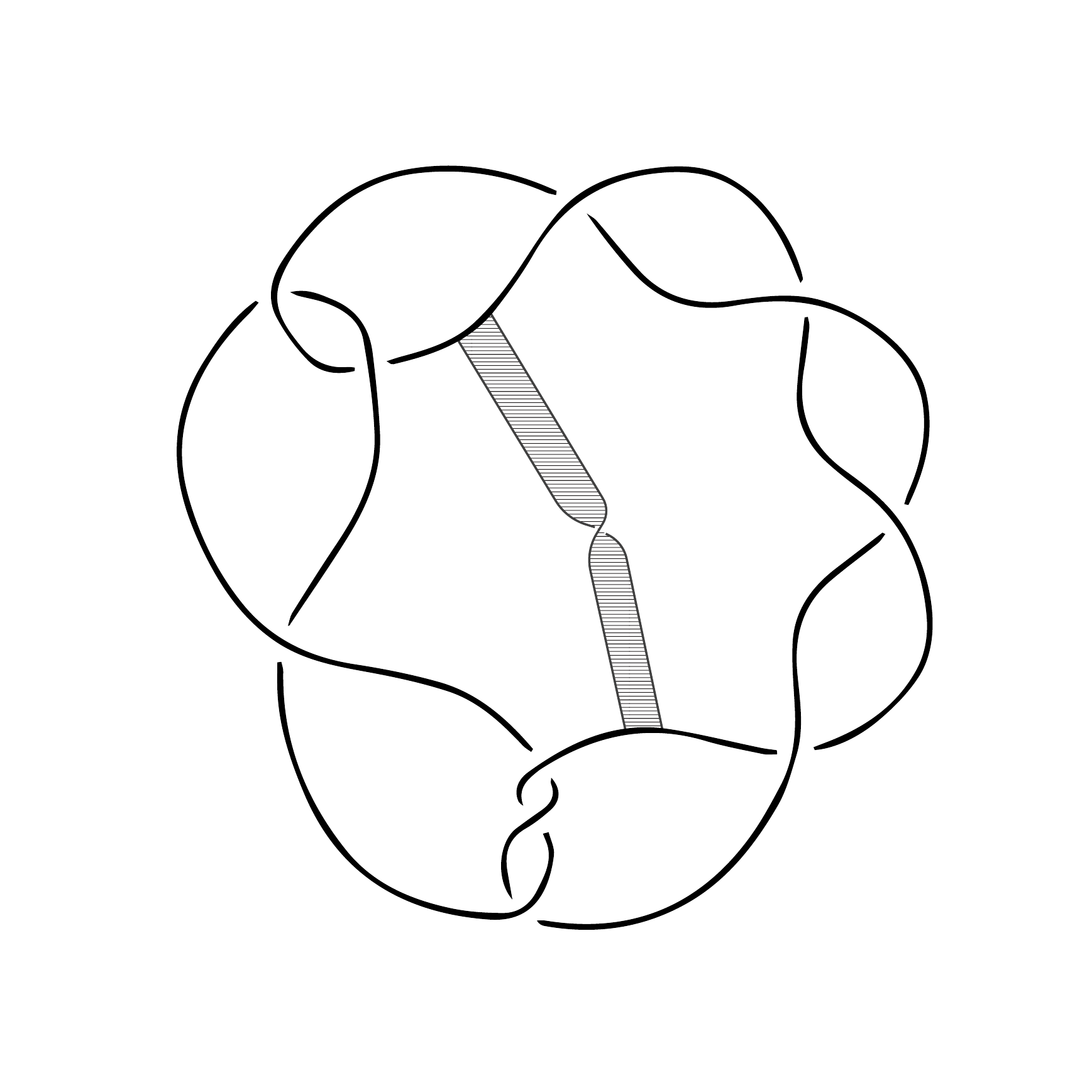}
		\caption{$10_{20}\stackrel{1}{\longrightarrow} 8_{20}$}
		\label{FigureFor10-20}
	\end{subfigure}
	~
	\begin{subfigure}[b]{0.3\textwidth}
		\includegraphics[width=\textwidth]{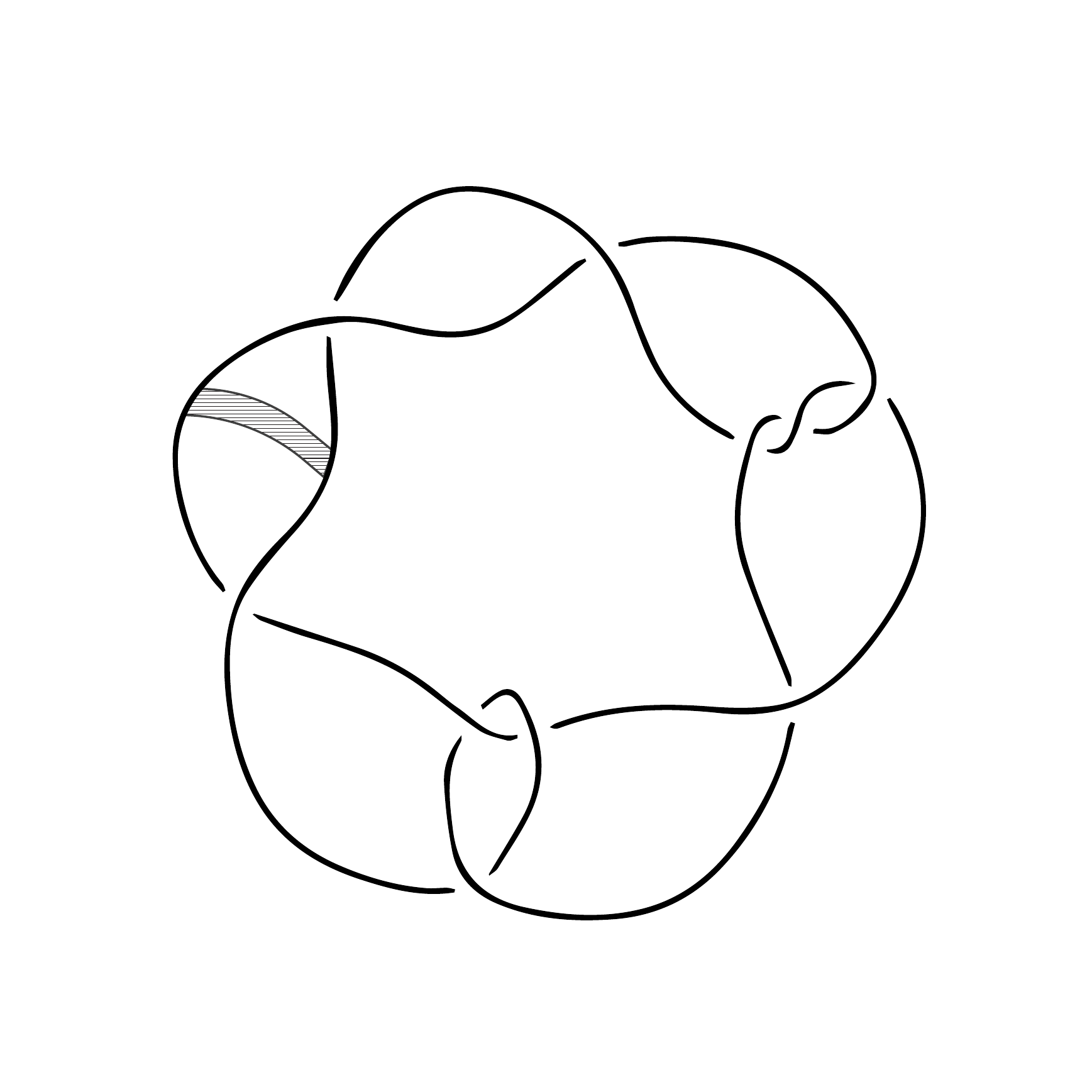}
		\caption{$10_{21}\stackrel{0}{\longrightarrow} 3_{1}\#-3_1$}
		\label{FigureFor10-21}
	\end{subfigure}
	\vskip3mm
	\begin{subfigure}[b]{0.3\textwidth}
		\includegraphics[width=\textwidth]{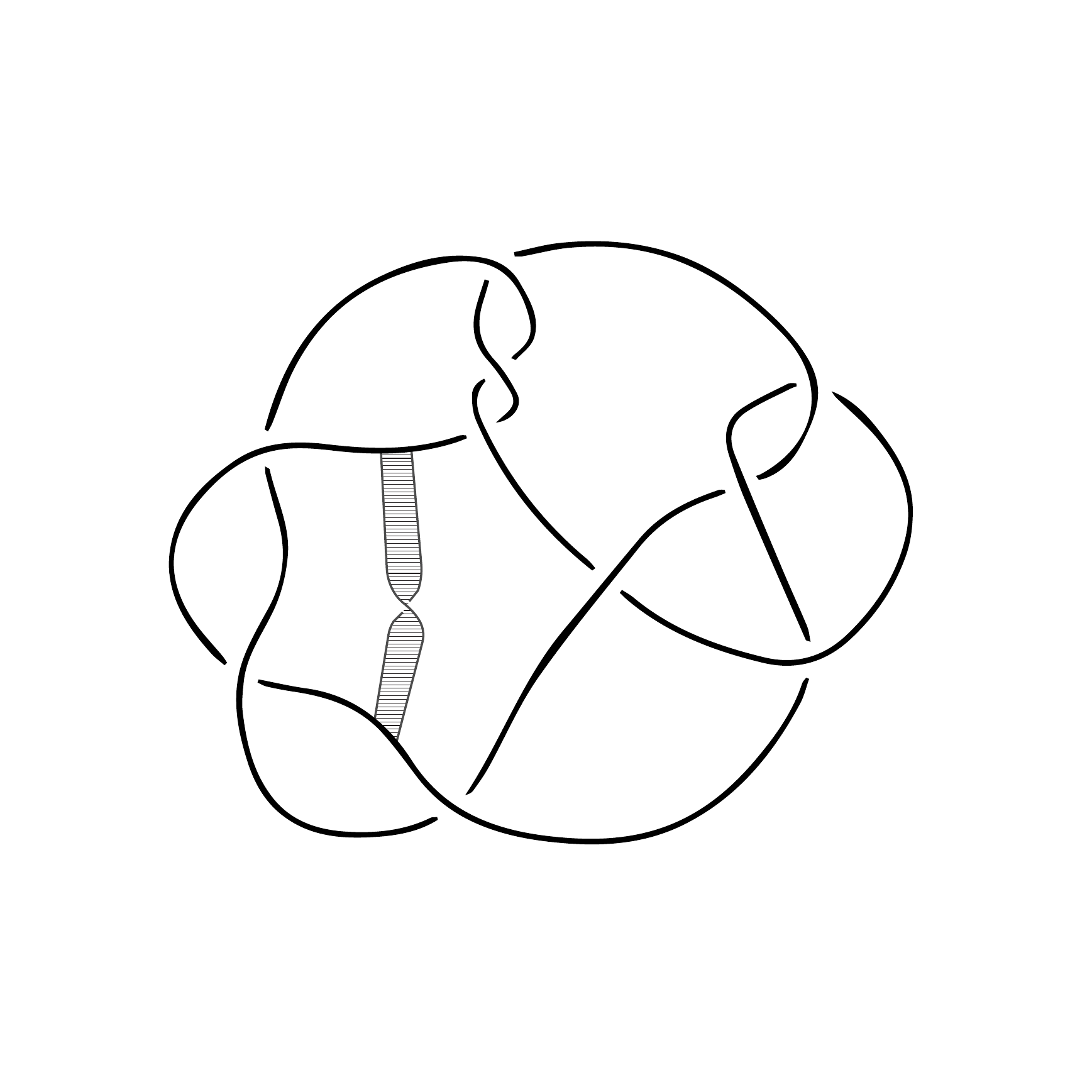}
		\caption{$10_{23}\stackrel{-1}{\longrightarrow} 0_{1}$}
		\label{FigureFor10-23}
	\end{subfigure}
	~
	\begin{subfigure}[b]{0.3\textwidth}
		\includegraphics[width=\textwidth]{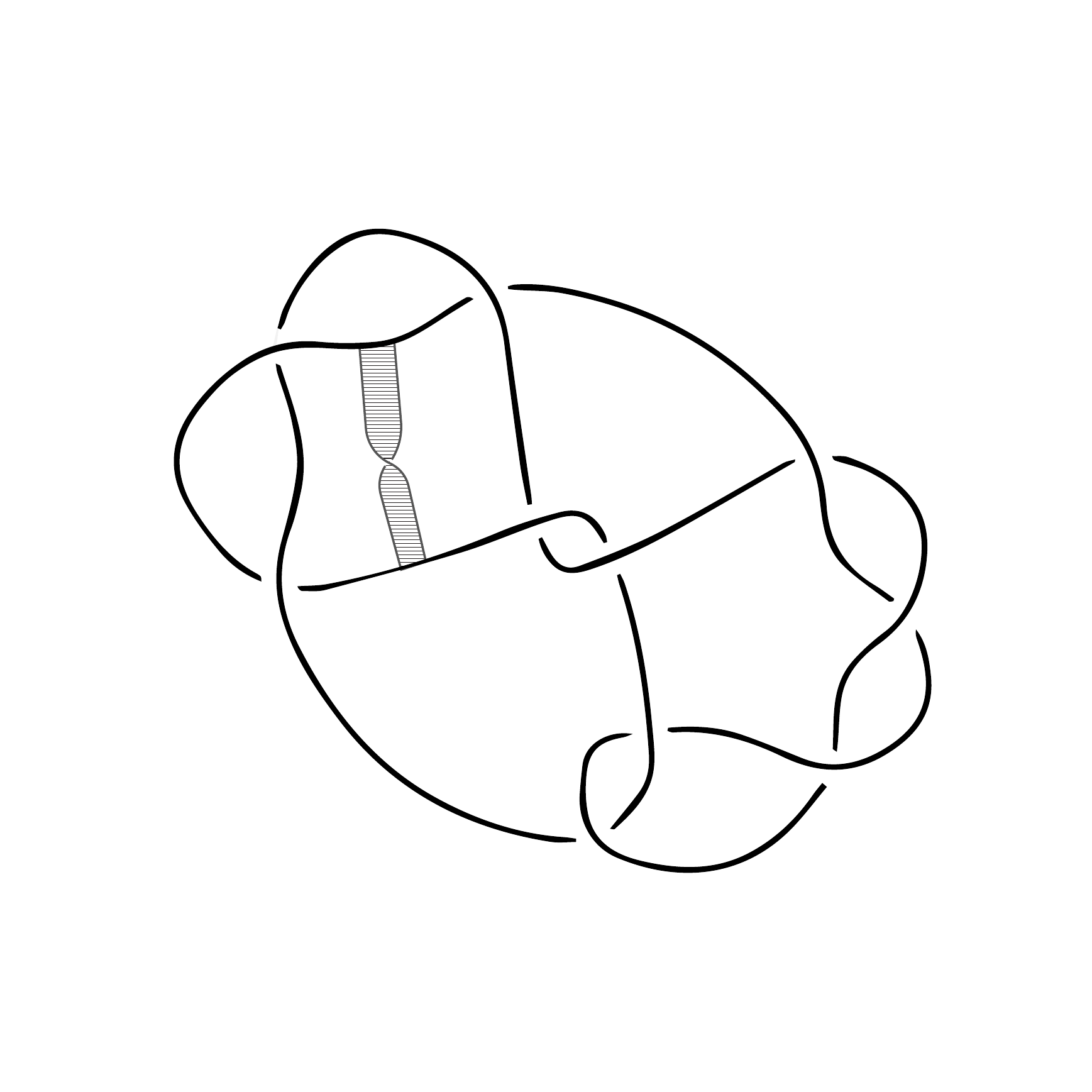}
		\caption{$10_{24}\stackrel{-1\phantom{i}}{\longrightarrow} 6_{1}$}
		\label{FigureFor10-24}
	\end{subfigure}
	~
	\begin{subfigure}[b]{0.3\textwidth}
		\includegraphics[width=\textwidth]{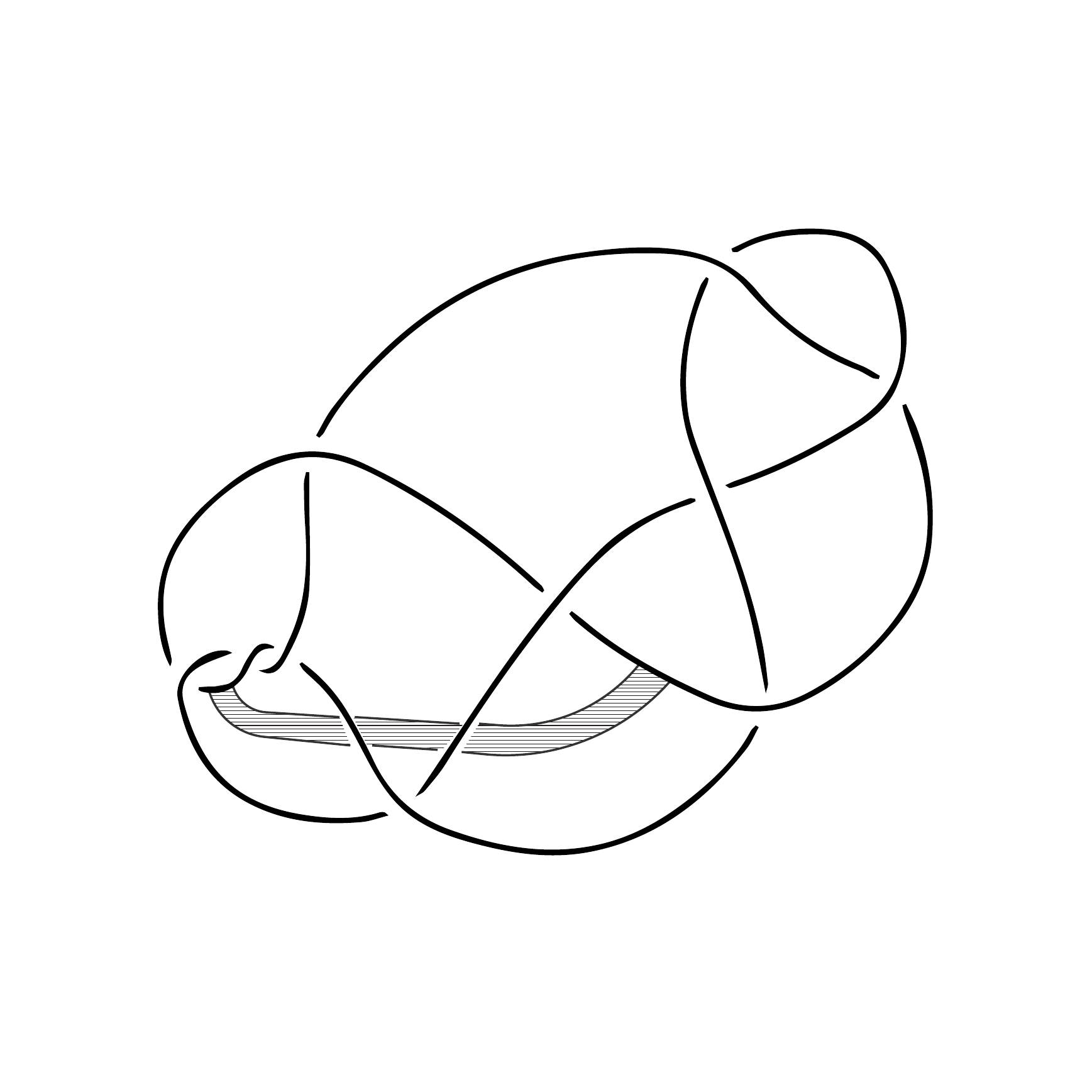}
		\caption{$10_{27}\stackrel{0}{\longrightarrow} 8_{8}$}
		\label{FigureFor10-27}
	\end{subfigure}
	~     
	\vskip3mm
	~
	\begin{subfigure}[b]{0.3\textwidth}
		\includegraphics[width=\textwidth]{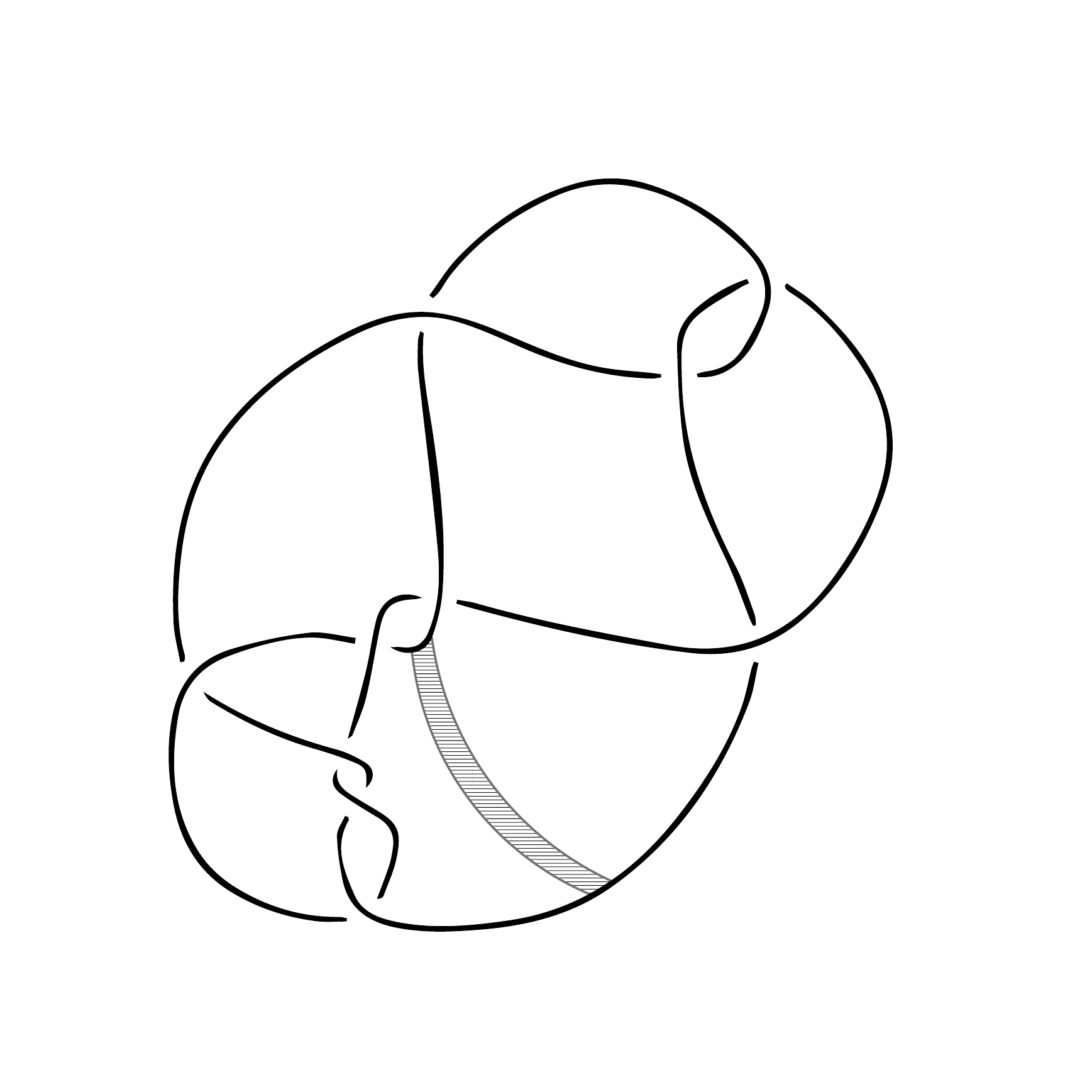}
		\caption{$10_{29}\stackrel{0}{\longrightarrow}5_2 \# -5_2$}
		\label{FigureFor10-29}
	\end{subfigure}
	~
	\begin{subfigure}[b]{0.3\textwidth}
		\includegraphics[width=\textwidth]{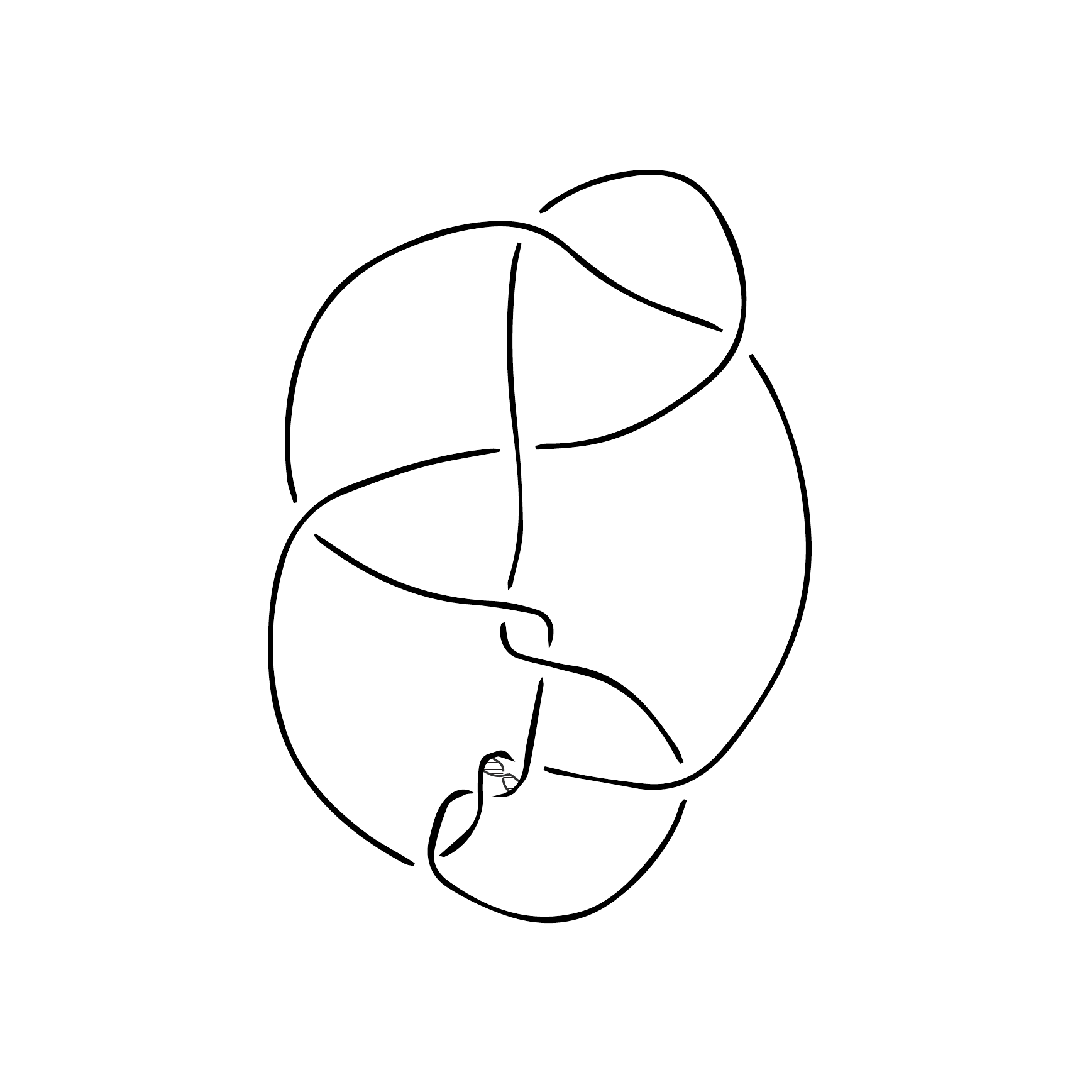}
		\caption{$10_{30}\stackrel{-1}{\longrightarrow} 9_{27}$}
		\label{FigureFor10-30}
	\end{subfigure}
	~
	\begin{subfigure}[b]{0.27\textwidth}
		\includegraphics[width=\textwidth]{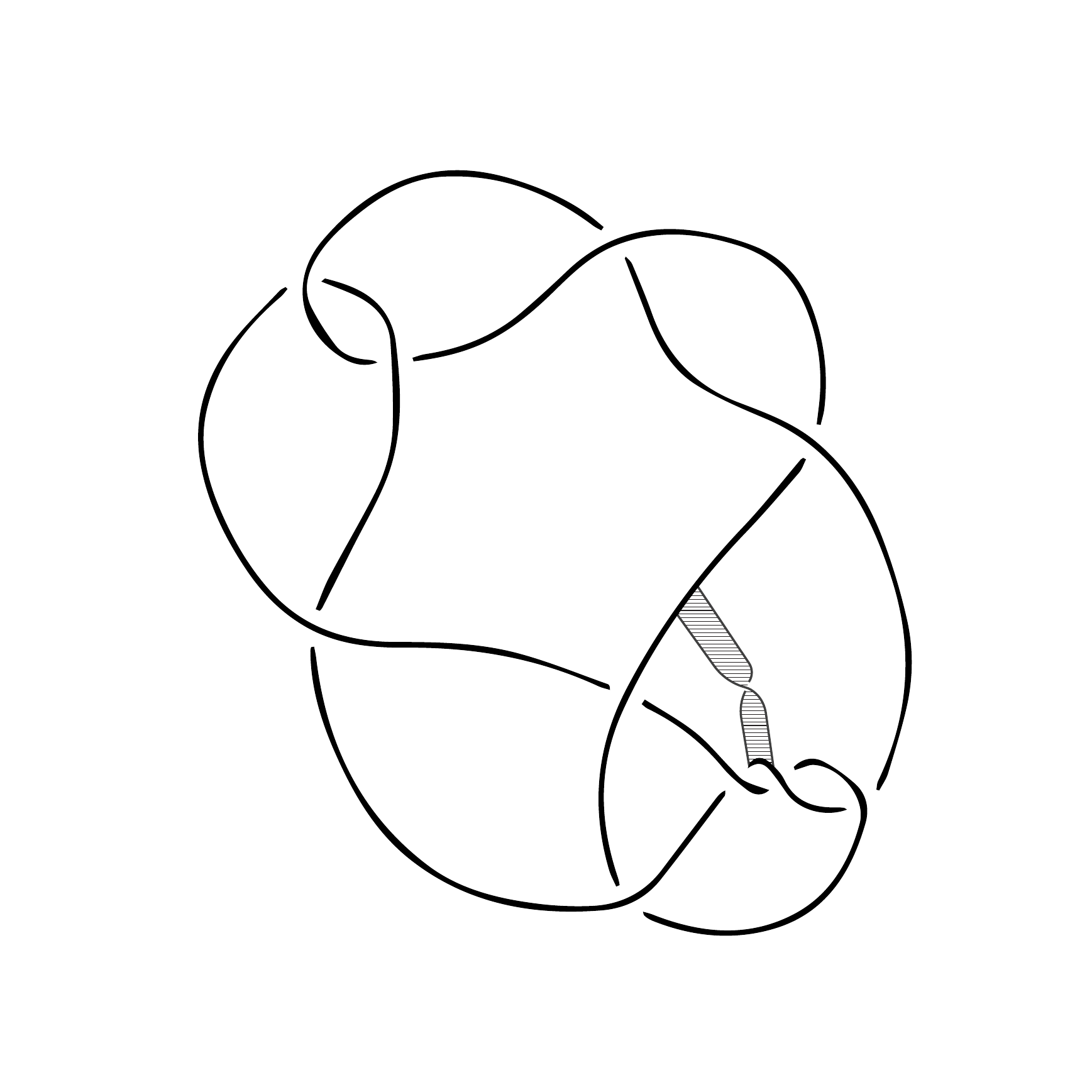}
		\caption{$10_{31}\stackrel{-1}{\longrightarrow} 0_{1}$}
		\label{FigureFor10-31}
	\end{subfigure}
	~
	\vskip3mm	

	\caption{Non-oriented band moves from the knots  $10_{17}$, $10_{20}$, $10_{21}$, $10_{23}$, $10_{24}$, $10_{27}$, $10_{29}$, $10_{30}$, $10_{31}$ to slice knots}
	\label{slice2}
\end{figure}
\begin{figure}[h]
	\centering
			\begin{subfigure}[b]{0.27\textwidth}
		\includegraphics[width=\textwidth]{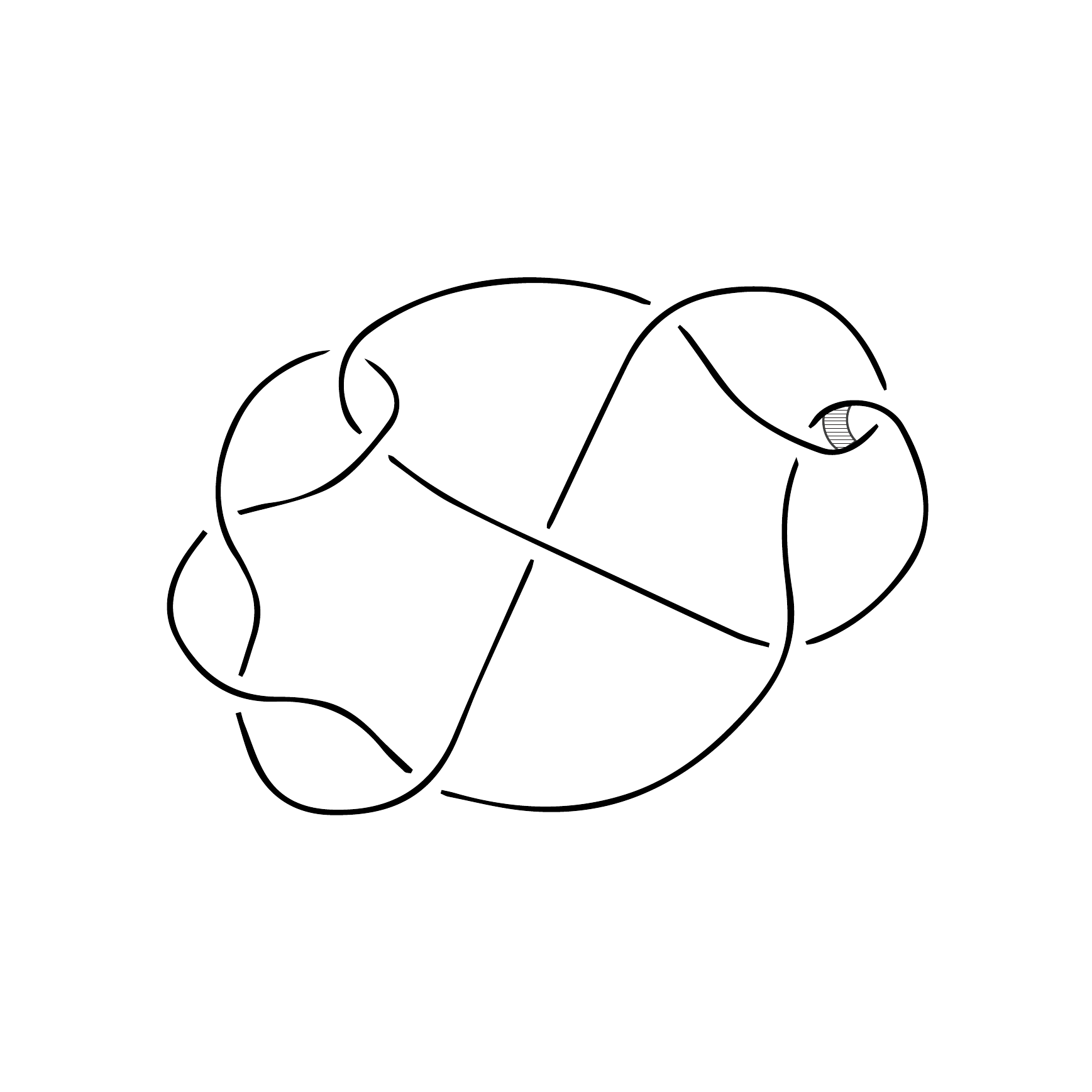}
		\caption{$10_{38}\stackrel{0}{\longrightarrow} 8_{8}$}
		\label{FigureFor10-38}
	\end{subfigure}
	~
	\begin{subfigure}[b]{0.3\textwidth}
		\includegraphics[width=\textwidth]{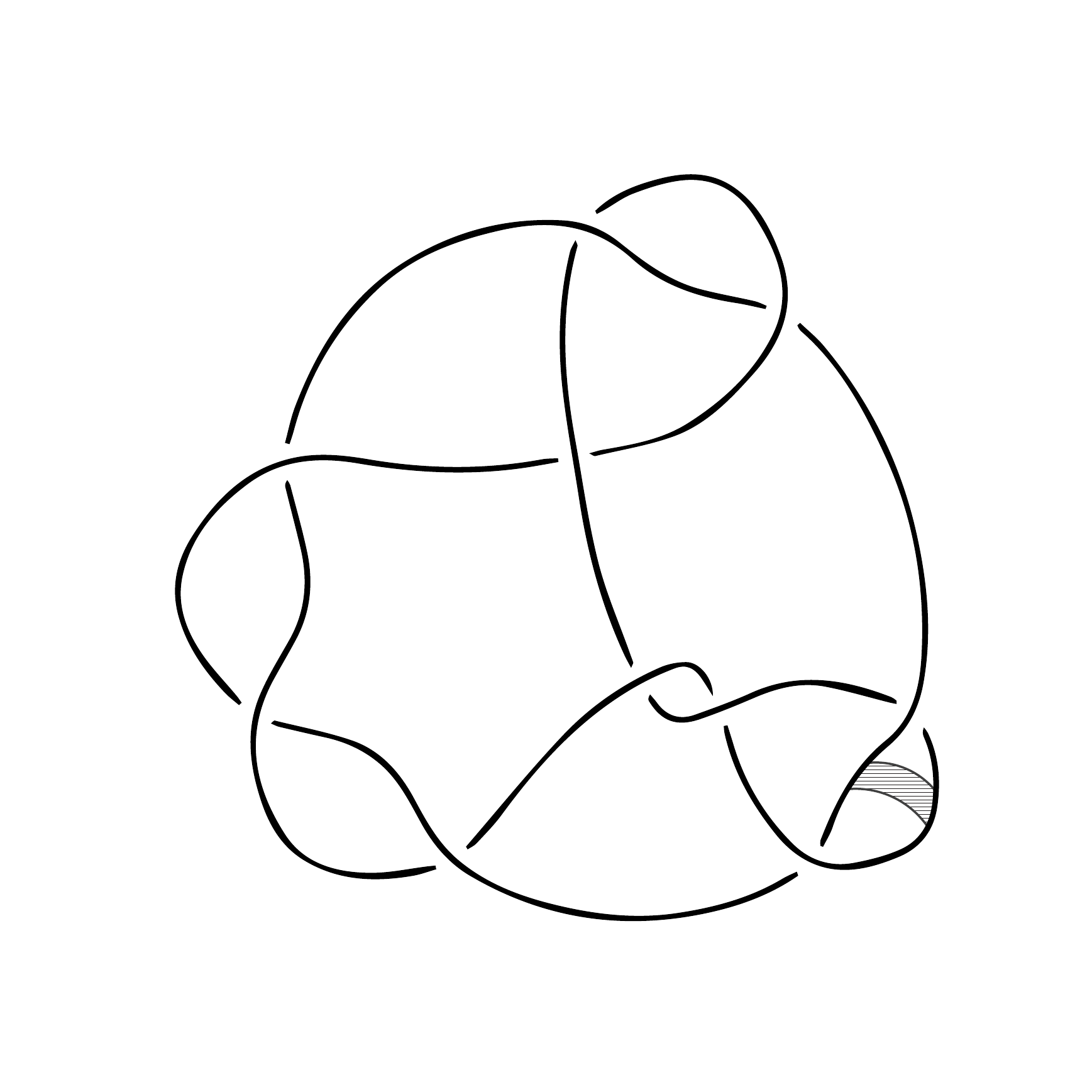}
		\caption{$10_{39}\stackrel{0}{\longrightarrow} 8_{8}$}
		\label{FigureFor10-39}
	\end{subfigure}
	~
	\begin{subfigure}[b]{0.27\textwidth}
		\includegraphics[width=\textwidth]{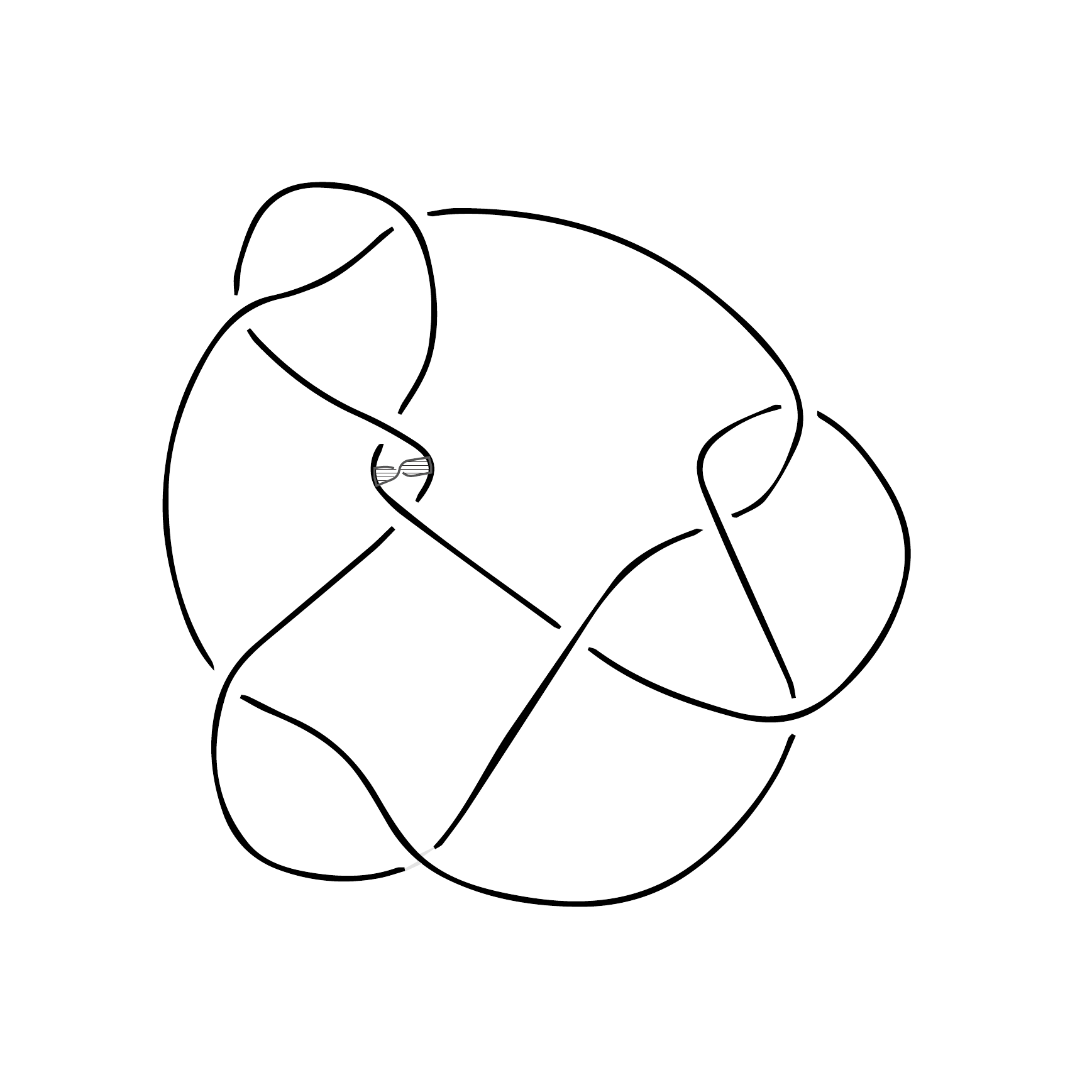}
		\caption{$10_{40}\stackrel{-1}{\longrightarrow} 9_{27}$}
		\label{FigureFor10-40}
	\end{subfigure}
	\vskip3mm
	\begin{subfigure}[b]{0.27\textwidth}
		\includegraphics[width=\textwidth]{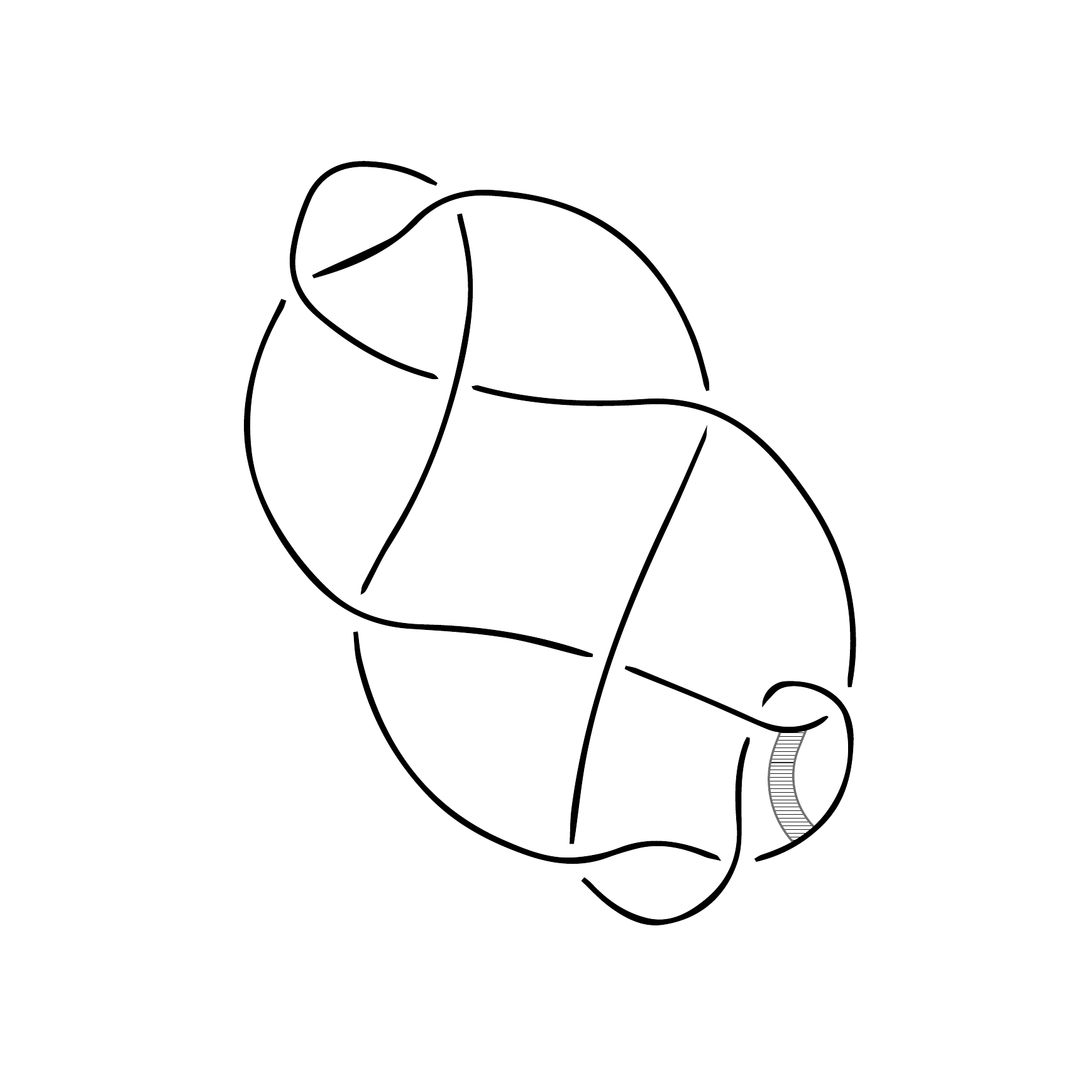}
		\caption{$10_{41}\stackrel{-1}{\longrightarrow} 9_{27}$}
		\label{FigureFor10-41}
	\end{subfigure}
	~
	\begin{subfigure}[b]{0.27\textwidth}
		\includegraphics[width=\textwidth]{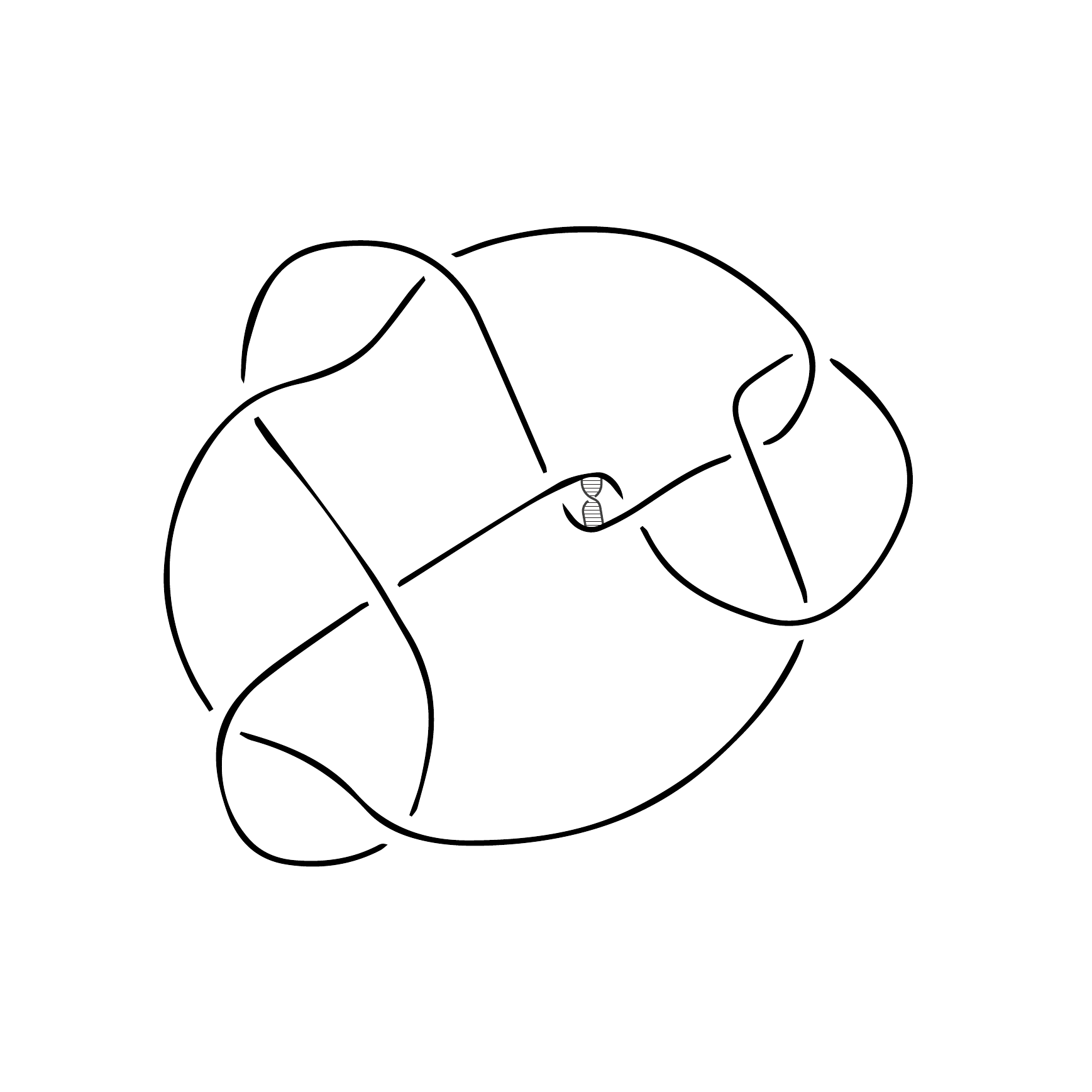}
		\caption{$10_{43}\stackrel{-1}{\longrightarrow} 9_{27}$}
		\label{FigureFor10-43}
	\end{subfigure}
	~
	\begin{subfigure}[b]{0.27\textwidth}
		\includegraphics[width=\textwidth]{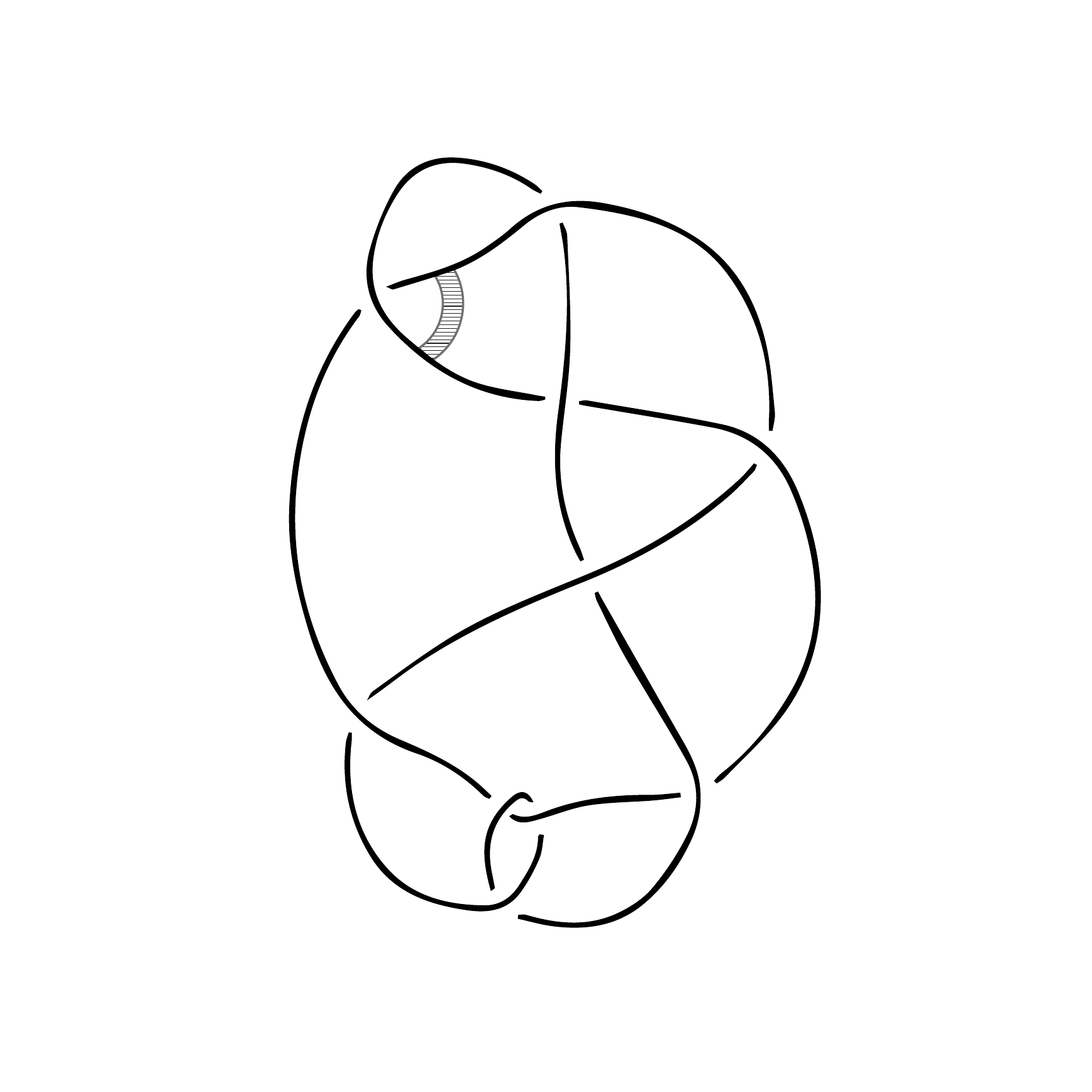}
		\caption{$10_{44}\stackrel{0\phantom{i}}{\longrightarrow} 9_{27}$}
		\label{FigureFor10-44}
	\end{subfigure}
	~     
	\vskip3mm
	~
	\begin{subfigure}[b]{0.3\textwidth}
		\includegraphics[width=\textwidth]{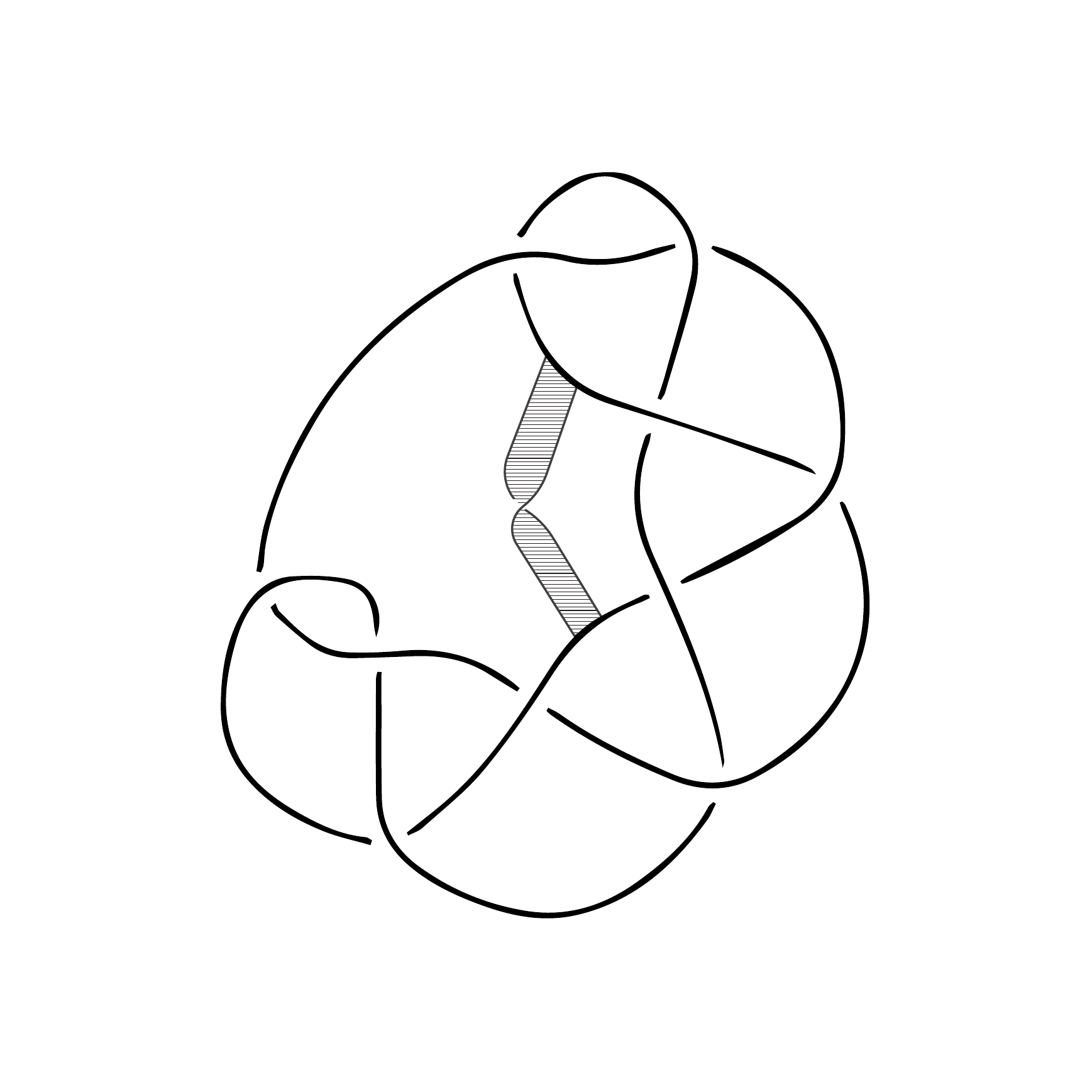}
		\caption{$10_{45}\stackrel{1}{\longrightarrow}8_8$}
		\label{FigureFor10-45}
	\end{subfigure}
	~
	\begin{subfigure}[b]{0.3\textwidth}
		\includegraphics[width=\textwidth]{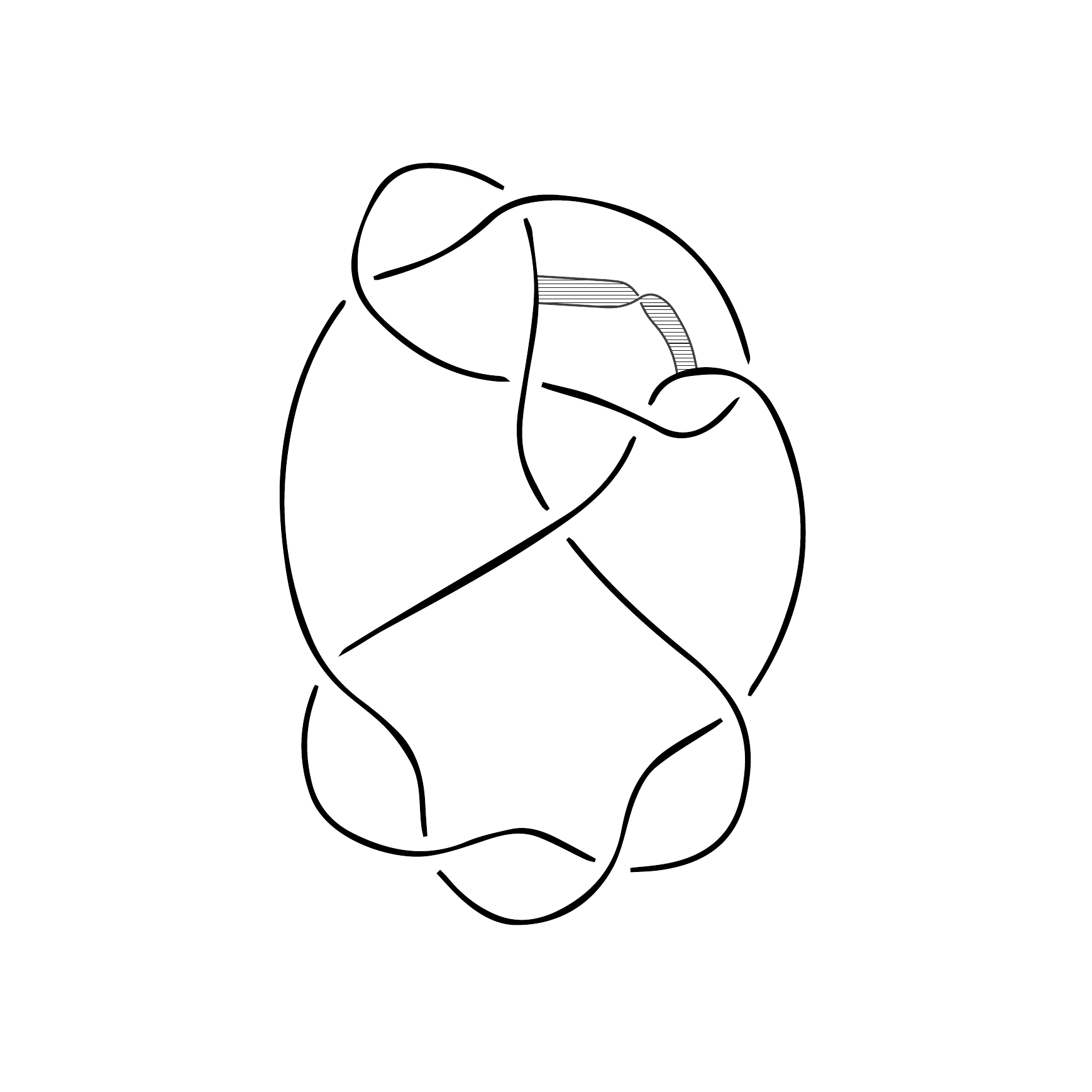}
		\caption{$10_{49}\stackrel{-1}{\longrightarrow} 6_{1}$}
		\label{FigureFor10-49}
	\end{subfigure}
	~
	\begin{subfigure}[b]{0.27\textwidth}
		\includegraphics[width=\textwidth]{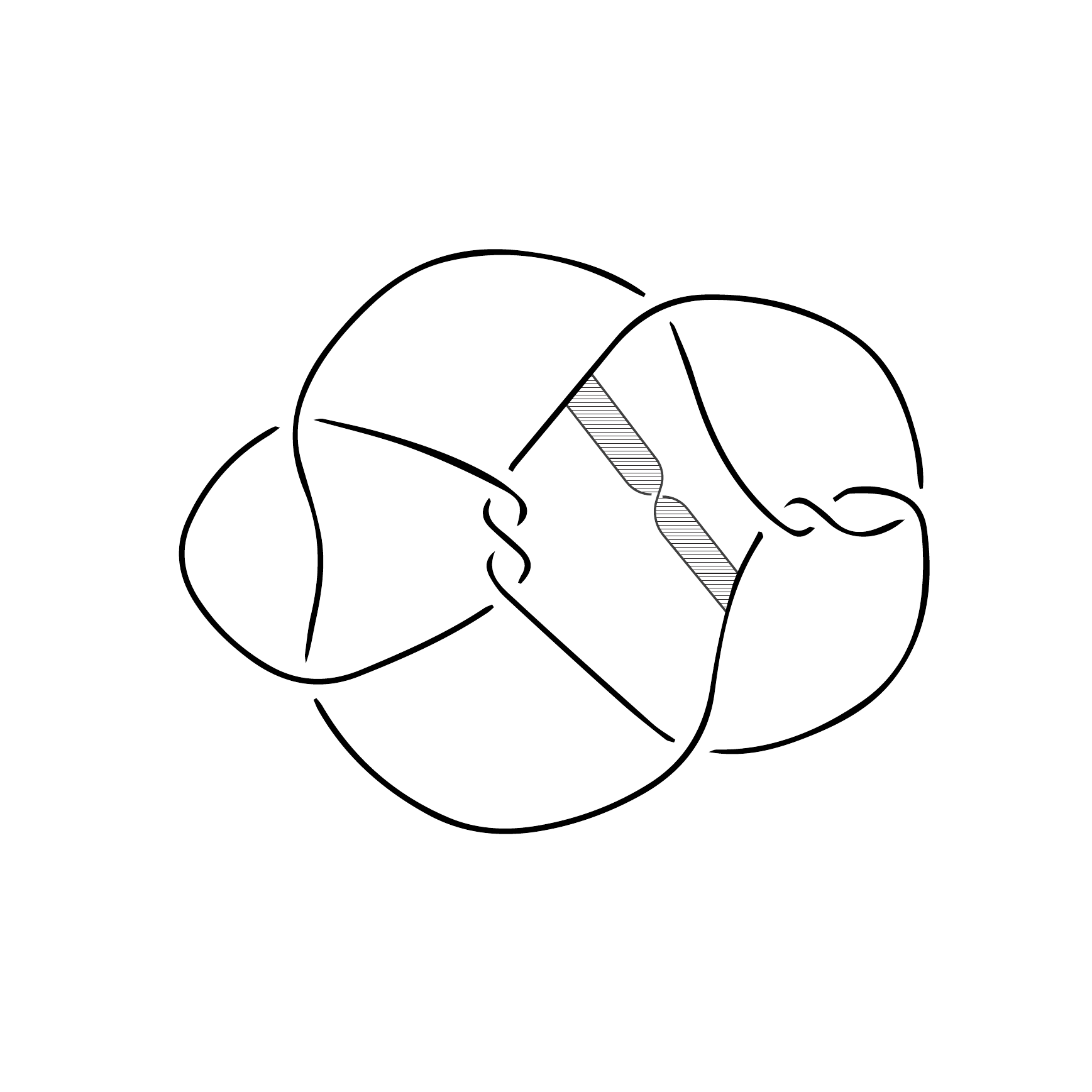}
		\caption{$10_{50}\stackrel{1}{\longrightarrow} 8_{20}$}
		\label{FigureFor10-50}
	\end{subfigure}
	~
	\vskip3mm	

	\caption{Non-oriented band moves from the knots  $10_{38}$, $10_{39}$, $10_{40}$, $10_{41}$, $10_{43}$, $10_{44}$, $10_{45}$, $10_{49}$, $10_{50}$ to slice knots}\label{slice3}
\end{figure}
\begin{figure}[h]
	\centering
			\begin{subfigure}[b]{0.27\textwidth}
		\includegraphics[width=\textwidth]{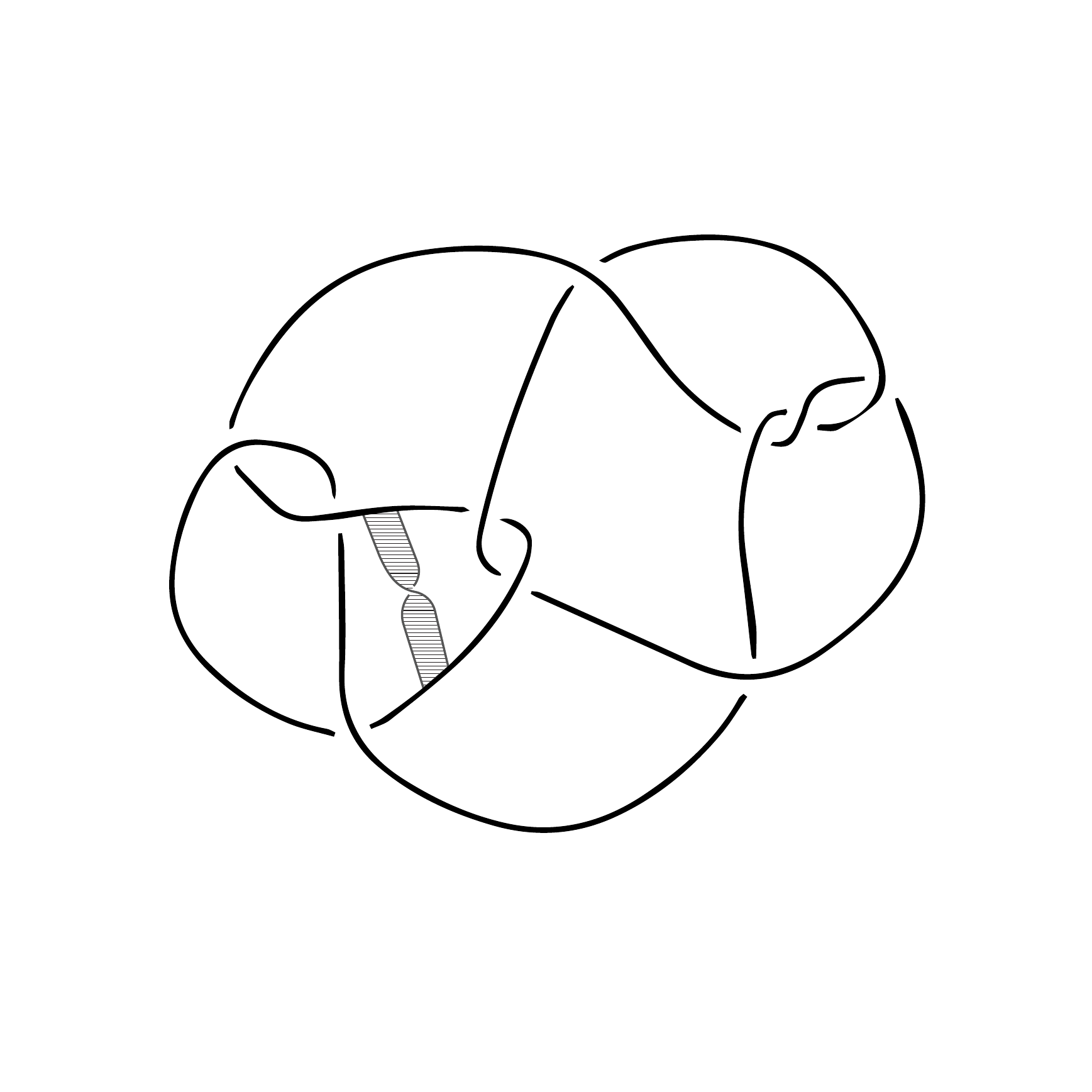}
		\caption{$10_{51}\stackrel{-1}{\longrightarrow} 10_{129}$}
		\label{FigureFor10-51}
	\end{subfigure}
	~
	\begin{subfigure}[b]{0.3\textwidth}
		\includegraphics[width=\textwidth]{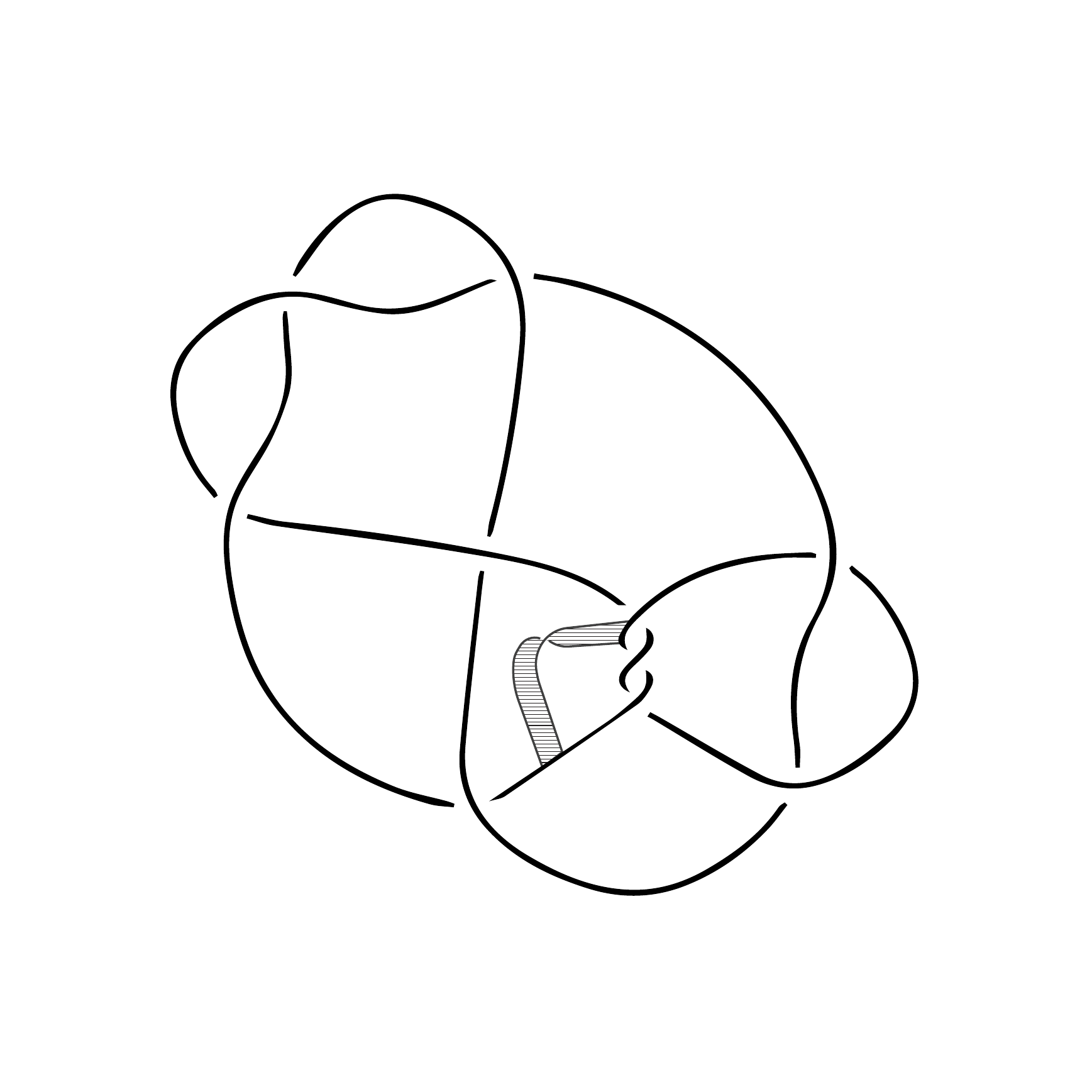}
		\caption{$10_{52}\stackrel{-1}{\longrightarrow} 0_{1}$}
		\label{FigureFor10-52}
	\end{subfigure}
	~
	\begin{subfigure}[b]{0.27\textwidth}
		\includegraphics[width=\textwidth]{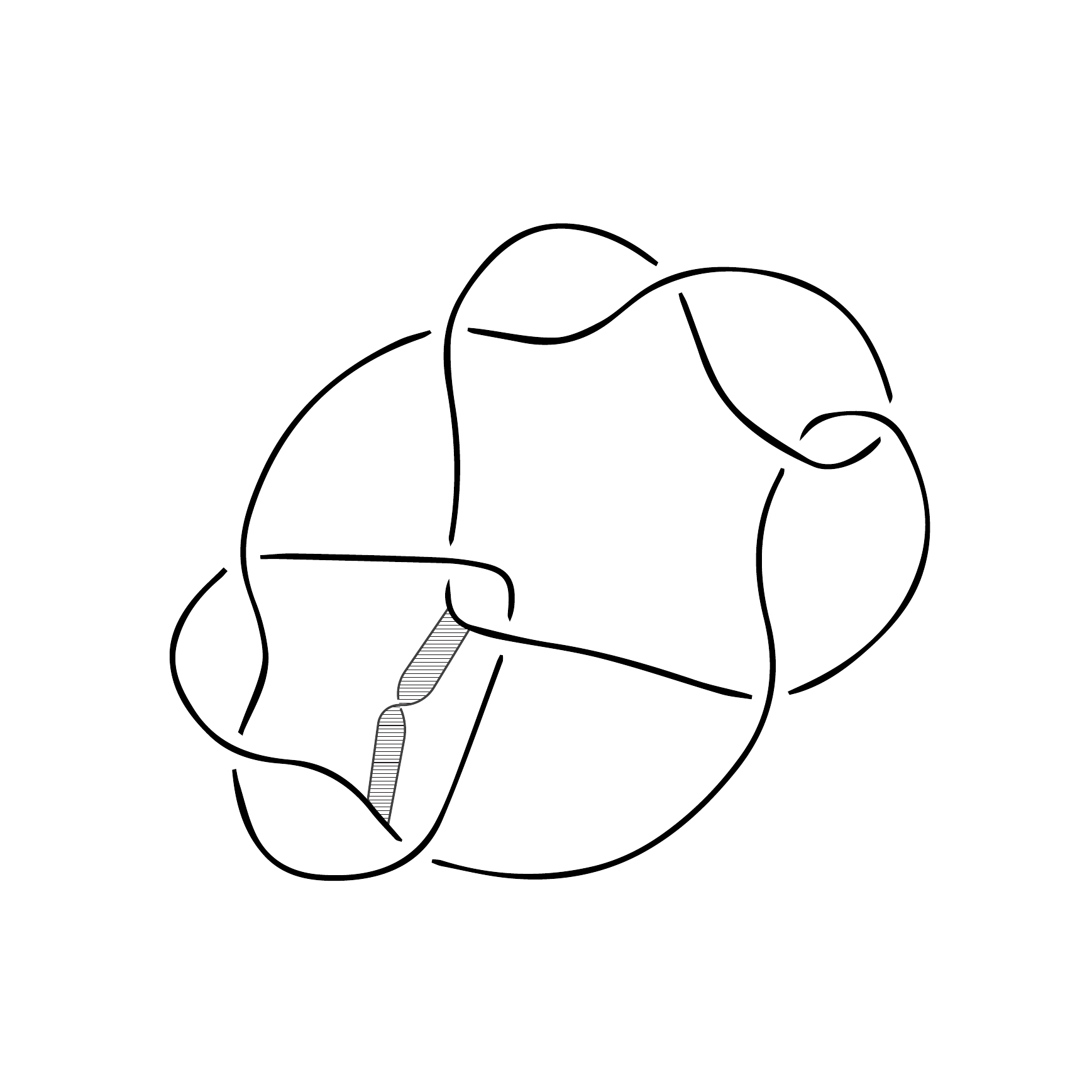}
		\caption{$10_{54}\stackrel{1}{\longrightarrow} 6_{1}$}
		\label{FigureFor10-54}
	\end{subfigure}
	\vskip3mm
	\begin{subfigure}[b]{0.27\textwidth}
		\includegraphics[width=\textwidth]{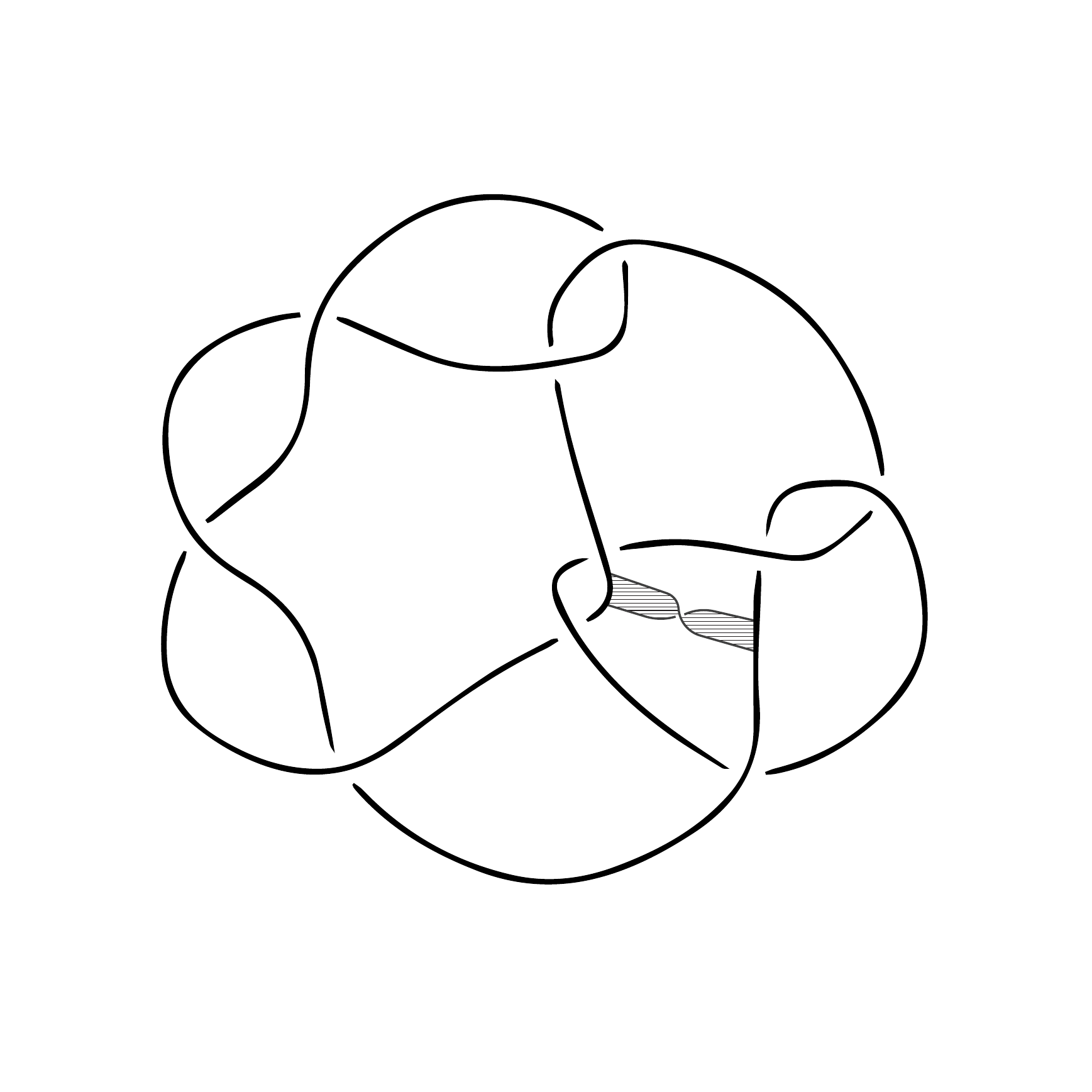}
		\caption{$10_{55}\stackrel{1}{\longrightarrow} 6_{1}$}
		\label{FigureFor10-55}
	\end{subfigure}
	~
	\begin{subfigure}[b]{0.27\textwidth}
		\includegraphics[width=\textwidth]{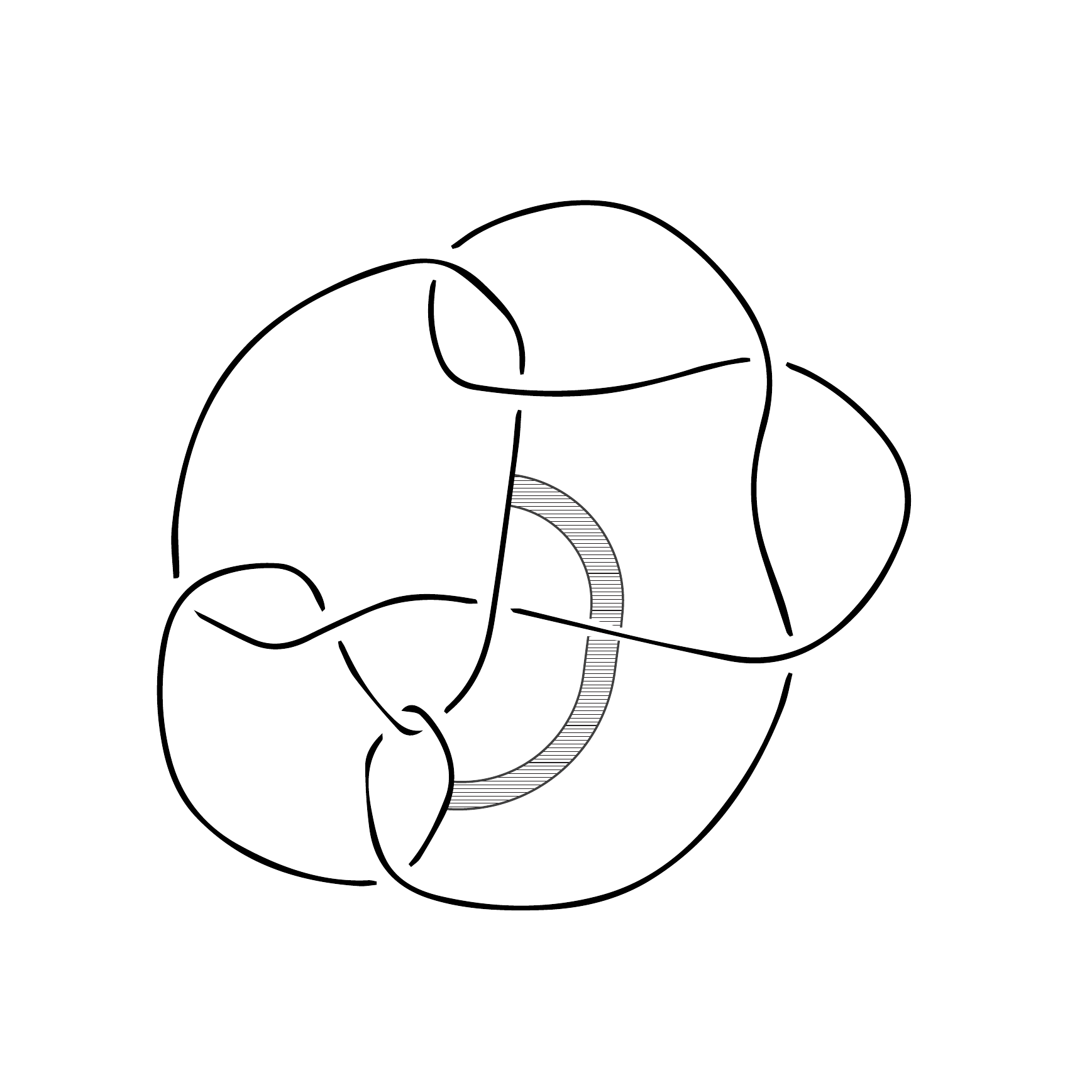}
		\caption{$10_{57}\stackrel{0}{\longrightarrow} 6_{1}$}
		\label{FigureFor10-57}
	\end{subfigure}
	~
	\begin{subfigure}[b]{0.27\textwidth}
		\includegraphics[width=\textwidth]{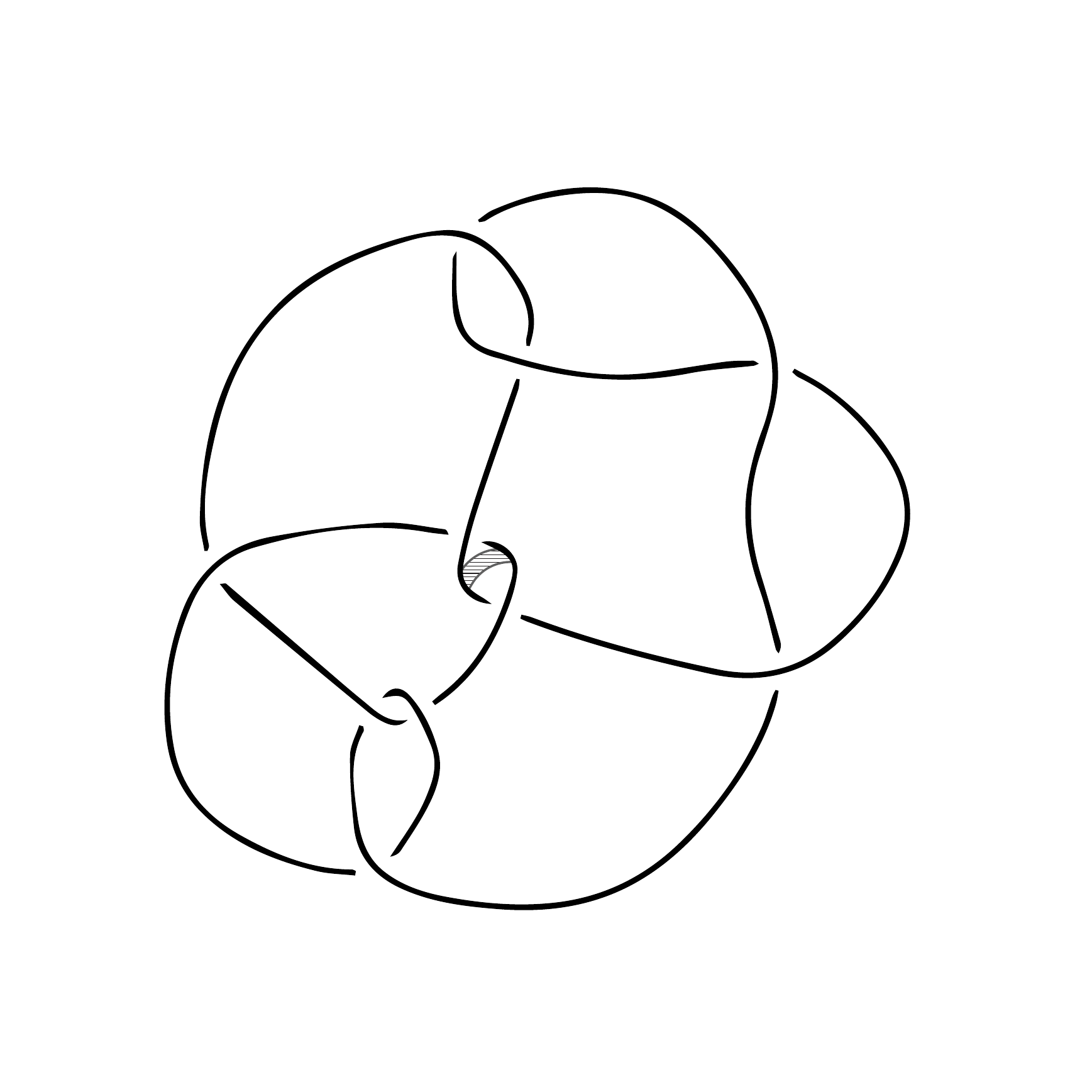}
		\caption{$10_{59}\stackrel{0\phantom{i}}{\longrightarrow} 8_{8}$}
		\label{FigureFor10-59}
	\end{subfigure}
	~      
	\vskip3mm
	~
	\begin{subfigure}[b]{0.27\textwidth}
		\includegraphics[width=\textwidth]{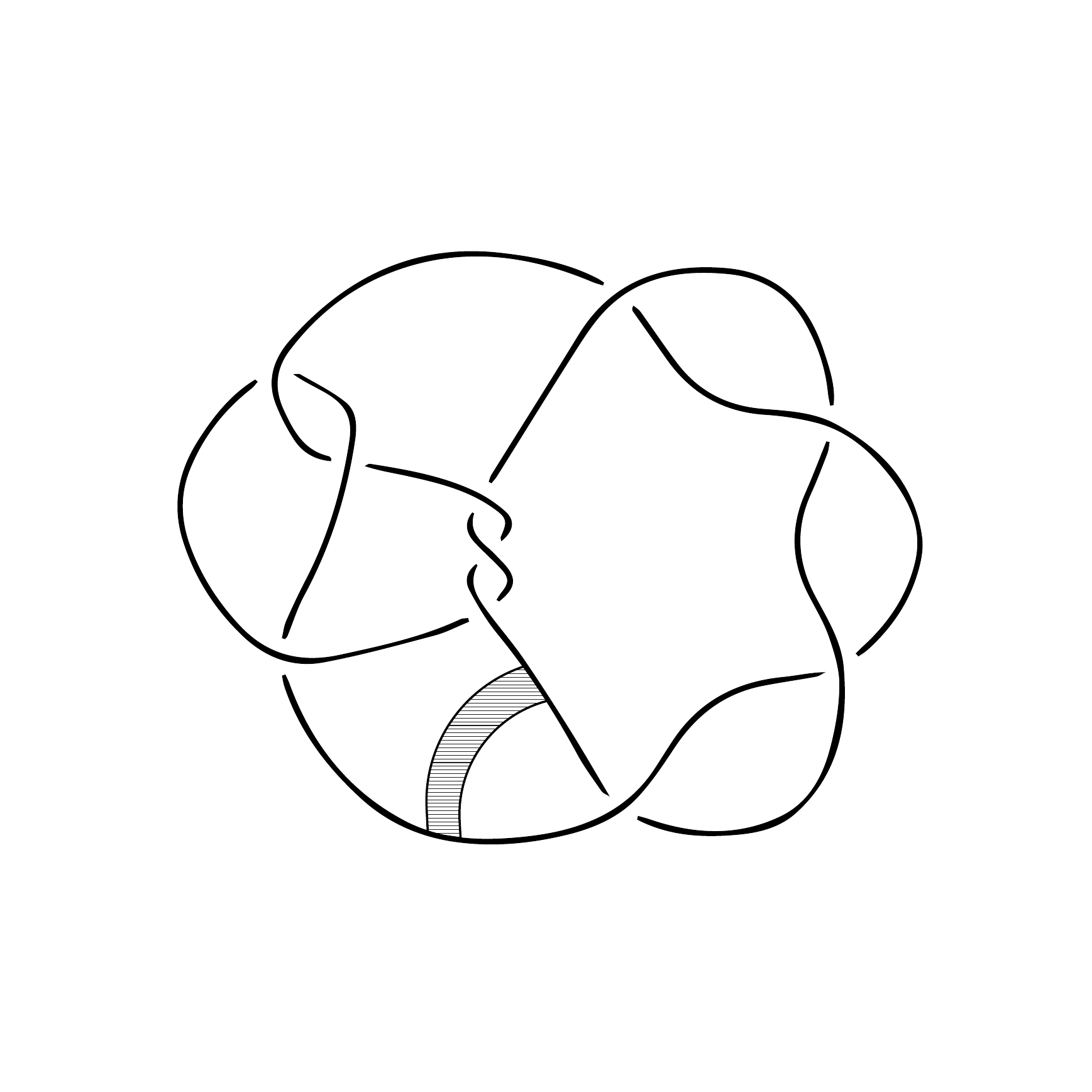}
		\caption{$10_{62}\stackrel{0\phantom{i}}{\longrightarrow} 6_{1}$}
		\label{FigureFor10-62}
	\end{subfigure}
~
	\begin{subfigure}[b]{0.3\textwidth}
		\includegraphics[width=\textwidth]{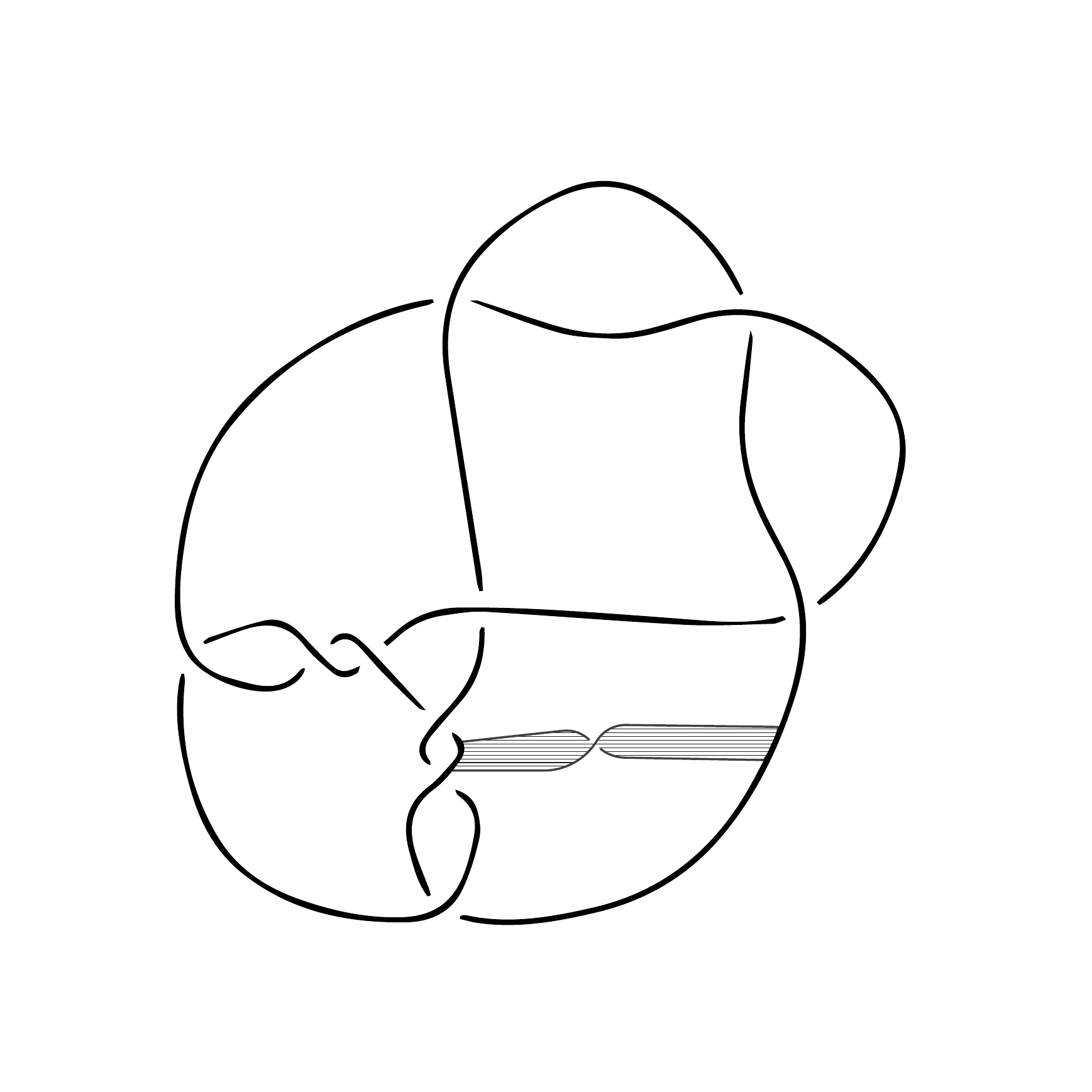}
		\caption{$10_{64}\stackrel{-1}{\longrightarrow}0_1$}
		\label{FigureFor10-64}
	\end{subfigure}
	~
	\begin{subfigure}[b]{0.3\textwidth}
		\includegraphics[width=\textwidth]{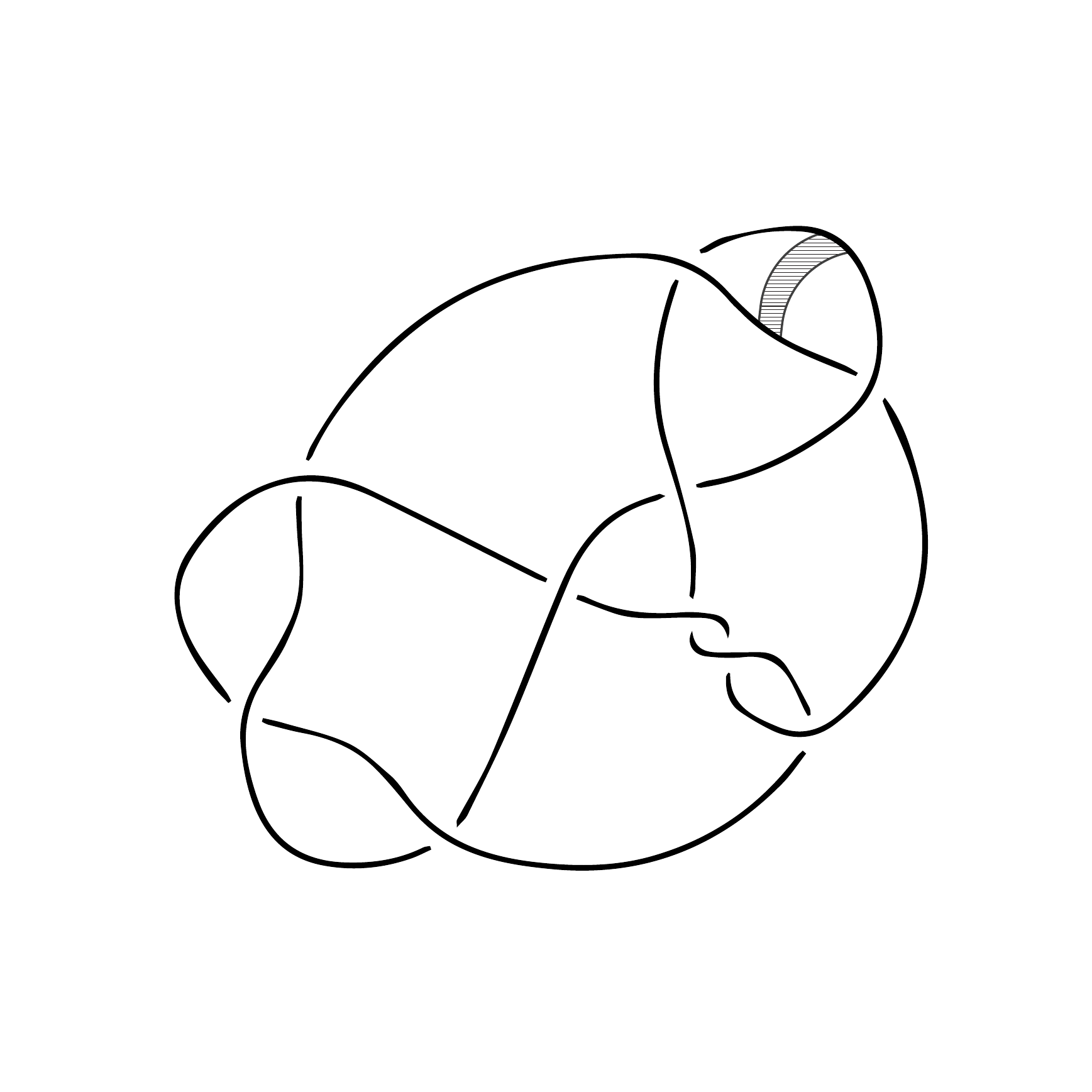}
		\caption{$10_{65}\stackrel{0}{\longrightarrow} 8_9$}
		\label{FigureFor10-65}
	\end{subfigure}
	~
	\vskip3mm

	\caption{Non-oriented band moves from the knots $10_{51}$, $10_{52}$, $10_{54}$, $10_{55}$, $10_{57}$, $10_{59}$, $10_{62}$, $10_{64}$, $10_{65}$ to slice knots}\label{slice4}
\end{figure}
\begin{figure}[h]
	\centering
	~
			\begin{subfigure}[b]{0.27\textwidth}
		\includegraphics[width=\textwidth]{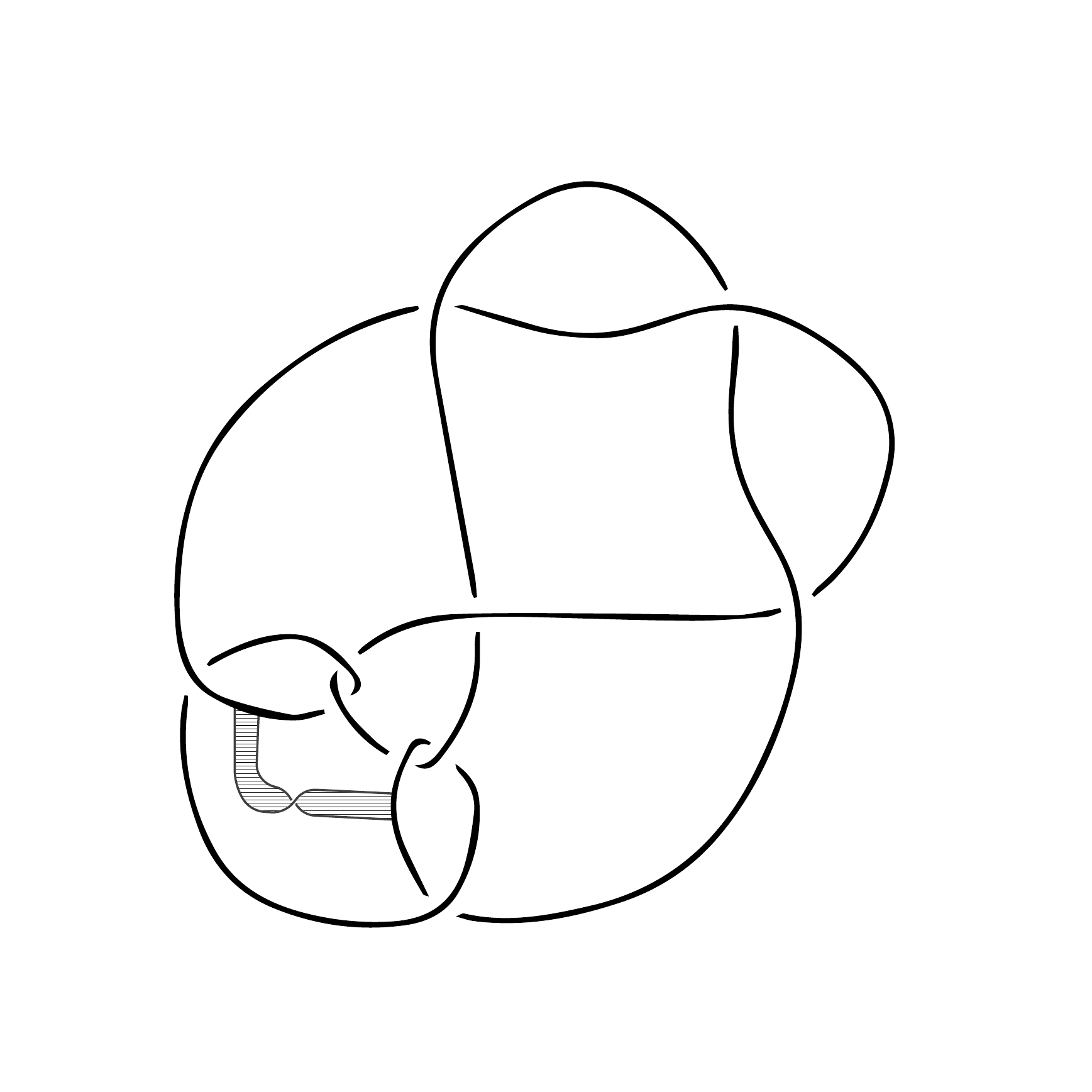}
		\caption{$10_{66}\stackrel{-1}{\longrightarrow} 8_9$}
		\label{FigureFor10-66}
	\end{subfigure}
		~
		\begin{subfigure}[b]{0.27\textwidth}
		\includegraphics[width=\textwidth]{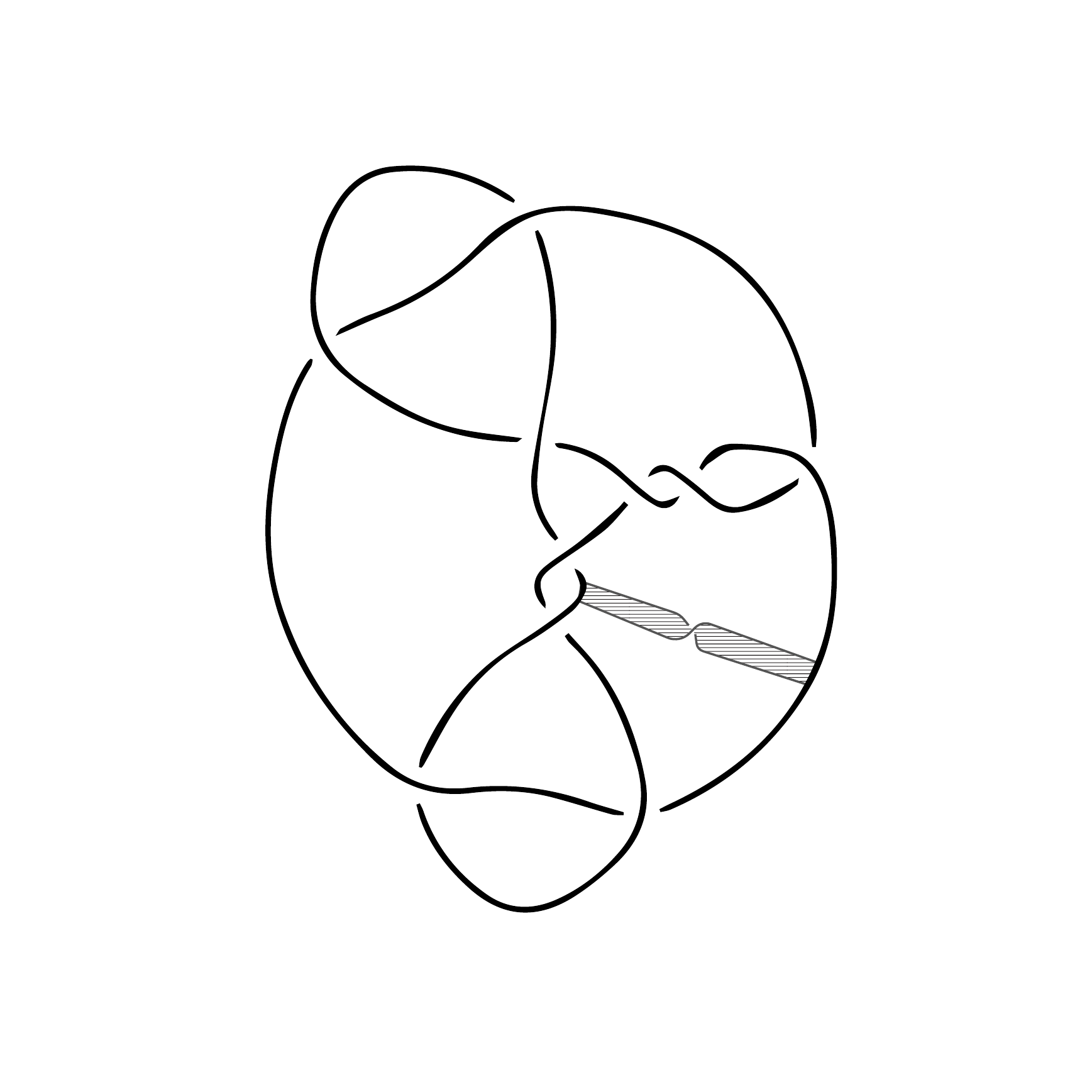}
		\caption{$10_{67}\stackrel{-1}{\longrightarrow} 8_{20}$}
		\label{FigureFor10-67}
	\end{subfigure}
	~
		\begin{subfigure}[b]{0.3\textwidth}
		\includegraphics[width=\textwidth]{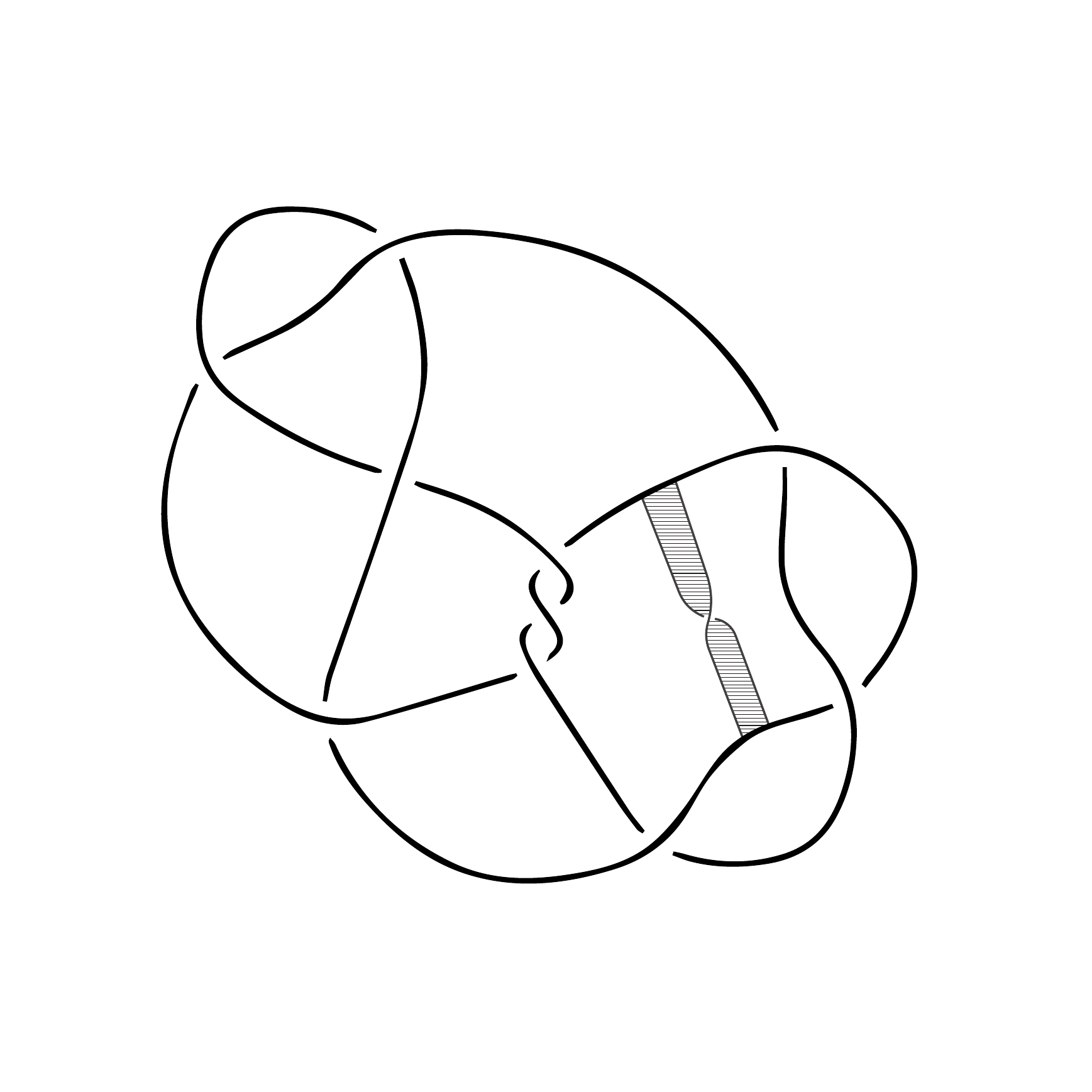}
		\caption{$10_{68}\stackrel{1}{\longrightarrow} 0_{1}$}
		\label{FigureFor10-68}
	\end{subfigure}
	\vskip3mm
	\begin{subfigure}[b]{0.27\textwidth}
		\includegraphics[width=\textwidth]{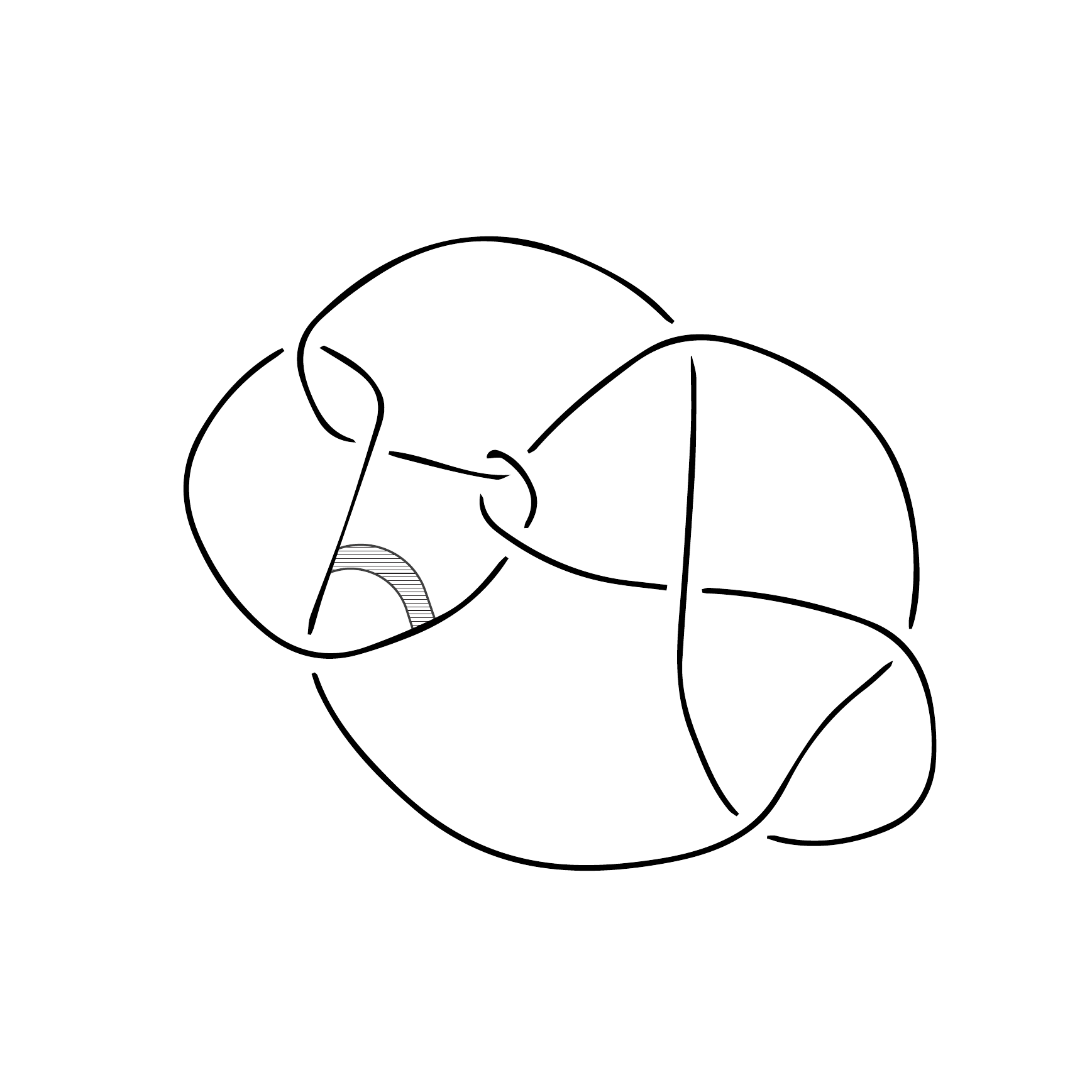}
		\caption{$10_{69}\stackrel{0}{\longrightarrow} 9_{27}$}
		\label{FigureFor10-69}
	\end{subfigure}
	~
	\begin{subfigure}[b]{0.27\textwidth}
		\includegraphics[width=\textwidth]{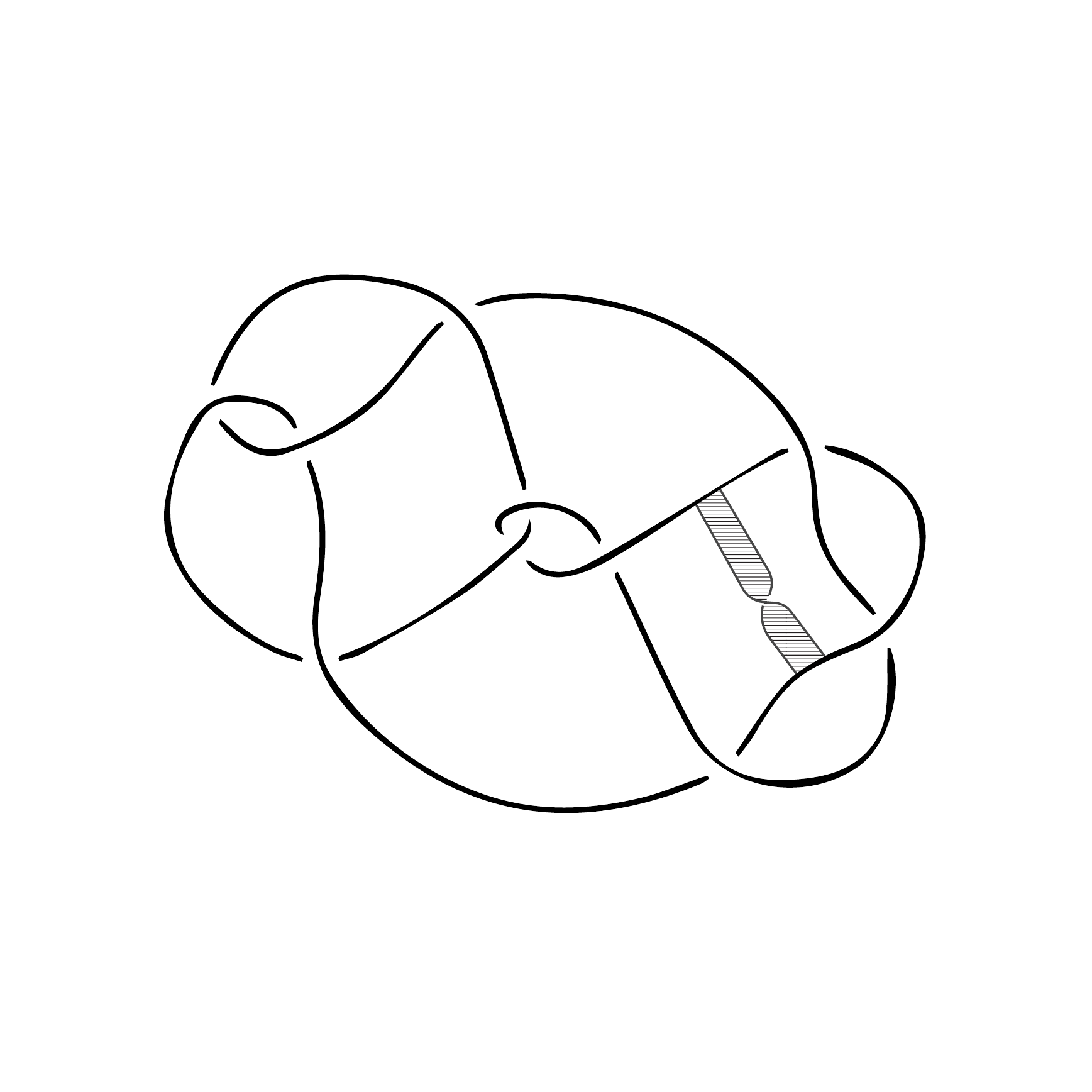}
		\caption{$10_{70}\stackrel{-1}{\longrightarrow} 6_{1}$}
		\label{FigureFor10-70}
	\end{subfigure}
	~
	\begin{subfigure}[b]{0.27\textwidth}
		\includegraphics[width=\textwidth]{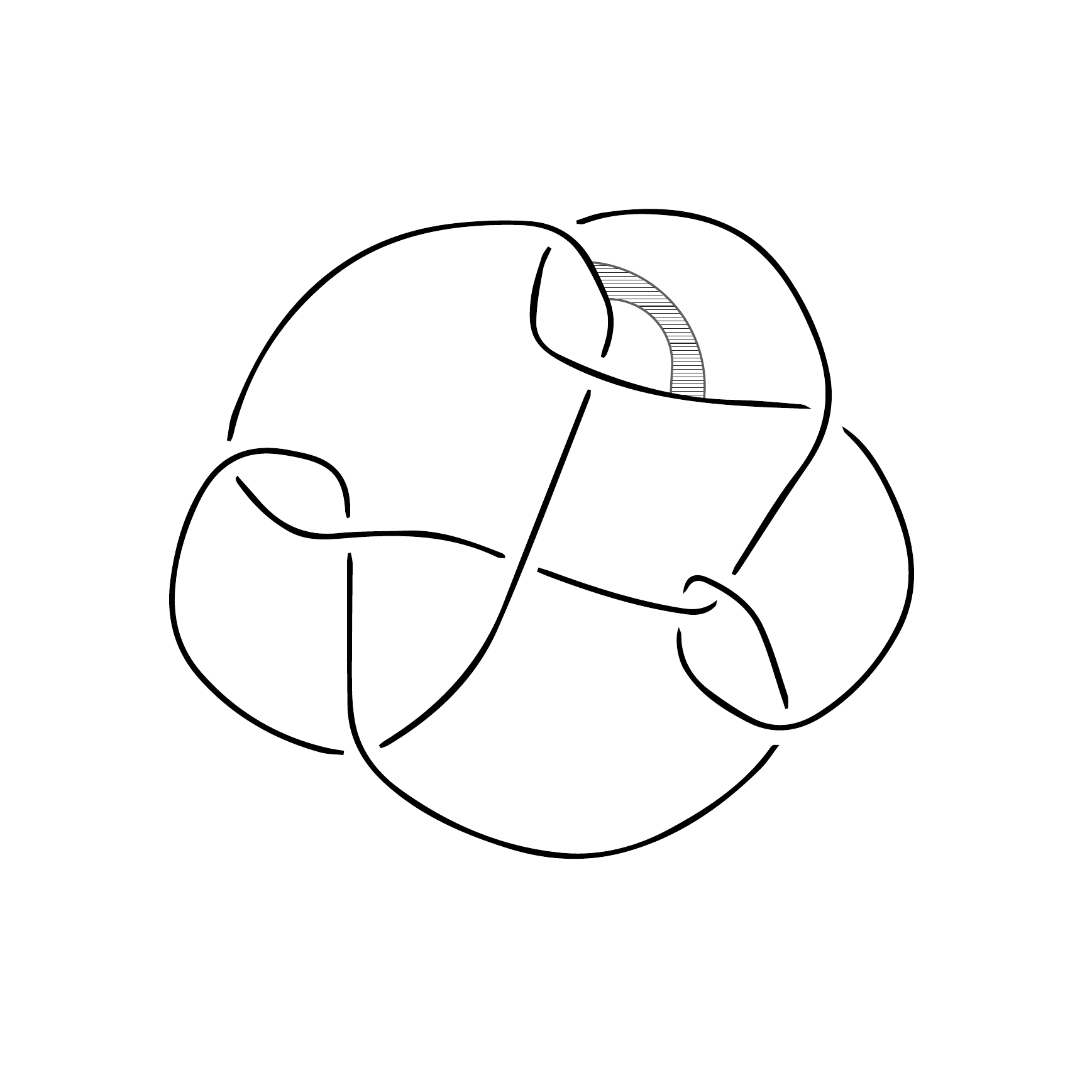}
		\caption{$10_{73}\stackrel{0}{\longrightarrow} 9_{27}$}
		\label{FigureFor10-73}
	\end{subfigure}
	~     
	\vskip3mm
	~
	\begin{subfigure}[b]{0.27\textwidth}
		\includegraphics[width=\textwidth]{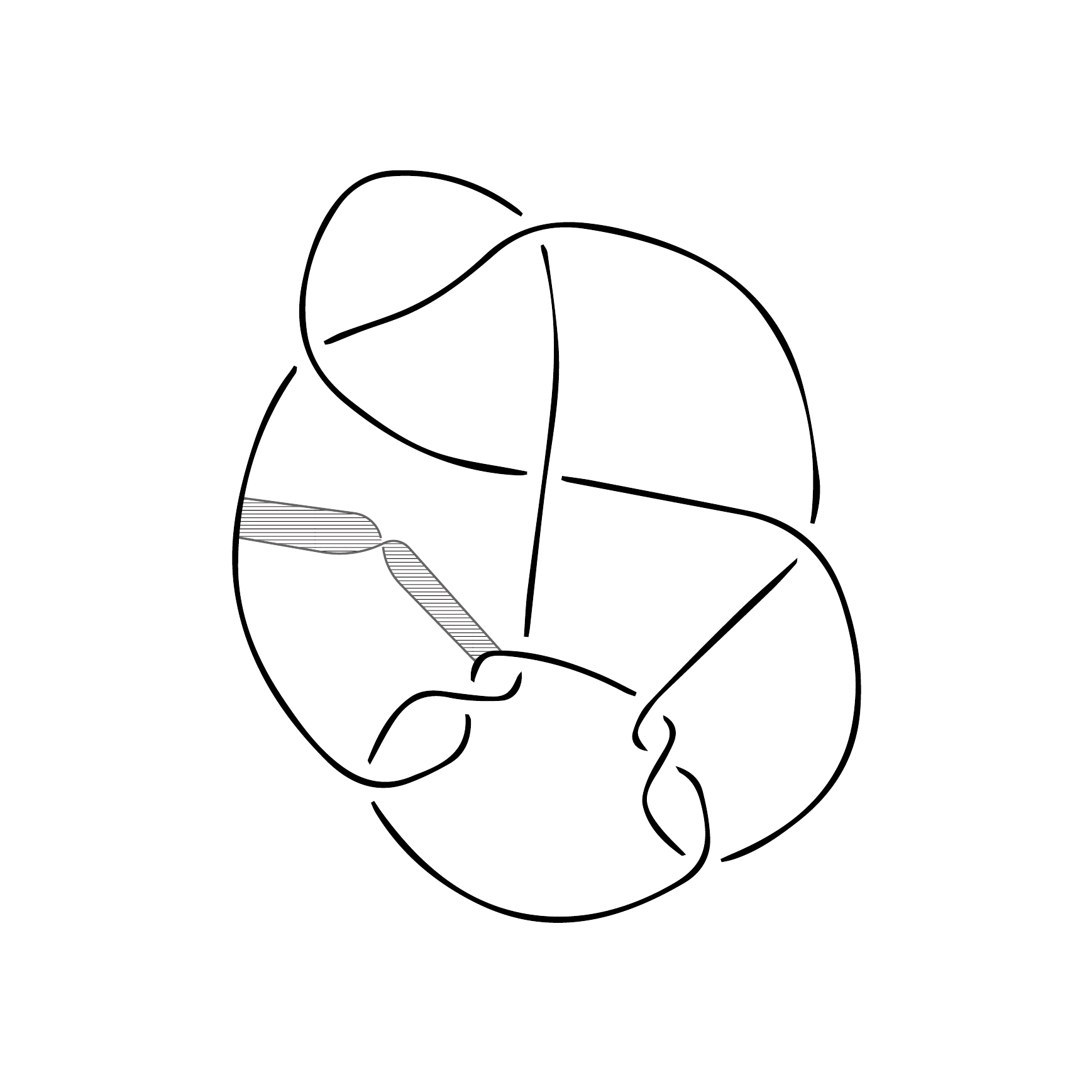}
		\caption{$10_{74}\stackrel{-1\phantom{i}}{\longrightarrow} 6_{1}$}
		\label{FigureFor10-74}
	\end{subfigure}
	~
		\begin{subfigure}[b]{0.3\textwidth}
		\includegraphics[width=\textwidth]{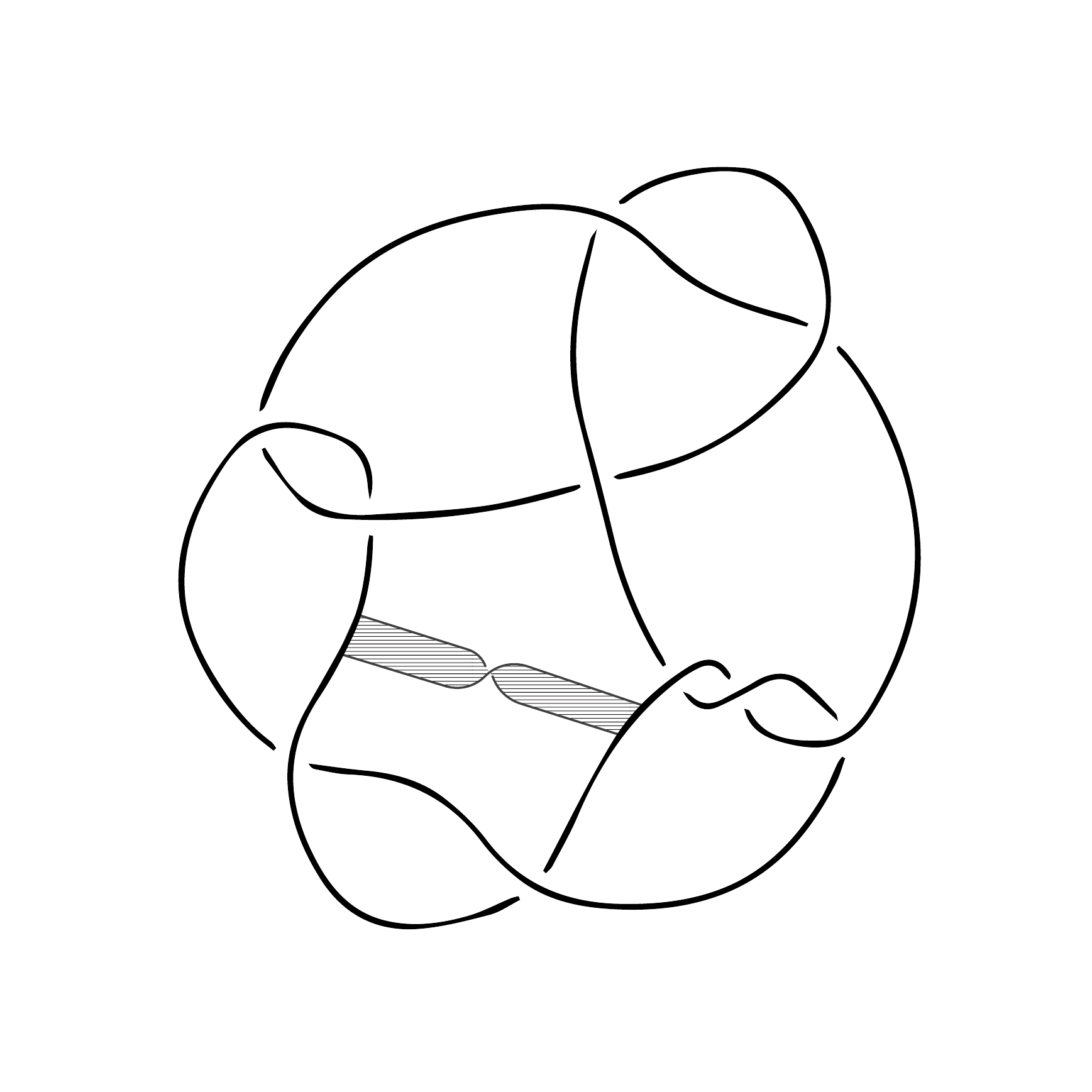}
		\caption{$10_{77}\stackrel{-1}{\longrightarrow}8_{20}$}
		\label{FigureFor10-77}
	\end{subfigure}
	~
	\begin{subfigure}[b]{0.3\textwidth}
		\includegraphics[width=\textwidth]{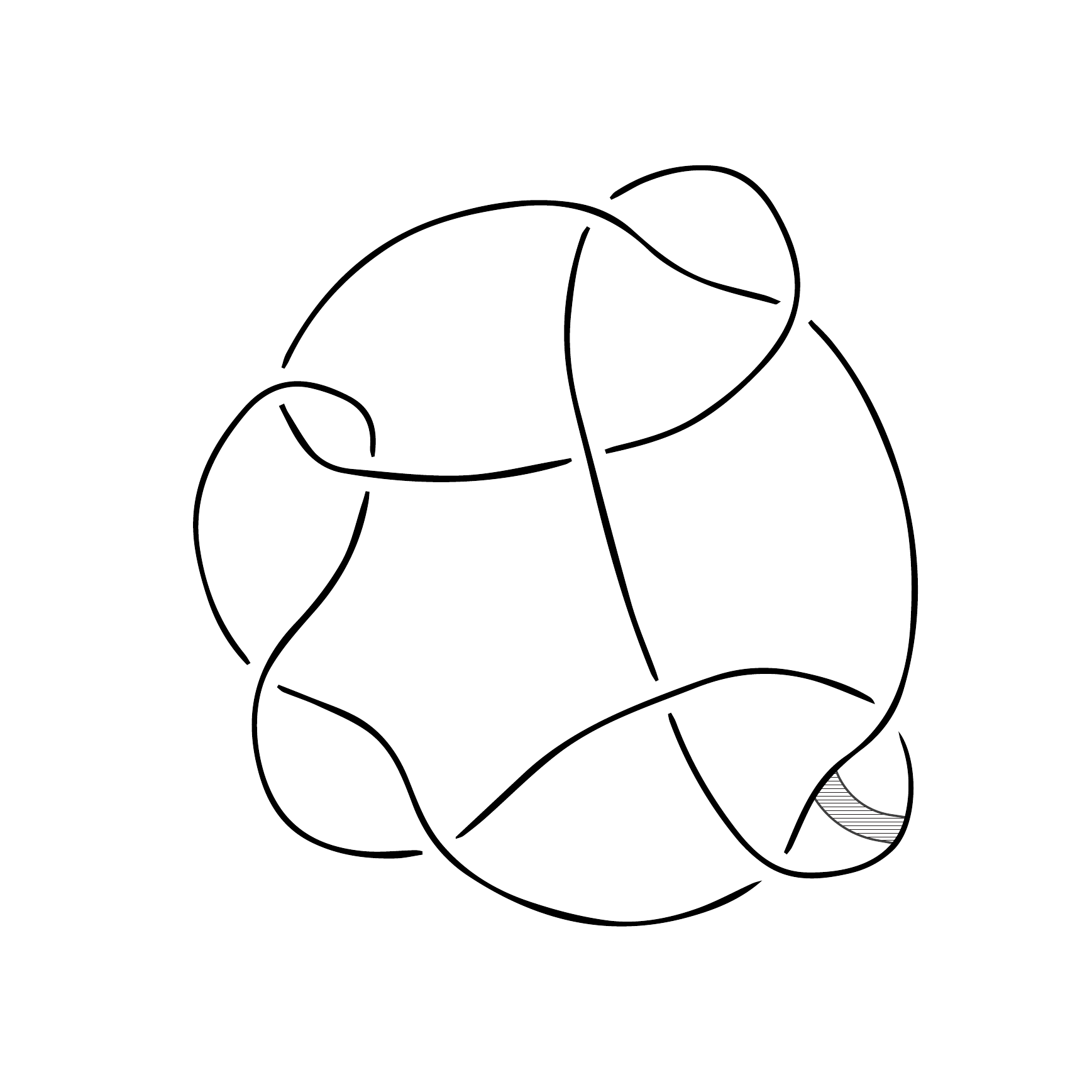}
		\caption{$10_{78}\stackrel{0}{\longrightarrow} 8_{8}$}
		\label{FigureFor10-78}
	\end{subfigure}
	~
	\vskip3mm

	\caption{Non-oriented band moves from the knots $10_{66}$, $10_{67}$, $10_{68}$, $10_{69}$, $10_{70}$, $10_{73}$, $10_{74}$, $10_{77}$, $10_{78}$ to slice knots}\label{slice5}
\end{figure}
\begin{figure}[h]
	\centering
			\begin{subfigure}[b]{0.27\textwidth}
		\includegraphics[width=\textwidth]{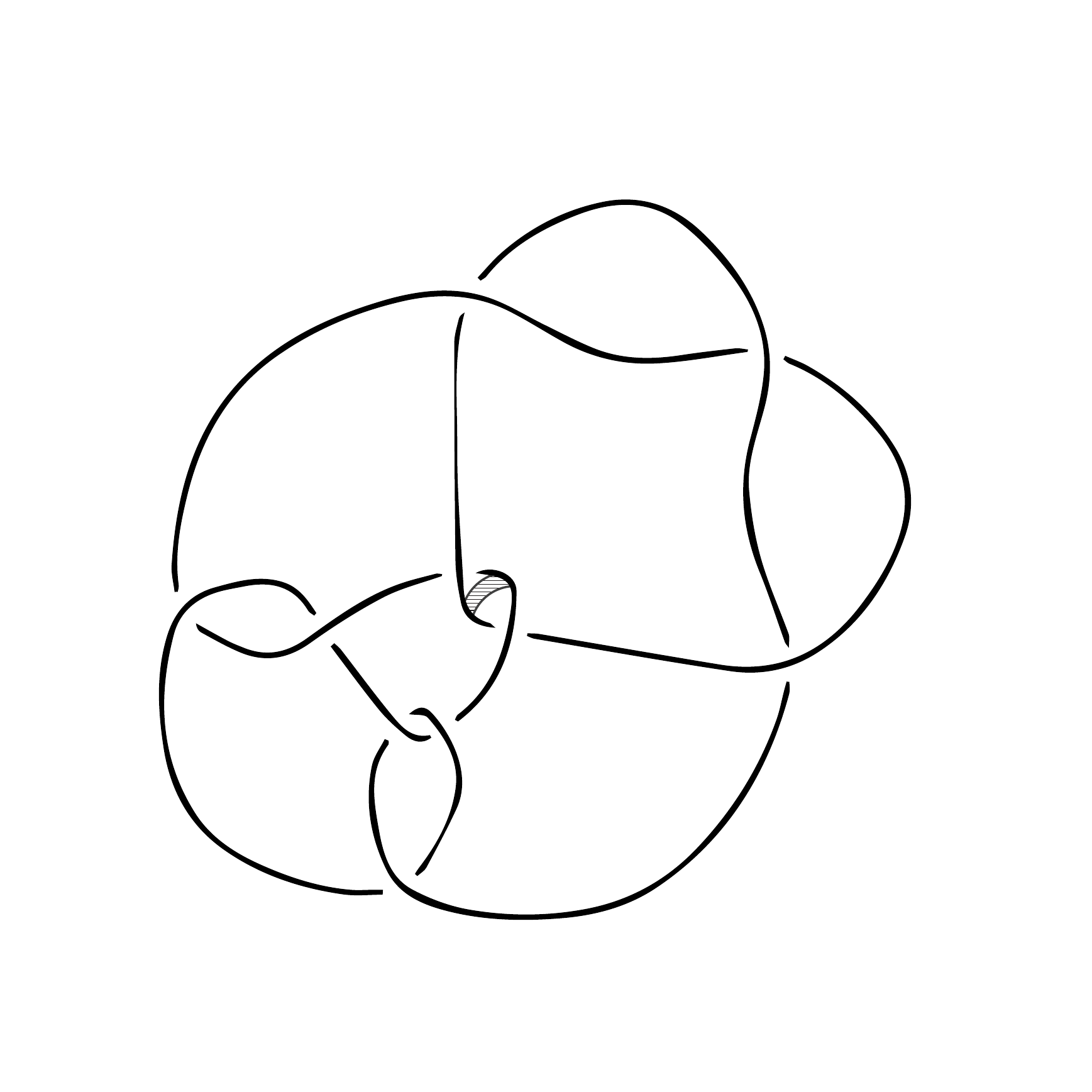}
		\caption{$10_{80}\stackrel{0}{\longrightarrow} 8_{8}$}
		\label{FigureFor10-80}
	\end{subfigure}
		~
		\begin{subfigure}[b]{0.27\textwidth}
		\includegraphics[width=\textwidth]{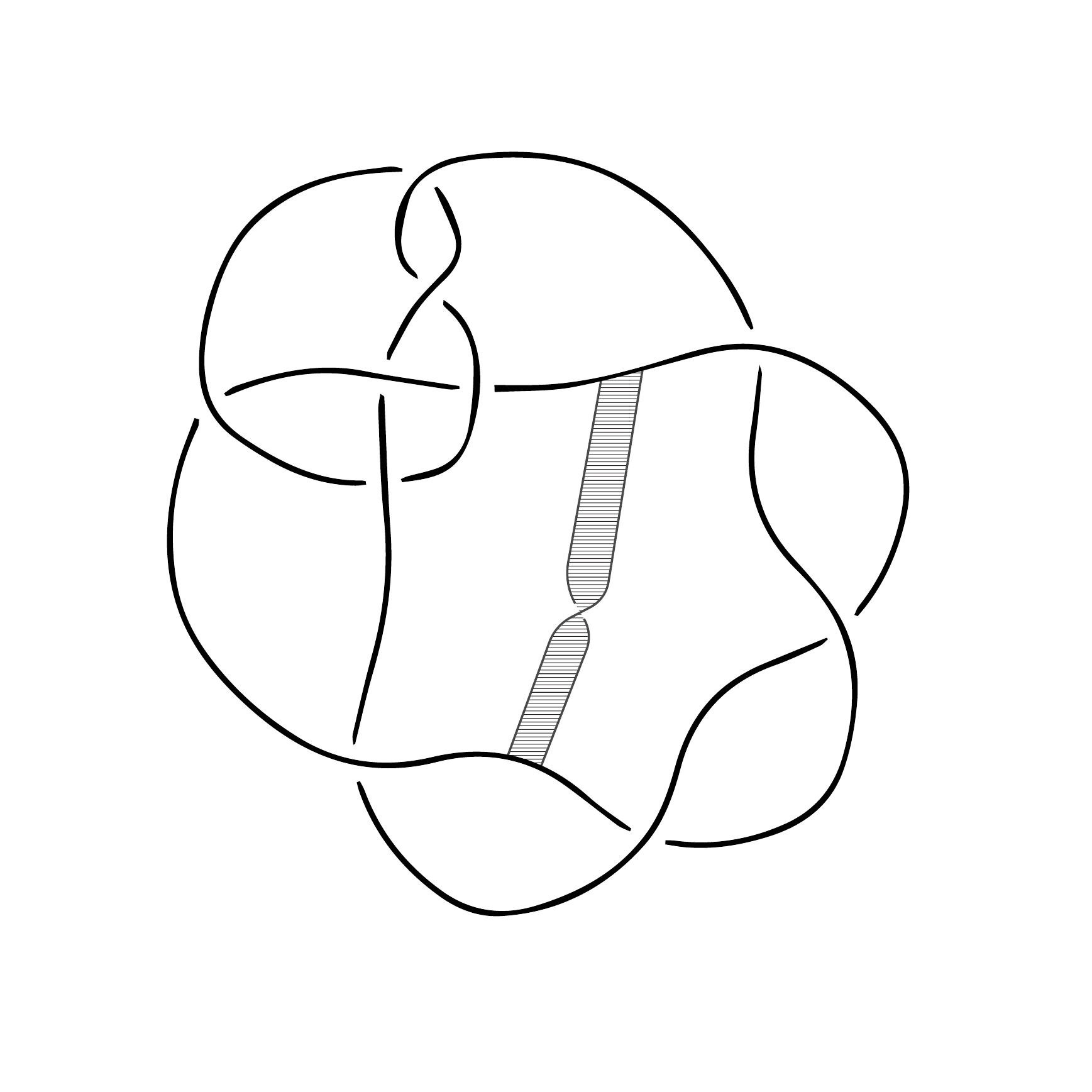}
		\caption{$10_{82}\stackrel{1}{\longrightarrow} 8_{20}$}
		\label{FigureFor10-82}
	\end{subfigure}
	~	
		\begin{subfigure}[b]{0.27\textwidth}
		\includegraphics[width=\textwidth]{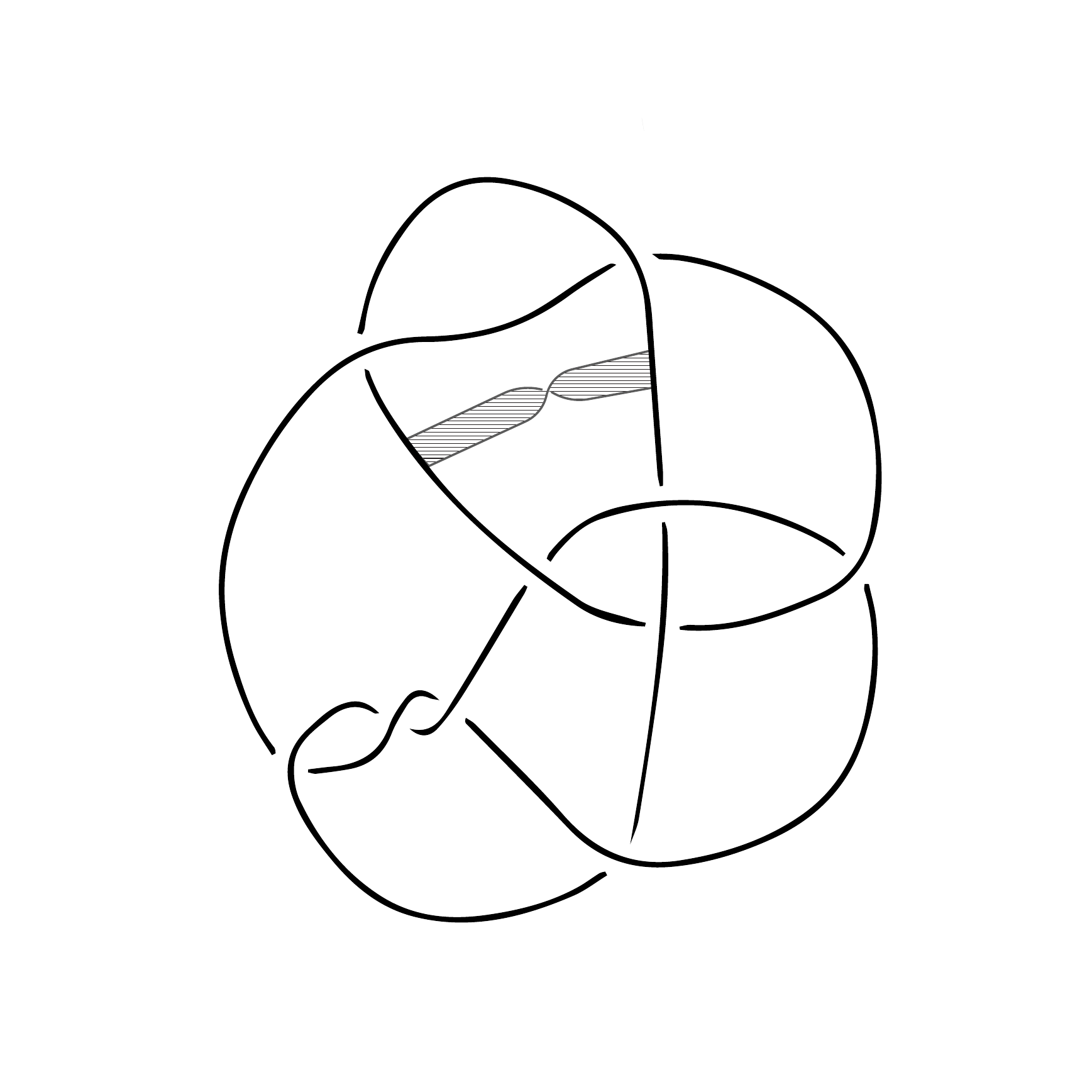}
		\caption{$10_{83}\stackrel{-1}{\longrightarrow} 10_{129}$}
		\label{FigureFor10-83}
	\end{subfigure}
	\vskip3mm
	~
	\begin{subfigure}[b]{0.27\textwidth}
		\includegraphics[width=\textwidth]{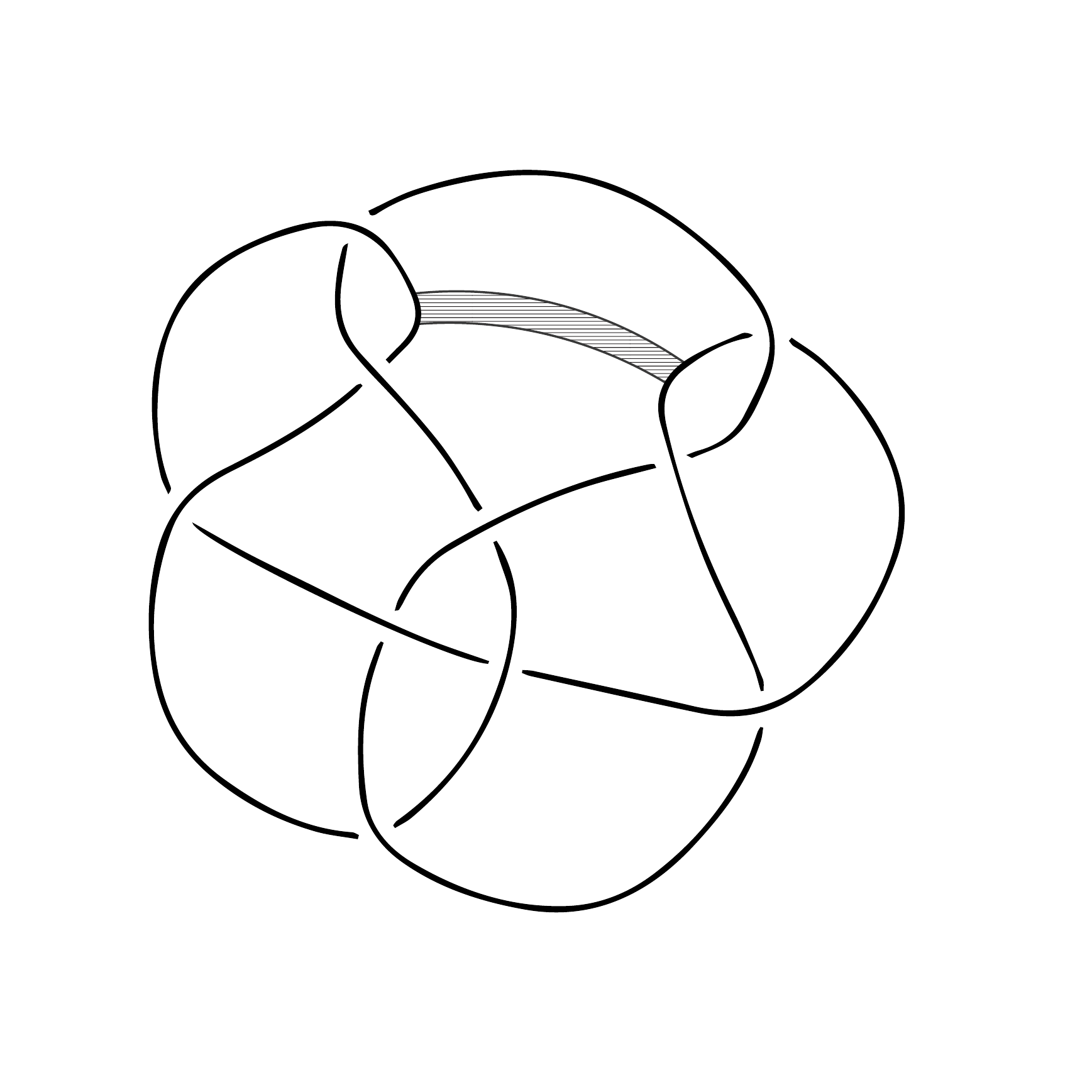}
		\caption{$10_{89}\stackrel{0}{\longrightarrow} 10_{87}$}
		\label{FigureFor10-89}
	\end{subfigure}
~
		\begin{subfigure}[b]{0.3\textwidth}
		\includegraphics[width=\textwidth]{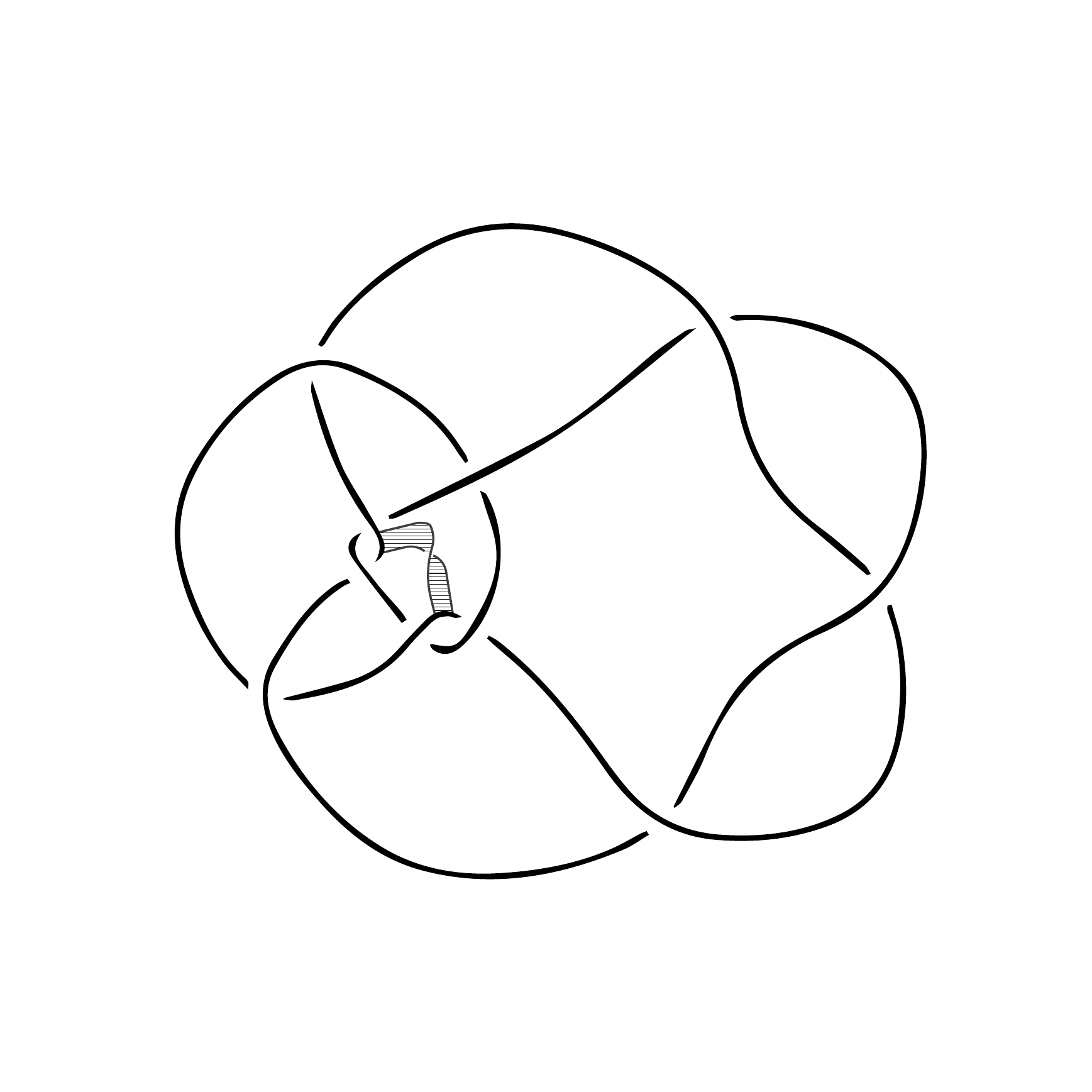}
		\caption{$10_{91}\stackrel{1}{\longrightarrow} 6_{1}$}
		\label{FigureFor10-91}
	\end{subfigure}
	~
		\begin{subfigure}[b]{0.27\textwidth}
		\includegraphics[width=\textwidth]{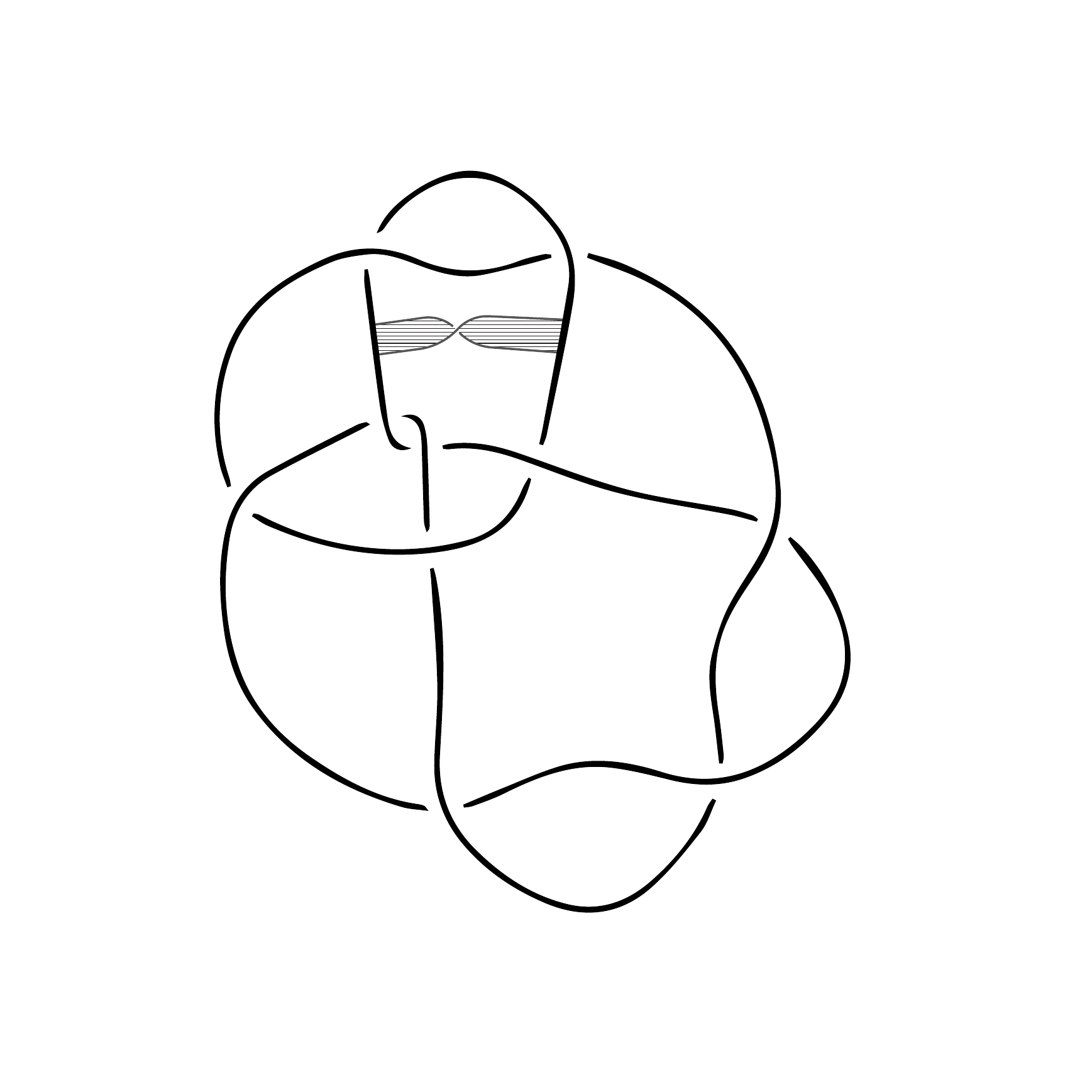}
		\caption{$10_{93}\stackrel{-1}{\longrightarrow} 10_{140}$}
		\label{FigureFor10-93}
	\end{subfigure}
	\vskip3mm
	~
	\begin{subfigure}[b]{0.27\textwidth}
		\includegraphics[width=\textwidth]{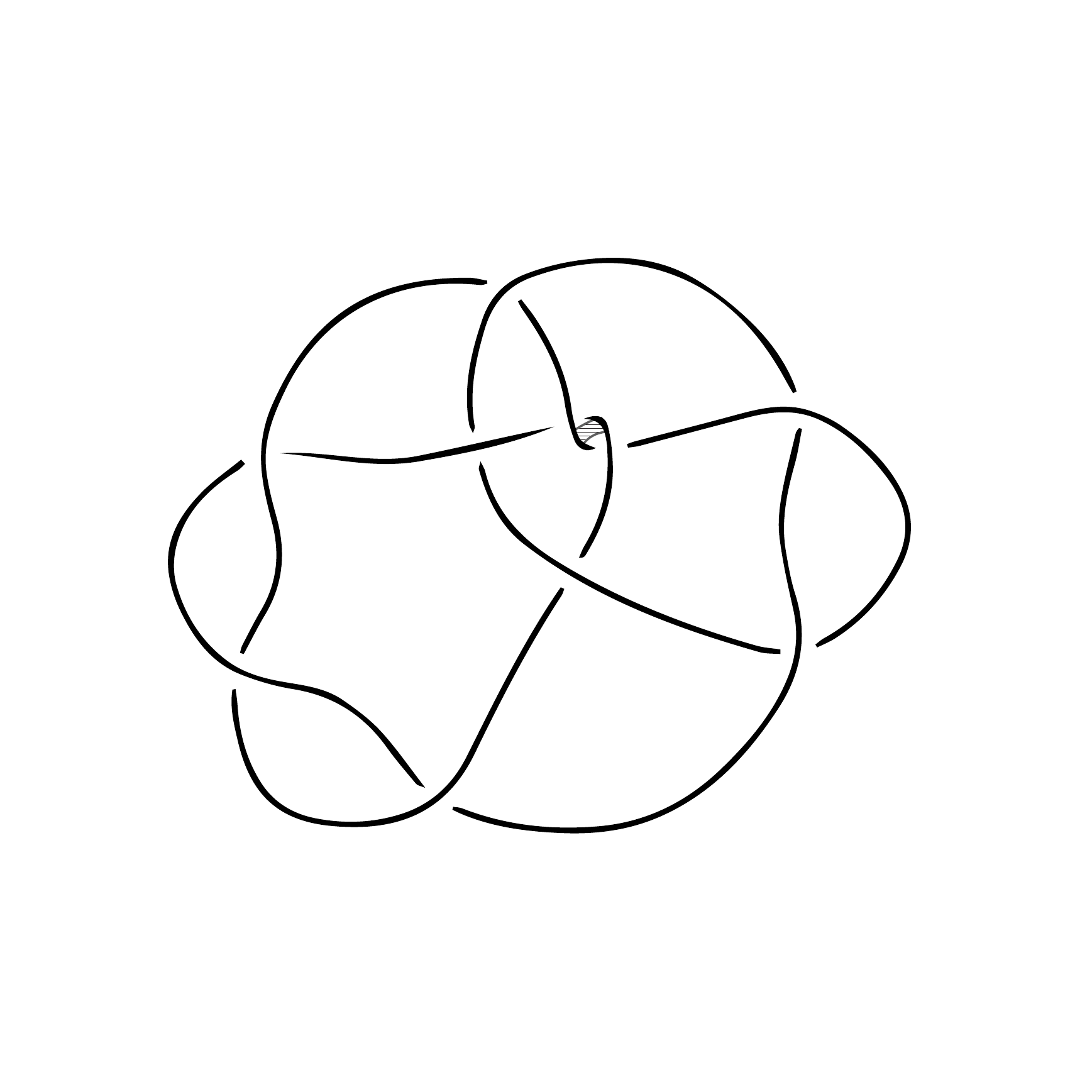}
		\caption{$10_{94}\stackrel{0}{\longrightarrow} 8_{8}$}
		\label{FigureFor10-94}
	\end{subfigure}
	~	      
	\begin{subfigure}[b]{0.27\textwidth}
		\includegraphics[width=\textwidth]{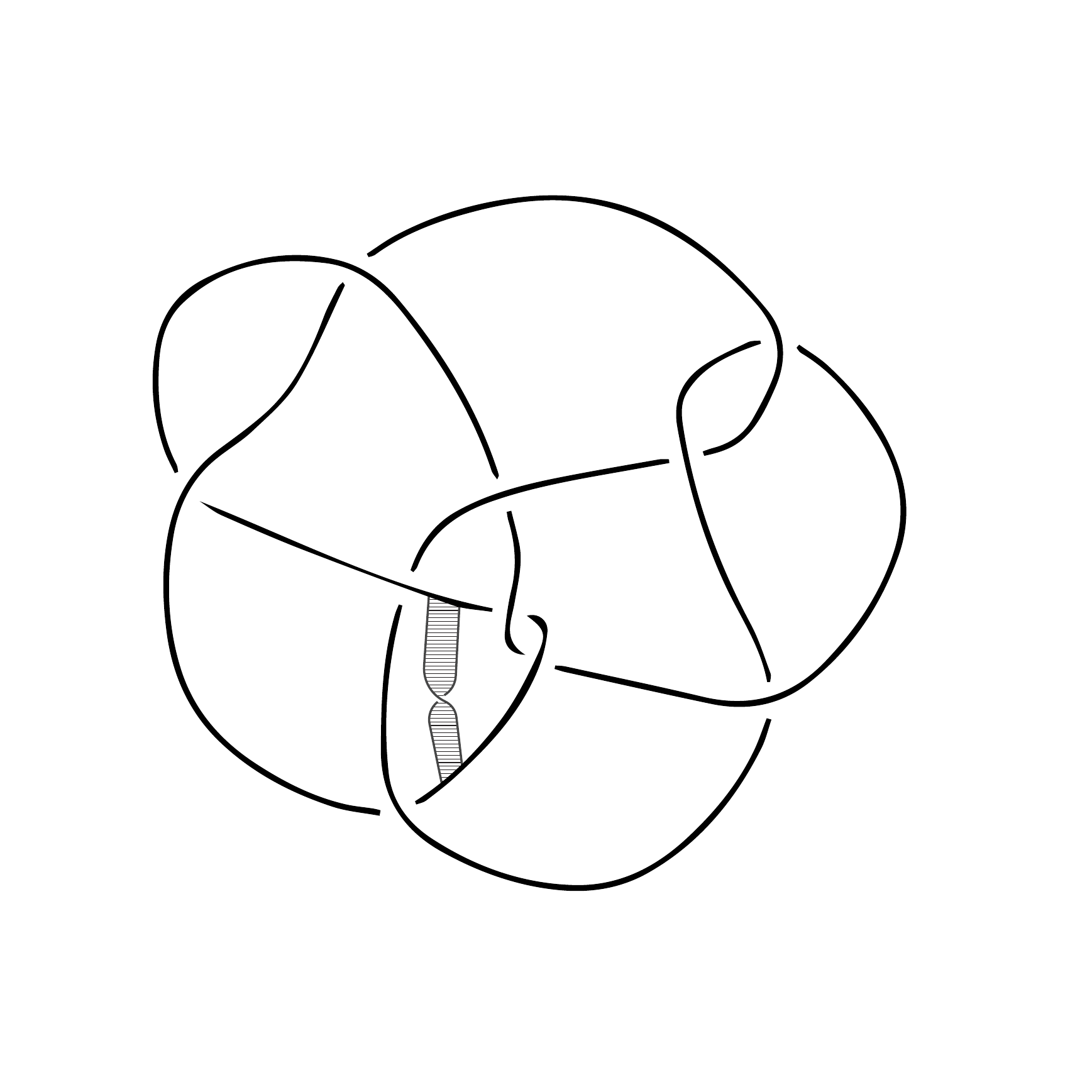}
		\caption{$10_{97}\stackrel{-1}{\longrightarrow} 10_{137}$}
		\label{FigureFor10-97}
	\end{subfigure}
	~
	\begin{subfigure}[b]{0.27\textwidth}
		\includegraphics[width=\textwidth]{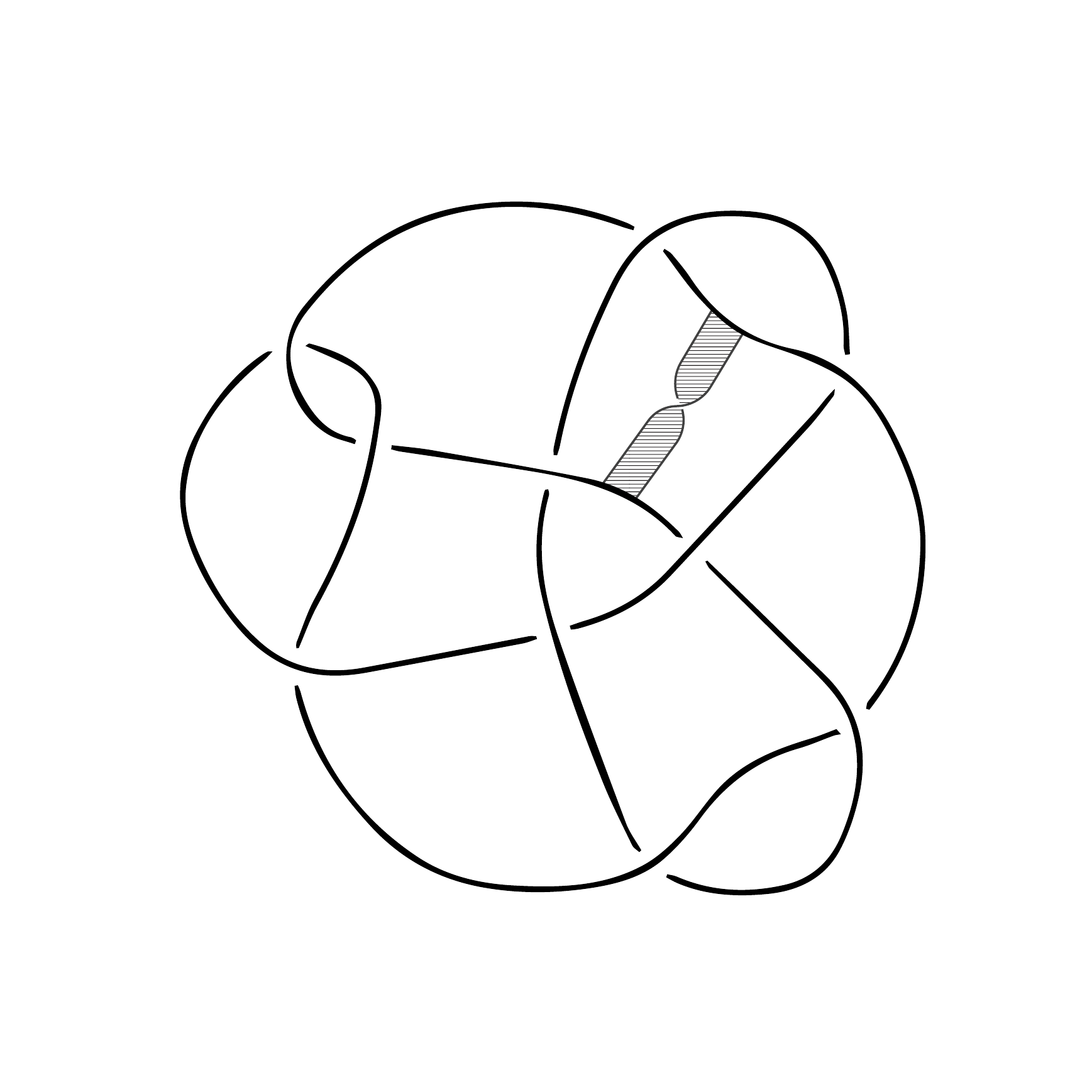}
		\caption{$10_{101}\stackrel{1\phantom{i}}{\longrightarrow} 6_{1}$}
		\label{FigureFor10-101}
	\end{subfigure}
	\vskip3mm
	~
	
	\caption{Non-oriented band moves from the knots $10_{80}$, $10_{82}$, $10_{83}$, $10_{89}$, $10_{91}$, $10_{93}$, $10_{94}$, $10_{97}$, $10_{101}$, to slice knots}\label{slice6}
\end{figure}
\begin{figure}[h]
	\centering
	\begin{subfigure}[b]{0.3\textwidth}
		\includegraphics[width=\textwidth]{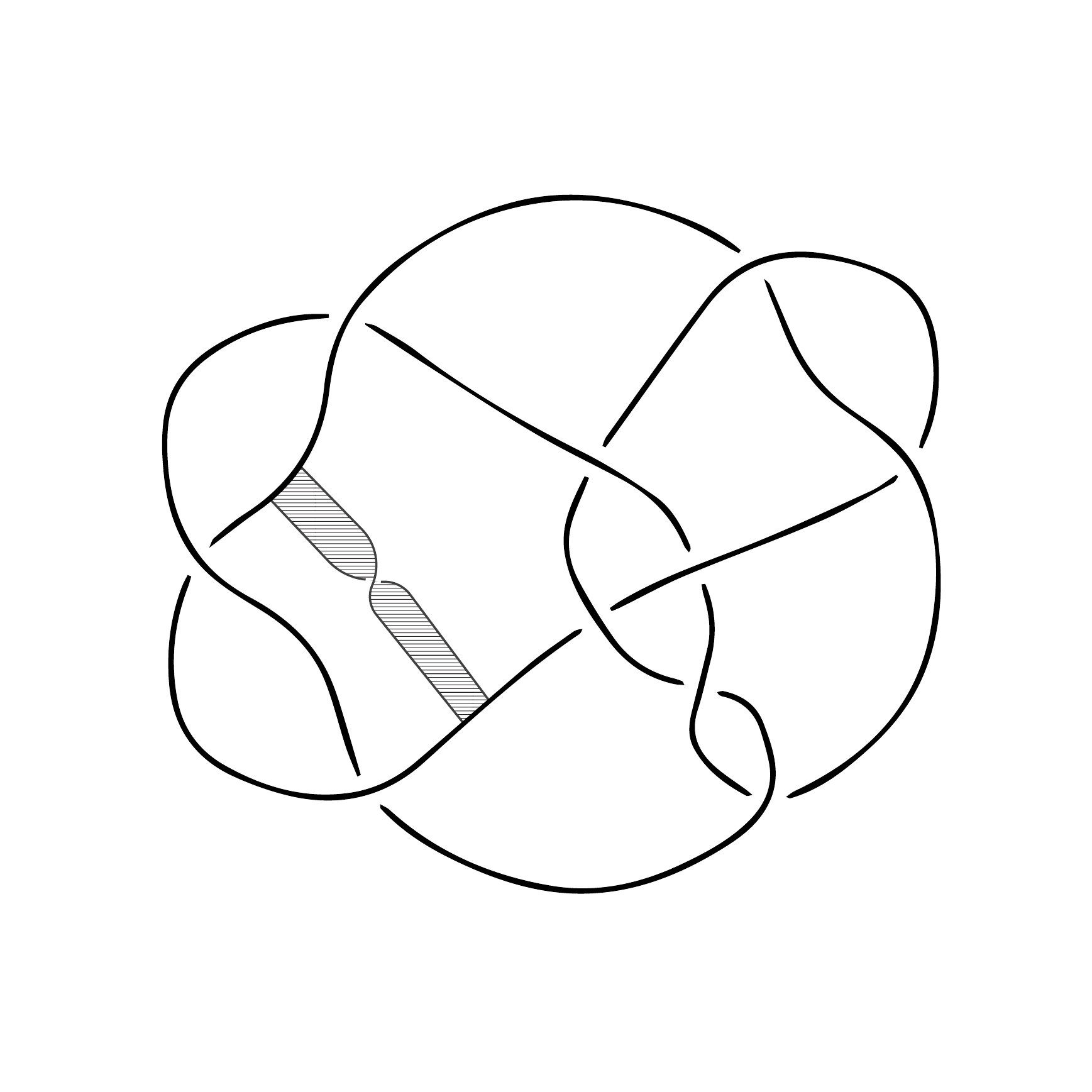}
		\caption{$10_{102}\stackrel{1}{\longrightarrow}0_1$}
		\label{FigureFor10-102}
	\end{subfigure}
~
		\begin{subfigure}[b]{0.3\textwidth}
		\includegraphics[width=\textwidth]{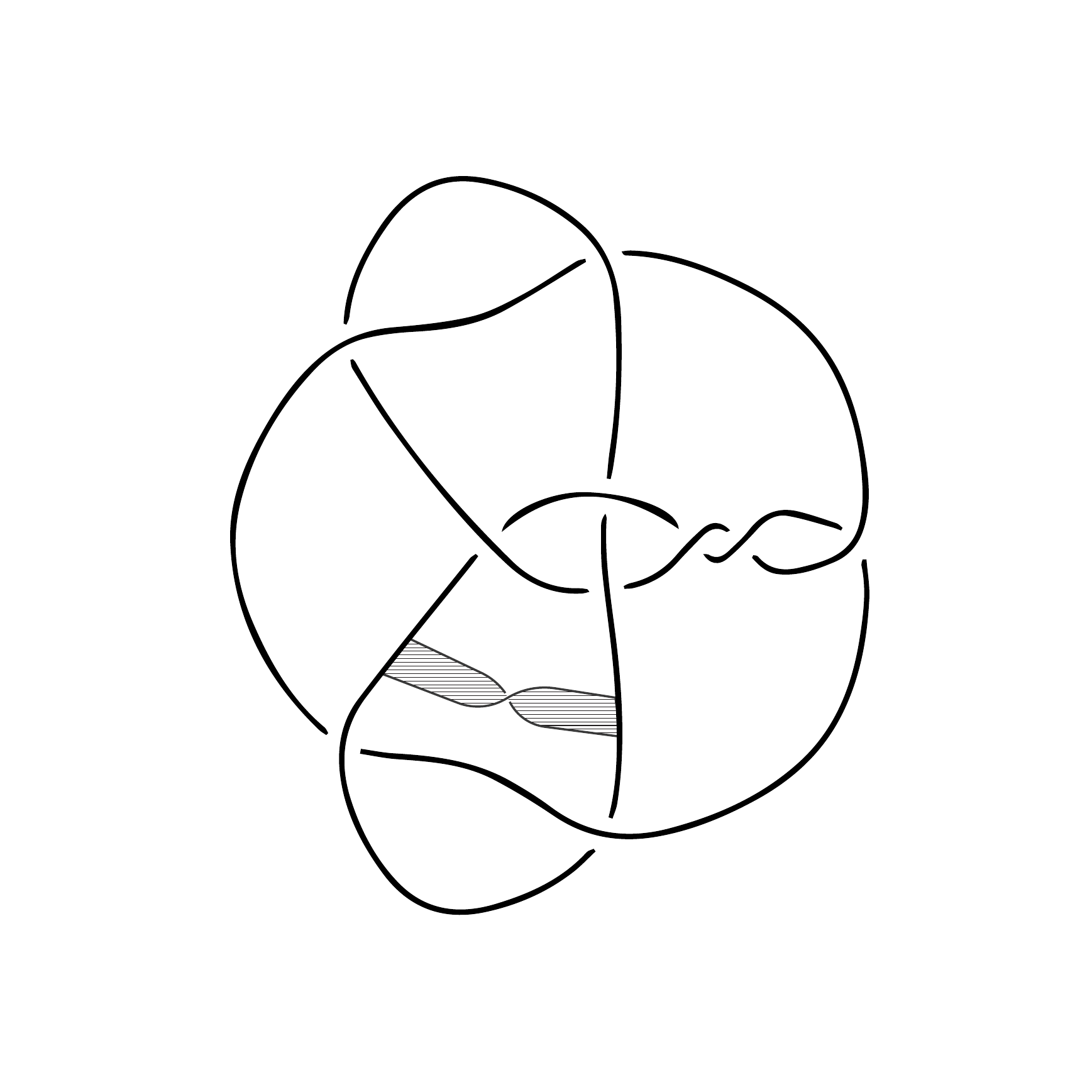}
		\caption{$10_{103}\stackrel{-1}{\longrightarrow} 10_{129}$}
		\label{FigureFor10-103}
	\end{subfigure}
	~
	\begin{subfigure}[b]{0.3\textwidth}
		\includegraphics[width=\textwidth]{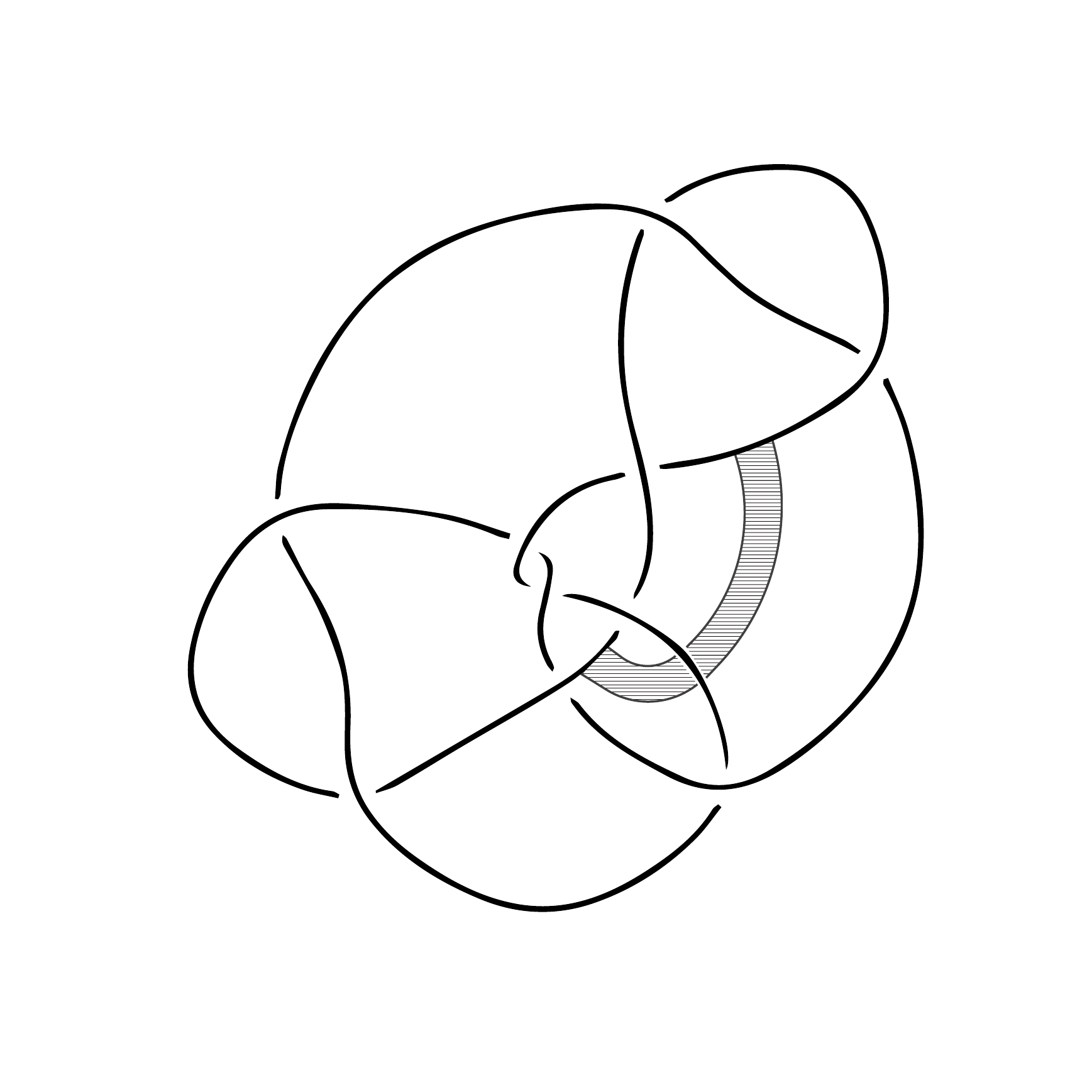}
		\caption{$10_{105}\stackrel{0}{\longrightarrow} 9_{41}$}
		\label{FigureFor10-105}
	\end{subfigure}
	~
	\vskip3mm
	~
			\begin{subfigure}[b]{0.27\textwidth}
		\includegraphics[width=\textwidth]{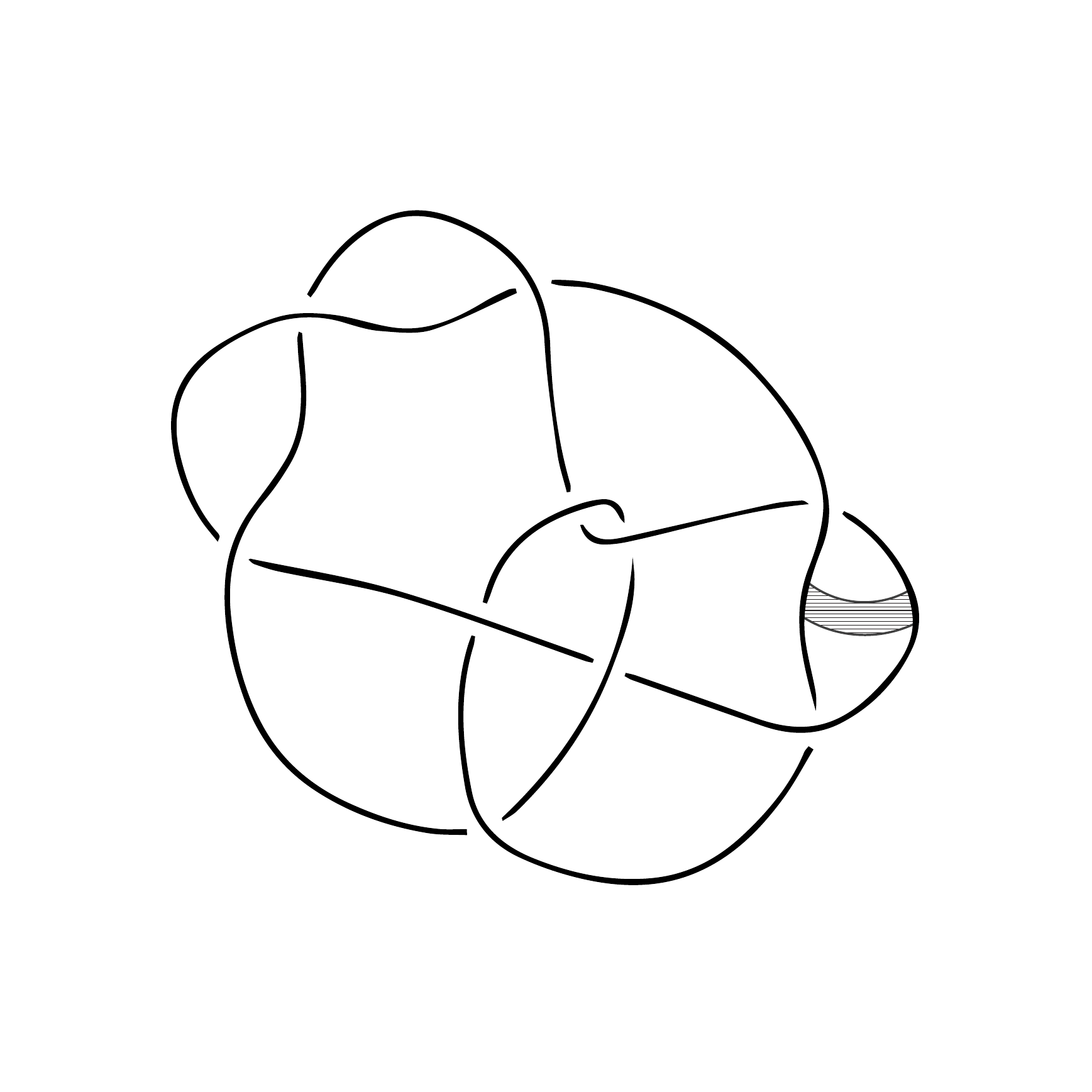}
		\caption{$10_{106}\stackrel{0}{\longrightarrow} 8_{8}$}
		\label{FigureFor10-106}
	\end{subfigure}
		~
		\begin{subfigure}[b]{0.27\textwidth}
		\includegraphics[width=\textwidth]{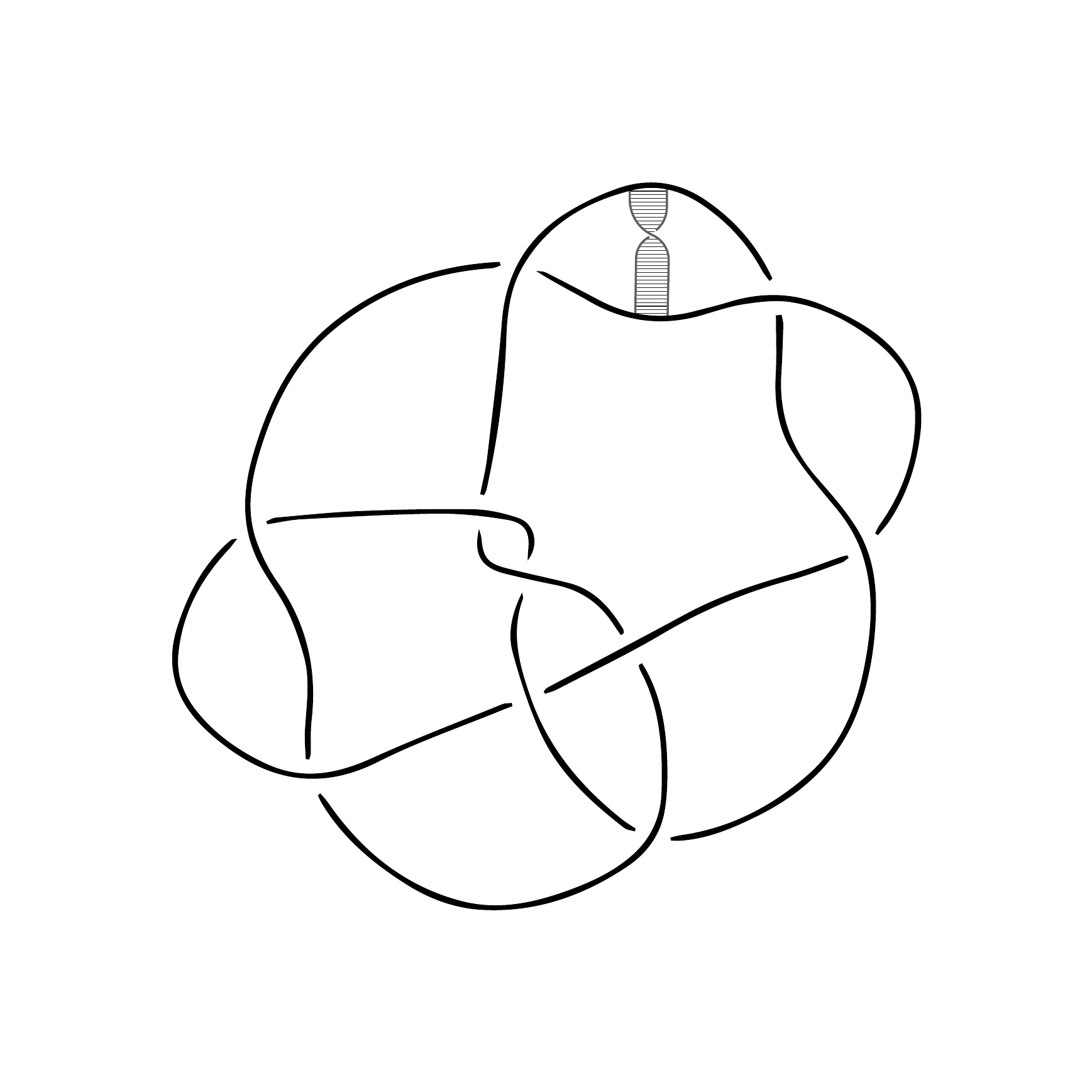}
		\caption{$10_{108}\stackrel{-1}{\longrightarrow} 9_{41}$}
		\label{FigureFor10-108}
	\end{subfigure}
	~
	\begin{subfigure}[b]{0.3\textwidth}
		\includegraphics[width=\textwidth]{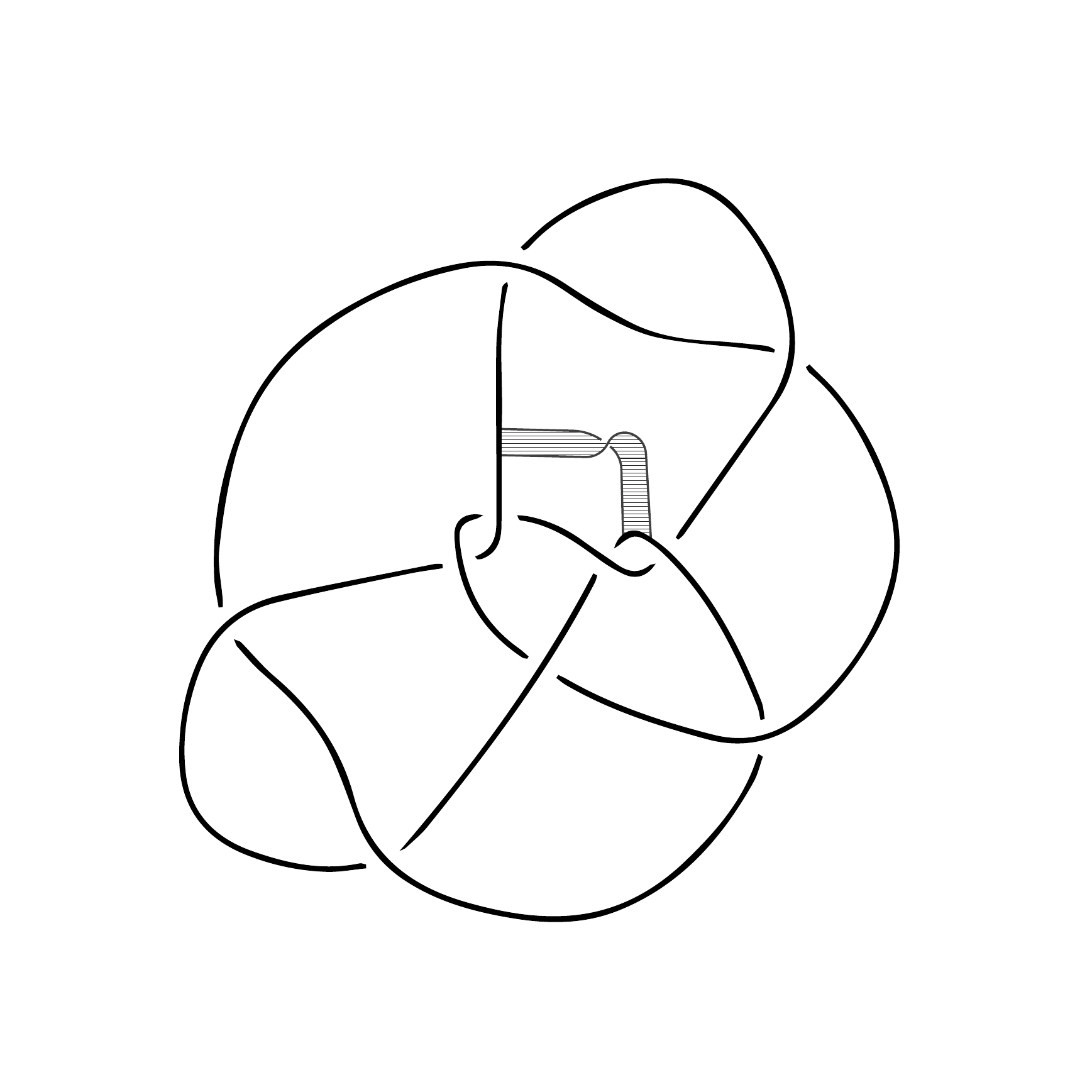}
		\caption{$10_{110}\stackrel{-1}{\longrightarrow} 8_{20}$}
		\label{FigureFor10-110}
	\end{subfigure}
	~
	\vskip3mm
	~
	\begin{subfigure}[b]{0.27\textwidth}
		\includegraphics[width=\textwidth]{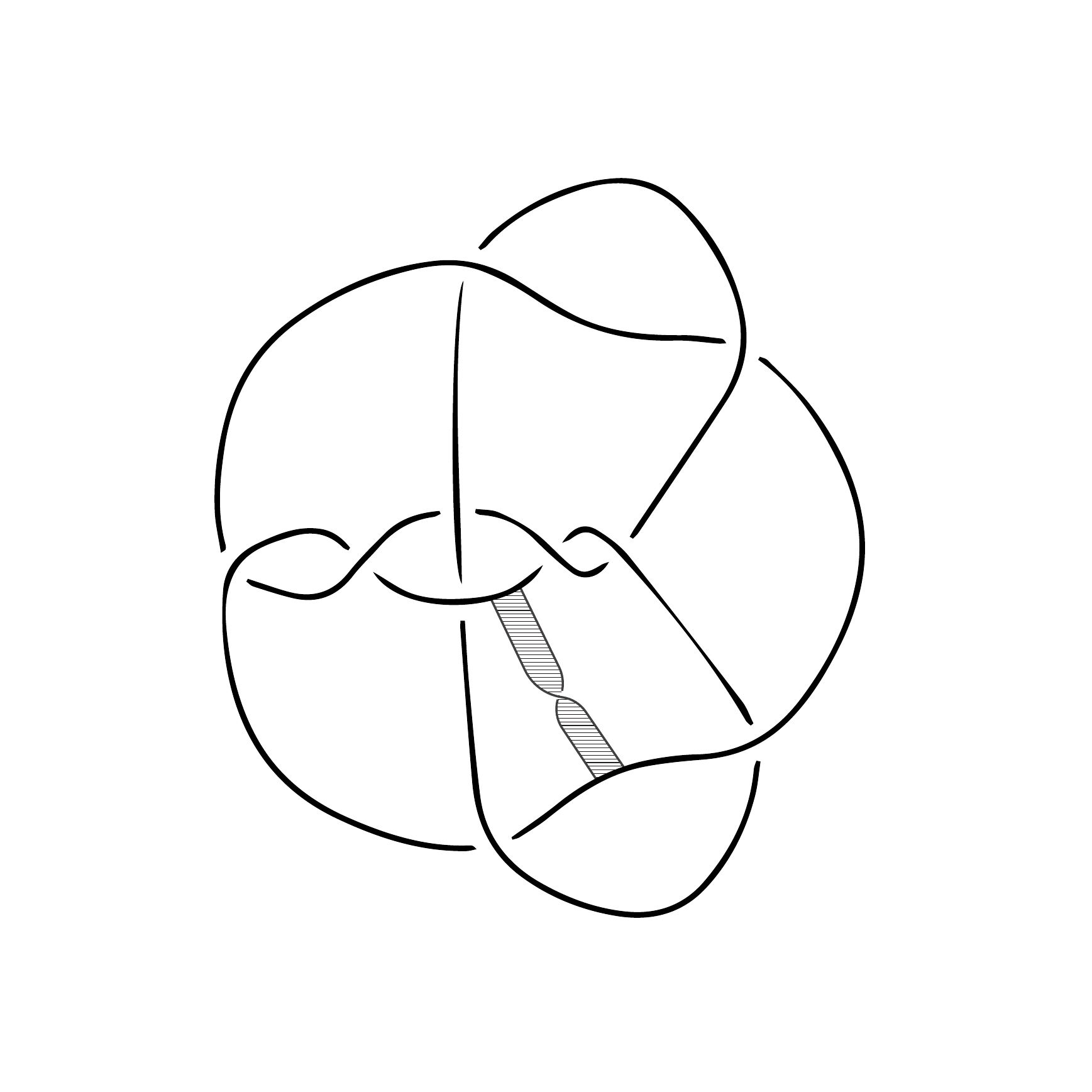}
		\caption{$10_{111}\stackrel{-1}{\longrightarrow} 0_{1}$}
		\label{FigureFor10-111}
	\end{subfigure}
	~
		\begin{subfigure}[b]{0.27\textwidth}
		\includegraphics[width=\textwidth]{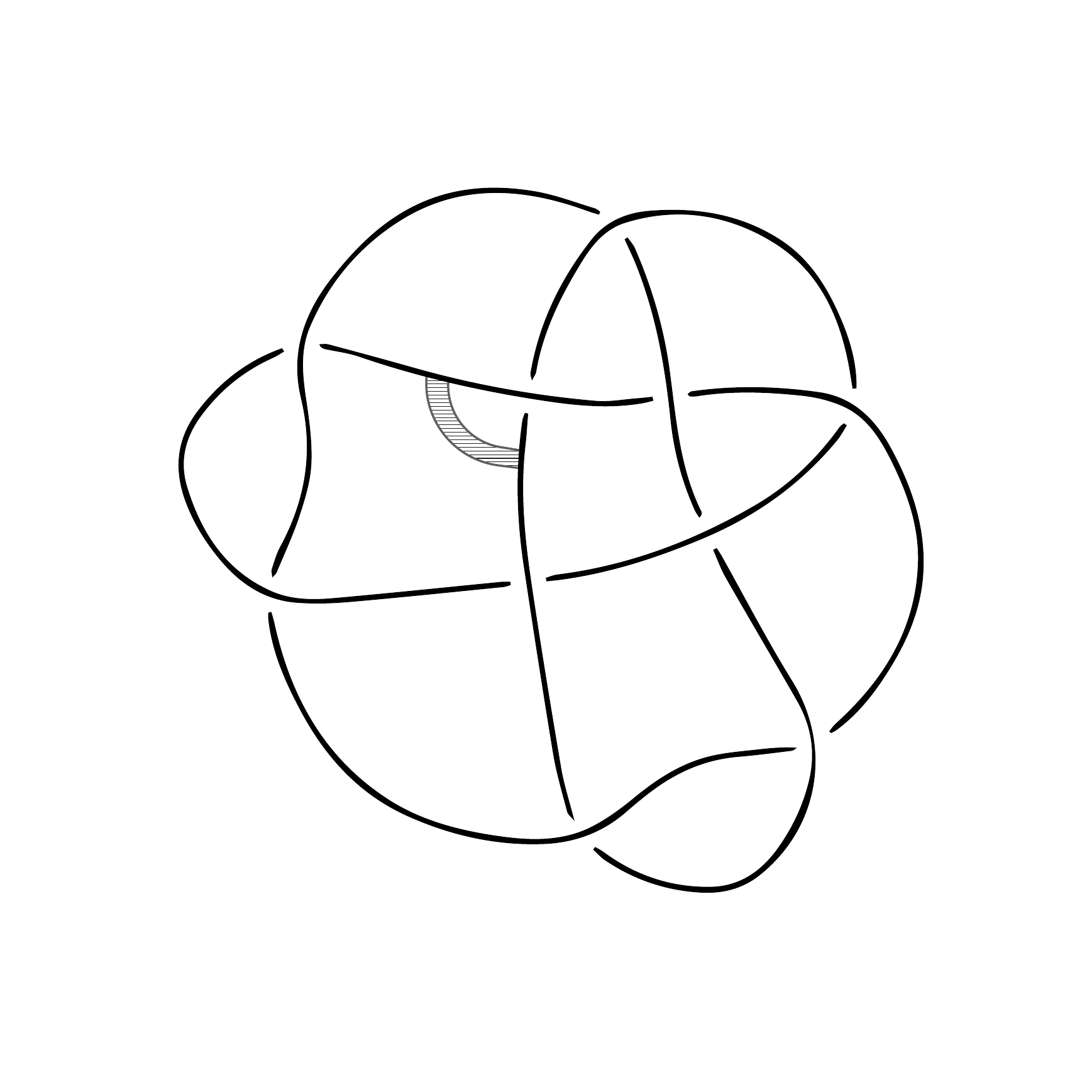}
		\caption{$10_{116}\stackrel{0}{\longrightarrow} 9_{27}$}
		\label{FigureFor10-116}
	\end{subfigure}
	~
	\begin{subfigure}[b]{0.3\textwidth}
		\includegraphics[width=\textwidth]{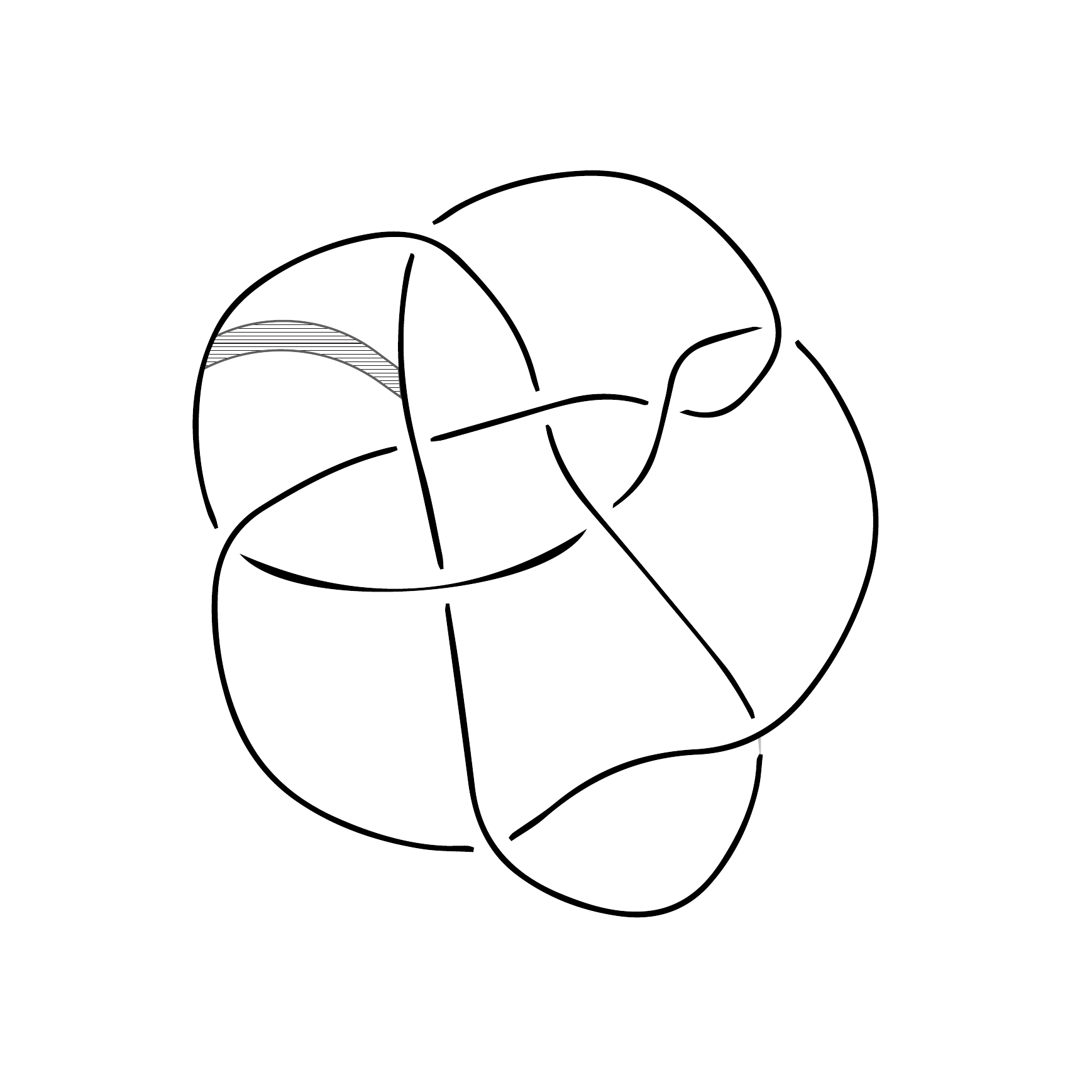}
		\caption{$10_{117}\stackrel{0}{\longrightarrow} 9_{27}$}
		\label{FigureFor10-117}
	\end{subfigure}
	~
	\vskip3mm
	~

	\caption{Non-oriented band moves from the knots $10_{102}$, $10_{103}$, $10_{105}$, $10_{106}$, $10_{108}$, $10_{110}$, $10_{111}$, $10_{116}$, $10_{117}$ to slice knots}\label{slice7}
		\end{figure}
\begin{figure}[h]
	\centering
		~
		\begin{subfigure}[b]{0.27\textwidth}
		\includegraphics[width=\textwidth]{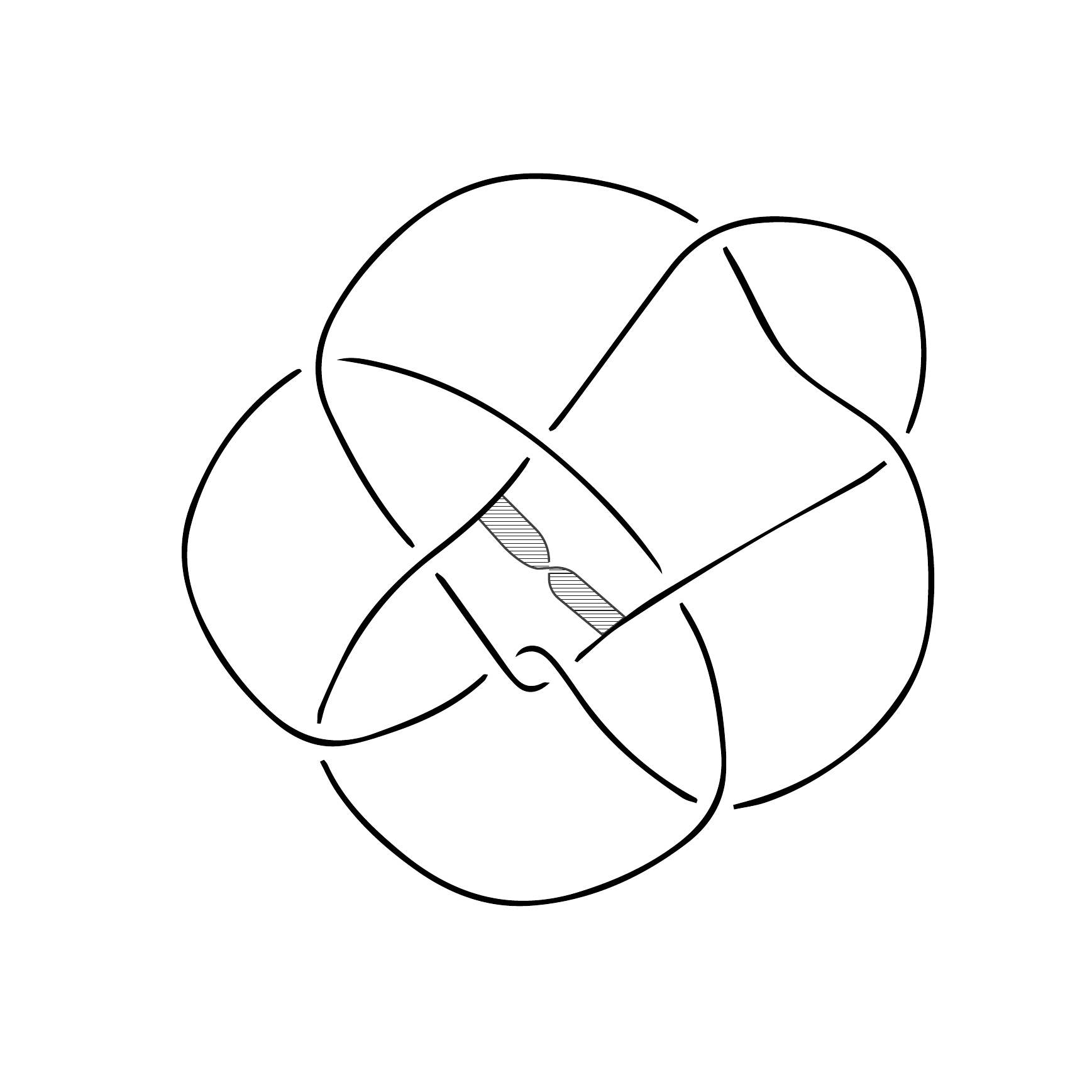}
		\caption{$10_{118}\stackrel{-1}{\longrightarrow} 10_{129}$}
		\label{FigureFor10-118}
	\end{subfigure}
~
	\begin{subfigure}[b]{0.27\textwidth}
		\includegraphics[width=\textwidth]{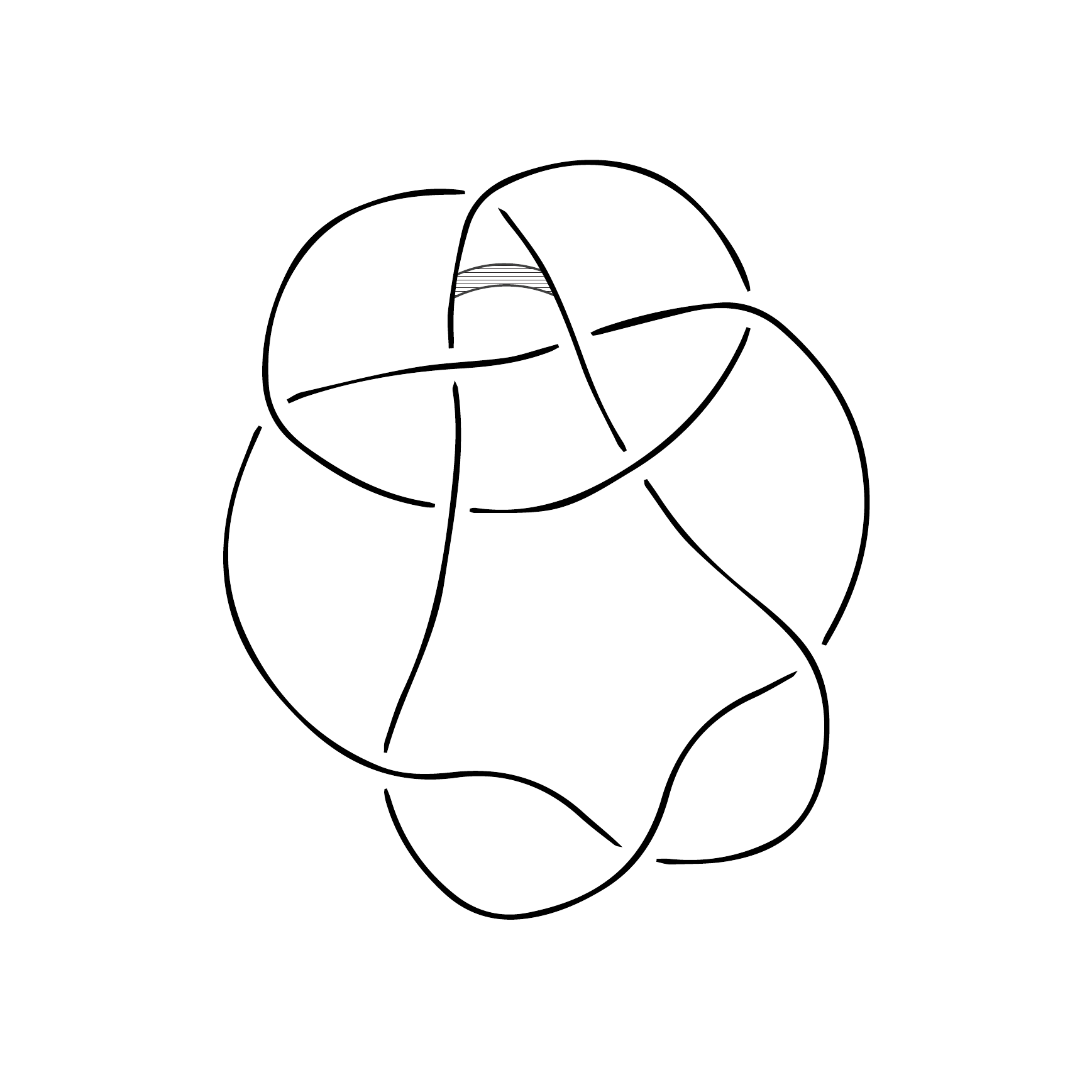}
		\caption{$10_{121}\stackrel{0}{\longrightarrow} 9_{27}$}
		\label{FigureFor10-121}
	\end{subfigure}
	~
	\begin{subfigure}[b]{0.27\textwidth}
		\includegraphics[width=\textwidth]{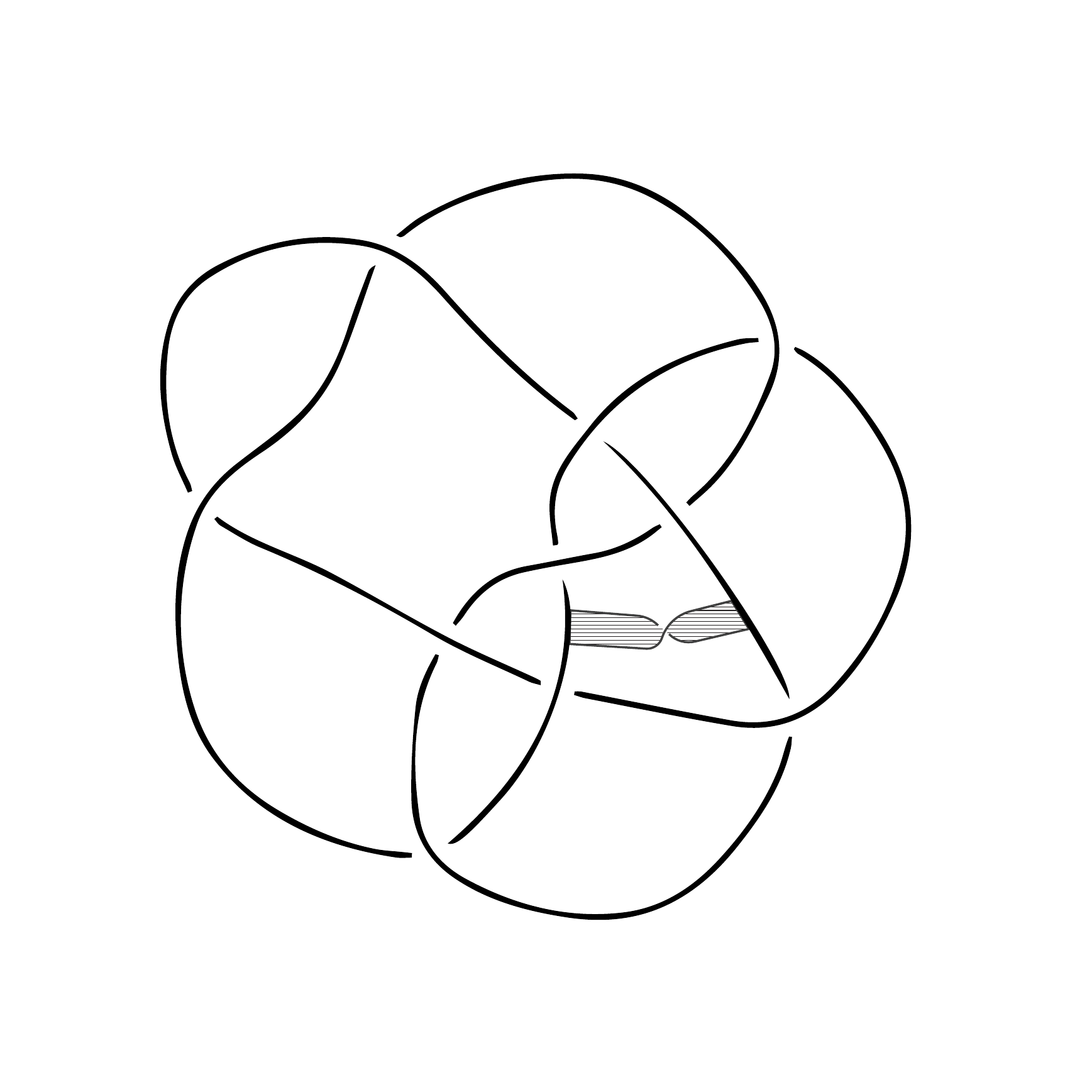}
		\caption{$10_{122}\stackrel{-1\phantom{i}}{\longrightarrow} 10_{155}$}
		\label{FigureFor10-122}
	\end{subfigure}	
	~
	\vskip3mm
	~
	\begin{subfigure}[b]{0.3\textwidth}
		\includegraphics[width=\textwidth]{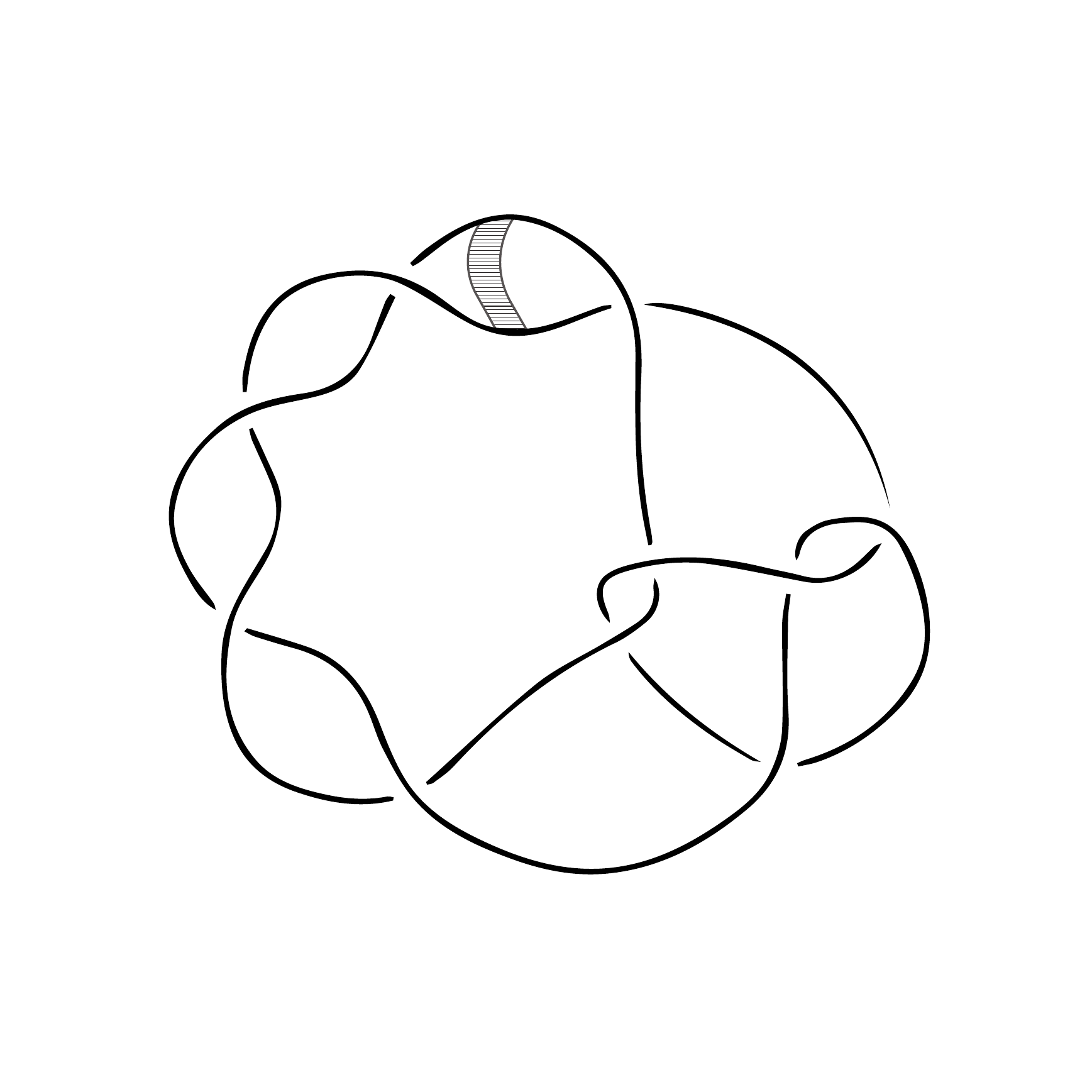}
		\caption{$10_{124}\stackrel{0}{\longrightarrow}0_1$}
		\label{FigureFor10-124}
	\end{subfigure}
	~
	\begin{subfigure}[b]{0.3\textwidth}
		\includegraphics[width=\textwidth]{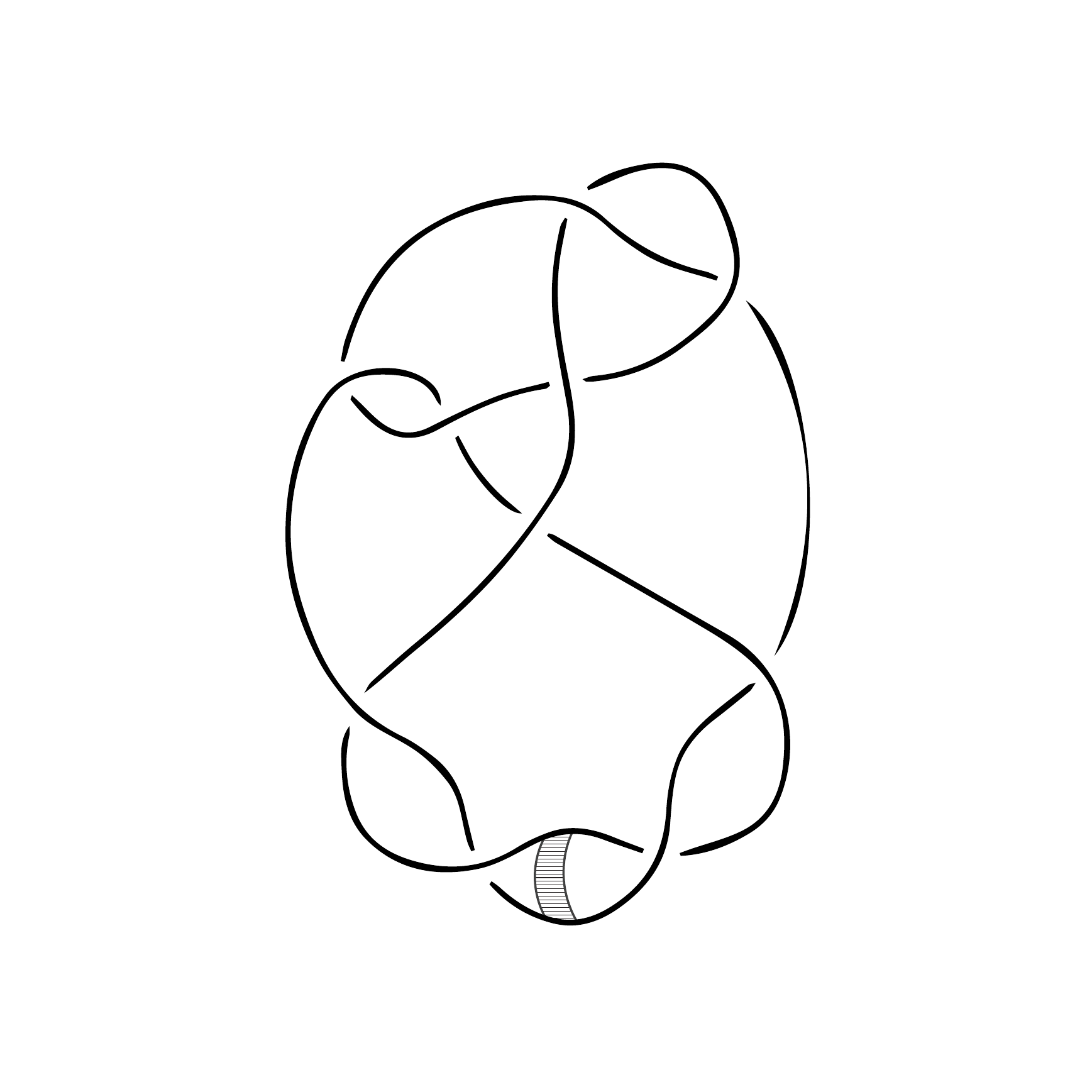}
		\caption{$10_{125}\stackrel{0}{\longrightarrow} 0_{1}$}
		\label{FigureFor10-125}
	\end{subfigure}
	~
	\begin{subfigure}[b]{0.27\textwidth}
		\includegraphics[width=\textwidth]{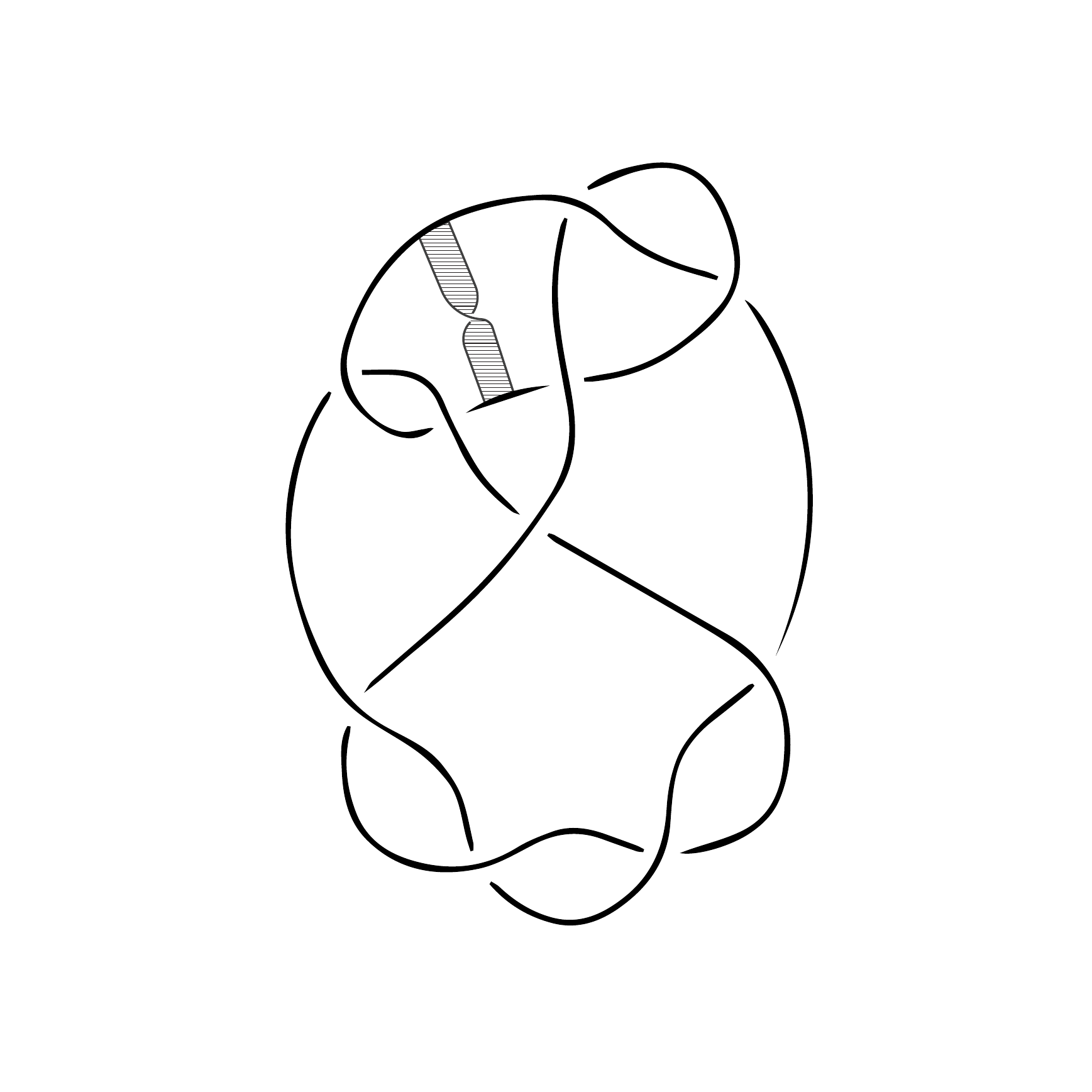}
		\caption{$10_{126}\stackrel{-1}{\longrightarrow} 10_{148}$}
		\label{FigureFor10-126}
	\end{subfigure}
	~
	\vskip3mm
	~	
		\begin{subfigure}[b]{0.27\textwidth}
		\includegraphics[width=\textwidth]{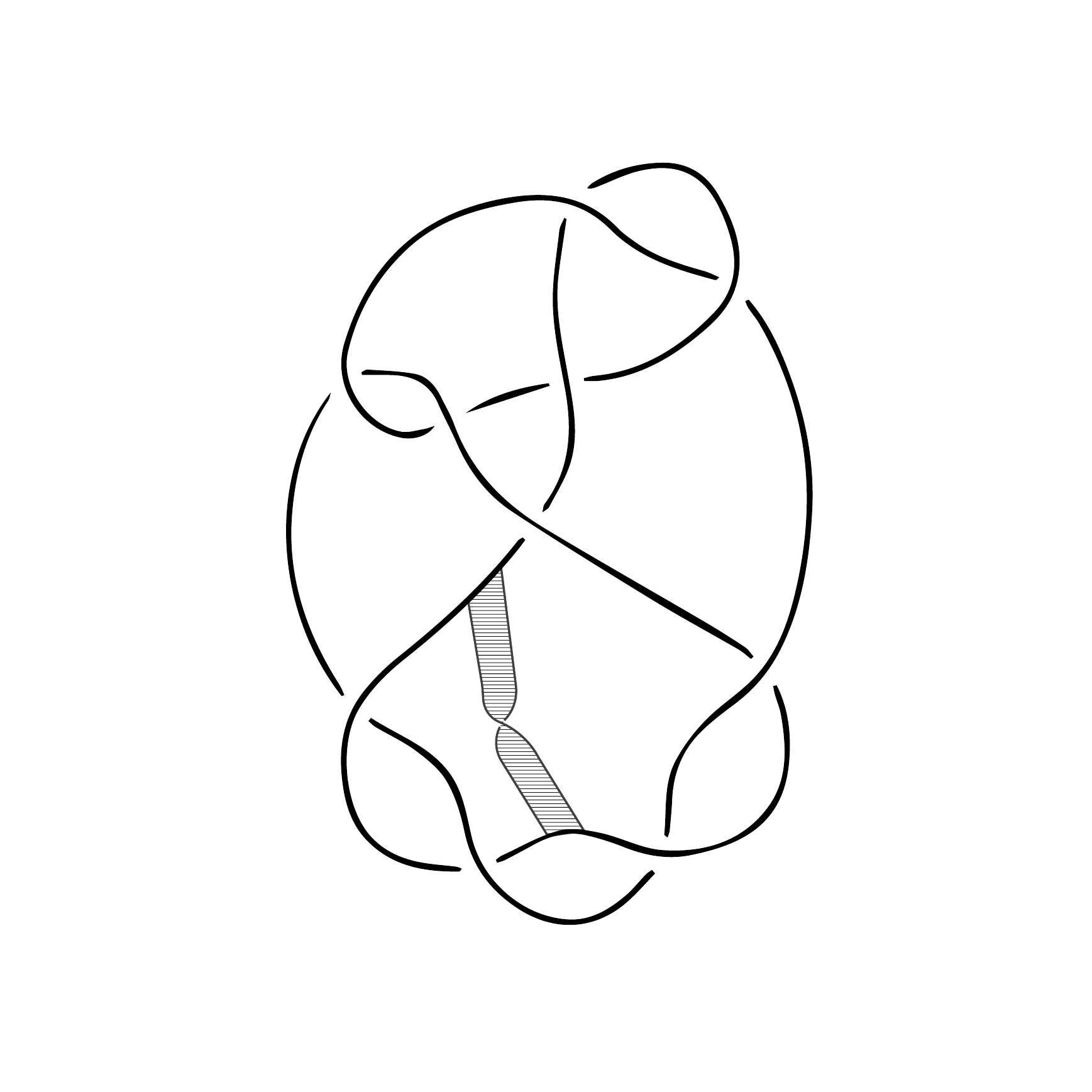}
		\caption{$10_{127}\stackrel{-1}{\longrightarrow} 0_{1}$}
		\label{FigureFor10-127}
	\end{subfigure}
	~
		\begin{subfigure}[b]{0.3\textwidth}
		\includegraphics[width=\textwidth]{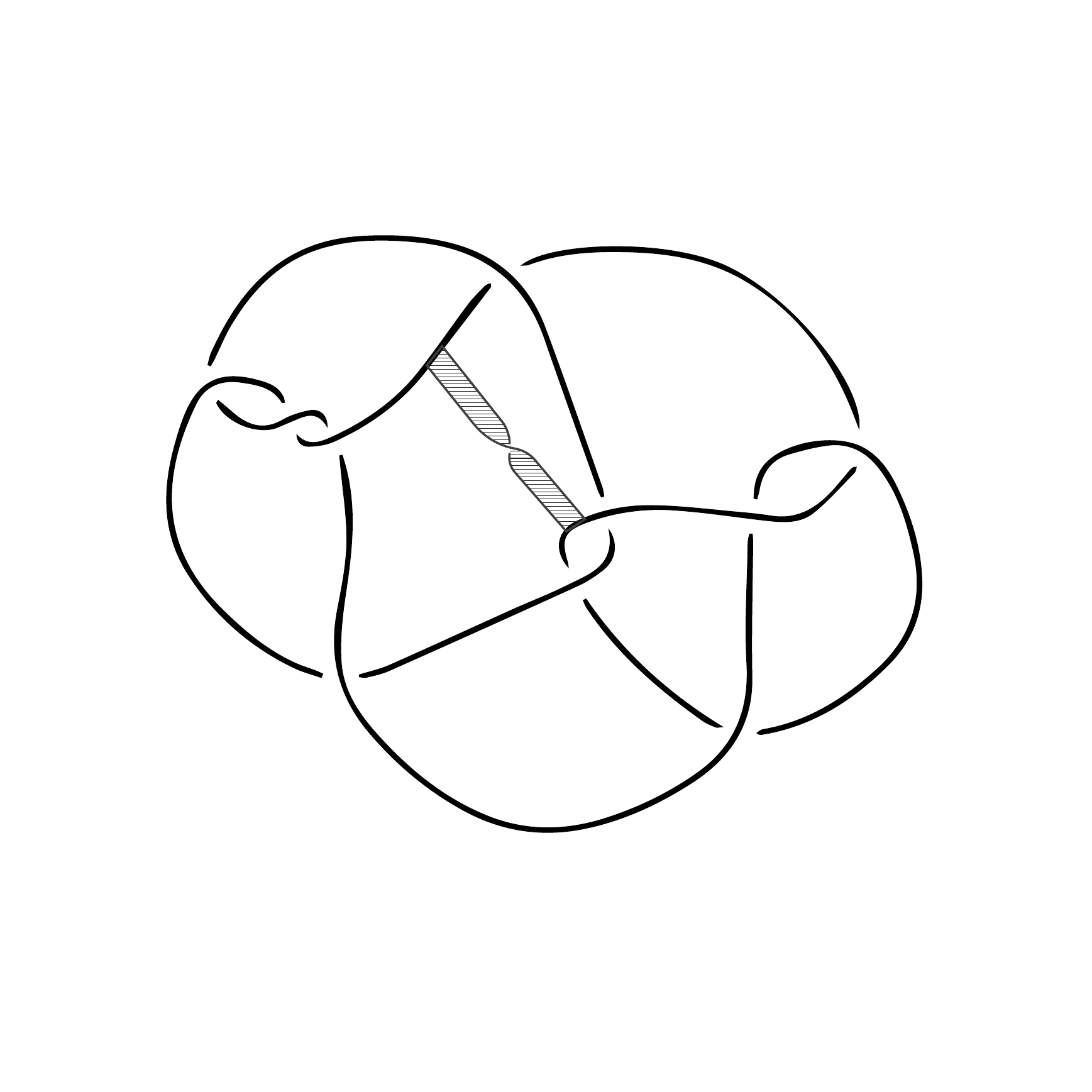}
		\caption{$10_{128}\stackrel{-1}{\longrightarrow} 0_{1}$}
		\label{FigureFor10-128}
	\end{subfigure}
	~
	\begin{subfigure}[b]{0.27\textwidth}
		\includegraphics[width=\textwidth]{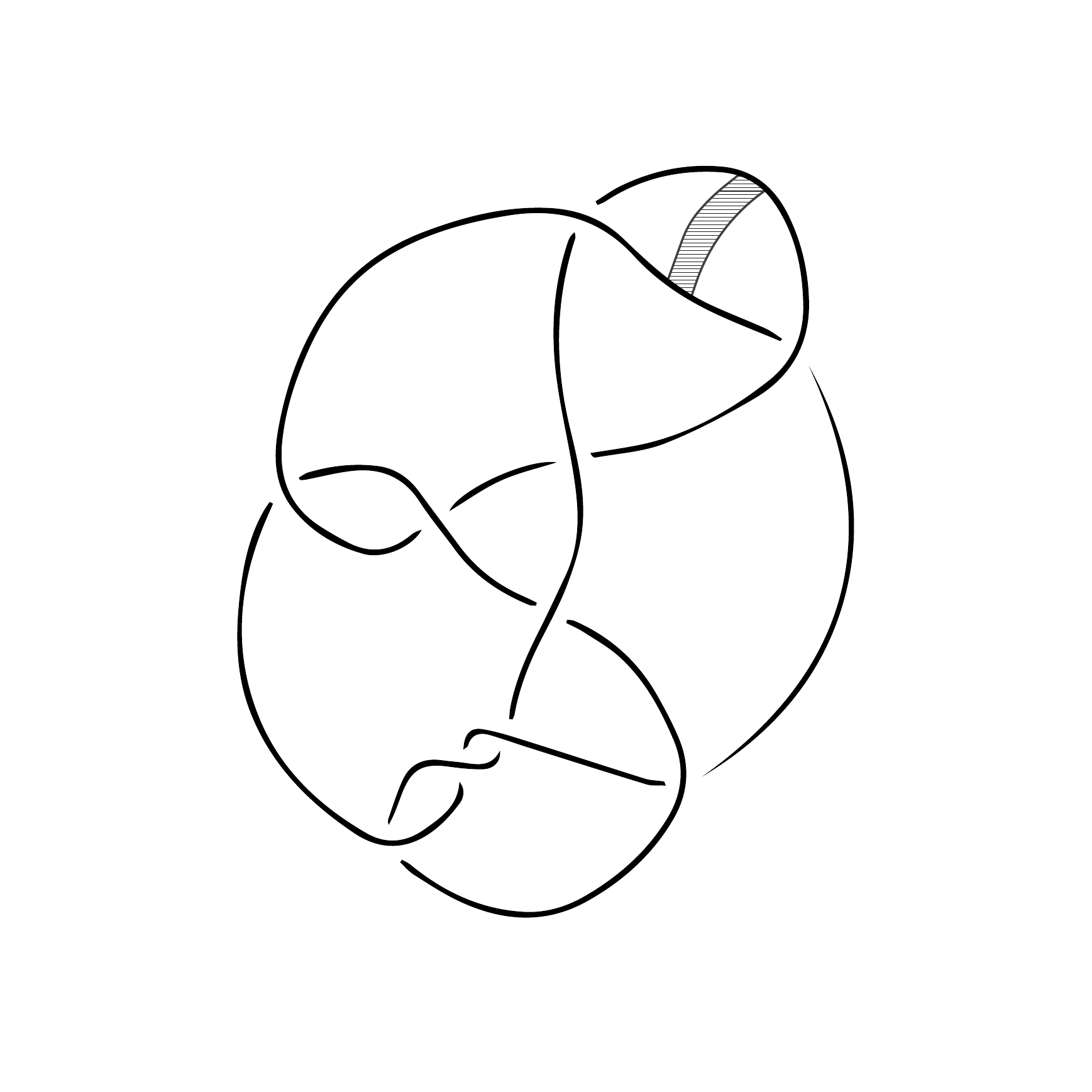}
		\caption{$10_{130}\stackrel{0}{\longrightarrow} 0_{1}$}
		\label{FigureFor10-130}
	\end{subfigure}
	~
	\vskip3mm
	~

	\caption{Non-oriented band moves from the knots $10_{118}$, $10_{121}$, $10_{122}$, $10_{124}$, $10_{125}$, $10_{126}$, $10_{127}$, $10_{128}$, $10_{130}$ to slice knots}\label{slice8}
\end{figure}
\begin{figure}[h]
	\centering
	\begin{subfigure}[b]{0.27\textwidth}
		\includegraphics[width=\textwidth]{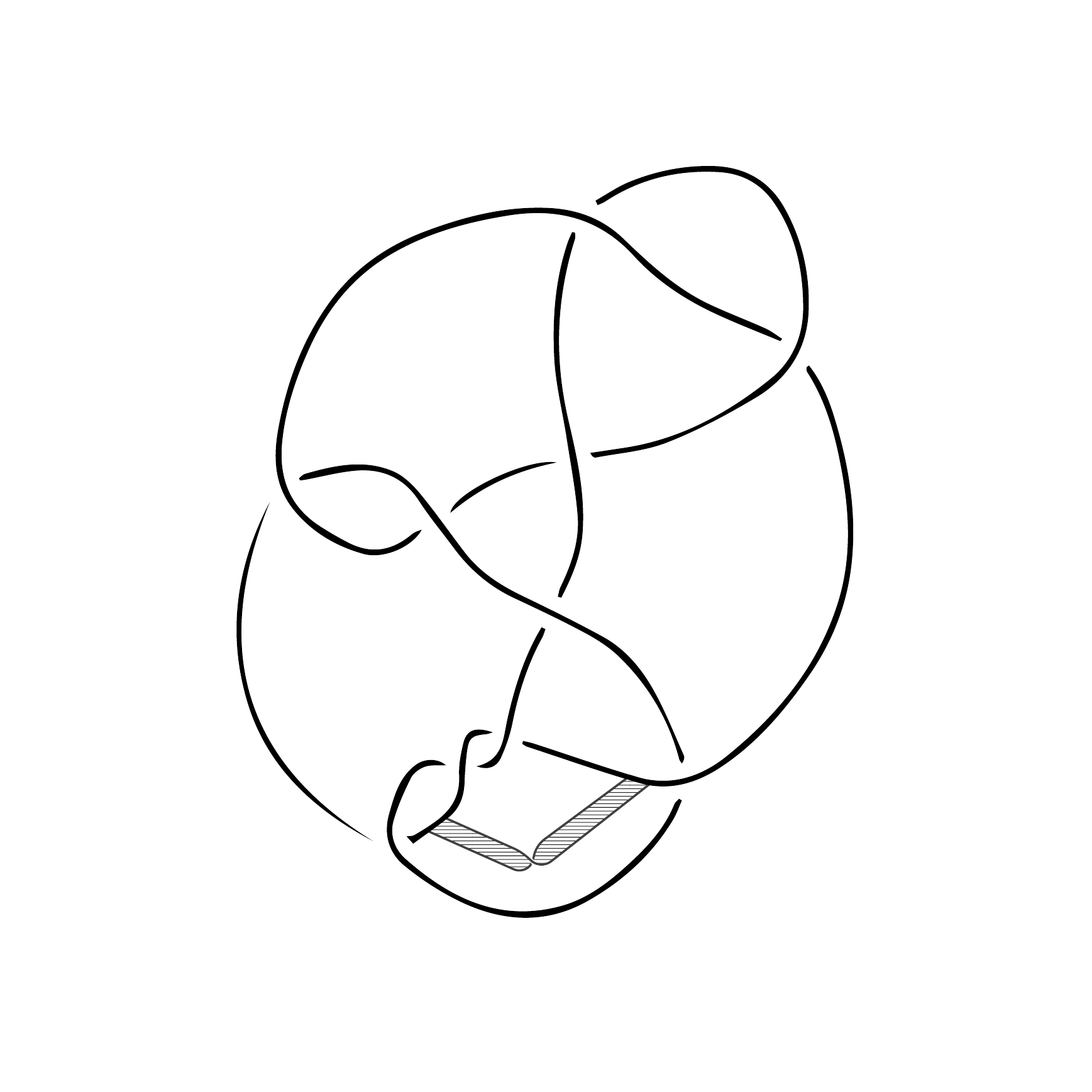}
		\caption{$10_{131}\stackrel{1}{\longrightarrow} 0_{1}$}
		\label{FigureFor10-131}
	\end{subfigure}
		~
	\begin{subfigure}[b]{0.27\textwidth}
		\includegraphics[width=\textwidth]{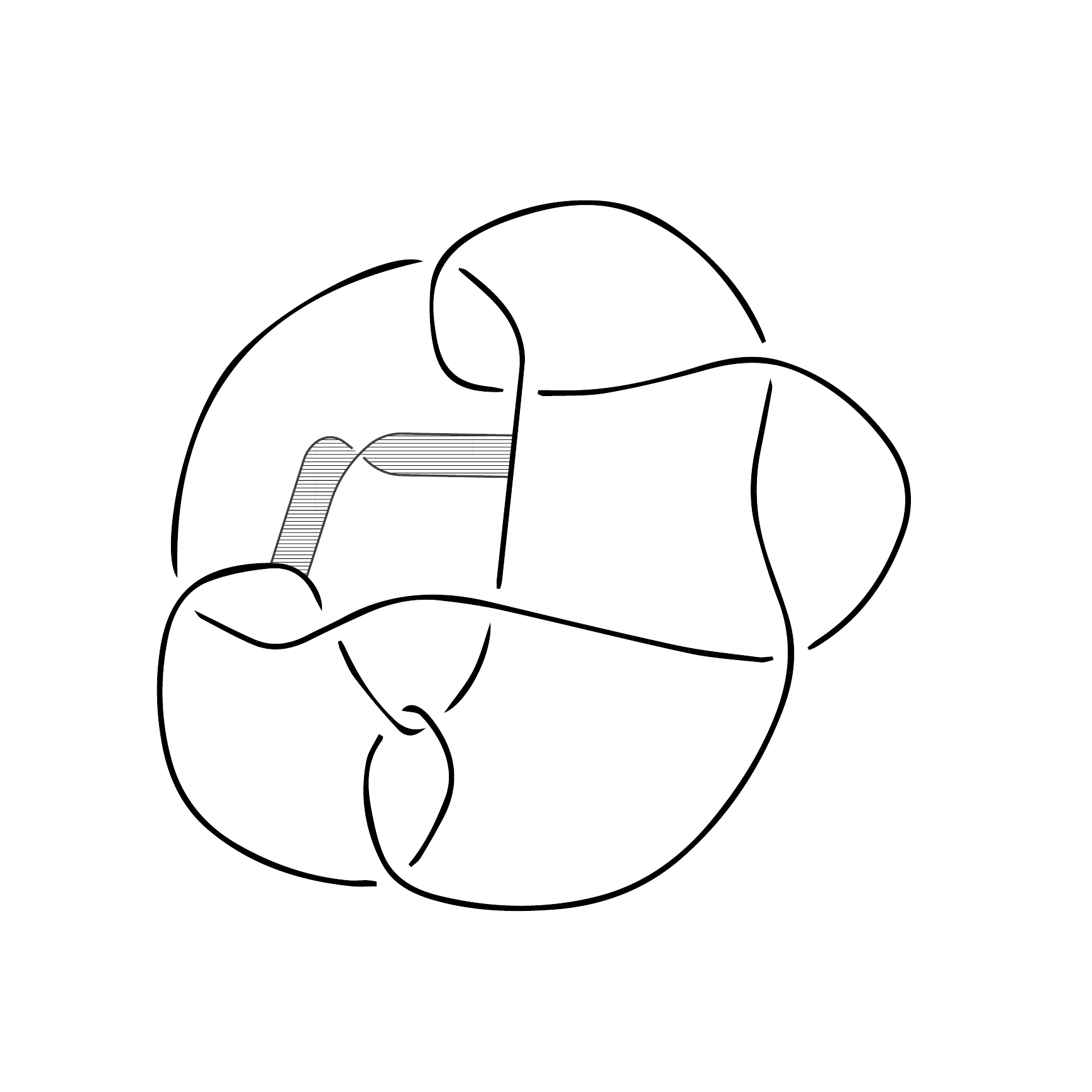}
		\caption{$10_{133}\stackrel{-1}{\longrightarrow} 8_{8}$}
		\label{FigureFor10-133}
	\end{subfigure}
	~
	\begin{subfigure}[b]{0.27\textwidth}
		\includegraphics[width=\textwidth]{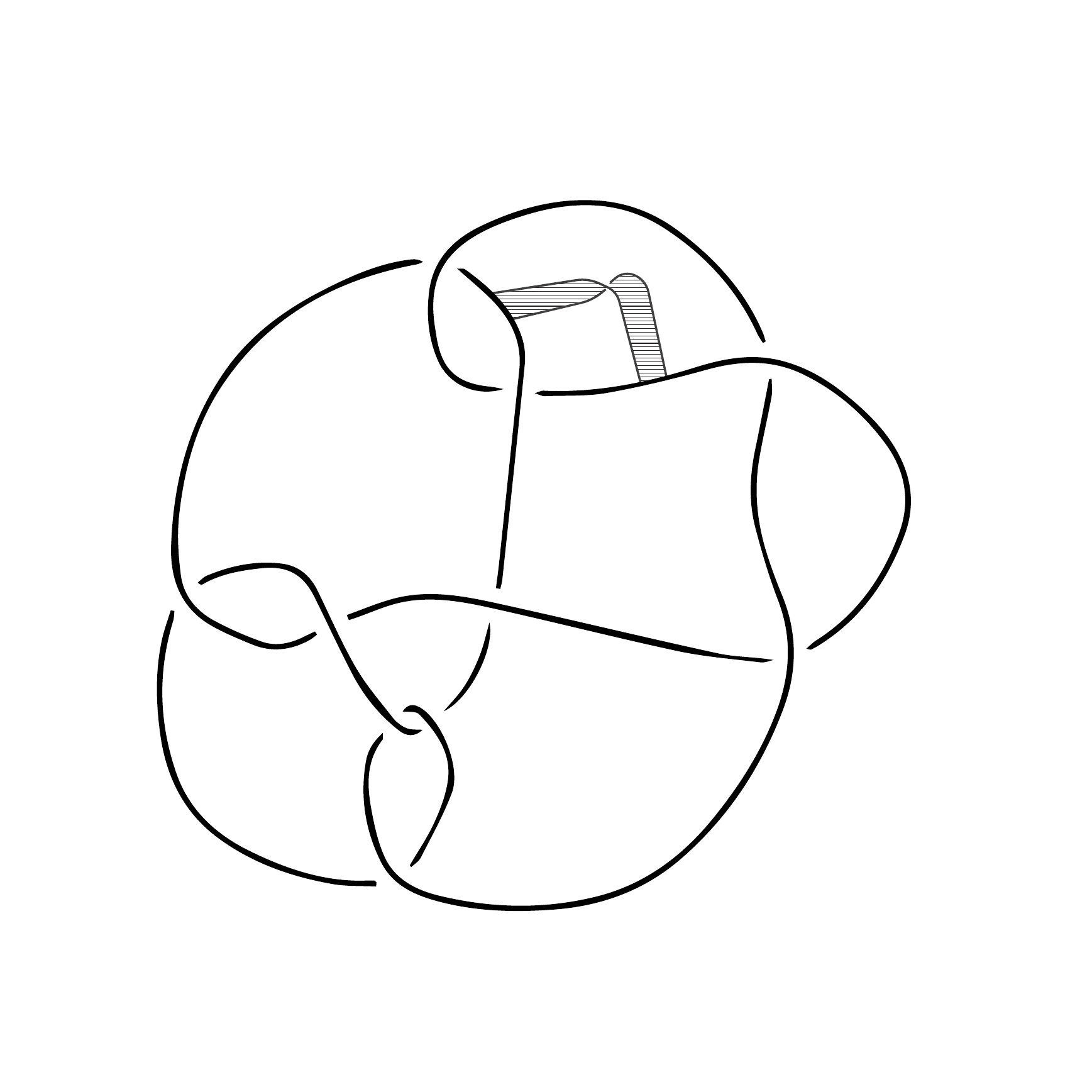}
		\caption{$10_{134}\stackrel{-1\phantom{i}}{\longrightarrow} 8_{20}$}
		\label{FigureFor10-134}
	\end{subfigure}
	~
	\vskip3mm
	~
	\begin{subfigure}[b]{0.3\textwidth}
		\includegraphics[width=\textwidth]{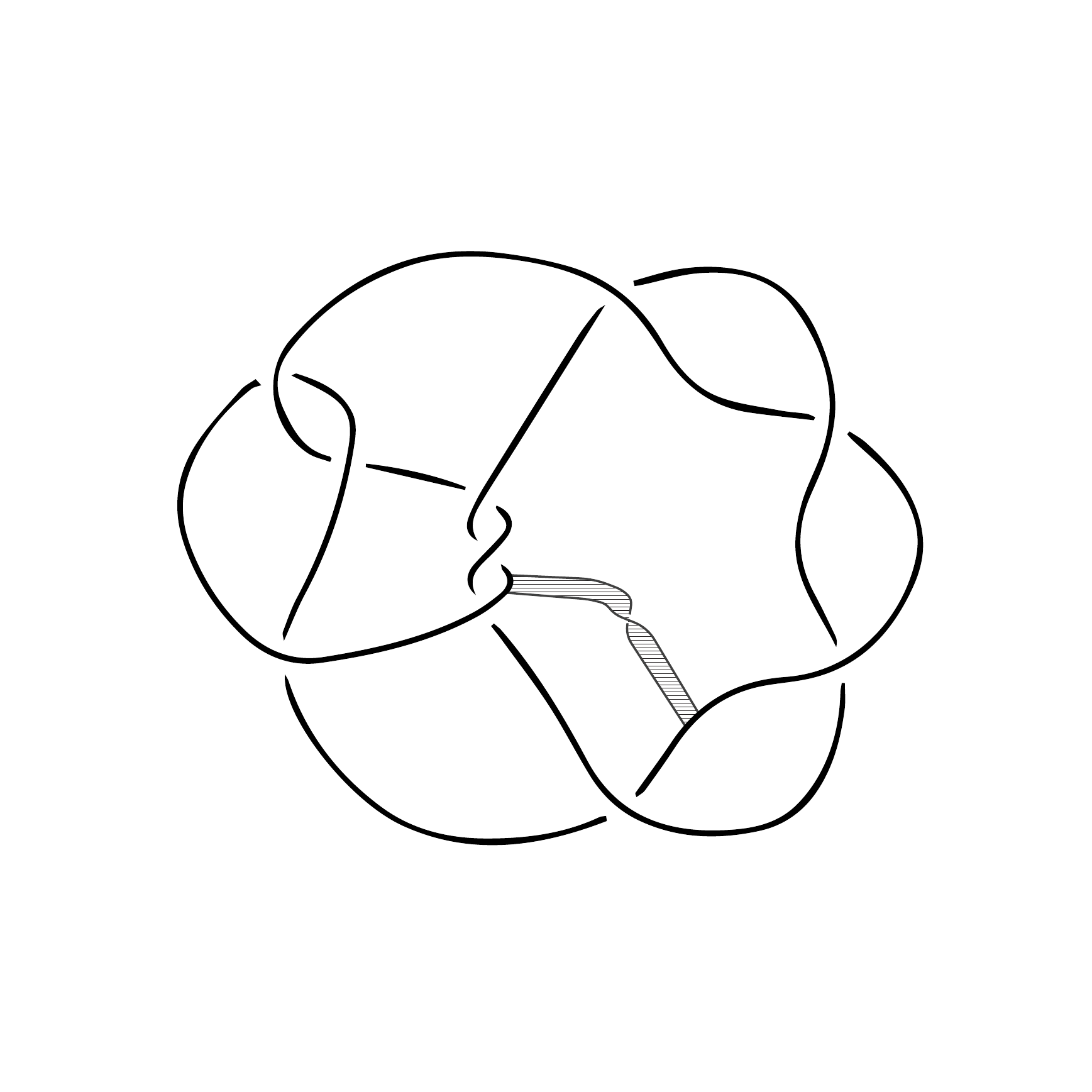}
		\caption{$10_{139}\stackrel{-1}{\longrightarrow}0_1$}
		\label{FigureFor10-139}
	\end{subfigure}
	~
	\begin{subfigure}[b]{0.3\textwidth}
		\includegraphics[width=\textwidth]{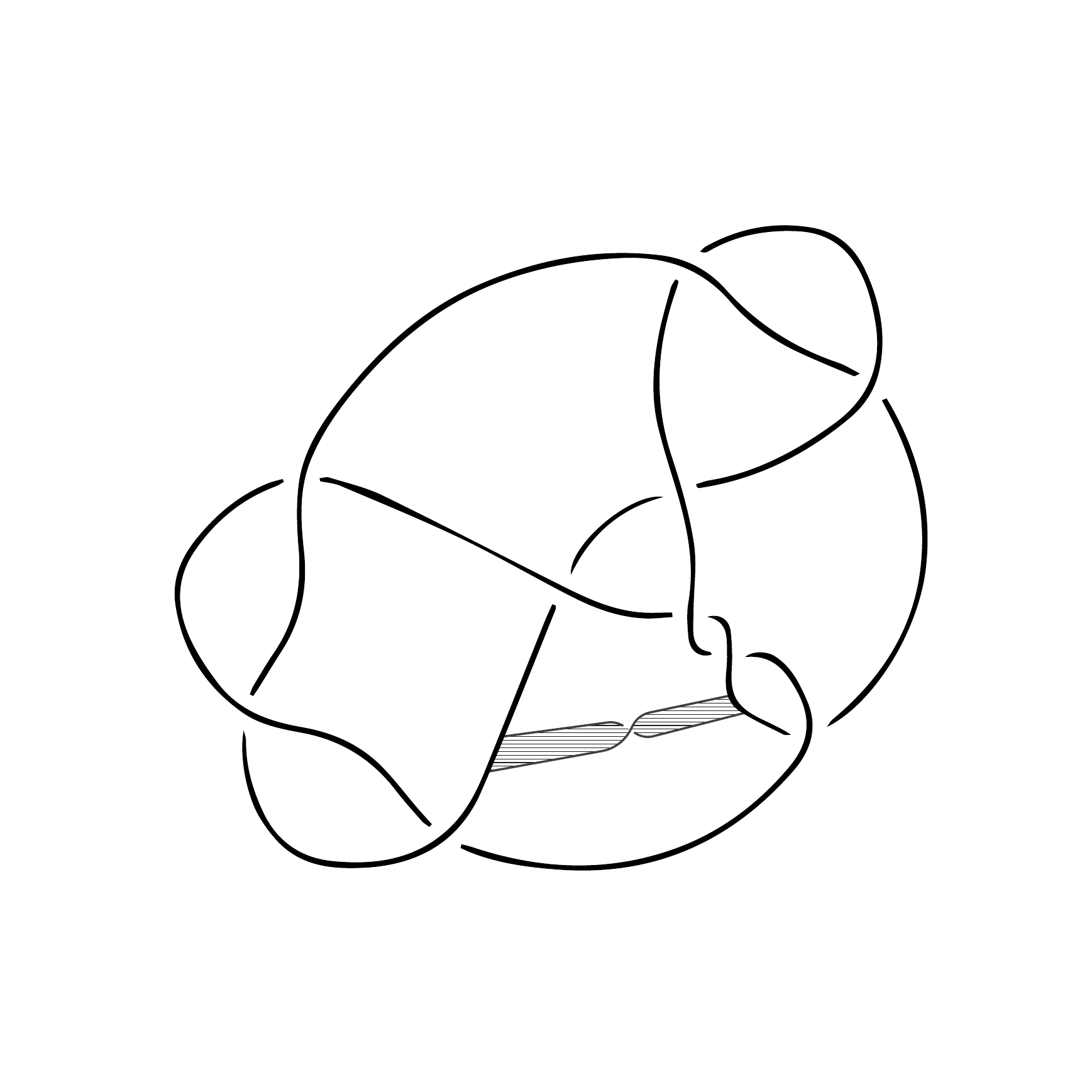}
		\caption{$10_{142}\stackrel{-1}{\longrightarrow} 0_1$}
		\label{FigureFor10-142}
	\end{subfigure}
	~
\begin{subfigure}[b]{0.27\textwidth}
		\includegraphics[width=\textwidth]{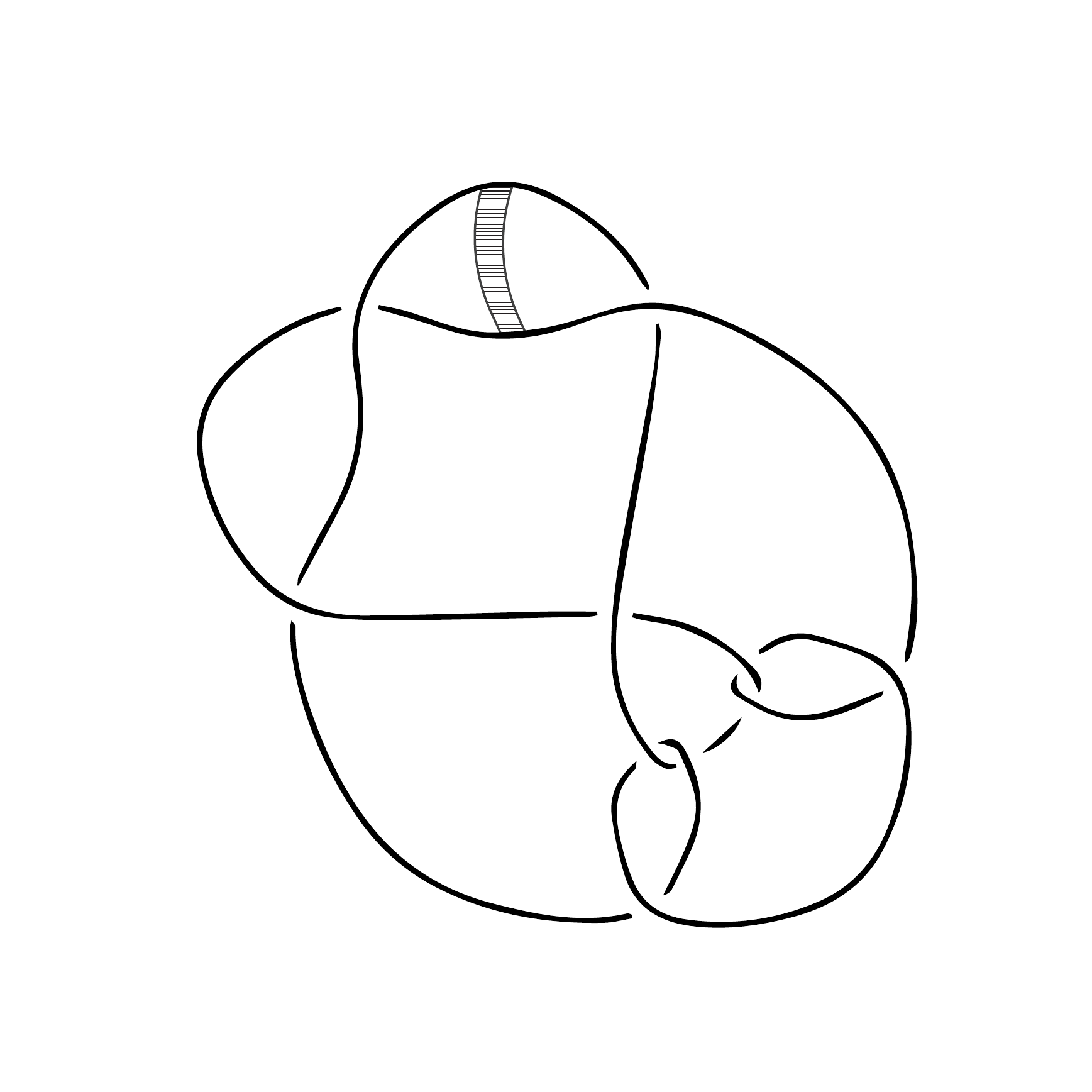}
		\caption{$10_{143}\stackrel{0}{\longrightarrow} 6_{1}$}
		\label{FigureFor10-143}
	\end{subfigure}
	~
	\vskip3mm
	~
		\begin{subfigure}[b]{0.27\textwidth}
		\includegraphics[width=\textwidth]{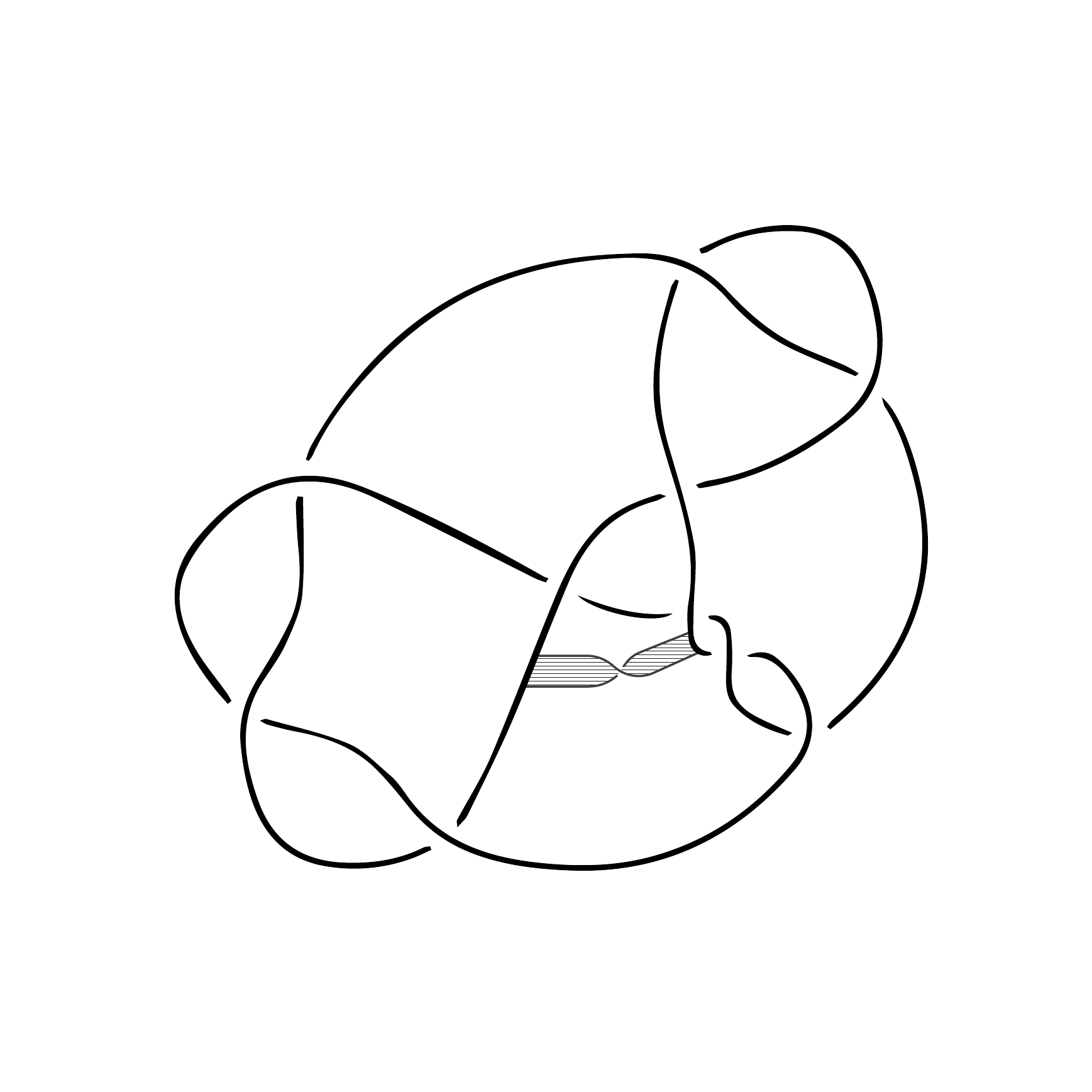}
		\caption{$10_{144}\stackrel{1}{\longrightarrow} 11n_{83}$}
		\label{FigureFor10-144}
	\end{subfigure}
	~    
			\begin{subfigure}[b]{0.3\textwidth}
		\includegraphics[width=\textwidth]{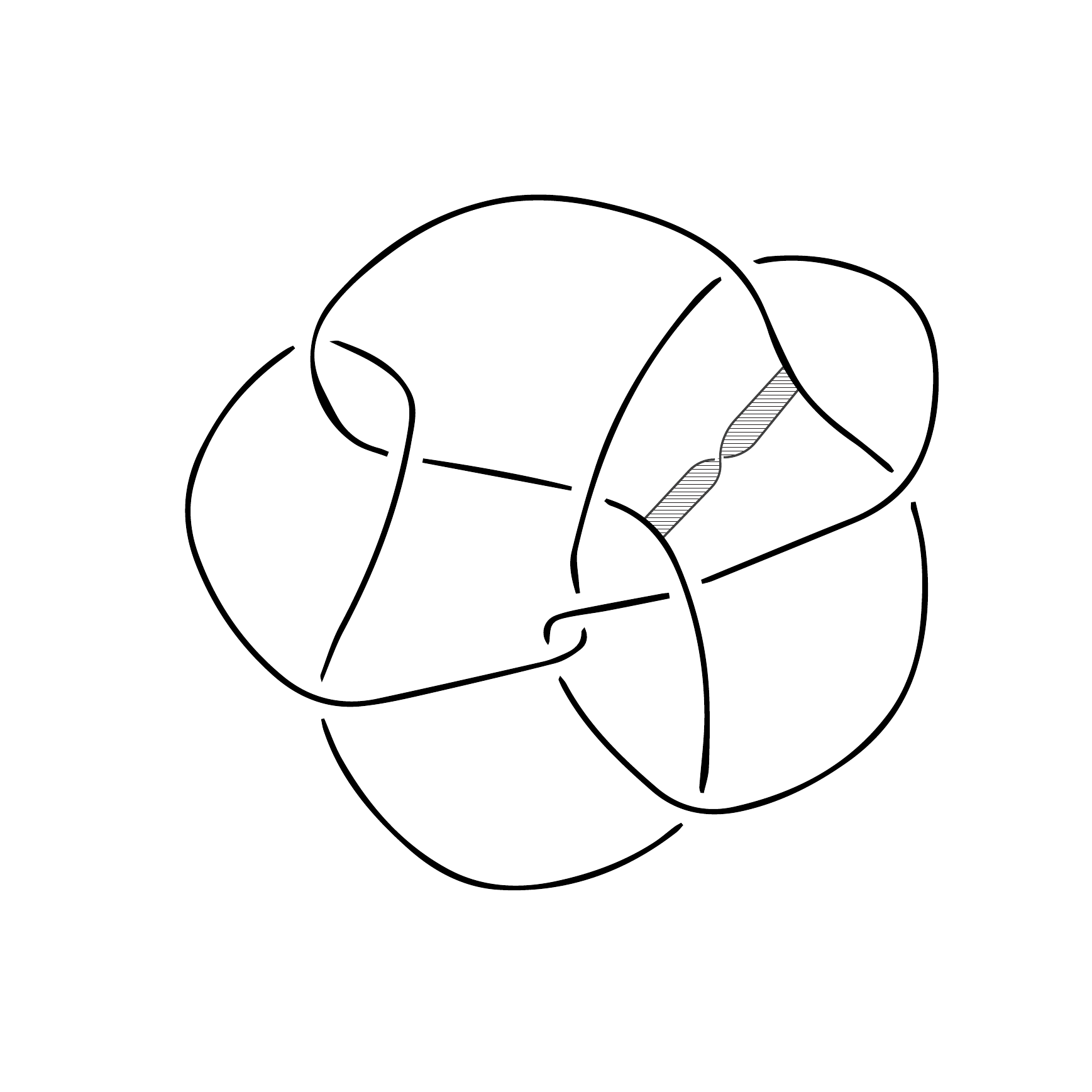}
		\caption{$10_{145}\stackrel{-1}{\longrightarrow} 0_{1}$}
		\label{FigureFor10-145}
	\end{subfigure}
	~
	\begin{subfigure}[b]{0.27\textwidth}
		\includegraphics[width=\textwidth]{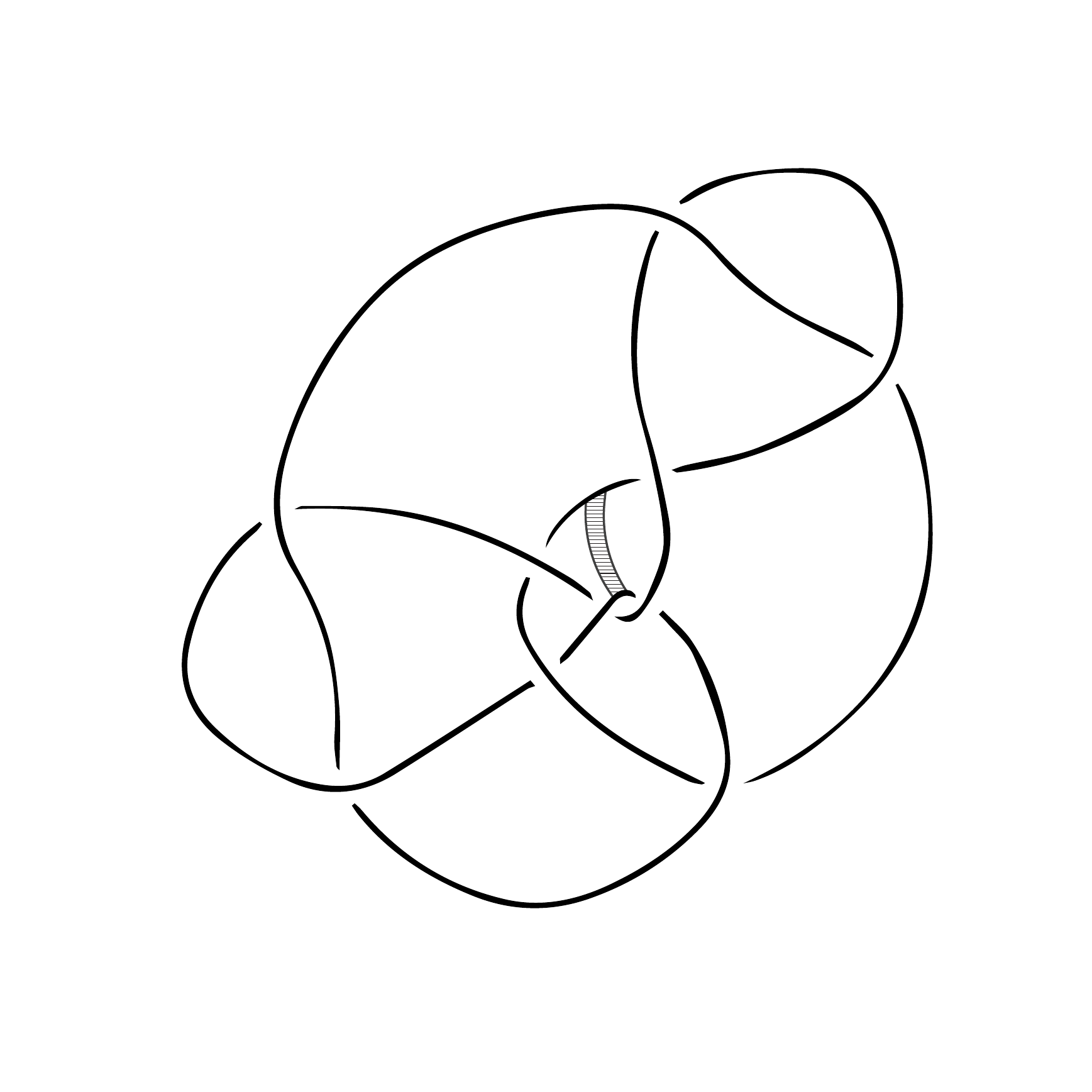}
		\caption{$10_{146}\stackrel{0}{\longrightarrow} 0_{1}$}
		\label{FigureFor10-146}
	\end{subfigure}
	~
			\vskip3mm

	\caption{Non-oriented band moves from the knots $10_{131}$, $10_{133}$, $10_{134}$, $10_{139}$, $10_{142}$, $10_{143}$, $10_{144}$, $10_{145}$, $10_{146}$ to slice knots}\label{slice9}
\end{figure}
\begin{figure}[h]
	\centering
	\begin{subfigure}[b]{0.27\textwidth}
		\includegraphics[width=\textwidth]{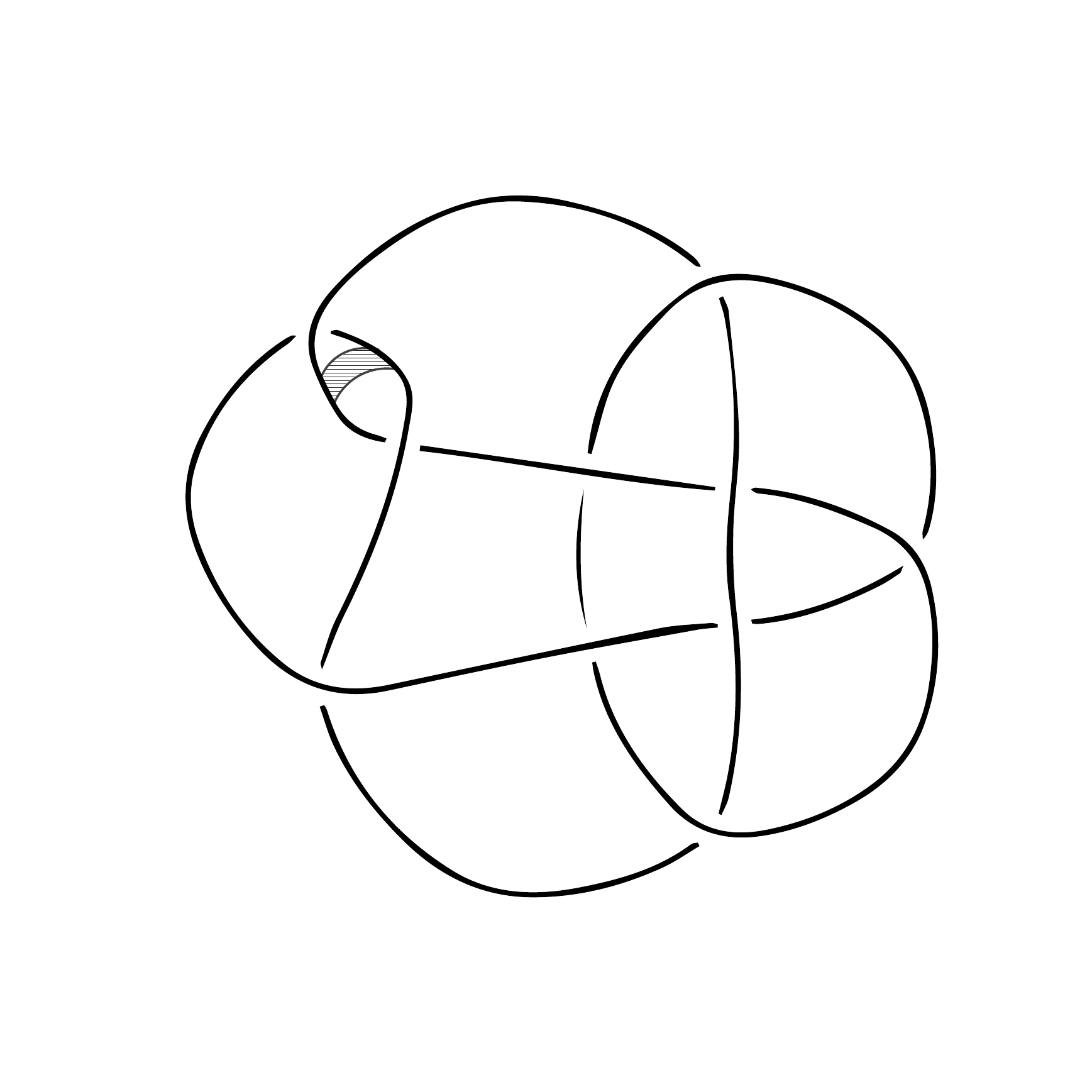}
		\caption{$10_{147}\stackrel{0}{\longrightarrow} 8_{20}$}
		\label{FigureFor10-147}
	\end{subfigure}
		~
	\begin{subfigure}[b]{0.27\textwidth}
		\includegraphics[width=\textwidth]{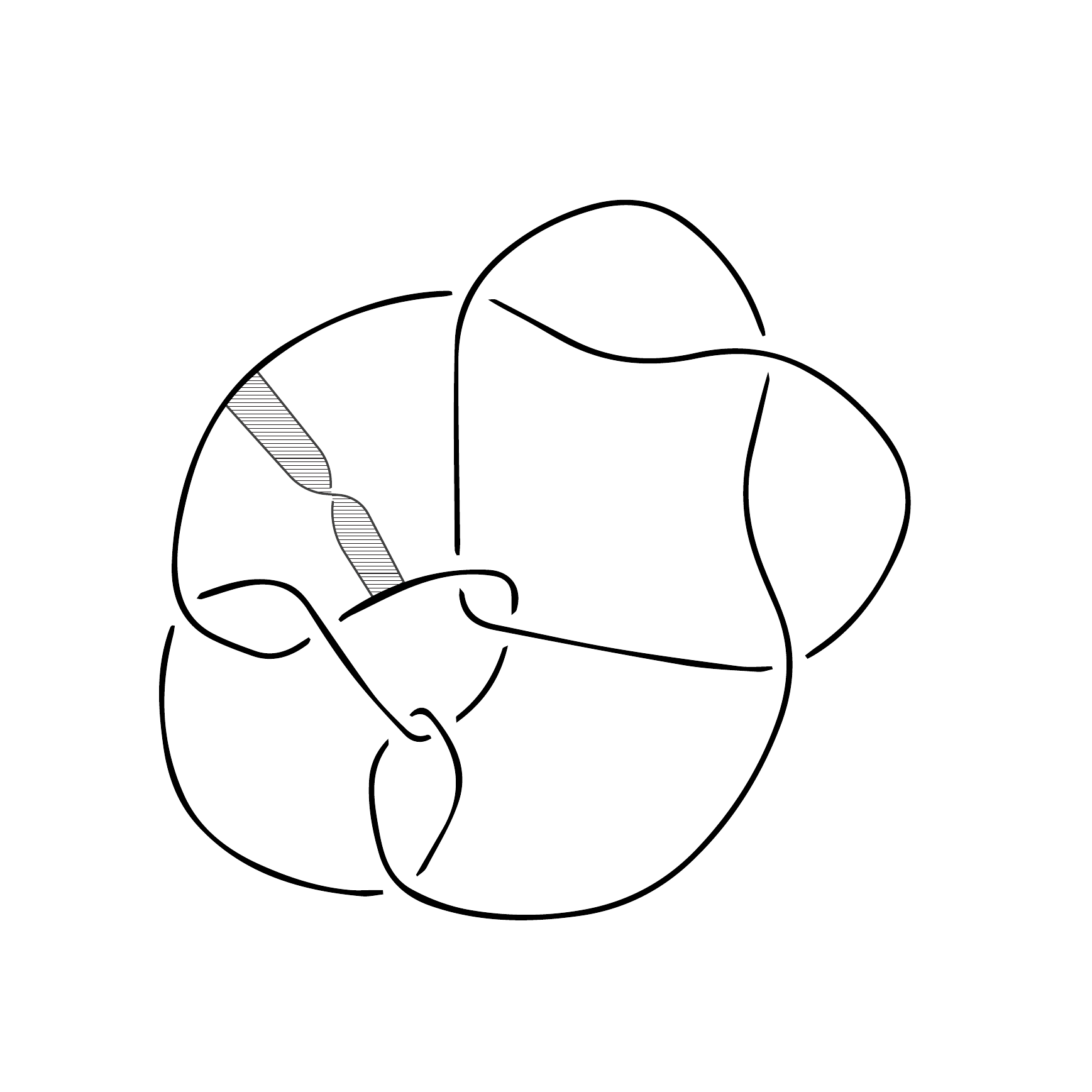}
		\caption{$10_{148}\stackrel{-1}{\longrightarrow} 10_{153}$}
		\label{FigureFor10-148}
	\end{subfigure}
	~
	\begin{subfigure}[b]{0.27\textwidth}
		\includegraphics[width=\textwidth]{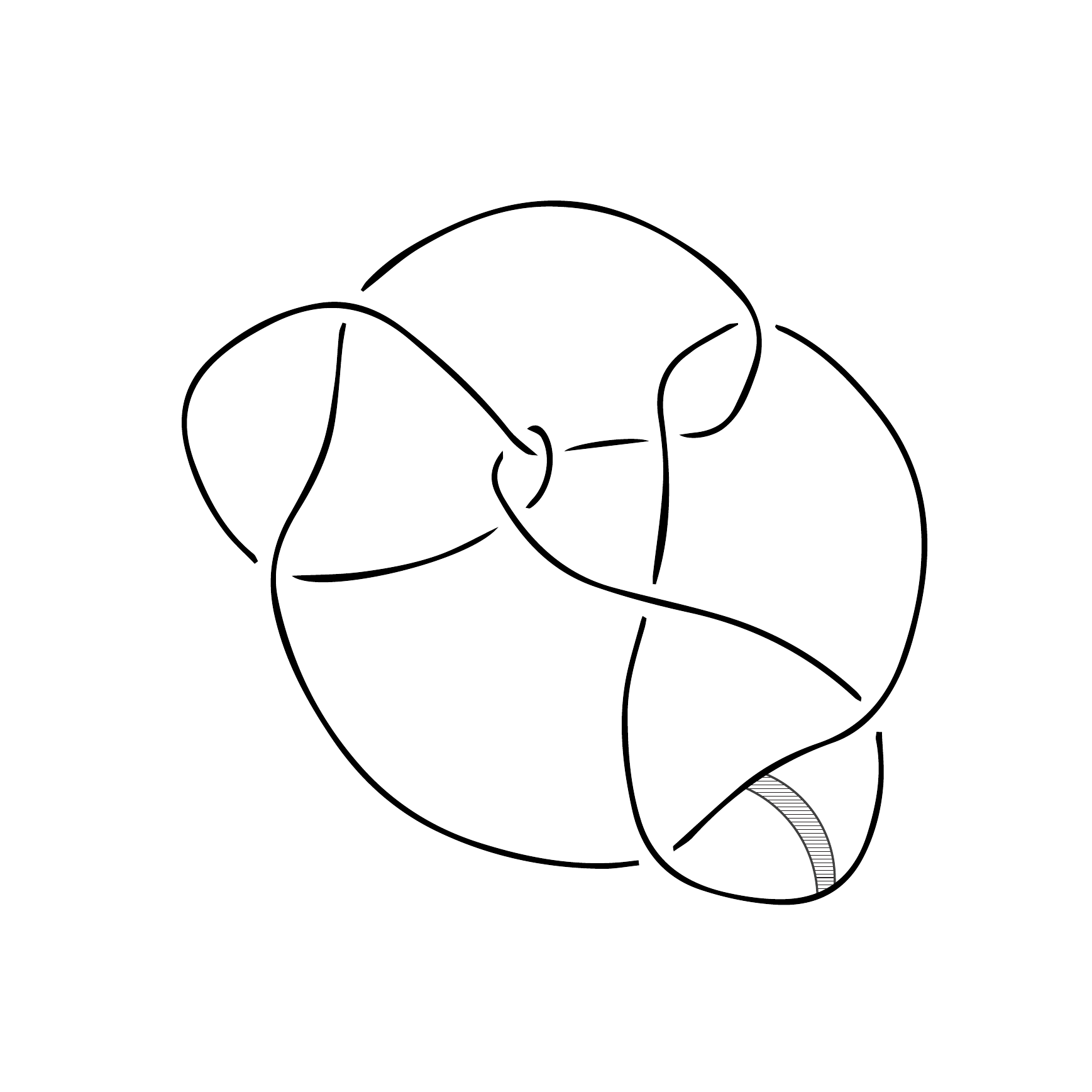}
		\caption{$10_{150}\stackrel{0\phantom{i}}{\longrightarrow} 8_{20}$}
		\label{FigureFor10-150}
	\end{subfigure}
	~
	\vskip3mm
	~
	\begin{subfigure}[b]{0.27\textwidth}
		\includegraphics[width=\textwidth]{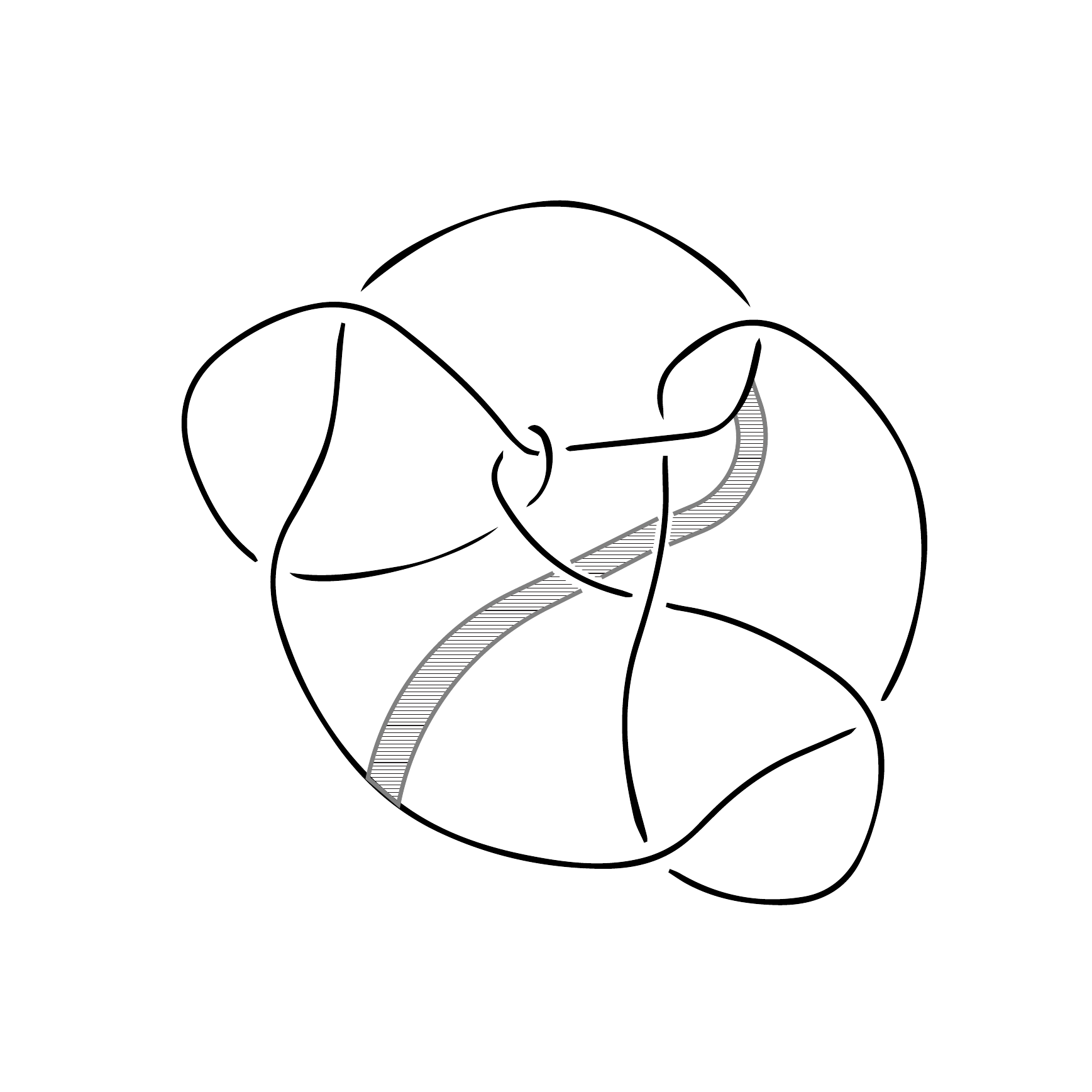}
		\caption{$10_{151}\stackrel{0\phantom{i}}{\longrightarrow} 10_{153}$}
		\label{FigureFor10-151}
	\end{subfigure}
		~
		\begin{subfigure}[b]{0.3\textwidth}
		\includegraphics[width=\textwidth]{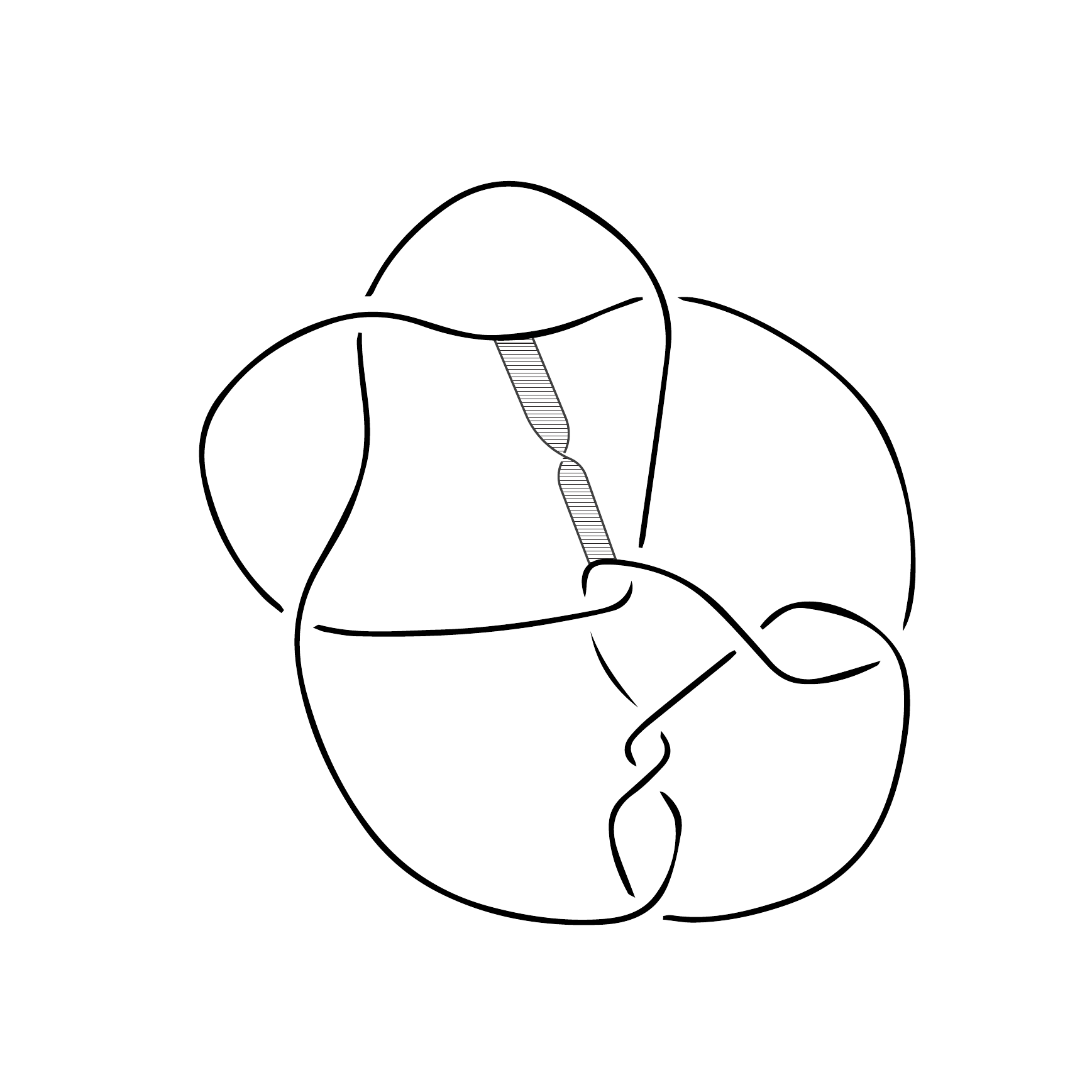}
		\caption{$10_{152}\stackrel{-1}{\longrightarrow}0_1$}
		\label{FigureFor10-152}
	\end{subfigure}
	~
	\begin{subfigure}[b]{0.3\textwidth}
		\includegraphics[width=\textwidth]{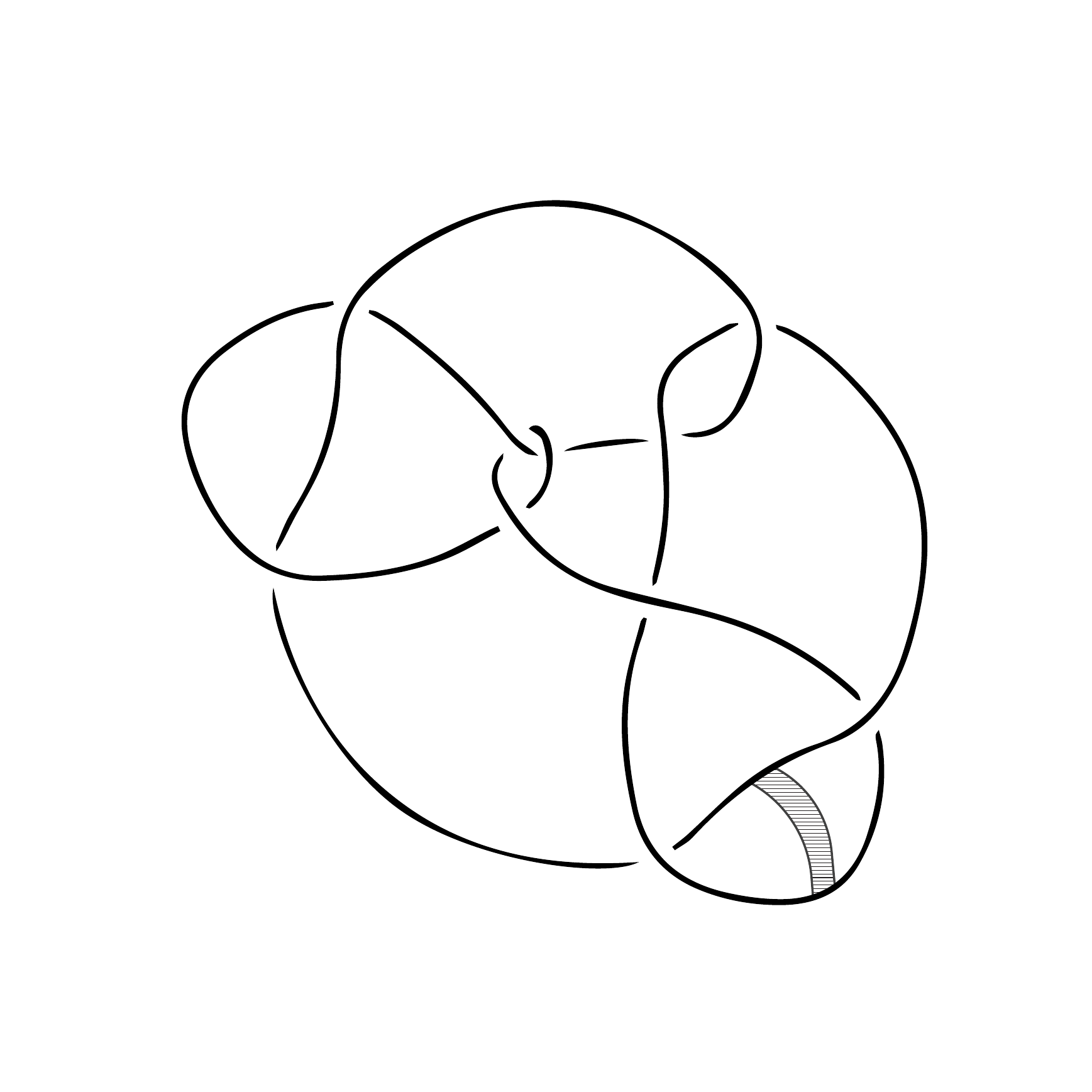}
		\caption{$10_{154}\stackrel{0}{\longrightarrow} 8_{20}$}
		\label{FigureFor10-154}
	\end{subfigure}
	~
	\vskip3mm
	~
	\begin{subfigure}[b]{0.27\textwidth}
		\includegraphics[width=\textwidth]{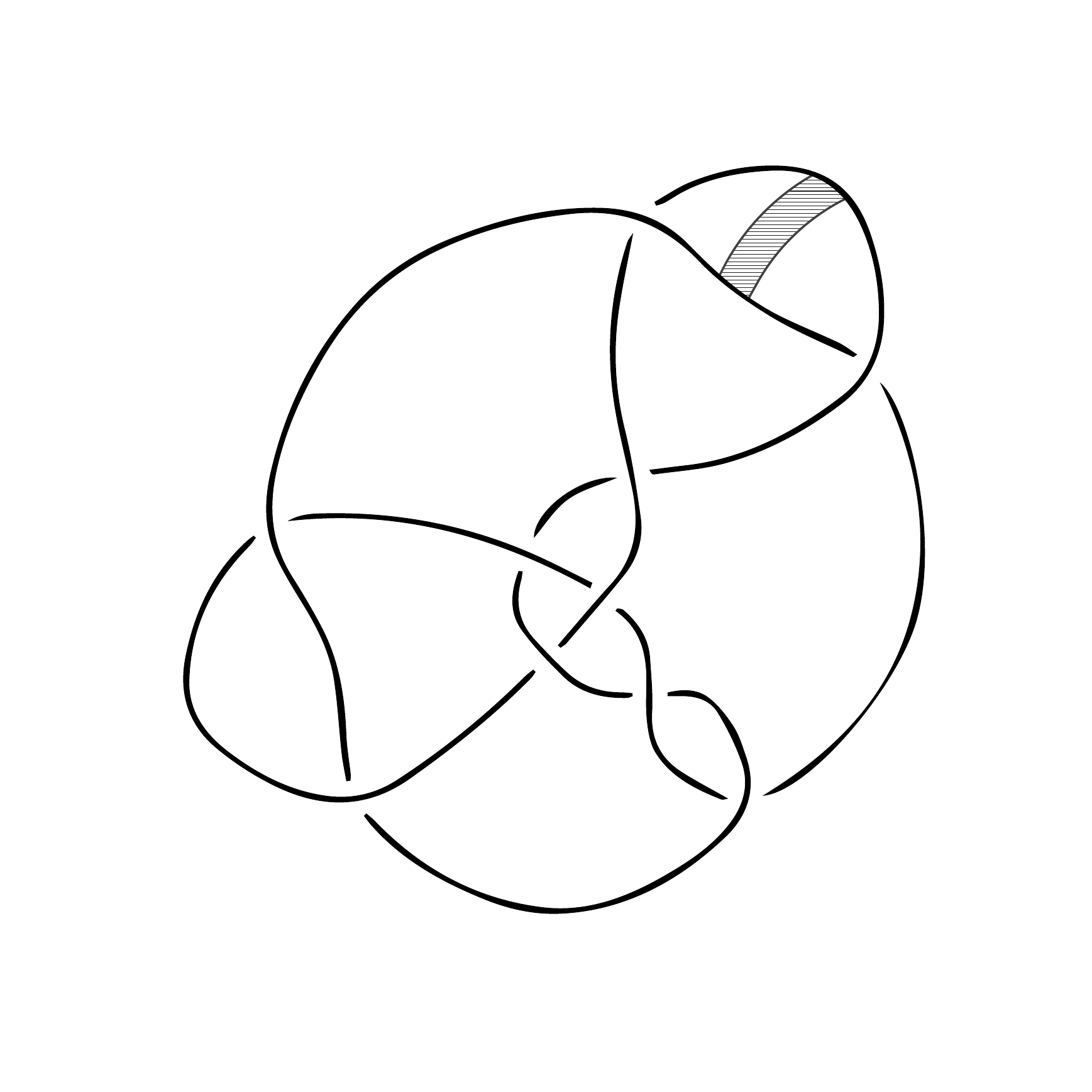}
		\caption{$10_{160}\stackrel{0\phantom{i}}{\longrightarrow} 0_{1}$}
		\label{FigureFor10-160}
	\end{subfigure}
	~
	\begin{subfigure}[b]{0.3\textwidth}
		\includegraphics[width=\textwidth]{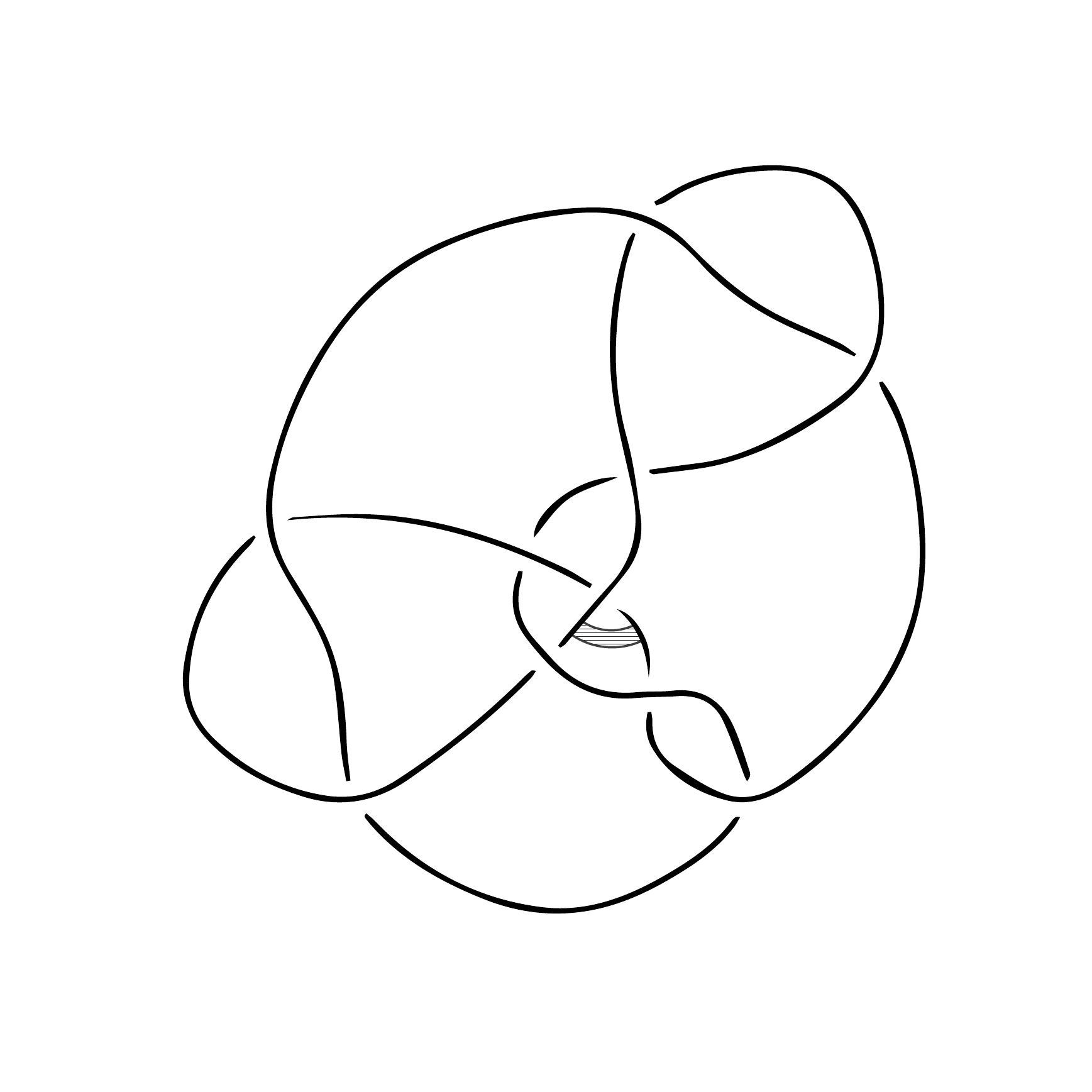}
		\caption{$10_{161}\stackrel{0}{\longrightarrow}0_1$}
		\label{FigureFor10-161}
	\end{subfigure}
~
	\begin{subfigure}[b]{0.3\textwidth}
		\includegraphics[width=\textwidth]{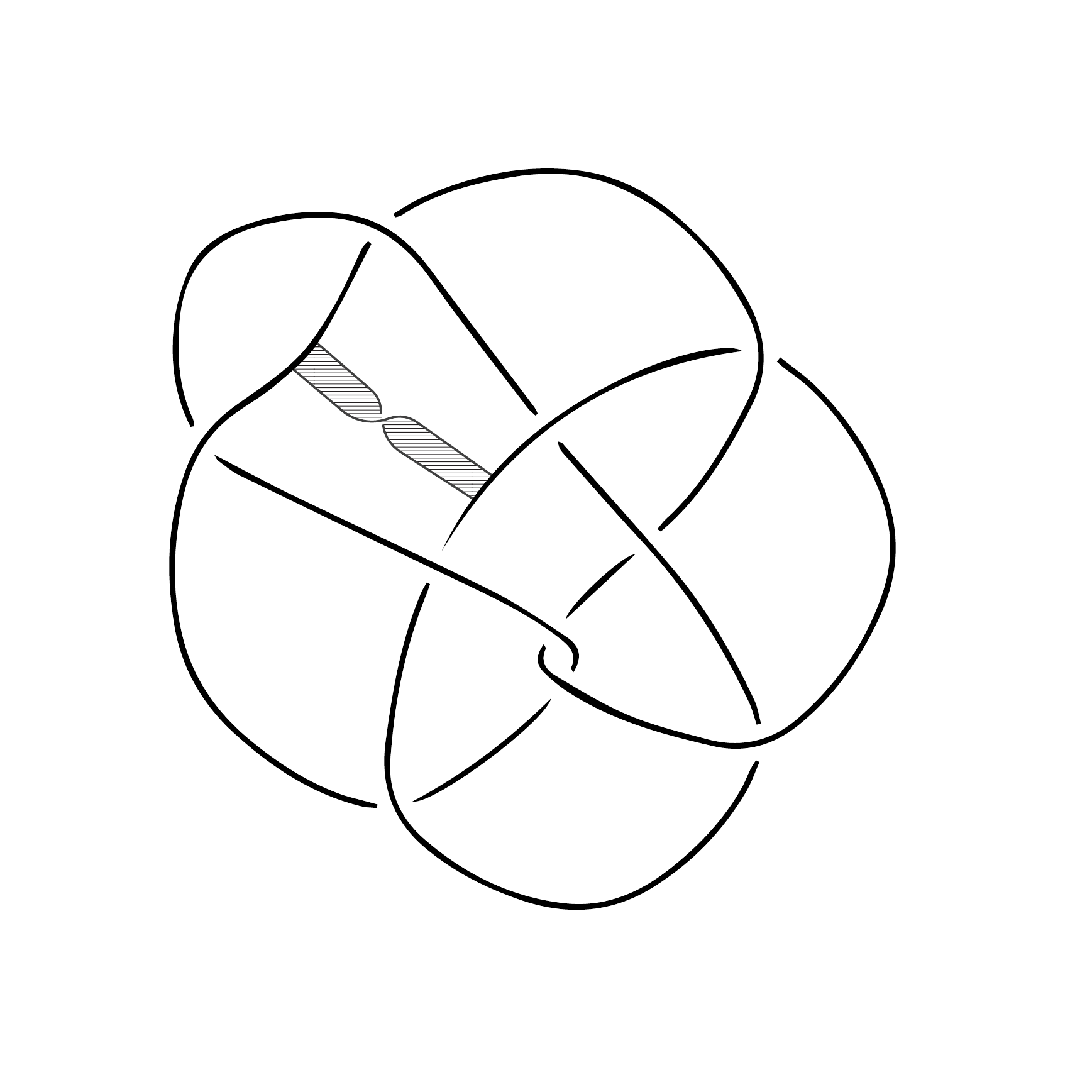}
		\caption{$10_{165}\stackrel{-1}{\longrightarrow} 0_{1}$}
		\label{FigureFor10-165}
	\end{subfigure}	
	~
\vskip3mm
	~

		\caption{Non-oriented band moves from the $10_{147}$, $10_{148}$, $10_{150}$, $10_{152}$, $10_{154}$, $10_{160}$, $10_{161}$, $10_{165}$ to slice knots}\label{slice10}
\end{figure}

 There are exactly 90 knots for which we were able to find a band move to a slice knot. According to Proposition \ref{band move}, $\gamma_4$ of these 90 knots are also equal to 1. These knots are addressed in Figures \ref{slice1}, \ref{slice2}, \ref{slice3}, \ref{slice4}, \ref{slice5}, \ref{slice6}, \ref{slice7}, \ref{slice8}, \ref{slice9} and \ref{slice10} at the end of this paper.
\subsection {Non-alternating knots}\label{NAK} 
\begin{figure}[h]
	\centering
		\begin{subfigure}[b]{0.3\textwidth}
		\includegraphics[width=\textwidth]{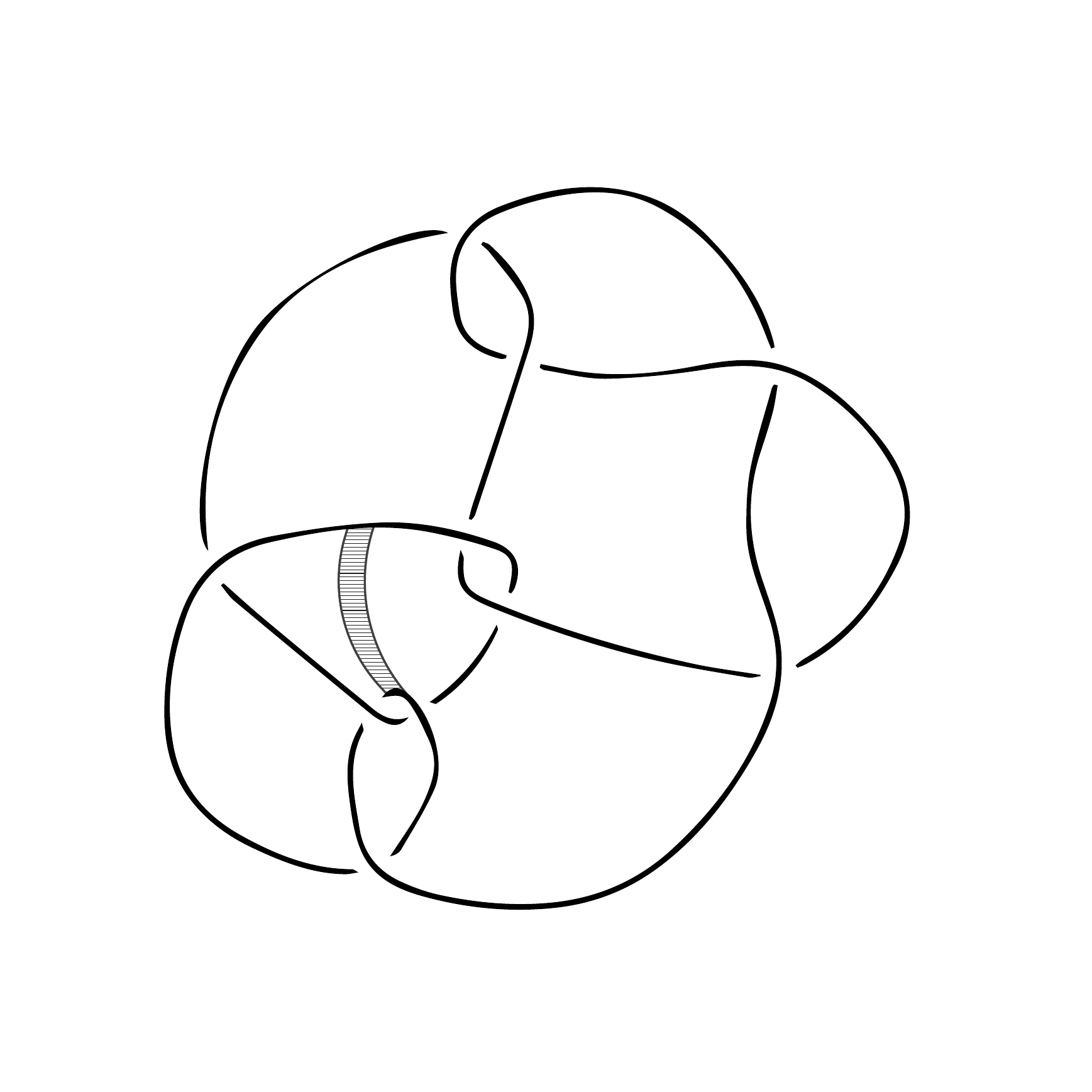}
		\caption{$10_{136}\stackrel{0}{\longrightarrow} 9_{45}$}
		\label{FigureFor10-136}
	\end{subfigure}
	~
	\begin{subfigure}[b]{0.3\textwidth}
		\includegraphics[width=\textwidth]{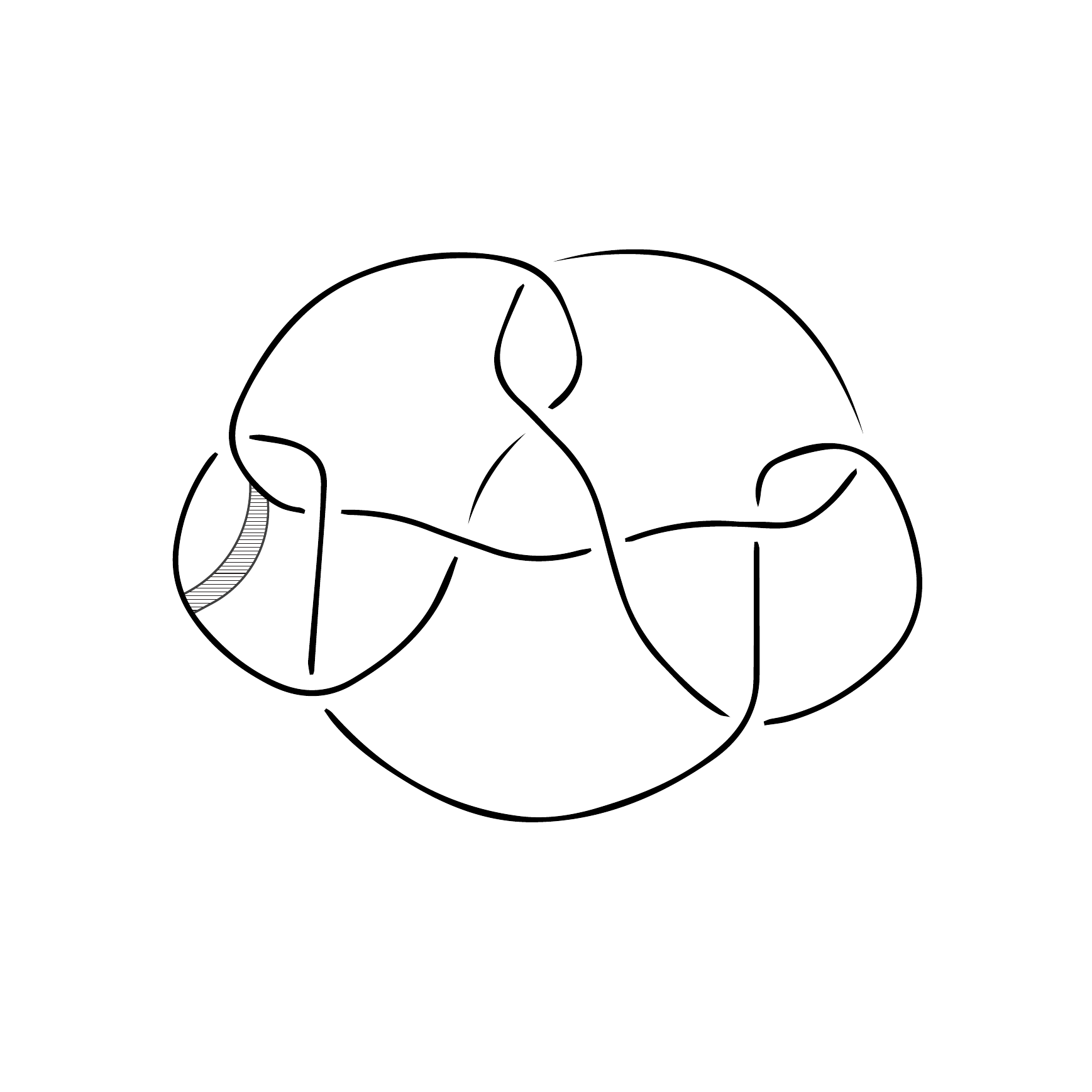}
		\caption{$10_{138}\stackrel{0}{\longrightarrow} 9_{45}$}
		\label{FigureFor10-138}
	\end{subfigure}	
	~
	\begin{subfigure}[b]{0.3\textwidth}
		\includegraphics[width=\textwidth]{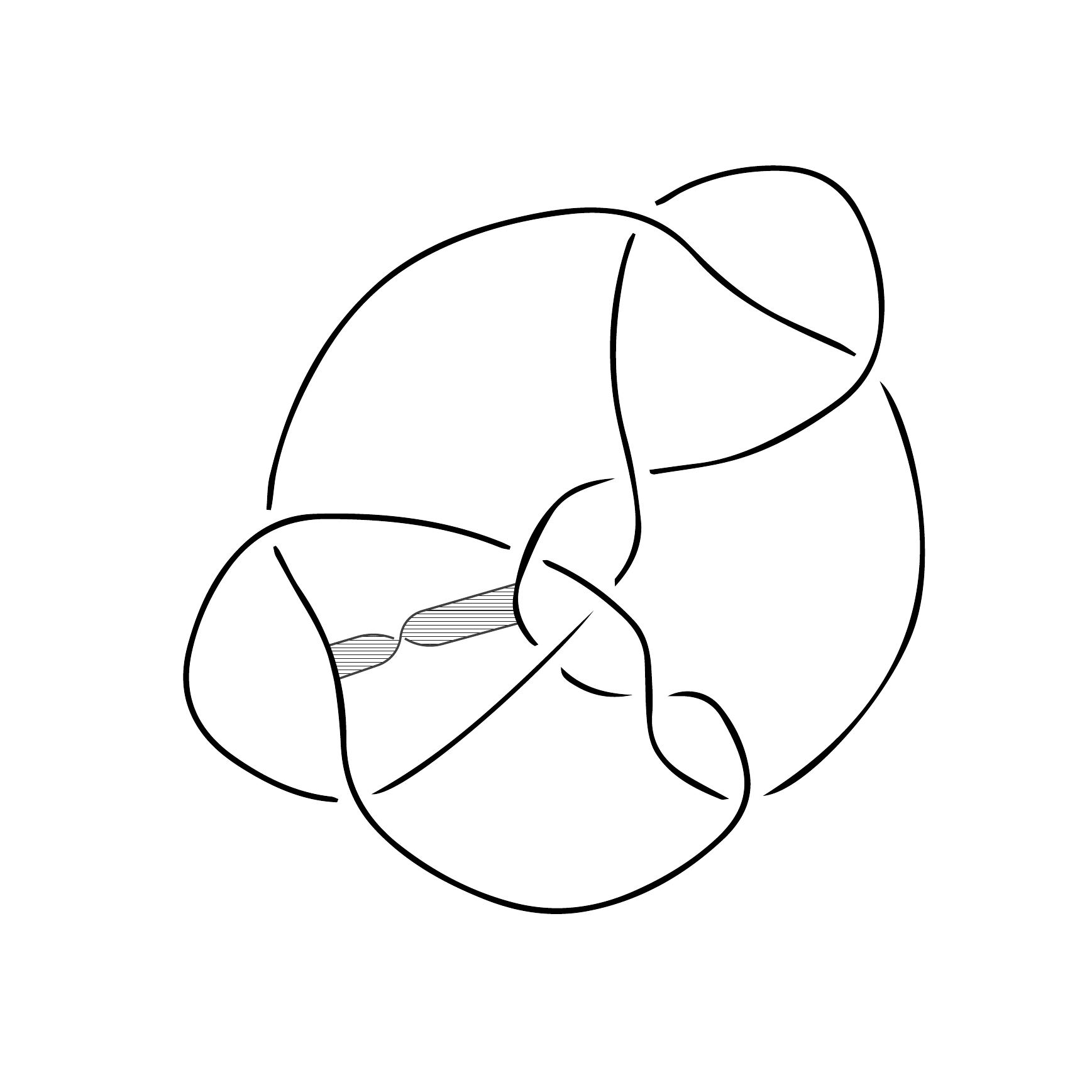}
		\caption{$10_{156}\stackrel{-1}{\longrightarrow} 5_{2}$}
		\label{FigureFor10-156}
	\end{subfigure}
	~
	\vskip3mm	
		 	    \begin{subfigure}[b]{0.3\textwidth}
		\includegraphics[width=\textwidth]{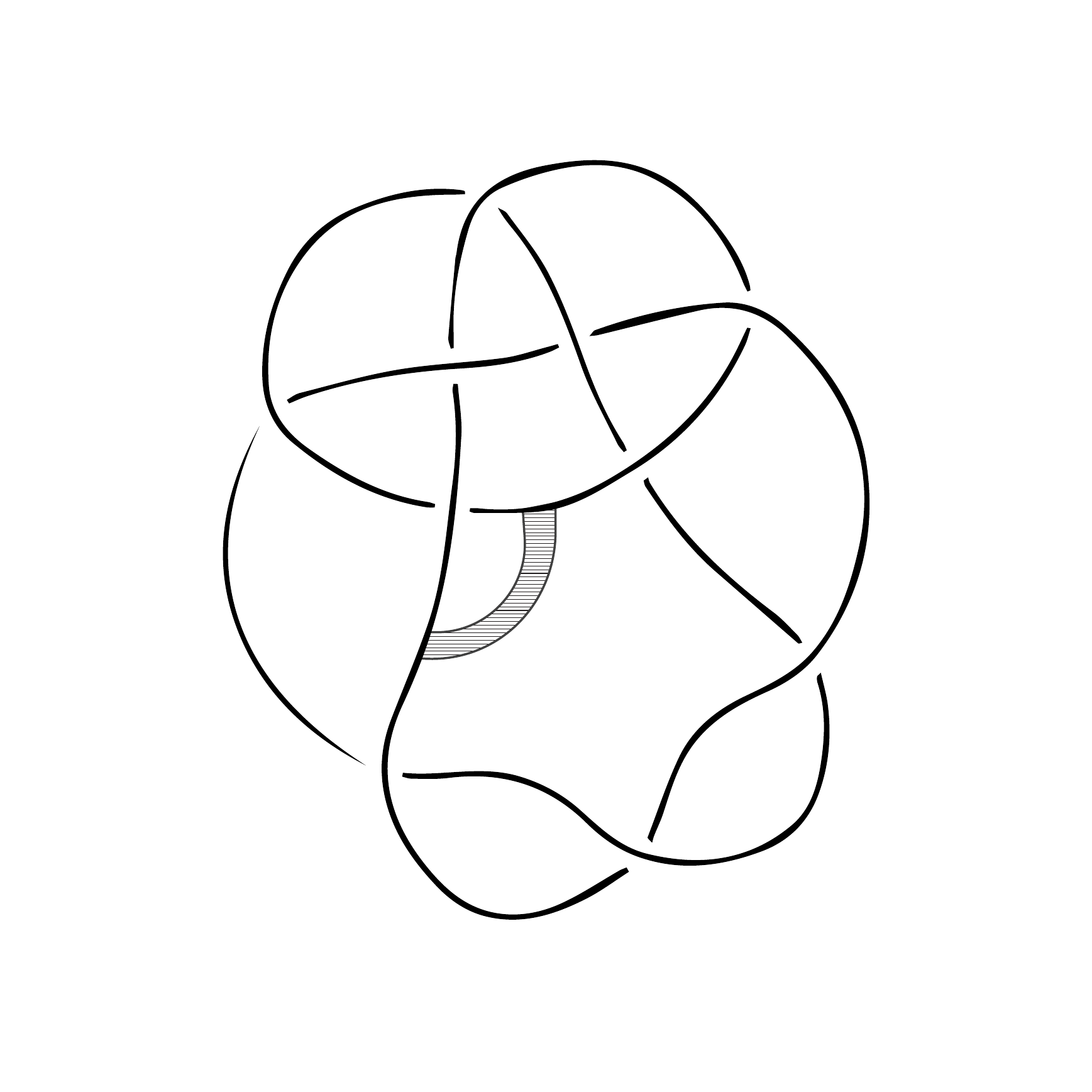}
		\caption{$10_{159}\stackrel{0}{\longrightarrow} 8_{14}$}
		\label{FigureFor10-159}
	\end{subfigure}
	~
	\begin{subfigure}[b]{0.3\textwidth}
		\includegraphics[width=\textwidth]{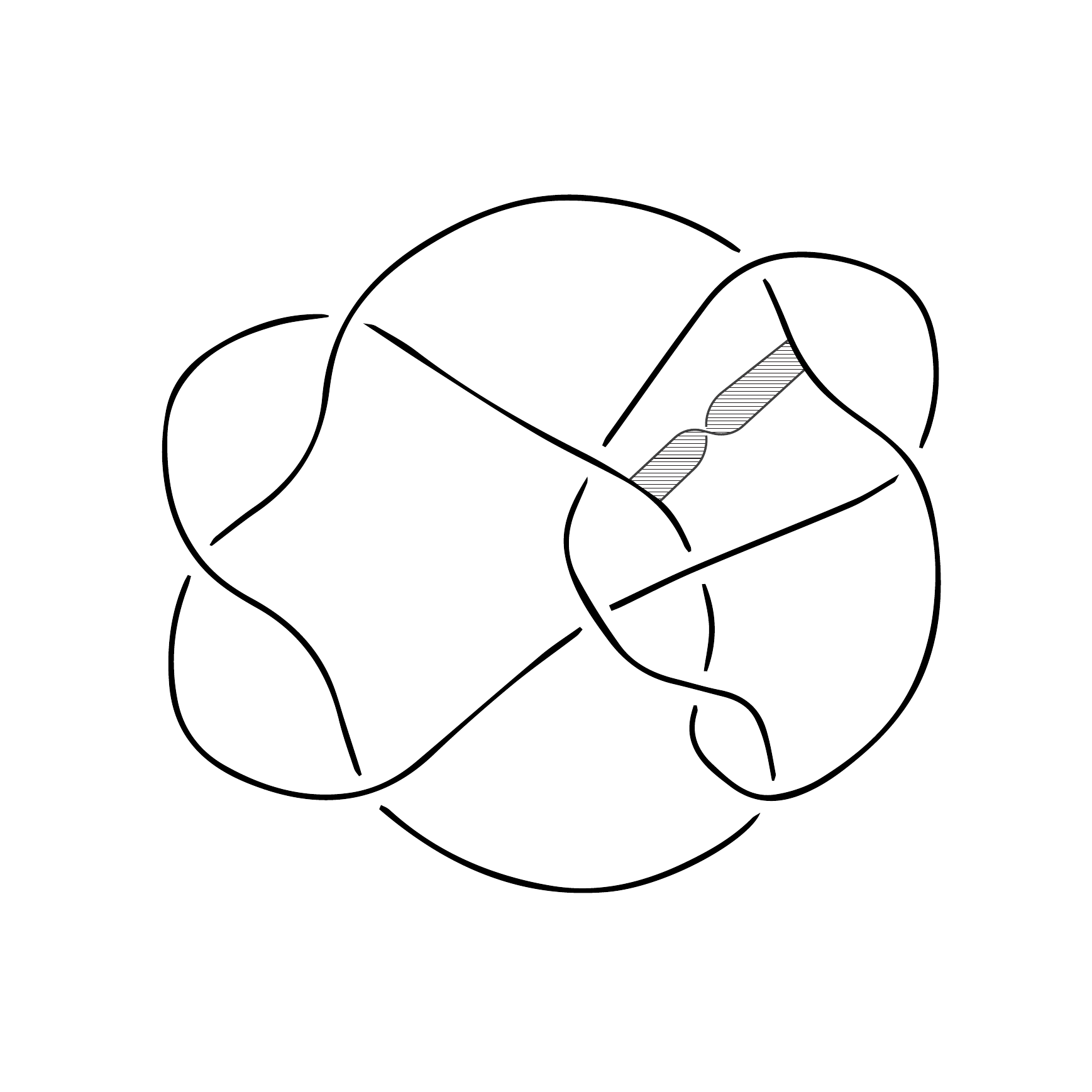}
		\caption{$10_{162}\stackrel{1}{\longrightarrow} 5_{2}$}
		\label{FigureFor10-162}
	\end{subfigure}
		~
	\begin{subfigure}[b]{0.3\textwidth}
		\includegraphics[width=\textwidth]{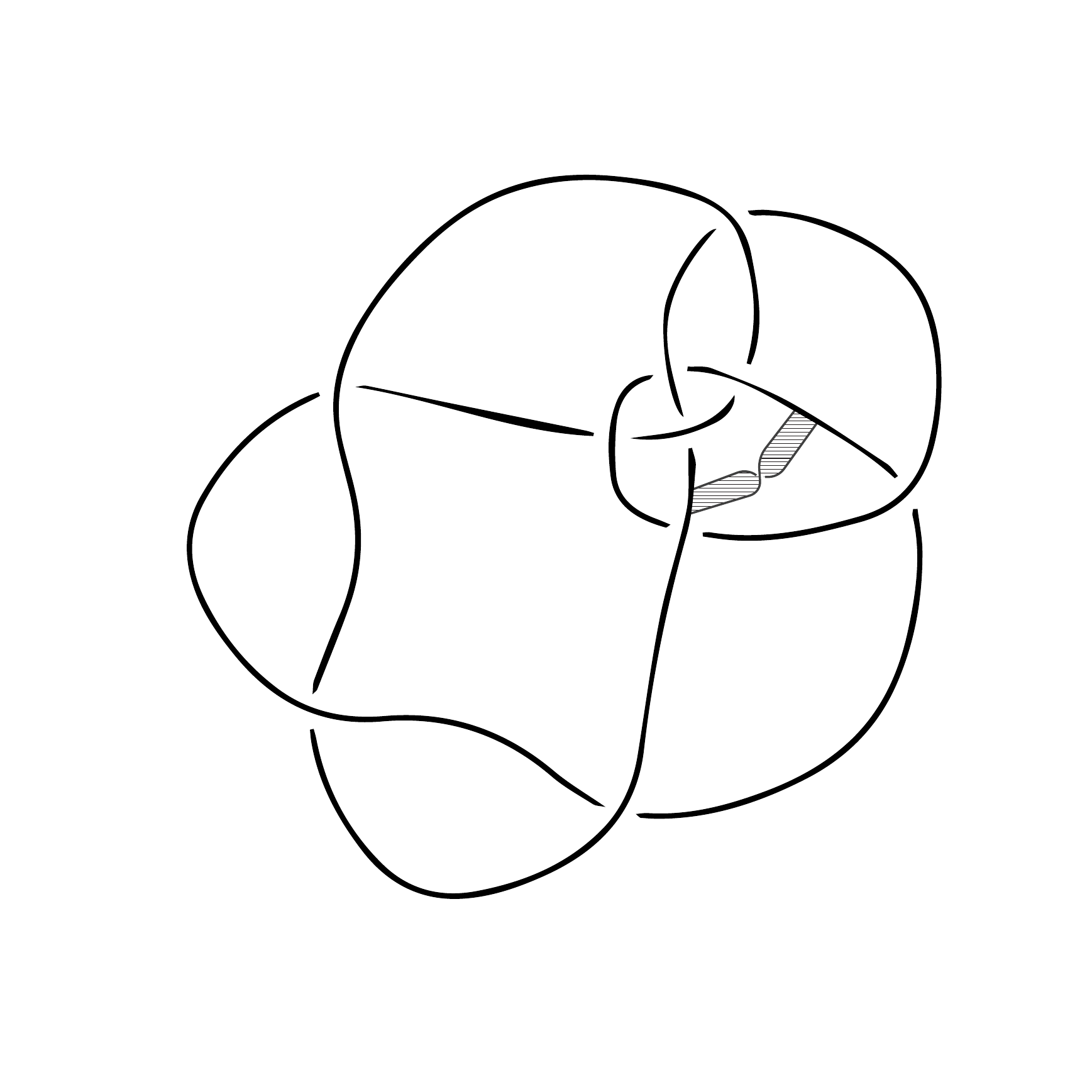}
		\caption{$10_{163}\stackrel{-1}{\longrightarrow} 9_{44}$}
		\label{FigureFor10-163}
	\end{subfigure}
	\vskip3mm
	\caption{Non-oriented band moves from the knots $10_{136}$, $10_{138}$ ,$10_{156}$, $10_{159}$, $10_{162}$, $10_{163}$ to  knots with $\gamma_4=1$.}\label{6knots}
\end{figure}

According to Proposition \ref{band move} if we find a band move for any knot $K$ to a slice knot, we can immediately conclude that $\gamma_4(K)=1$. Among all 165 knots with crossing number 10, we have exactly 42 non-alternating knots. We have shown in Section \ref{sliceknots} that these 29 knots
 $$10_{124}, \, 10_{125}, \, 10_{126}, \, 10_{127}, \,10_{128}, \, 10_{129}, \, 10_{130}, \, 10_{131}$$
$$ 10_{133}, \, 10_{134}, \, 10_{137}, \,10_{139}, \,10_{140}, \,10_{142}, \,10_{143}, \, 10_{144}, \, 10_{145}$$
$$10_{146}, \, 10_{147}, \, 10_{148}, \, 10_{150}, \, 10_{151}, \,10_{152}, \, 10_{153}, \, 10_{154}, \ 10_{155}, \,10_{160}, \, 10_{161}, 10_{165} $$
have $\gamma_4$ equals to 1.
It is shown in Section \ref{first group} that these 7 knots
 $$10_{132}, \,10_{135}, \,10_{141}, \, 10_{149}, \,10_{157}, \, 10_{158}, \,10_{164}$$
 have $\gamma_4$ equals to $2$.
 \\
 In this section we will show that the remaining non-alternating knots 
  $$10_{136}, \,10_{138}, \,10_{156}, \, 10_{159}, \,10_{162}, \, 10_{163},$$
also have $\gamma_4=2$. Figure \ref{6knots} shows that these 6 knots have $\gamma_4\leq 2$. Now to find a lower bound for these 6 non-alternating knots, as we mentioned before, we can not use the Georitz Matrix. Therefore by using the following theorem from Gilmer and Livingston \cite{GilmerLivingston} we will show that for these 6 knots $\gamma_4\geq{2}$.
 \begin{theorem}{\label{lk}}(Corollary 3 in \cite{GilmerLivingston})
 Suppose that $H_1(M(K))=Z_n$ where $n$ is the product of primes, all with odd exponenet. Then if $K$ bounds a Mobius band in $B^4$, there is a generator $a\in H_1(M(K))$ such that $\ell k(a,a)=\pm\frac{1}{n}$.
 \end{theorem}
 Here $(H_1(M(K)),\ell k)$ denotes the linking form on $M(K)$,  the 2-fold cover of $S^3$ branched over $K$. A direct (computer aided) calculation shows that for the non-alternating knots 
$$K=10_{136}, \, 10_{138}, \, 10_{156}, \, 10_{159}, \, 10_{162}, \, 10_{163},$$
the first homology groups of their 2-fold branched covers are isomorphic to 
$$  \mathbb Z_{15}, \, \mathbb Z_{35}, \, \mathbb Z_{35}, \, \mathbb Z_{39}, \, \mathbb Z_{35}, \, \mathbb Z_{51}, \, $$
and their respective linking forms are given by multiplication by 
$$\frac{8}{15}, \,\frac{12}{35}, \,\frac{32}{35}, \,\frac{19}{39}, \, \frac{8}{35}, \, \frac{20}{51}.$$ 
By easy inspection, none of these pass the obstruction to $\gamma_4=1$ from Theorem \ref{lk} and thus $\gamma_4(K)\geq 2$ for all 6 of these knots. 
\subsection { Knots with $\sigma(K)+4\cdot \text{Arf}(K)\equiv{4} \pmod 8$ } \label{first group} There are exactly 43 knots which satisfy the congruance $\sigma(K)+4\cdot \text{Arf}(K)\equiv{4} \pmod 8$. According to Proposition \ref{lowerbound}, for each knot in this group $\gamma_4(K)\geq{2}$. Therefore, according to Proposition \ref{band move} it is enough to find a band move to a knot with $\gamma_4=1$ to prove that $\gamma_4(K)$ is equal to 2 for all 43 knots in this group. Such band moves are shown in Figures \ref{gamma4,1}, \ref{gamma4,2}, \ref{gamma4,3}, \ref{gamma4,4} and \ref{gamma4,5}.
\subsection { Knots with $\sigma(K)+4\cdot \text{Arf}(K)\equiv{-2} \pmod 8$ } \label {second group} There are exactly 42 knots  with crossing number 10 which satisfy the congruence $\sigma(K)+4\cdot \text{Arf}(K)\equiv{-2} \pmod 8$ as follows
$$10_4, \, 10_{9}, \, 10_{15}, \, 10_{16}, \, 10_{18}, \, 10_{23}, \, 10_{24}, \, 10_{29}, \, 10_{40}, \, 10_{41}, \, 10_{44}, \, 10_{49}, \, 10_{51}, \, 10_{66}, \, 10_{67}, \,10_{70}$$
$$10_{73}, \, 10_{74}, \, 10_{82}, \, 10_{83}, \, 10_{84}, \, 10_{89}, \, 10_{93}, \, 10_{94}, \, 10_{97}, \, 10_{103}, \, 10_{108}, \, 10_{113}, \, 10_{125}, \, 10_{126}, \, 10_{128}$$
$$10_{131}, \,10_{139}, \,10_{143}, \, 10_{144}, \, 10_{145}, \, 10_{151}, \, 10_{152}, \, 10_{156}, \, 10_{156}, \, 10_{159}, \, 10_{162}$$
Among these 42 knots we have only 5 knots with $\gamma_4$ equal to 2, which are
$$10_9, \, 10_{18}, \,10_{84},\,10_{95}, \,10_{113}$$
By Theorem \ref{group2,3}, we show that $\gamma_4(K)\geq{2}$ and then according to Proposition \ref{band move}, we find a band move to a knot with $\gamma_4$ equal to 1. Thus we prove that $\gamma_4(K)$ for these five knots equals 2. Note that in this group of knots as we mentioned before, since we are looking for the negative-definite Goeritz matrix, we work with the mirror of knot $K$ and it is denoted by $-K$. We also can consider knot $K$ and find the positive-definite Goeritz matrix. Both ways are completely correct.
\begin{figure}[h]
\centering
\begin{subfigure}[b]{0.30\textwidth}
\includegraphics[width=\textwidth]{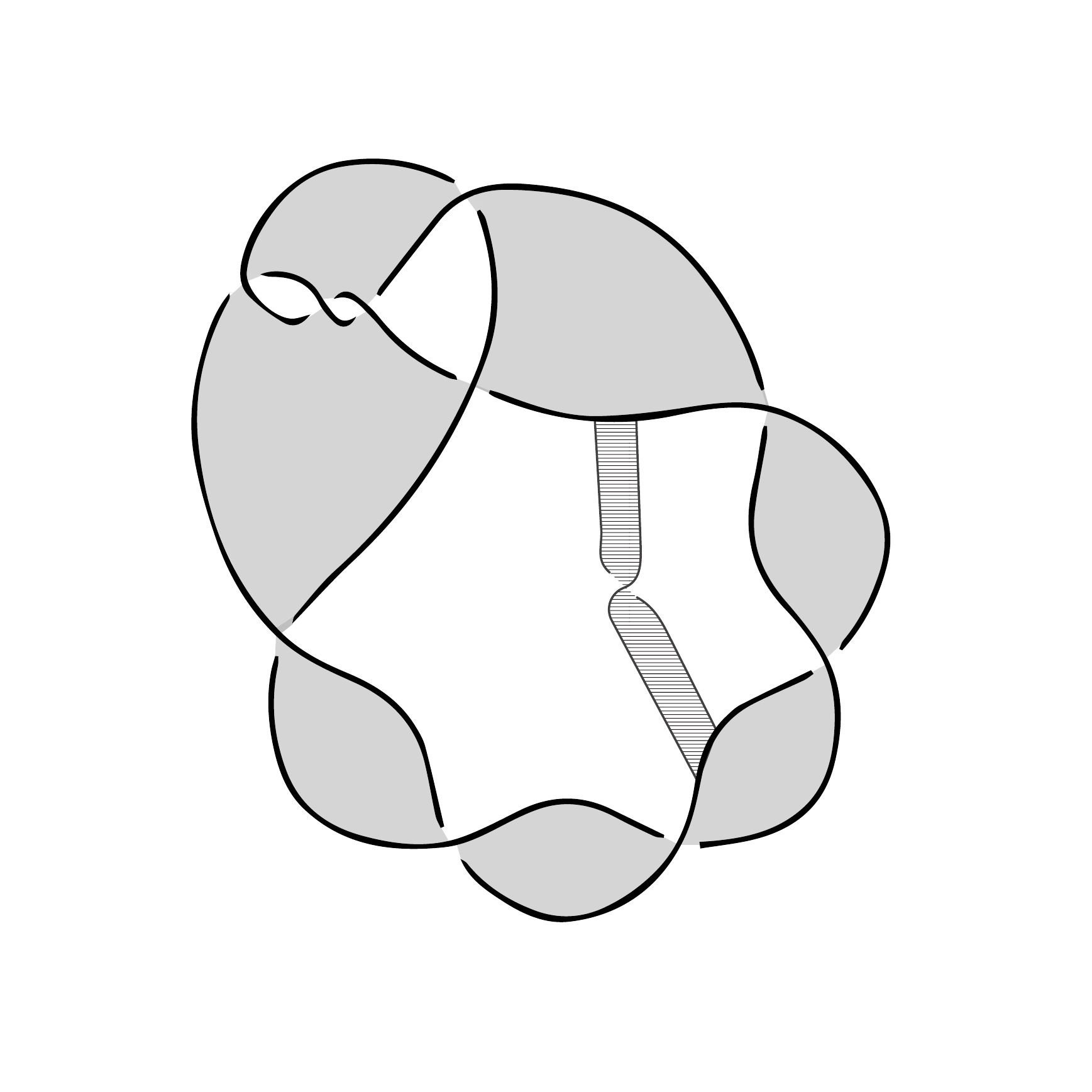}
\caption{Checkerboard diagram for $K=-10_{9}\stackrel{1}{\longrightarrow} 6_{2}$} \label{10-9a}
\end{subfigure}
\qquad \qquad \qquad \qquad
\begin{subfigure}[b]{0.30\textwidth}
\includegraphics[width=\textwidth]{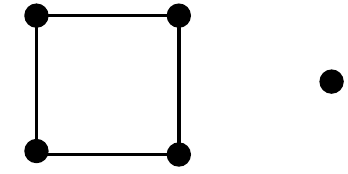}
\put(-135,5){$e_1$}
\put(-130, -10){$-7$}
\put(-80,-10){$-3$}
\put(-130,68){$-2$}
\put(-135,55){$e_2$}
\put(-80, 70){$-2$}
\put(-55,5){$e_4$}
\put(-55, 52){$e_3$}
\put(-10,23){$e_5$}
\put(-12,42){$-39$}
\caption{Incidence graph.}
\label{10-9b}
\end{subfigure}
\caption{Case of $K=-10_{9}$}
\end{figure}
\par
\textbf{Case of $K=-10_{9}$} by using the checkerboard coloring method as shown in Figure \ref{10-9a} we can find the negative definite Goeritz matrix G such as :
\\
\\
$$G = \begin{bmatrix}
-7& 1 & 0 & 1 \\
1 & -2 & 1 & 0\\
0 & 1 & -2 & 1\\
1 & 0 & 1 & -3 &
\end{bmatrix}.$$
\\
\\
\\	
Figure \ref{10-9b}. is the geometric representation of the Goertiz Matrix. Since  $\det (-10_{9}) = 39$   is square-free, we seek an embedding $\phi : (\mathbb{Z}^5,G\oplus [-39])\hookrightarrow(\mathbb{Z}^5, -Id)$. If such a $\phi$ existed, we assume that  $\phi(e_2)=f_1+f_2$ then up to isomorphism, we have only one possibility for $\phi(e_3)=-f_2+f_3$. Now to find $\phi(e_4)$ we have 2 cases :
\begin{enumerate}
\item $\phi(e_4)=-f_1+f_2+f_4$
\\
Let $\phi(e_5)=\sum_{i=1}^{5}\lambda_if_i$, where $\lambda_i's$ are integers. Then for $i=2,3,4$ , the relation $e_5.e_i=0$ gives us the following equations:
\begin{align*}
-\lambda_1-\lambda_2&=0 \cr
\lambda_2-\lambda_3&=0 \cr
\lambda_1-\lambda_2-\lambda_4&=0 \cr
\end{align*}
Solving for $\lambda_1$, $\lambda_3$ and $\lambda_4$ in terms of $\lambda_2$, we obtain, $\lambda_1=-\lambda_2$, $ \lambda_3=\lambda_2$ and $\lambda_4=-2\lambda_2$.
Since $\lambda_1^2+\lambda_2^2+\lambda_3^2+\lambda_4^2+\lambda_5^2=39$ we have $7\lambda_2^2+\lambda_5^2=39$.Thus $|\lambda _2| \le 2$ and neither of the 5 possibilities of $\lambda _1\in \{0, \pm 1, \pm2\}$ leads to an integral solution of $\lambda _5$, showing that  the case of $\phi(e_4)=-f_1+f_2+f_4$ cannot occur.
\item $\phi(e_4)=-f_3+f_4+f_5$
\\
Let $\phi(e_5)=\sum_{i=1}^{5}\lambda_if_i$, where $\lambda_i's$ are integers. Then for $i=2,3,4$ , the relation $e_5.e_i=0$ gives us the following equations:
\begin{align*}
-\lambda_1-\lambda_2&=0 \cr
\lambda_2-\lambda_3&=0 \cr
\lambda_3-\lambda_4-\lambda_5&=0
\end{align*}
By solving the above equations we have, $\lambda_1=-\lambda_2$ , $\lambda_3=\lambda_2$ and $\lambda_5=\lambda_2-\lambda_4$.Since $\lambda_1^2+\lambda_2^2+\lambda_3^2+\lambda_4^2+\lambda_5^2=39$ we have $4\lambda_2^2+2\lambda_4^2-2\lambda_2\lambda_4=39$ and this is contradiction since the left side of the equation is an even number, and the right side is an odd number.Thus the case of $\phi(e_4)=-f_3+f_4+f_5$ cannot occur.
\end{enumerate}
These two cases show that the embedding $\phi$ does not exist and therefore that $\gamma_4(-10_{9})\geq2$. As we have shown in Figure \ref{10-9a}, there is a non-oriented band from $-10_{9}$ to $6_2$, and since $\gamma_4(6_2)=1$, we can conclude that $\gamma_4(-10_{9})=2$.\

	\begin{figure}[h]
\centering
\begin{subfigure}[b]{0.30\textwidth}
\includegraphics[width=\textwidth]{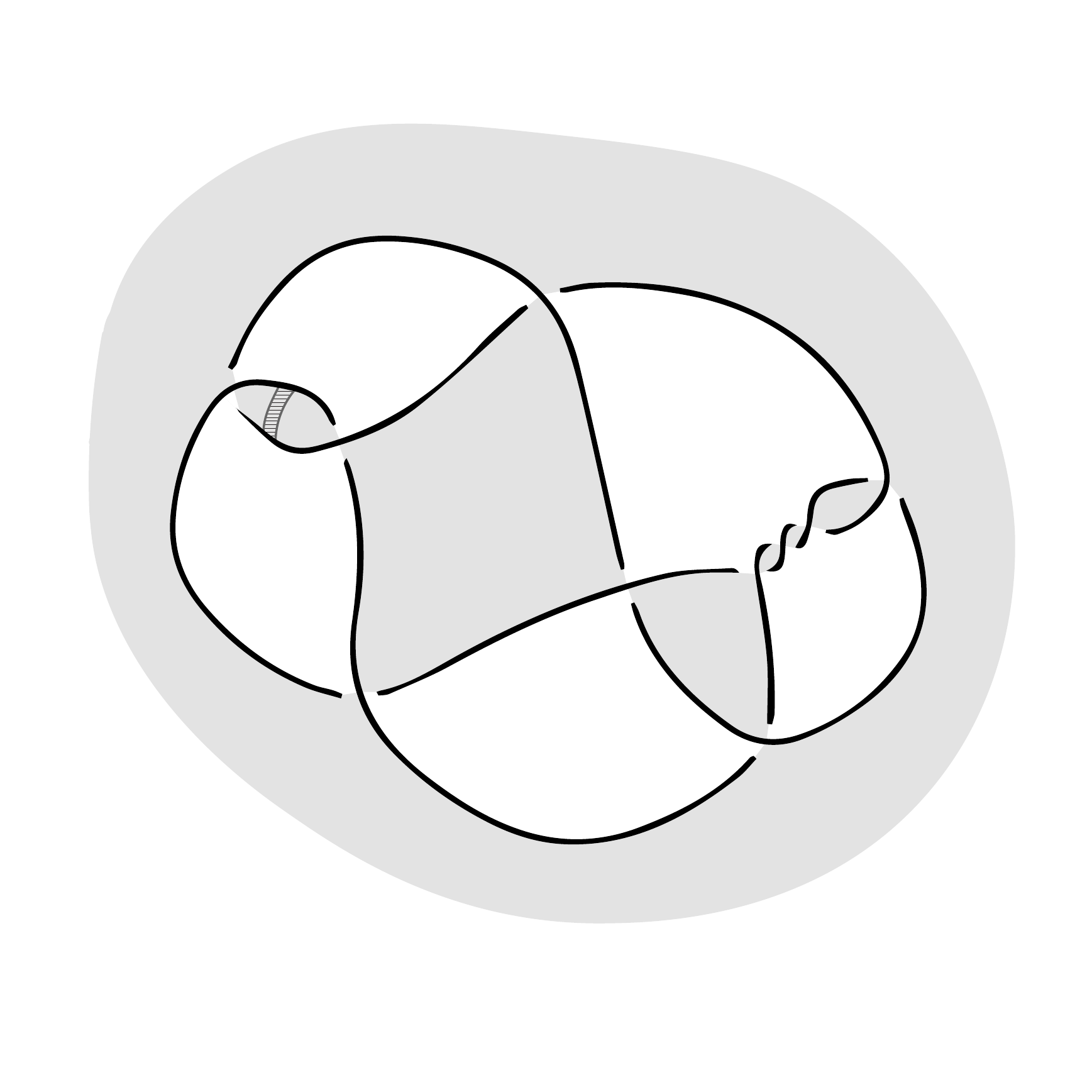}
\caption{Checkerboard diagram for $K=-10_{18}\stackrel{0}{\longrightarrow} 8_{7}$} \label{10-18a}
\end{subfigure}
\qquad \qquad \qquad \qquad
\begin{subfigure}[b]{0.30\textwidth}
\includegraphics[width=\textwidth]{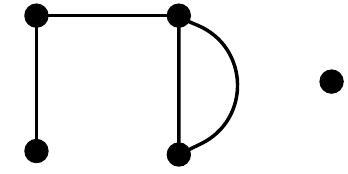}
\put(-135,5){$e_1$}
\put(-125, -10){$-5$}
\put(-65,-10){$-3$}
\put(-135,68){$-3$}
\put(-135,55){$e_2$}
\put(-65, 70){$-3$}
\put(-80,5){$e_4$}
\put(-80, 52){$e_3$}
\put(-10,23){$e_5$}
\put(-12,42){$-55$}
\caption{Incidence graph.}
\label{10-18b}
\end{subfigure}
\caption{Case of $K=-10_{18}$}
\end{figure}
\par	
\textbf{Case of $K=-10_{18}$} by using the checkerboard coloring method as shown in Figure \ref{10-18a} we can find the negative definite Goeritz matrix G such as :
\\
\\
$$G = \begin{bmatrix}
-5 & 1 & 0 & 0\\
1 & -3 & 1 & 0 \\
0 & 1 & -3 & 2 \\
0 & 0 & 2 & -3 \\
\end{bmatrix}.$$
\\
Figure \ref{10-18b}. is the geometric representation of the Goertiz Matrix. Since  $\det (-10_{18}) =55$   is square-free, we seek an embedding $\phi : (\mathbb{Z}^5,G\oplus [-55])\hookrightarrow(\mathbb{Z}^5, -Id)$. If such a $\phi$ existed, we assume that $\phi(e_4)=f_1+f_2+f_3$ then up to isomorphism we have only one possibility for $\phi(e_3)=-f_1-f_2+f_4$. To find $\phi(e_2)$ we have two cases:
\begin{enumerate}
\item $\phi(e_2)=f_1-f_3+f_5$
\\
 Let $\phi(e_1)=\sum_{i=1}^{5}\lambda_if_i$, where $\lambda_i's$ are integers. Then $e_1.e_3=e_1.e_4=0$ and $e_1.e_2=1$ give us the following equations:
\begin{align*}
-\lambda_1-\lambda_2-\lambda_3&=0 \cr
\lambda_1+\lambda_2-\lambda_4&=0 \cr
-\lambda_1+\lambda_3-\lambda_5&=1 \cr
\end{align*}
Note that it is not possible that $|\lambda_i|=1$ for $i=1,...,5$ (using all five elements of the basis) because $e_1$ and $e_4$ share no edge, while they share three elements of the basis. Therefore $\phi(e_1)=\lambda_mf_m+\lambda_nf_n$ for some $1\leq m,n\leq5$, where $|\lambda_m|=2$ and $|\lambda_n|=1$. Now by adding the first and second equations, we obtain $\lambda_3=\lambda_4$. But, they both must be zero because of the form of the $\phi(e_1)$. Now using the first equation, $\lambda_1=-\lambda_2=0$ by the same justification. This is a contradiction since only $\lambda_5$ can be nonzero. it shows that $\phi ( e_2) = f_1-f_3+f_5$ cannot occur. 
\item $\phi(e_2)=f_1-f_2-f_4$
\\
Let $\phi(e_5)=\sum_{i=1}^{5}\lambda_if_i$, where $\lambda_i's$ are integers. Then for $i=2,3,4$ , $e_5.e_i=0$ gives us the following equations:

\begin{align*}
-\lambda_1-\lambda_2-\lambda_3&=0 \cr
\lambda_1+\lambda_2-\lambda_4&=0 \cr
-\lambda_1+\lambda_2+\lambda_4&=0 \cr
\end{align*}
Solving $\lambda_1$, $\lambda_2$, $\lambda_4$ and $\lambda_5	$ in terms of $\lambda_3$ we obtain
$$\lambda_1=-\lambda_3, \quad \lambda_2=0, \quad \lambda_4=-\lambda_3$$
Since $\lambda_1^2+\lambda_2^2+\lambda_3^2+\lambda_4^2+\lambda_5^2=55$, we find that $3\lambda_3^2+\lambda_5^2=55$ and it is contradiction because for all possibilities of $|\lambda_3|\leq4$ there is no integral solution for $\lambda_5$. Thus the case of $\phi(e_2)=f_1-f_2-f_4$ cannot occur.
\end{enumerate}
These two cases show that the embedding $\phi$ does not exist and $\gamma_4(-10_{18})\geq2$. As we have shown in Figure \ref{10-18a}, there is a non-oriented band from $-10_{18}$ to $8_7$, and since $\gamma_4(8_7)=1$, we can conclude that $\gamma_4(-10_{18})=2$.\
\\
\\
\begin{figure}[h]
\centering
\begin{subfigure}[b]{0.30\textwidth}
\includegraphics[width=\textwidth]{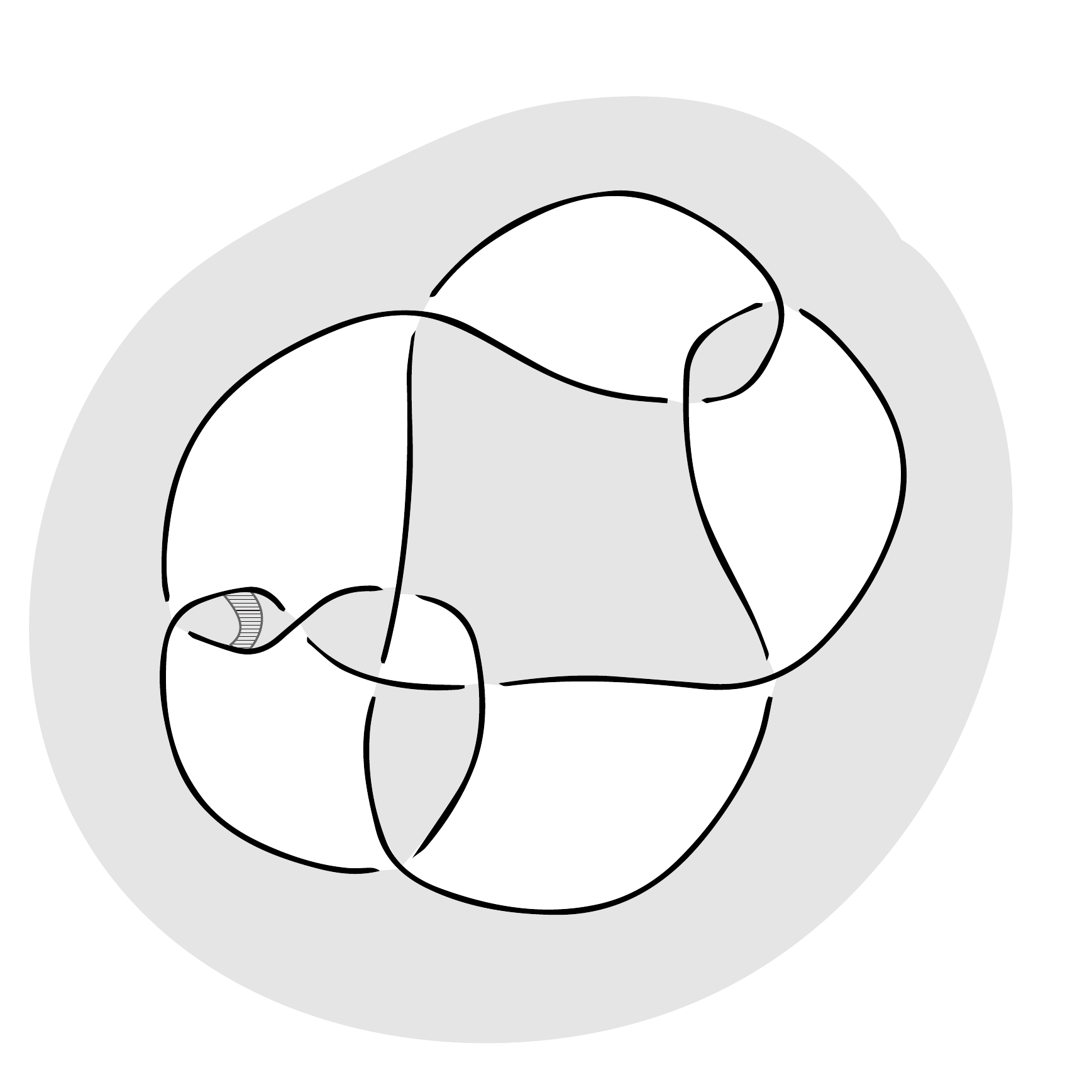}
\caption{Checkerboard diagram for $K=-10_{84}\stackrel{0}{\longrightarrow} 8_{14}$} \label{-10-84a}
\end{subfigure}
\qquad \qquad \qquad \qquad
\begin{subfigure}[b]{0.30\textwidth}
\includegraphics[width=\textwidth]{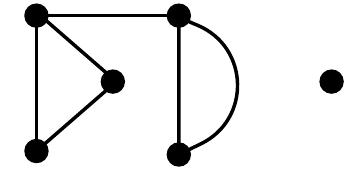}
\put(-135,5){$e_5$}
\put(-125, -10){$-4$}
\put(-65,-10){$-3$}
\put(-135,68){$-3$}
\put(-135,55){$e_3$}
\put(-65, 70){$-3$}
\put(-80,5){$e_1$}
\put(-80, 52){$e_2$}
\put(-10,23){$e_6$}
\put(-90,23){$e_4$}
\put(-12,42){$-87$}
\caption{Incidence graph.}
\label{-10-84b}
\end{subfigure}
\caption{Case of $K=-10_{84}$}
\end{figure}
\textbf{Case of $K=-10_{84}$} by using the checkerboard coloring method as shown in Figure \ref{-10-84a} we can find the negative definite Goeritz matrix G such as :
\\
$$G = \begin{bmatrix}
-3 & 2 & 0 & 0 & 0\\
2 & -3 & 1 & 0 & 0\\
0 & 1 & -3 & 1 & 1\\
0 & 0 & 1 & -3 & 1\\
0 & 0 & 1 & 1 & -4
\end{bmatrix}.$$
\\	
Figure \ref{-10-84b}. is the geometric representation of the Goertiz Matrix. Since  $\det (-10_{84}) = 87$   is   square-free, we seek an embedding,  $\phi : (\mathbb{Z}^6,G\oplus [-87])\hookrightarrow(\mathbb{Z}^6, -Id)$, .If such a $\phi$ existed, suppose $\phi(e_5)=f_1+f_2+f_3+f_4$ then we have two possibilities for $\phi(e_4)$:
\begin{enumerate}
\item $\phi(e_4)=-f_1-f_2+f_3$ then we need to check two cases for $\phi(e_3)$:
\begin{enumerate}
\item If $\phi(e_3) =-f_3+f_5+f_6$ then $\phi(e_2)=-f_5+f_1-f_2$.
\\
Let $\phi(e_1)=\sum_{i=1}^{6}\lambda_if_i$, where $\lambda_i's$ are integers. Then for $i=2,...,5$ , $e_1.e_i$ gives us the following equations:
\begin{align*}
\lambda_5-\lambda_1+\lambda_2 & =2\cr
\lambda_3-\lambda_5-\lambda_6&=0 \cr
\lambda_1+\lambda_2-\lambda_3&=0 \cr
-\lambda_1-\lambda_2-\lambda_3-\lambda_4&=0
\end{align*}
Solving for $\lambda_3$,  $\lambda_4$, $\lambda_5$ and $\lambda_6$  in terms of $\lambda_1$ and $\lambda_2$, we obtain
$$ \lambda_3=\lambda_1+\lambda_2, \quad  \lambda_4=-2\lambda_1-2\lambda_2, \quad \lambda_5=\lambda_1-\lambda_2+2, \quad \lambda_6=2\lambda_2-2$$
Since $\lambda_1^2+\lambda_2^2+\lambda_3^2+\lambda_4^2+\lambda_5^2+\lambda_6^2=3 $ , we find that $\lambda_1^2+\lambda_2^2+5(\lambda_1+\lambda_2)^2+(\lambda_1-\lambda_2+2)^2+4(\lambda_2-1)^2=3$. $\lambda_2$ must equal $1$, otherwise the left side is greater than 3. By the same justification, $\lambda_1$ must equal $-1$ which leads to a contradiction.
\item If $\phi(e_3)=f_1-f_2-f_3$, then let $\phi(e_2)=\sum_{i=1}^{6}\lambda_if_i$, where $\lambda_i's$ are integers.Then $e_2.e_3=1$ and $e_2.e_4=0$ lead to the following equations
\begin{align*}
\lambda_1+\lambda_2-\lambda_3&=0 \cr
-\lambda_1+\lambda_2+\lambda_3&=1 \cr
\end{align*}
By adding these equations we have $\lambda_2=\frac{1}{2}$, which is impossible and the case of $\phi(e_3)=f_1-f_2-f_3$ is not acceptable.
\end{enumerate}
\item $\phi(e_4)=-f_1+f_5+f_6$ then we have two possibilities for $\phi(e_3)$:
\begin{enumerate}
\item If $\phi(e_3)=f_1-f_2-f_3$ then $\phi(e_2)=-f_1-f_5+f_4$
Let $\phi(e_1)=\sum_{i=1}^{6}\lambda_if_i$, where $\lambda_i's$ are integers. Then for $i=3,4,5$ , $e_1.e_i=0$ and $e_1.e_2=2$ gives us the following equations
\begin{align*}
-\lambda_1-\lambda_2-\lambda_3-\lambda_4&=0 \cr
\lambda_1-\lambda_5-\lambda_6&=0 \cr
-\lambda_1+\lambda_2+\lambda_3&=0 \cr
\lambda_1+\lambda_5-\lambda_4&=2 \cr
\end{align*}
By adding the first and third equations we have $\lambda_4=-2\lambda_1$, so $\lambda_4$ is an even number and must be zero because $\lambda_1^2+...+\lambda_6^2=3$. Then $\lambda_1=0$ and the last equation gives us $\lambda_5=2$ which is a contradiction.
 This shows that the case of $\phi (e_4) =-f_1+f_5+f_6$ cannot occur.
 \item If $\phi(e_3)=-f_1-f_5-f_6$ then let $\phi(e_2)=\sum_{i=1}^{6}\lambda_if_i$, where $\lambda_i's$ are integers. Then $e_2.e_3=1$ and $e_2.e_4=0$ lead to the following equations:
 \begin{align*}
 \lambda_1-\lambda_5-\lambda_6&=0 \cr
 \lambda_1+\lambda_5+\lambda_6&=1 \cr
 \end{align*}
 By adding these equations we get $\lambda_1=\frac{1}{2}$ which is impossible because $\lambda_1$ is an integer. This shows that the case of $\phi(e_3)=-f_1-f_5-f_6$ cannot occur.
 \end{enumerate}
\end{enumerate}
This shows that the embedding $\phi$ does not exist and therefore that $\gamma_4(-10_{84})\geq2$. As we have shown in Figure \ref{-10-84a}, there is a non-oriented band from $-10_{84}$ to $8_{14}$, and since $\gamma_4(8_{14})=1$, we can conclude that $\gamma_4(-10_{84})=2$.\
\begin{figure}[h]
\centering
\begin{subfigure}[b]{0.30\textwidth}
\includegraphics[width=\textwidth]{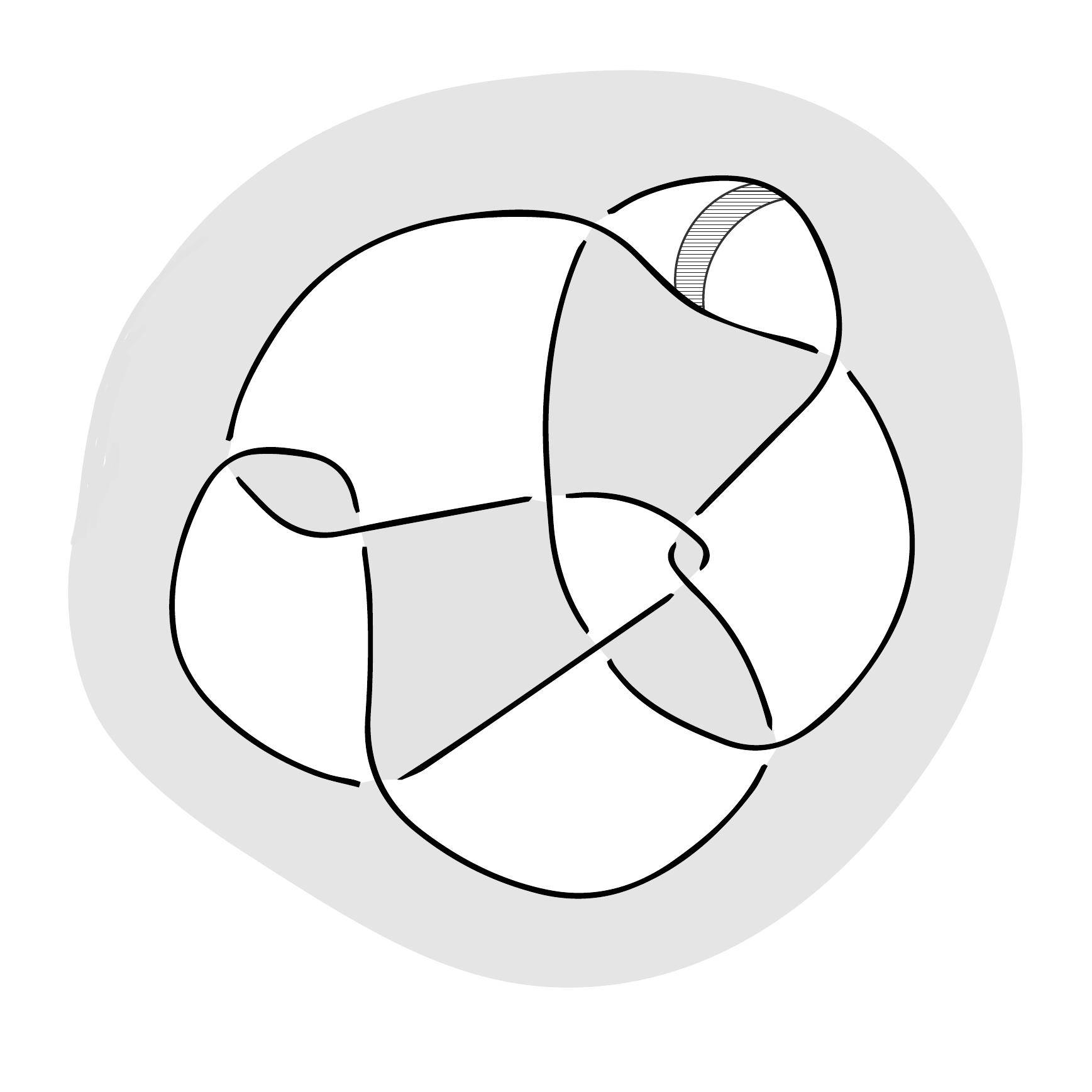}
\caption{Checkerboard diagram for $K=-10_{95}\stackrel{0}{\longrightarrow} 8_{14}$} \label{-10_95a}
\end{subfigure}
\qquad \qquad \qquad \qquad
\begin{subfigure}[b]{0.30\textwidth}
\includegraphics[width=\textwidth]{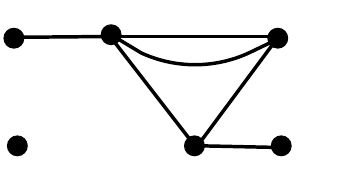}
\put(-145,10){$e_6$}
\put(-140,-5){$-91$}
\put(-80,10){$e_4$}
\put(-75,-5){$-3$}
\put(-15,10){$e_5$}
\put(-35,-5){$-3$}
\put(-145,50){$e_1$}
\put(-140,35){$-2$}
\put(-100,40){$e_2$}
\put(-100,60){$-4$}
\put(-15,50){$e_3$}
\put(-35,60){$-4$}
\caption{Incidence graph.}
\label{-10_95b}
\end{subfigure}
\caption{Case of $K=-10_{95}$}
\end{figure}
\\
\par
\textbf{Case of $K=-10_{95}$} by using the checkerboard coloring method as shown in Figure \ref{-10_95a} we can find the negative definite Goeritz matrix G such as :
\\
$$G = \begin{bmatrix}
-2 & 1 & 0 & 0 & 0\\
1 & -4 & 2 & 1 & 0\\
0 & 2 & -4 & 1 & 0\\
0 & 1 & 1 & -3 & 1\\
0 & 0 & 0 & 1 & -3
\end{bmatrix}.$$
\\
Figure \ref{-10_95b}. is the geometric representation of the Goertiz Matrix. Since  $\det (-10_{95}) = 91$   is   square-free, we seek the embedding,  $\phi : (\mathbb{Z}^6,G\oplus [-91])\hookrightarrow(\mathbb{Z}^6, -Id)$, .If such a $\phi$ existed, suppose $\phi(e_1)=f_1+f_2$ up to isomorphism the only possibilities for $\phi(e_2)$, $\phi(e_3)$ and $\phi(e_4)$ are as follows:
\begin{align*}
\phi(e_2)&=-f_2+f_3+f_4+f_5 \cr
\phi(e_3)&=f_2-f_1-f_3+f_6 \cr
\phi(e_4)&=f_3-f_4-f_5 \cr
\end{align*}
Let $\phi(e_5)=\sum_{i=1}^{6}\lambda_if_i$, where for $i=1,2,3$, $e_5.e_i=0$ and $e_5.e_4=1$ give us the following equations
\begin{align*}
-\lambda_1-\lambda_2&=0 \cr
\lambda_2-\lambda_3-\lambda_4-\lambda_5&=0 \cr
-\lambda_2+\lambda_1+\lambda_3-\lambda_6&=0 \cr
-\lambda_3+\lambda_4+\lambda_5&=1 \cr
\end{align*}
By adding the second and fourth equations we have, $\lambda_2=2\lambda_3+1$ and from the first equation we know that $\lambda_1=-\lambda_2$, so $\lambda_1=-2\lambda_3-1$. By plugging $\lambda_1$ and $\lambda_2$ into the third equation we obtain, $\lambda_6=5\lambda_3+2$. For all possibilities of $\lambda_3=0,\pm1$ there is no possible solution for $\lambda_6$ since $\lambda_1^2+...+\lambda_6^2=3$.It shows that $\phi$ does not exist and therefore that $\gamma_4(-10_{95})\geq2$. As we have shown in Figure \ref{-10_95a}, there is a non-oriented band from $-10_{95}$ to $8_{14}$, and since $\gamma_4(8_{14})=1$, we can conclude that $\gamma_4(-10_{95})=2$.\
\begin{figure}[h]
\centering
\begin{subfigure}[b]{0.30\textwidth}
\includegraphics[width=\textwidth]{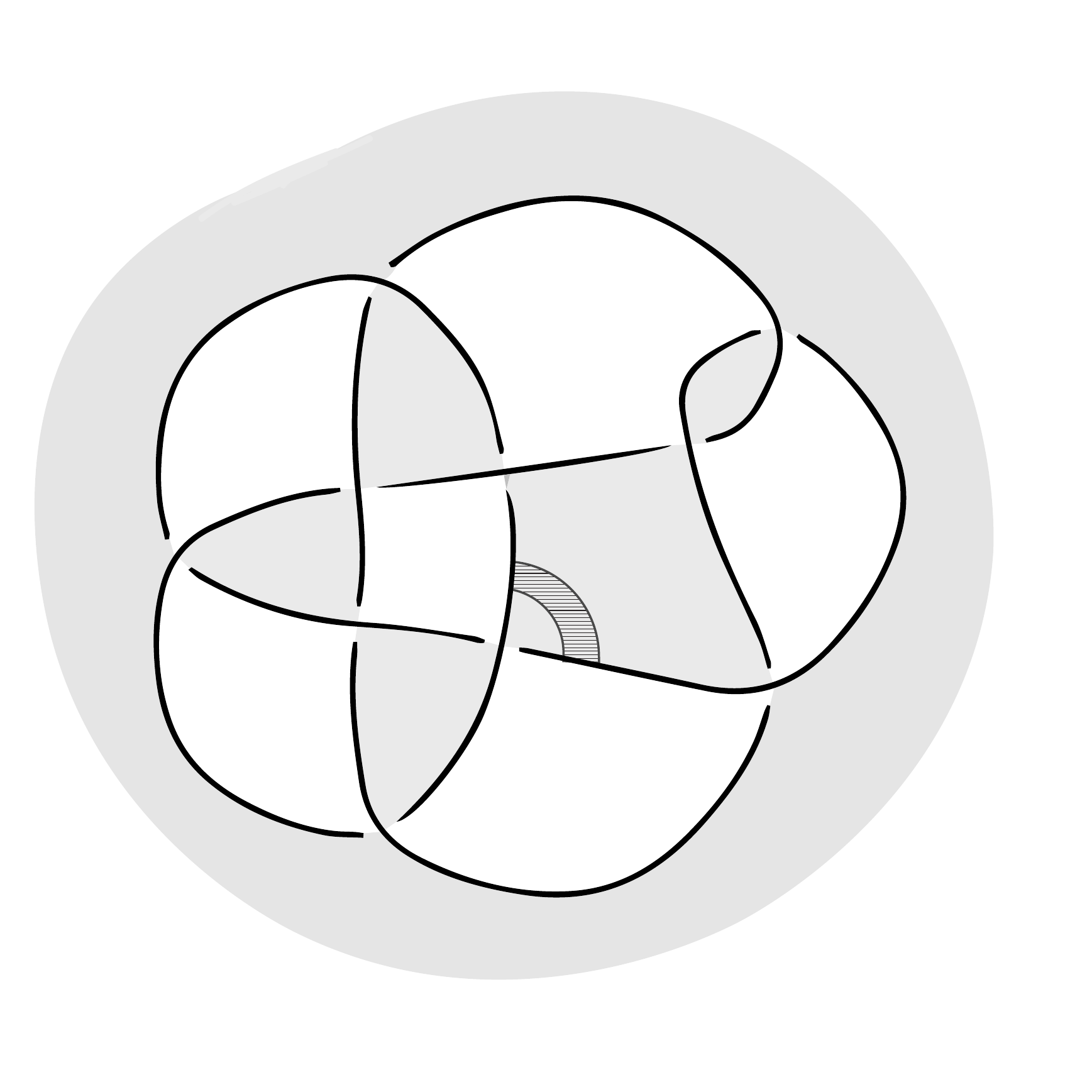}
\caption{Checkerboard diagram for $K=-10_{113}\stackrel{0}{\longrightarrow} 9_{31}$} \label{-10_113a}
\end{subfigure}
\qquad \qquad \qquad \qquad
\begin{subfigure}[b]{0.30\textwidth}
\includegraphics[width=\textwidth]{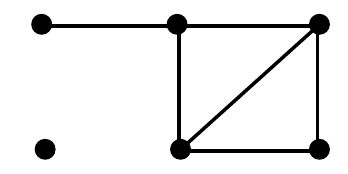}
\put(-135,10){$e_6$}
\put(-130,-5){$-111$}
\put(-80,10){$e_4$}
\put(-75,-5){$-3$}
\put(0,10){$e_5$}
\put(-25,-5){$-3$}
\put(-135,50){$e_1$}
\put(-130,64){$-3$}
\put(-80,45){$e_2$}
\put(-80,64){$-3$}
\put(0,50){$e_3$}
\put(-25,64){$-4$}
\caption{Incidence graph.}
\label{-10_113b}
\end{subfigure}
\caption{Case of $K=-10_{113}$}
\end{figure}
\\
\textbf{Case of $K=-10_{113}$} by using the checkerboard coloring method as shown in Figure \ref{-10_113a} we can find the negative definite Goeritz matrix G such as :
\\
$$G = \begin{bmatrix}
-3 & 1 & 0 & 0 & 0\\
1 & -3 & 1 & 1 & 0\\
0 & 1 & -4 & 1 & 1\\
0 & 1 & 1 & -3 & 1\\
0 & 0 & 1 & 1 & -3
\end{bmatrix}.$$
\\
Figure \ref{-10_113b}. is the geometric representation of the Goertiz Matrix. Since  $\det (-10_{113}) = 111$   is   square-free, we seek the embedding,  $\phi : (\mathbb{Z}^6,G\oplus [-111])\hookrightarrow(\mathbb{Z}^6, -Id)$, .If such a $\phi$ existed, suppose $\phi(e_3)=f_1+f_2+f_3+f_4$, then there are two cases for $\phi(e_5)$,
\begin{enumerate}
\item $\phi(e_5)=-f_4+f_5+f_6$ then we have two possibilities for $\phi(e_4)$ as follows
\begin{enumerate}
 \item If $\phi(e_4)=-f_1-f_2+f_4$.
\\
Let $\phi(e_1)=\sum_{i=1}^{6}\lambda_if_i$, where $\lambda_i's$ are integers and $e_1.e_3=e_1.e_4=e_1.e_5=0$ leads to the following equations
\begin{align*}
-\lambda_1-\lambda_2-\lambda_3-\lambda_4&=0 \cr
\lambda_4-\lambda_5-\lambda_6&=0 \cr
\lambda_1+\lambda_2-\lambda_4&=0 \cr
\end{align*}
By adding the first and third equations $\lambda_3=-2\lambda_4$ which is an even number, but since $\lambda_1^2+...+\lambda_6^2=3$, $\lambda_3$ must be zero and then $\lambda_4=0$. Thus $\lambda_1=-\lambda_2$ and $\lambda_5=-\lambda_6$ which is not possible because just three of the $\lambda_i's$ can be nonzero. Therefore, the case of $\phi(e_5)=-f_4+f_5+f_6$ and $\phi(e_4)=-f_1-f_2+f_4$ cannot occur.
\item If $\phi(e_4)=-f_4-f_5-f_6$
\\
Let $\phi(e_2)=\sum_{i=1}^{6}\lambda_if_i$, where $\lambda_i's$ are integers and $e_2.e_4=1$ and $e_2.e_5=0$, lead to the following equations
\begin{align*}
\lambda_4-\lambda_5-\lambda_6&=0 \cr
\lambda_4+\lambda_5+\lambda_5&=1
\end{align*}
By adding these two equations we find that $\lambda_4=\frac{1}{2}$ which contradicts $\lambda_4$ being an integer. Thus the case of $\phi(e_5)=-f_4+f_5+f_6$ and $\phi(e_4)=-f_4-f_5-f_6$ cannot occur.
\end{enumerate}
\item $\phi(e_5)=-f_1-f_2+f_3$, then
\begin{align*}
\phi(e_1)&=-f_1+f_2+f_5 \cr
\phi(e_2)&=-f_2-f_3+f_4 \cr
\end{align*}
Let $\phi(e_4)=\sum_{i=1}^{6}\lambda_if_i$, where $\lambda_i's$ are integers and $e_4.e_5=e_4.e_3=e_4.e_2=1$ and $e_4.e_1=0$ give us the following equations,
\begin{align*}
-\lambda_1-\lambda_2-\lambda_3-\lambda_4&=1 \cr
\lambda_1+\lambda_2-\lambda_3&=1 \cr
\lambda_1-\lambda_2-\lambda_5&=0 \cr
\lambda_2+\lambda_3-\lambda_4&=1 \cr
\end{align*}
By adding the first and forth equations we obtain $\lambda_1=-2\lambda_4-2$ and because $\lambda_1$ is  an even number, the only possibility for $\lambda_1$ is zero, since $\lambda_1^2+...+\lambda_6^2=3$. It implies that $\lambda_4=-1$.Now by adding the first two equations, we find that $\lambda_3=-\frac{1}{2}$ which is a contradiction because $\lambda_3$ is an integer. Thus the case of $\phi(e_5)=-f_1-f_2+f_3$ cannot occur.
\end{enumerate}

These two cases show that the embedding $\phi$ does not exist and therefore that $\gamma_4(-10_{113})\geq2$. As we have shown in Figure \ref{-10_113a}, there is a non-oriented band from $-10_{113}$ to $9_{31}$, and since $\gamma_4(9_{31})=1$, we can conclude that $\gamma_4(-10_{113})=2$.\ 
\\
\subsection { Knots with $\sigma(K) + 4\cdot \text{Arf}(K) \equiv{2} \pmod 8$ }\label{third group}
Among all 10-crossing knots there exists exactly 36 knots which satisfy the congruance $\sigma(K) + 4.\text{Arf}(K) \equiv{2} \pmod 8$ as follows
$$10_2, \, 10_7, \,10_{11}, \, 10_{12}, \, 10_{19}, \, 10_{20}, \, 10_{27}, \, 10_{30}, \, 10_{36}, \, 10_{38}, \, 10_{46}, \, 10_{52}, \, 10_{54}, \, 10_{57}, \, 10_{59}, \, 10_{64}, \, 10_{65}$$
$$10_{69}, \, 10_{77}, \, 10_{80}, \,10_{105}, \, 10_{106}, \, 10_{110}, \, 10_{112}, \, 10_{116}, \, 10_{117}, \, 10_{121}, \, 10_{133}$$
$$10_{134}, \, 10_{136}, \, 10_{138}, \, 10_{142}, \, 10_{147}, \, 10_{148}, \, 10_{163}, \, 10_{165}$$ 
Among these 36 knots we have exactly 5 knots with $\gamma_4$ equal 2:
$$10_2, \,10_{19}, \,10_{36}, \, 10_{46},\,10_{112}$$
\begin{figure}[h]
\centering
\begin{subfigure}[b]{0.30\textwidth}
\includegraphics[width=\textwidth]{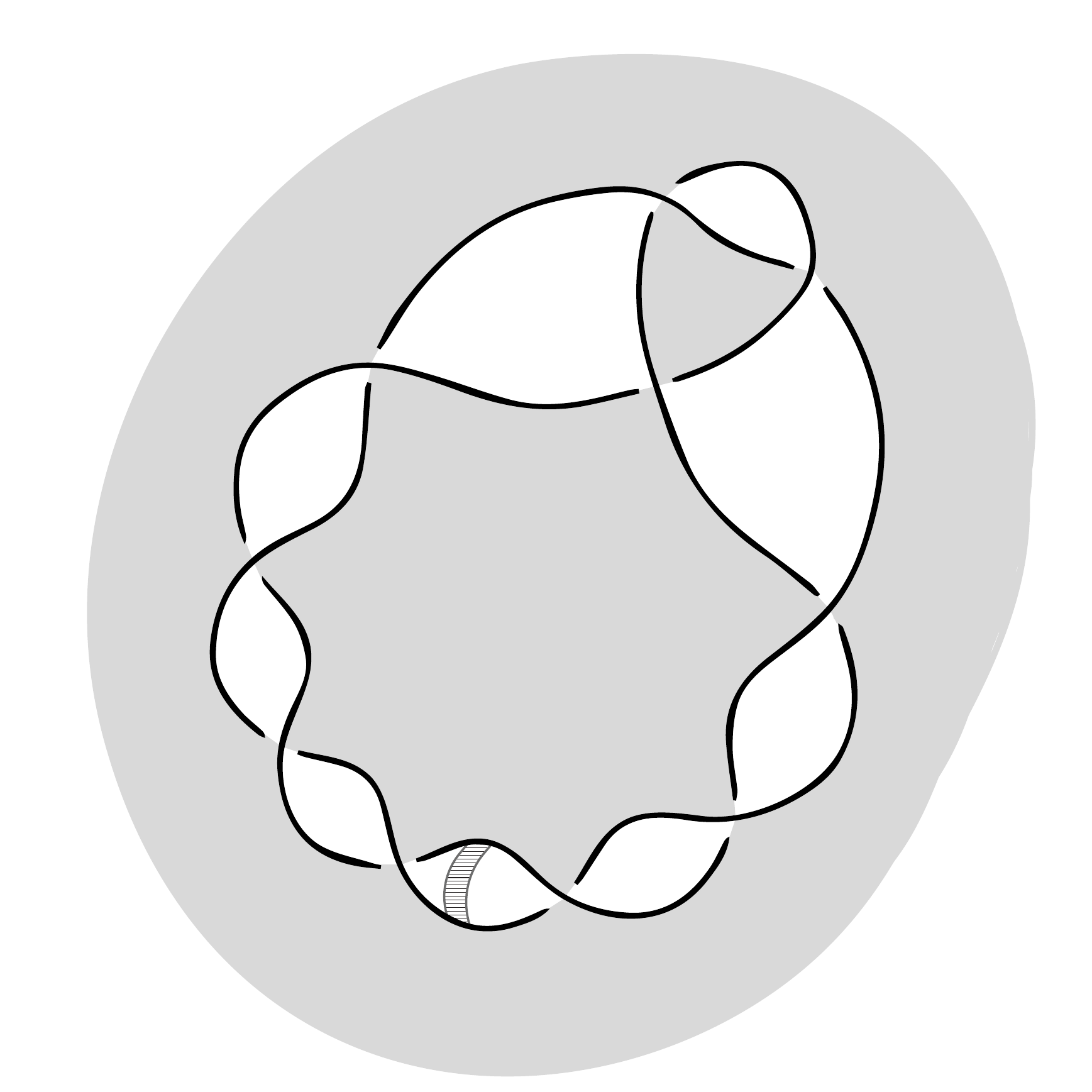}
\put(95,15){$e_5$}
\put(100,0){$-2$}
\put(145,22){$e_4$}
\put(145,0){$-2$}
\put(175,22){$e_3$}
\put(175,0){$-2$}
\put(225,15){$e_2$}
\put(210,0){$-2$}
\put(90,45){$e_6$}
\put(100,50){$-2$}
\put(140,32){$e_7$}
\put(140,50){$-3$}
\put(185,32){$e_8$}
\put(185,50){$-3$}
\put(230,40){$e_1$}
\put(215,51){$-2$}
\caption{Checkerboard diagram for $K=10_{2}\stackrel{0}{\longrightarrow} 3_{1}$} \label{10-2a}
\end{subfigure}
\qquad \qquad \qquad \qquad
\begin{subfigure}[b]{0.30\textwidth}
\includegraphics[width=\textwidth]{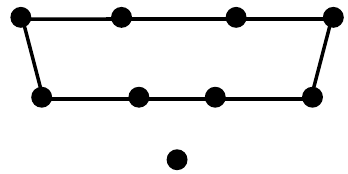}
\put(-80,2){$e_9$}
\put(-60,2){$-23$}
\caption{Incidence graph.}
\label{10-2b}
\end{subfigure}
\caption{Case of $K=10_{2}$}
\end{figure}
\par
\textbf{Case of $K=10_{2}$} by using the checkerboard coloring method as shown in Figure \ref{10-2a} we can find the negative definite Goeritz matrix G such as :
\\
$$G = \begin{bmatrix}
-3 & 1 & 0 & 0 & 0&0&0&1\\
1 & -3 & 1 & 0 & 0 &0&0&0\\
0 & 1 & -2 & 1 & 0&0&0&0\\
0 & 0 & 1 & -2 & 1&0&0&0\\
0 & 0 & 0 & 1 & -2&1&0&0\\
0& 0 &  0 & 0 & 1&-2&1&0\\
0 & 0 & 0& 0 & 0& 1 & -2&1\\
1 & 0 & 0 & 0& 0& 0 &1&-2
\end{bmatrix}.$$
Figure \ref{10-2b}. is the geometric representation of the Goertiz Matrix. Since  $\det (10_{2}) = 23$   is   square-free, we seek an embedding,  $\phi : (\mathbb{Z}^9,G\oplus [-23])\hookrightarrow(\mathbb{Z}^9, -Id)$, .If such a $\phi$ existed, suppose $\phi(e_1)=f_1-f_2$:
\begin{align*}
\phi(e_2) & =f_2-f_3 \cr
\phi(e_3) & =f_3-f_4 \cr
\phi(e_4) & =f_4-f_5\cr
\phi(e_5) & = f_5-f_6\cr
\phi(e_6) & = f_6-f_7\cr
\phi(e_7) & = f_7+f_8+f_9 
\end{align*}
Let $\phi(e_9)=\sum_{i=1}^{9}\lambda_if_i$, where $\lambda_i's$ are integers. Then for $i=1,...,7$ , $e_9.e_i=0$ gives us the following equations:
\begin{align*}
-\lambda_1+\lambda_2& =0\cr
-\lambda_2+\lambda_3&=0 \cr
-\lambda_3+\lambda_4&=0 \cr
-\lambda_4+\lambda_5&=0\cr
-\lambda_5+\lambda_6&=0\cr
-\lambda_6+\lambda_7&=0\cr
-\lambda_7-\lambda_8-\lambda_9&=0
\end{align*}
We obtain $$\lambda_1=\lambda_2=\lambda_3=\lambda_4=\lambda_5=\lambda_6=\lambda_7=0, \quad \lambda_8=-\lambda_9$$
Since $\lambda_1^2+\lambda_2^2+\lambda_3^2+\lambda_4^2+\lambda_5^2+\lambda_6^2+\lambda_7^2+\lambda_8^2+\lambda_9^2=23$ , we find that $2\lambda_8^2=23$ which is a contradiction, because $\lambda_8$ is an integer.
This shows that the embedding $\phi$ does not exist and therefore that $\gamma_4(10_{2})\geq2$. As we have shown in Figure \ref{10-2a}, there is a non-oriented band from $10_{2}$ to $3_{1}$, and since $\gamma_4(3_{1})=1$, we can conclude that $\gamma_4(10_{2})=2$.\ 
\begin{figure}[h]
\centering
\begin{subfigure}[b]{0.30\textwidth}
\includegraphics[width=\textwidth]{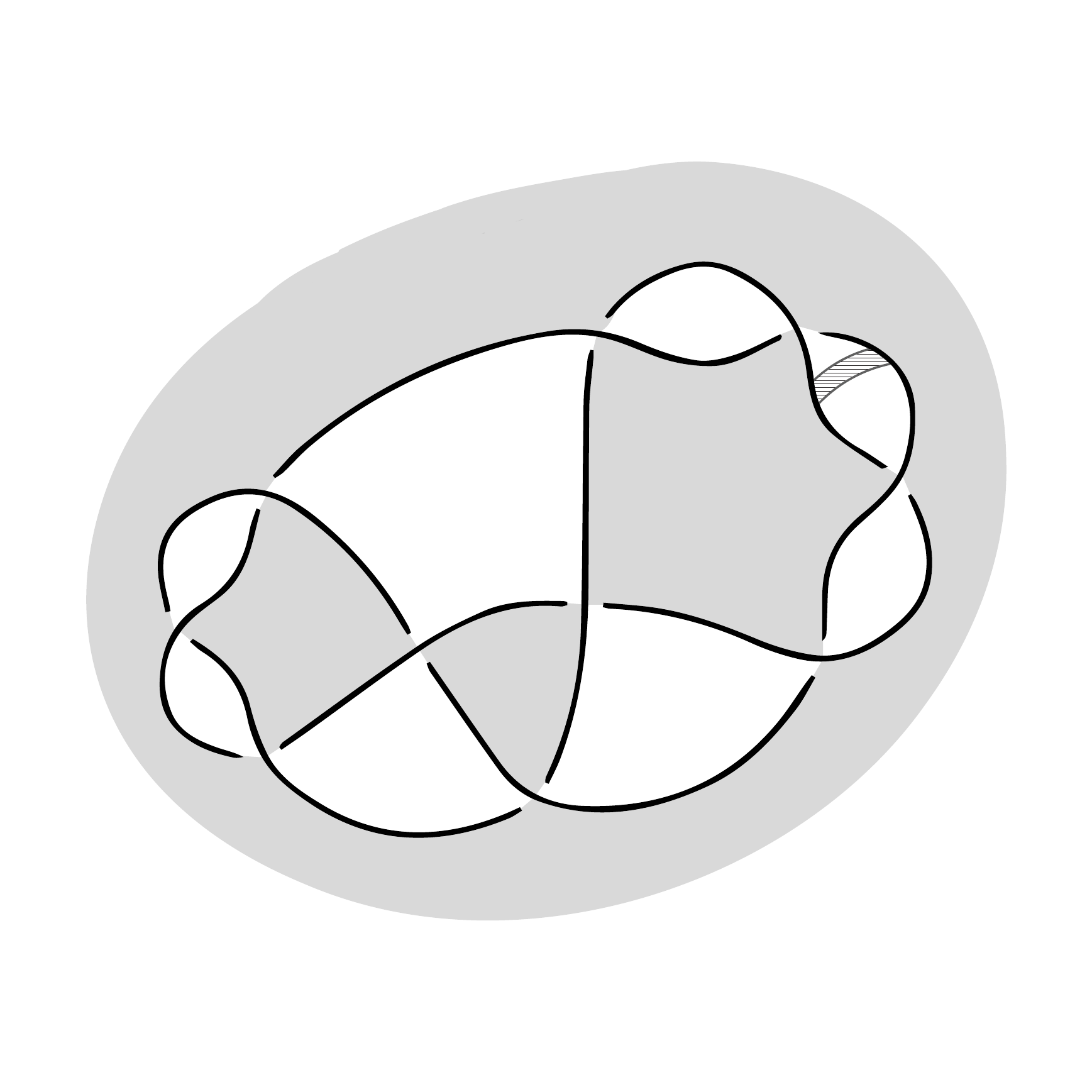}
\caption{Checkerboard diagram for $K=10_{19}\stackrel{0}{\longrightarrow} 6_{2}$} \label{10_19a}
\end{subfigure}
\qquad \qquad \qquad \qquad
\begin{subfigure}[b]{0.30\textwidth}
\includegraphics[width=\textwidth]{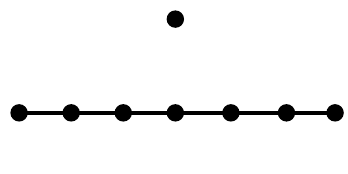}
\put(-10,10){$e_7$}
\put(-15,30){$-2$}
\put(-75,50){$e_8$}
\put(-80,65){$-51$}
\put(-30,10){$e_6$}
\put(-35,30){$-2$}
\put(-50,10){$e_5$}
\put(-55,30){$-3$}
\put(-130,10){$e_1$}
\put(-135,30){$-2$}
\put(-113,10){$e_2$}
\put(-115,30){$-2$}
\put(-91,10){$e_3$}
\put(-95,30){$-2$}
\put(-70,10){$e_4$}
\put(-75,30){$-3$}
\caption{Incidence graph.}
\label{10_19b}
\end{subfigure}
\caption{Case of $K=10_{19}$}
\end{figure}\
\newpage
\textbf{Case of $K=10_{19}$ }by using the checkerboard coloring method as shown in Figure \ref{10_19a} we can find the negative definite Goeritz matrix G such as :
\\
$$G = \begin{bmatrix}
-2 & 1 & 0 & 0 & 0 & 0 & 0 \\
1 & -2 & 1 & 0 & 0 & 0 & 0 \\
0 & 1 & -2 &1 & 0 & 0 & 0  \\
0 & 0 &1 & -3 & 1 & 0 &0 \\
0 & 0 & 0 & 1 & -3 & 1 & 0\\
0 & 0 & 0 & 0 & 1 & -2 & 1\\
0 & 0 & 1 & 1 & 0 & 1 & -2
\end{bmatrix}.$$
\\
Figure \ref{10_19b}. is the geometric representation of the Goertiz Matrix. Since  $\det (10_{19}) = 51$   is square-free, we seek an embedding $\phi : (\mathbb{Z}^8,G\oplus [-51])\hookrightarrow(\mathbb{Z}^8-Id)$ .If such a $\phi$ existed, assume that $\phi(e_1)=f_1+f_2$ then the only possibility for $\phi(e_i)$ for $i=2,3,4,5,6$ are as follows: 
\begin{align*}
\phi(e_2)&=-f_2+f_3 \cr
\phi(e_3)&=-f_3+f_4 \cr
\phi(e_4)&=-f_4+f_5+f_6 \cr
\phi(e_5)&=-f_5+f_7+f_8 \cr
\phi(e_6)&=f_5-f_6 \cr
\end{align*}
Let $\phi(e_7)=\sum_{i=1}^{8}\lambda_if_i$, where $\lambda_i's$ are integers. Then for $i=1,2,3,4,5$ , $e_7.e_i=0$ and $e_7.e_5=1$ give us the following equations:
\begin{align*}
-\lambda_1-\lambda_2&=0 \cr
\lambda_2-\lambda_3&=0 \cr
\lambda_4-\lambda_5-\lambda_6&=0 \cr
\lambda_5-\lambda_7-\lambda_8&=0 \cr
-\lambda_5+\lambda_6&=1 \cr
\end{align*}
From the first and second equations since $\lambda_1^2+...+\lambda_8^2=2$ we obtain that
$$-\lambda_1=\lambda_2=\lambda_3=\lambda_4=0$$
Therefore by using the forth equation we find that $\lambda_5=-\lambda_6$. However in the last equation we have $-\lambda_5+\lambda_6=1$. Thus $\lambda_6=\frac{1}{2}$ which is a contradiction, since $\lambda_6$ is an integer. Hence we conclude that the embedding $\phi$ does not exist, therefore that $\gamma_4(10_{19})\geq2$. As we have shown in Figure \ref{10_19a}, there is a non-oriented band from $10_{19}$ to $6_2$, and since $\gamma_4(6_2)=1$, we can conclude that $\gamma_4(10_{19})=2$.\
\begin{figure}[h]
\centering
\begin{subfigure}[b]{0.30\textwidth}
\includegraphics[width=\textwidth]{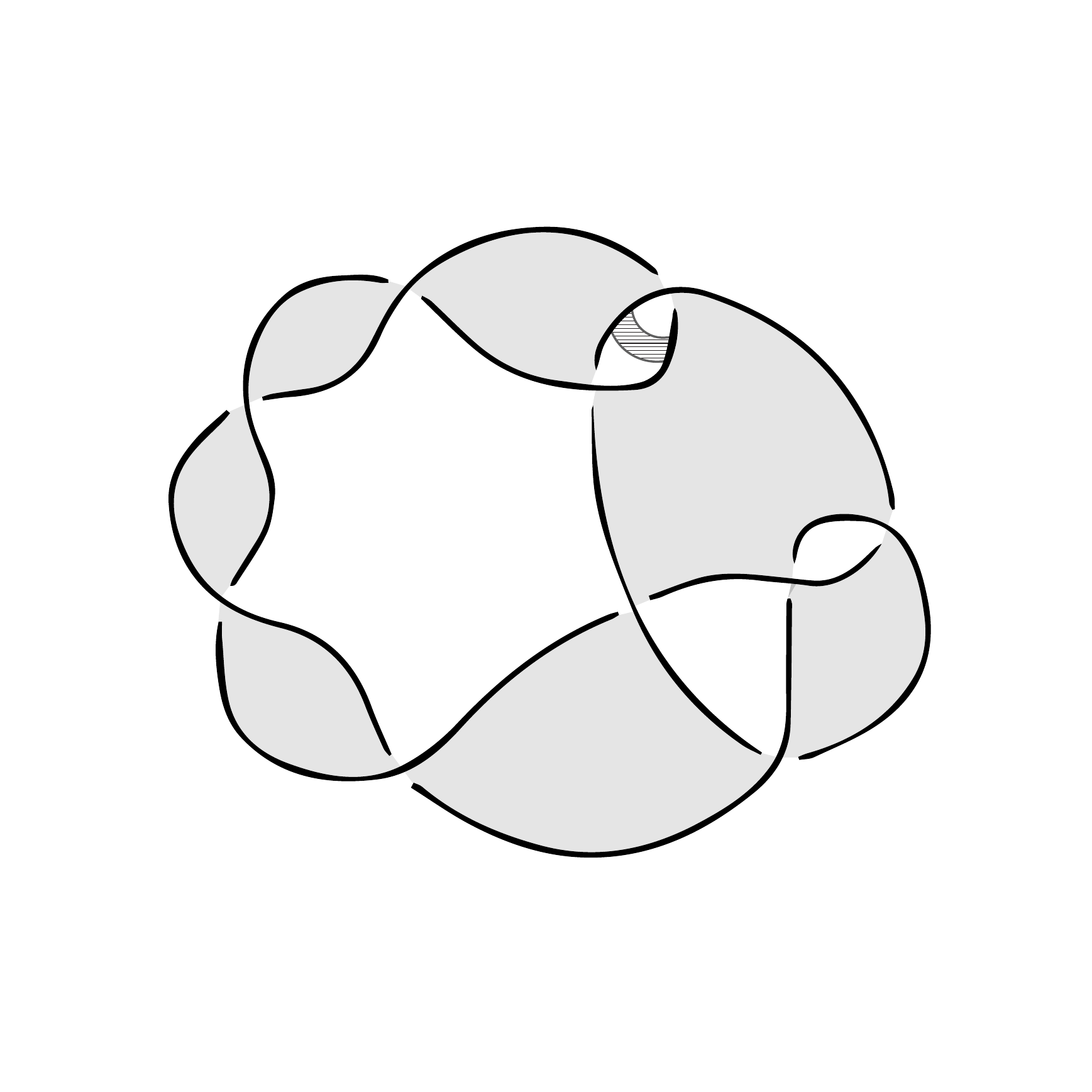}
\caption{Checkerboard diagram for $K=10_{36}\stackrel{0}{\longrightarrow} 8_{7}$} \label{10-36a}
\end{subfigure}
\qquad \qquad \qquad \qquad
\begin{subfigure}[b]{0.30\textwidth}
\includegraphics[width=\textwidth]{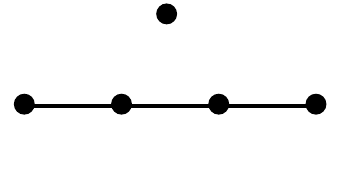}
\put(-130,10){$e_1$}
\put(-130,35){$-2$}
\put(-90,10){$e_4$}
\put(-95,35){$-6$}
\put(-50,10){$e_3$}
\put(-60,35){$-3$}
\put(-10,10){$e_2$}
\put(-20,35){$-2$}
\put(-70,50){$e_5$}
\put(-80,70){$-51$}
\caption{Incidence graph.}
\label{10-36b}
\end{subfigure}
\caption{Case of $K=10_{36}$}
\end{figure}
\\
\par
\textbf{Case of $K=10_{36}$} by using the checkerboard coloring method as shown in Figure \ref{10-36a} we can find the negative definite Goeritz matrix G such as :
\\
$$G = \begin{bmatrix}
-2 & 0 & 0 & 1 \\
0 & -2 & 1 & 0 \\
0 & 1 & -3& 1 \\
1 & 0 & 1 & -6\\
\end{bmatrix}.$$
\\
Figure \ref{10-36b}. is the geometric representation of the Goertiz Matrix. Since  $\det (10_{36}) = 51$   is square-free, we seek an embedding $\phi : (\mathbb{Z}^5,G\oplus [-51])\hookrightarrow(\mathbb{Z}^5, -Id)$. If such a $\phi$ existed, suppose $\phi(e_1)=f_1+f_2$, then up to isomorphism the only possibilities for $\phi(e_2)$ and $\phi(e_3)$ are as follows
\begin{align*}
\phi(e_2) & =f_3+f_4 \cr
\phi(e_3) & =f_1-f_2-f_3 \cr
\end{align*}
Let $\phi(e_5)=\sum_{i=1}^{5}\lambda_if_i$, where $\lambda_i's$ are integers. Then for $i=1,2,3$ , $e_5.e_i=0$ gives us the following equations
\begin{align*}
-\lambda_1-\lambda_2& =0\cr
-\lambda_3-\lambda_4&=0 \cr
-\lambda_1+\lambda_2+\lambda_3&=0 \cr
\end{align*}
Solving $\lambda_1$ , $\lambda_3$ and $\lambda_4$ in terms of $\lambda_2$, we obtain 
$$\lambda_1=-\lambda_2, \quad \lambda_3=-\lambda_4=2\lambda_2$$
Since $\lambda_1^2+\lambda_2^2+\lambda_3^2+\lambda_4^2+\lambda_5^2=51$ we find that $10\lambda_2^2+\lambda_5^2=51$. Thus $|\lambda _2| \le 2$ and neither of the 5 possibilities of $\lambda _1\in \{0, \pm 1, \pm 2\}$ leads to an integral solution of $\lambda _5$, showing that the case of $\phi ( e_1) = f_1+f_2+f_3$ cannot occur.
This shows that the embedding $\phi$ does not exist and therefore that $\gamma_4(10_{36})\geq2$. As we have shown in Figure \ref{10-36a}, there is a non-oriented band from $10_{36}$ to $8_{7}$, and since $\gamma_4(8_{7})=1$, we can conclude that $\gamma_4(10_{36})=2$.\ 
\begin{figure}[h]
\centering
\begin{subfigure}[b]{0.30\textwidth}
\includegraphics[width=\textwidth]{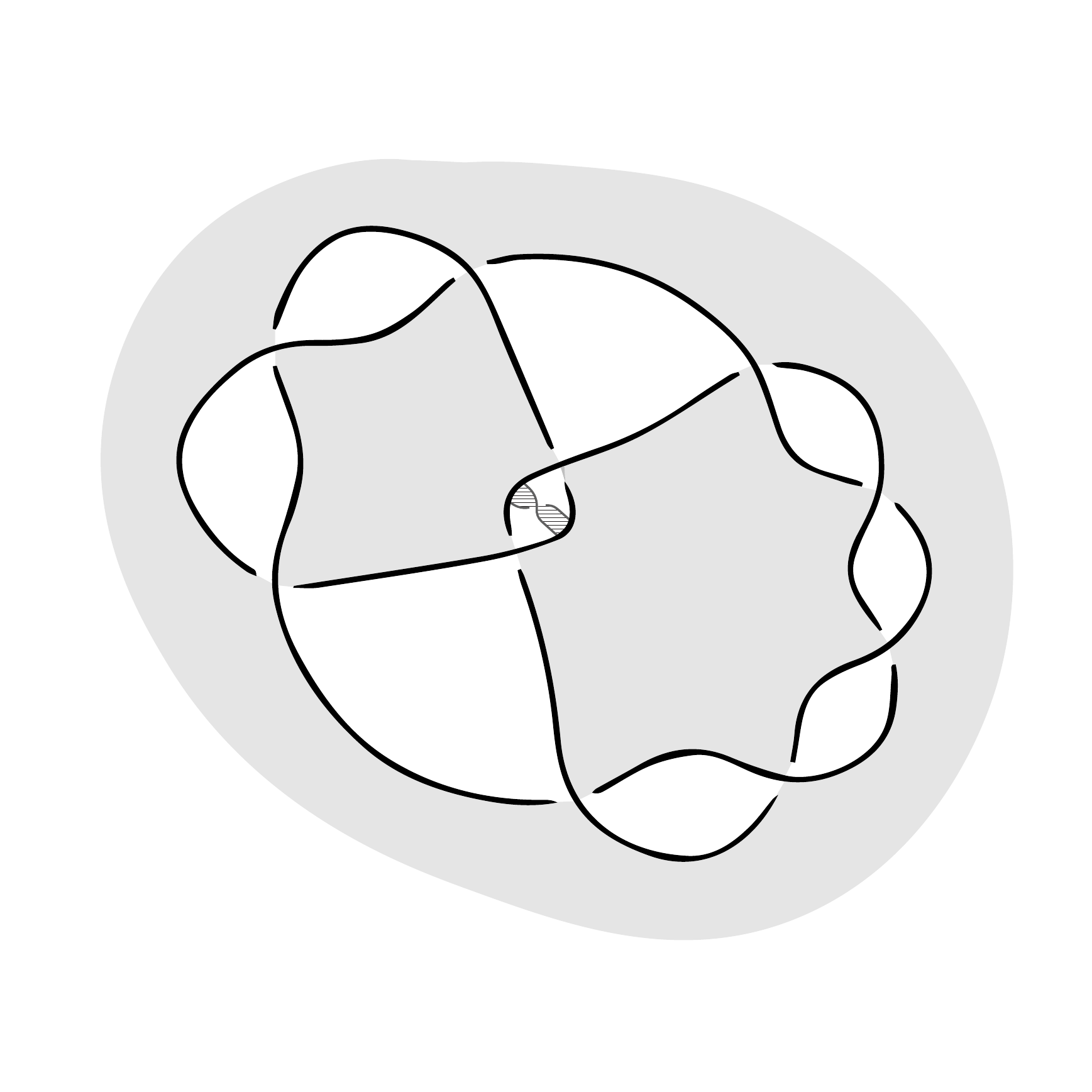}
\caption{Checkerboard diagram for $K=10_{46}\stackrel{1}{\longrightarrow} 9_{5}$} \label{10_46a}
\end{subfigure}
\qquad \qquad \qquad \qquad
\begin{subfigure}[b]{0.40\textwidth}
\includegraphics[width=\textwidth]{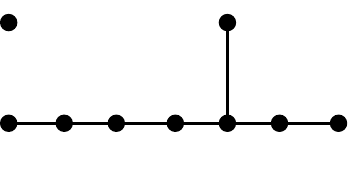}
\put(-180,12){$e_4$}
\put(-185,32){$-2$}
\put(-150,12){$e_3$}
\put(-160,32){$-2$}
\put(-120,12){$e_2$}
\put(-128,32){$-2$}
\put(-90,12){$e_1$}
\put(-98,32){$-2$}
\put(-60,12){$e_6$}
\put(-60,32){$-3$}
\put(-30,12){$e_7$}
\put(-34,32){$-2$}
\put(0,12){$e_8$}
\put(-5,32){$-2$}
\put(-80,75){$e_5$}
\put(-55,72){$-2$}
\put(-190,75){$e_9$}
\put(-165,72){$-31$}
\caption{Incidence graph.}
\label{10_46b}
\end{subfigure}
\caption{Case of $K=10_{46}$}
\end{figure}
\par
\textbf{Case of $K=10_{46}$} by using the checkerboard coloring method as shown in Figure \ref{10_46a} we can find the negative definite Goeritz matrix G such as :
\\
$$G = \begin{bmatrix}
-2 & 1 & 0 & 0 & 0 & 0 & 0 & 0\\
1 & -2 & 1 & 0 & 0 & 0 & 0 & 0 \\
0 & 1 & -2 & 1 & 0 & 0 & 0 & 0 \\
0 & 0 & 1 & -2 & 0 & 0 & 0 & 0 \\
0 & 0 & 0 & 0 & -2 & 1 & 0 & 0\\
0 & 0 & 0 & 0 & 1 & -3 & 1 & 0\\
0 & 0 & 0 & 0 & 0 & 1 & -2 & 1\\
0 & 0 & 0 & 0 &0 & 0 & 1 & -2\\
\end{bmatrix}.$$
\\
Figure \ref{10_46b}. is the geometric representation of the Goertiz Matrix. Since  $\det (10_{46}) = 31$   is square-free, we seek an embedding $\phi : (\mathbb{Z}^9,G\oplus [-87])\hookrightarrow(\mathbb{Z}^9, -Id)$. Suppose that $\phi(e_5)=f_1+f_2$ then up to isomorphism the only possibilities for $\phi(e_i), i=2,...,6$ are as follows
\begin{align*}
\phi(e_6)&=-f_1+f_3+f_4 \cr
\phi(e_1)&=-f_3+f_5 \cr
\phi(e_2)&=-f_5+f_6 \cr
\phi(e_3)&=-f_6+f_7 \cr
\phi(e_4)&=-f_7+f_8 \cr
\end{align*}
But we have two cases for $\phi(e_7)$
\begin{enumerate}
\item $\phi(e_7)=f_1-f_2$
\\
Let $\phi(e_8)=\sum_{i=1}^{9}\lambda_if_i$, where $\lambda_i's$ are integers. Then for $i=1,2,3,4,6,7$ , $e_8.e_i=0$ and $e_8.e_5=1$ lead to the following equations
\begin{align*}
-\lambda_1-\lambda_2&=1 \cr
\lambda_1-\lambda_3-\lambda_4&=0 \cr
\lambda_3-\lambda_5&=0 \cr
\lambda_5-\lambda_6&=0 \cr
\lambda_6-\lambda_7&=0 \cr
\lambda_7-\lambda_8&=0 \cr
-\lambda_1+\lambda_2&=0 \cr
\end{align*}
By adding the first and last equations we obtain$$\lambda_1=-\frac{1}{2}$$ Which is contradiction, since $\lambda_1$ is an integer. Thus the case of $\phi(e_7)=f_1-f_2$ cannot occur. 
\item $\phi(e_7)=-f_4+f_9$
\\
Let $\phi(e_8)=\sum_{i=1}^{9}\lambda_if_i$, where $\lambda_i's$ are integers. Then for $i=1,2,3,4,6,7$ , $e_8.e_i=0$ and $e_8.e_5=1$ lead to the following equations
\begin{align*}
-\lambda_1-\lambda_2&=1 \cr
\lambda_1-\lambda_3-\lambda_4&=0 \cr
\lambda_3-\lambda_5&=0 \cr
\lambda_5-\lambda_6&=0 \cr
\lambda_6-\lambda_7&=0 \cr
\lambda_7-\lambda_8&=0 \cr
\lambda_4-\lambda_9&=0 \cr
\end{align*}
By solving the equations we obtain 
$$\lambda_3=\lambda_5=\lambda_6=\lambda_7=\lambda_8=0$$
Since $\lambda_1^2+...+\lambda_9^2=2$.
And also
$$\lambda_1=-\lambda_2=\lambda_4, \quad \lambda_9=\lambda_4-1$$
Since $\lambda_1^2+...+\lambda_9^2=2$ we find that 
$$3\lambda_4^2+(\lambda_4-1)^2=2 \Longrightarrow 4\lambda_4^2+2\lambda_4-1=0$$
But this equation doesn't have any integral solution for $\lambda_4$ and the case $\phi(e_7)=-f_4+f_9$ cannot occur. 
\end{enumerate}
These two cases show that the embedding $\phi$ does not exist and therefore that $\gamma_4(10_{46})\geq2$. As we have shown in Figure \ref{10_46a}, there is a non-oriented band from $10_{46}$ to $9_{5}$, and since $\gamma_4(9_{5})=1$, we can conclude that $\gamma_4(10_{46})=2$.\ 
\begin{figure}[h]
\centering
\begin{subfigure}[b]{0.30\textwidth}
\includegraphics[width=\textwidth]{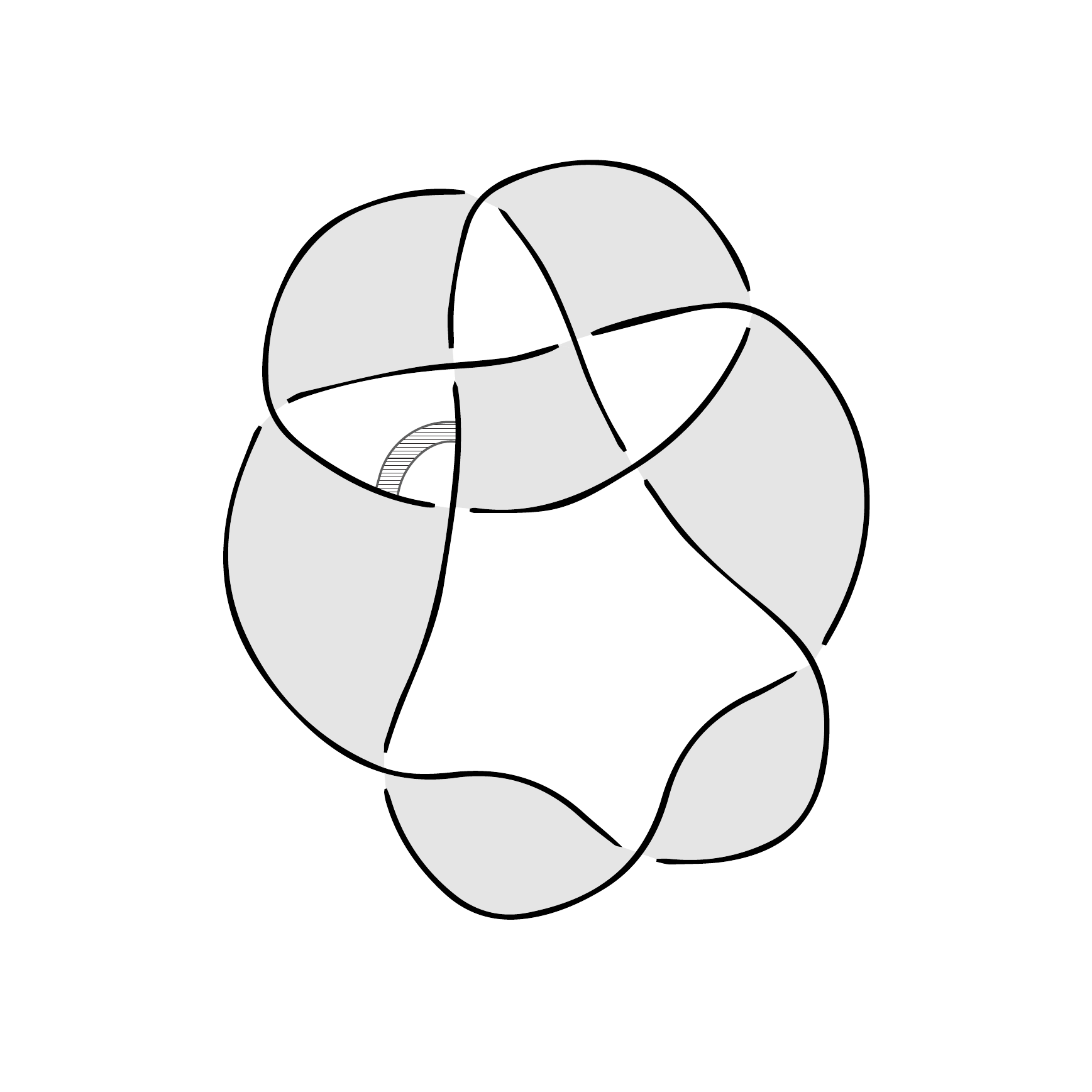}
\caption{Checkerboard diagram for $K=10_{112}\stackrel{0}{\longrightarrow} 9_{26}$} \label{10-112a}
\end{subfigure}
\qquad \qquad \qquad \qquad
\begin{subfigure}[b]{0.30\textwidth}
\includegraphics[width=\textwidth]{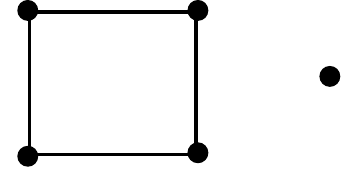}
\put(-60,-4){$e_2$}
\put(-50,5){$-3$}
\put(-125,-4){$e_1$}
\put(-145,5){$-3$}
\put(-125,70){$e_3$}
\put(-145,60){$-5$}
\put(-60,70){$e_4$}
\put(-50,60){$-3$}
\put(0,42){$e_5$}
\put(-4,25){$-87$}
\caption{Incidence graph.}
\label{10-112b}
\end{subfigure}
\caption{Case of $K=10_{112}$}
\end{figure}
\newpage
\textbf{Case of $K=10_{112}$} by using the checkerboard coloring method as shown in Figure \ref{10-112a} we can find the negative definite Goeritz matrix G such as :
\\
$$G = \begin{bmatrix}
-3 & 1 & 1 & 0 \\
1 & -3 & 0 & 1 \\
1 & 0 & -3& 1 \\
0 & 1 & 1 & -5\\
\end{bmatrix}.$$
\\
Figure \ref{10-112b}. is the geometric representation of the Goertiz Matrix. Since  $\det (10_{112}) = 87$   is square-free, we seek an embedding $\phi : (\mathbb{Z}^5,G\oplus [-87])\hookrightarrow(\mathbb{Z}^5, -Id)$. If such a $\phi$ existed, assume that $\phi(e_1)=f_1+f_2+f_3$ then there would be two cases for $\phi(e_2)$:
\begin{enumerate}
\item $\phi(e_2)=-f_3+f_4+f_5$ . 
\\
In this case the only possibility for $\phi(e_4)$ is $-f_1+f_2-f_4$.
Let $\phi(e_3)=\sum_{i=1}^{5}\lambda_if_i$, where $\lambda_i's$ are integers. Then for $i=1,4$ , $e_3.e_i=1$ and $e_3.e_2=0$ give us the following equations:
\begin{align*}
-\lambda_1-\lambda_2-\lambda_3&=1 \cr
\lambda_3-\lambda_4-\lambda_5&=0 \cr
\lambda_1-\lambda_2+\lambda_4&=1 \cr
\end{align*}
By adding all the equations we have
$$\lambda_5=-2-2\lambda_2 $$ It implies that $\lambda_5$ is an even number and since $e_3.e_3=3$, therefore $|\lambda _i| \le 1$ thus $\lambda_5=0$ and $\lambda_2=-1$ and we find that $\lambda_3=\lambda_4=-\lambda_1$.
Since $\lambda_1^2+\lambda_2^2+\lambda_3^2+\lambda_4^2+\lambda_5^2=5$ we find that $\lambda_1^2=\frac{2}{3}$ but this is a contradiction because $\lambda_1$ is an integer.It shows that the case of $\phi (e_3) =-f_3+f_4+f_5$ cannot occur. 
\item $\phi(e_2)= -f_1-f_2+f_3$.
\\
Let $\phi(e_4)=\sum_{i=1}^{5}\lambda_if_i$,where $\lambda_i's$ are integers. Then $e_1.e_4=0$ and $e_2.e_4=1$ give us the following equations:
\begin{align*}
-\lambda_1-\lambda_2-\lambda_3&=0 \cr
\lambda_1+\lambda_2-\lambda_3&=1 \cr
\end{align*}
Solving for $\lambda_3$, we obtain $\lambda_3=-\frac{1}{2}$ which is a contradiction because $\lambda_3$ is an integer. It shows that the case $\phi(e_2)=-f_1-f_2+f_3$ cannot occur.
\end{enumerate}
These two cases show that the embedding $\phi$ does not exists and therefore that $\gamma_4(10_{112})\geq2$. As we have shown in Figure \ref{10-112a}, there is a non-oriented band from $10_{112}$ to $9_{26}$, and since $\gamma_4(9_{26})=1$, we can conclude that $\gamma_4(10_{112})=2$.\
\subsection{ Knots with $\sigma(K)+4\cdot \text{Arf}(K)\equiv{0} \pmod 8$ } \label{fourth group} If knot $K$ satisfies the congruance $\sigma(K)+4\cdot \text{Arf}(K)\equiv{0} \pmod 8$ we need to consider $K$ and $-K$ (mirror of the knot $K$). Among all knots in this group we have only 2 knots with $\gamma_4$ equal 2 :
$$10_{33}, \,10_{58}$$
\begin{figure}[h]
\centering
\begin{subfigure}[b]{0.30\textwidth}
\includegraphics[width=\textwidth]{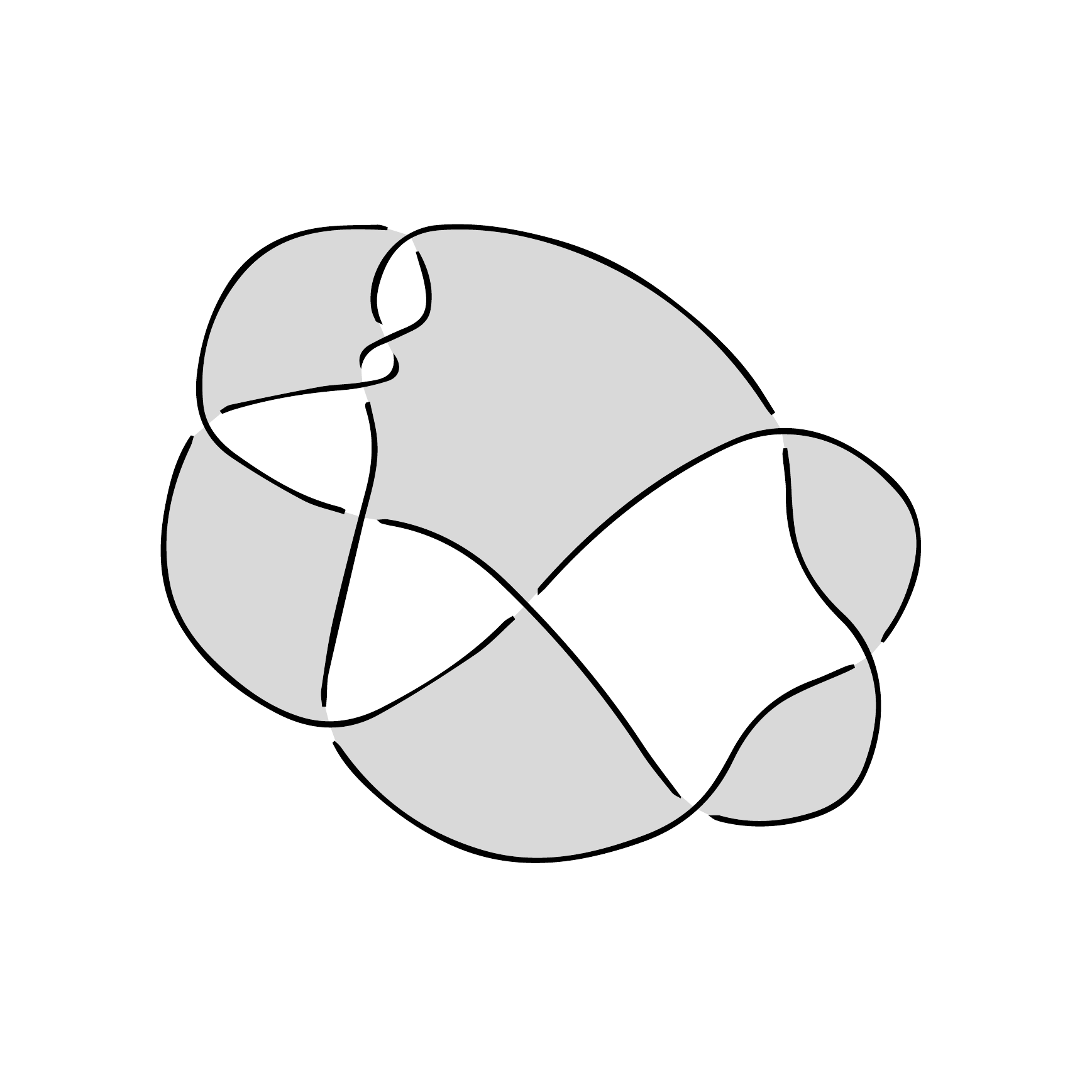}
\caption{Checkerboard diagram for $K=-10_{33}$} \label{-10-33a}
\end{subfigure}
\qquad \qquad \qquad \qquad
\begin{subfigure}[b]{0.40\textwidth}
\includegraphics[width=\textwidth]{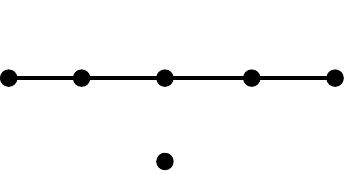}
\put(-110,5){$e_6$}
\put(-85,3){$-65$}
\put(-175,35){$e_1$}
\put(-180,60){$-2$}
\put(-140,35){$e_2$}
\put(-140,60){$-2$}
\put(-100,35){$e_3$}
\put(-100,60){$-3$}
\put(-55,35){$e_4$}
\put(-60,60){$-3$}
\put(-10,35){$e_5$}
\put(-15,60){$-4$}
\caption{Incidence graph.}
\label{-10-33b}
\end{subfigure}
\caption{Case of $K=-10_{33}$}
\end{figure}\
\par
\textbf{Case of $K=-10_{33}$} by using the checkerboard coloring method as shown in Figure \ref{-10-33a} we can find the negative definite Goeritz matrix G such as :
\\
$$G = \begin{bmatrix}
-2 & 1 & 0 & 0 &0\\
1 & -2 & 1 & 0 & 0 \\
0 & 1 & -3 & 1 &0 \\
0 & 0 & 1 & -3 &1 \\
0 & 0 & 0 & 1 & -4\\
\end{bmatrix}.$$
\\
Figure \ref{-10-33b}. is the geometric representation of the Goertiz Matrix. Since  $\det (-10_{33}) = 65$   is square-free, we seek an embedding $\phi : (\mathbb{Z}^5,G\oplus [65])\hookrightarrow(\mathbb{Z}^5, -Id)$. If such a $\phi$ existed, assume that $\phi(e_3)=f_1+f_2+f_3$ then up to isomorphism, we have $\phi(e_2)=-f_3+f_4$ then to write $\phi(e_1)$ we have two possibilities,
\begin{enumerate}
\item $\phi(e_1)=-f_1+f_3$ then $\phi(e_4)=-f_2+f_5+f_6$. Let $\phi(e_5)=\sum_{i=1}^{6}\lambda_if_i$, where $\lambda_i's$ are integers.Then $e_5.e_4=1$ and $e_5.e_i=0$ if $i=1,2,3$ lead to the following equations
\begin{align*}
-\lambda_1-\lambda_2-\lambda_3&=0 \cr
\lambda_3-\lambda_4&=0 \cr
\lambda_1-\lambda_3&=0 \cr
\lambda_2-\lambda_5-\lambda_6&=1 \cr
\end{align*}
From the second and third equations we have $\lambda_3=\lambda_4=\lambda_1$ and they must be nonzero, otherwise $\lambda_1^2+...+\lambda_6^2\leq3$ (as we have shown before $|\lambda_i's|\leq1$). From the first equation we have $\lambda_2=-2\lambda_1\neq0$, therefore $|\lambda_2|\geq2$ which is a contradiction.

\item $\phi(e_1)=-f_4+f_5$ then $\phi(e_4)=-f_3-f_4-f_5$.  Let $\phi(e_5)=\sum_{i=1}^{6}\lambda_if_i$, where $\lambda_i's$ are integers.Then $e_5.e_i=0$ if $i=1,2,3$ and $e_5.e_4=1$ lead to the following equations
\begin{align*}
-\lambda_1-\lambda_2-\lambda_3&=0 \cr
\lambda_3-\lambda_4&=0 \cr
\lambda_4-\lambda_5&=0 \cr
\lambda_3+\lambda_4+\lambda_5&=1 \cr
\end{align*}
From the second and third equation, $\lambda_3=\lambda_4=\lambda_5$. From the last equation $3\lambda_3=1$ which is not possible because $\lambda_3$ is an integer.
\end{enumerate}
These two cases show that the embedding $\phi$ does not exist and therefore that $\gamma_4(-10_{33})\geq2$. \par
Now we need to check if there exists an embedding for the knot $10_{33}$ \\
\begin{figure}[h]
\centering
\begin{subfigure}[b]{0.30\textwidth}
\includegraphics[width=\textwidth]{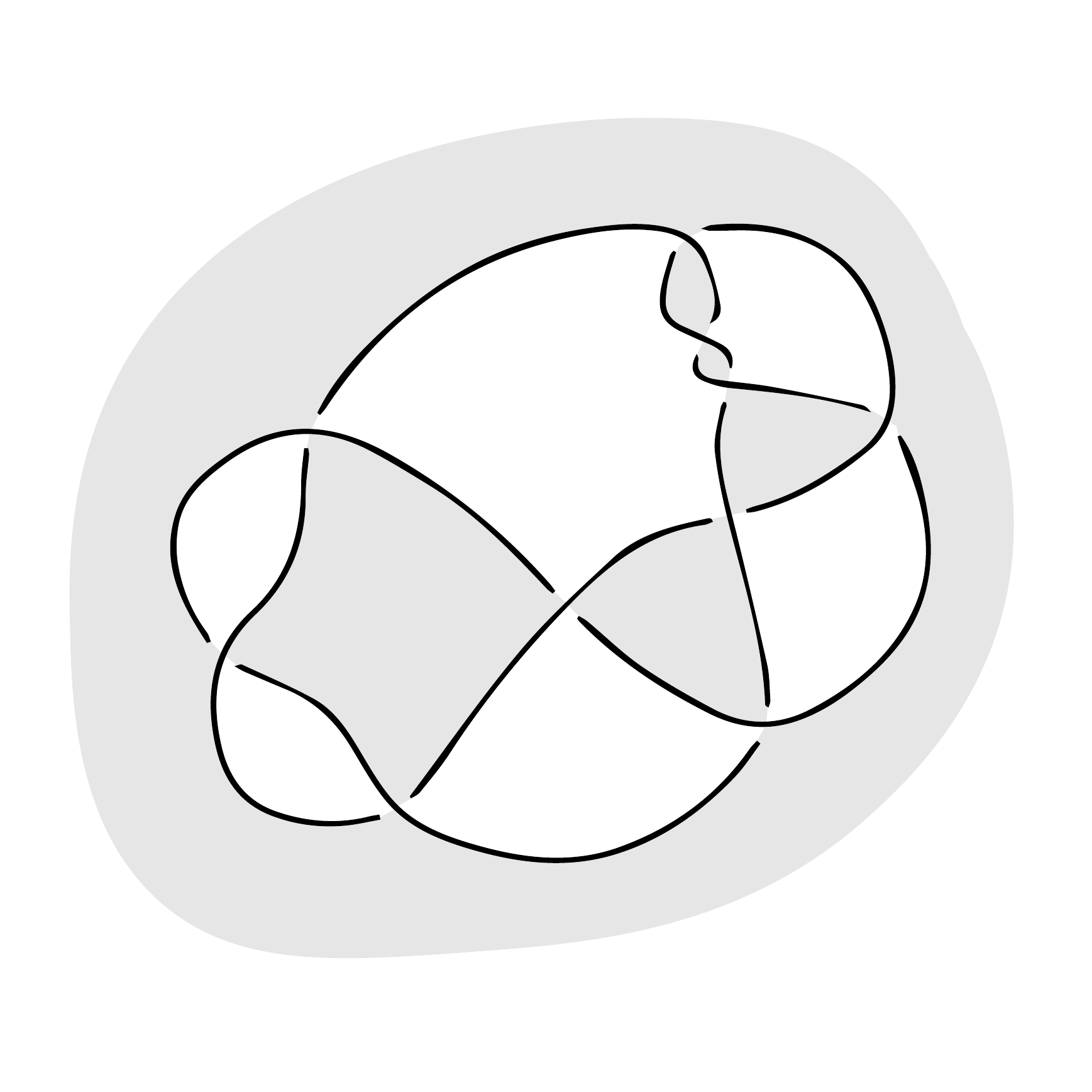}
\caption{Checkerboard diagram for $K=+10_{33}$} \label{+10-33a}
\end{subfigure}
\qquad \qquad \qquad \qquad
\begin{subfigure}[b]{0.35\textwidth}
\includegraphics[width=\textwidth]{./KnotDiagram/mirror10-33G}
\put(-100,5){$e_6$}
\put(-76,3){$-65$}
\put(-155,30){$e_1$}
\put(-160,50){$-2$}
\put(-120,30){$e_2$}
\put(-130,50){$-2$}
\put(-85,30){$e_3$}
\put(-95,50){$-3$}
\put(-45,30){$e_4$}
\put(-55,50){$-3$}
\put(-10,30){$e_5$}
\put(-20,50){$-4$}
\caption{Incidence graph.}
\label{+10-33b}
\end{subfigure}
\caption{Case of $K=+10_{33}$}
\end{figure}\
\par
\begin{figure}[h]
\centering
\includegraphics[width=70mm]{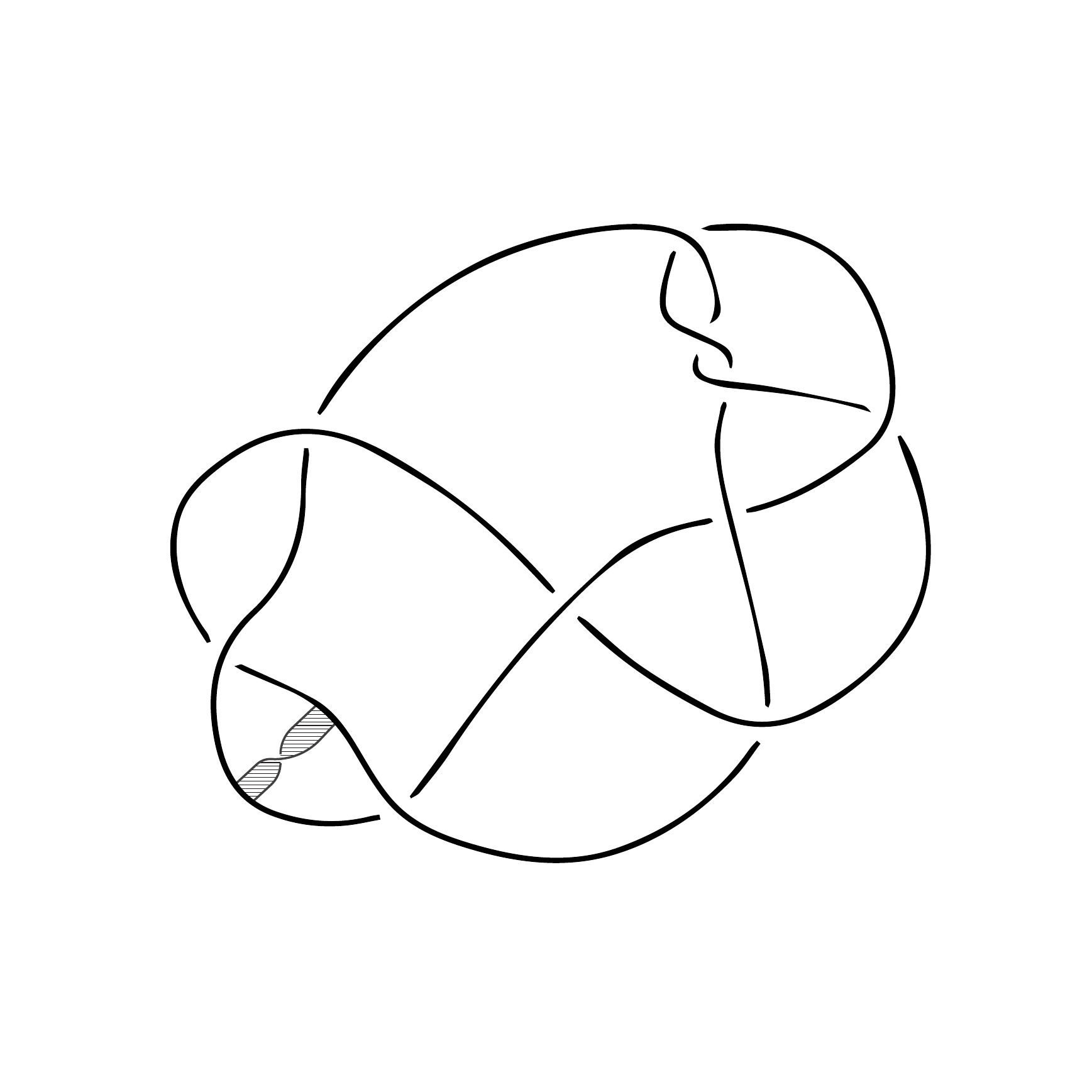}
\caption{$K=10_{33}\stackrel{1}{\longrightarrow} 9_{26}$} \label{10-33}
\end{figure}

\textbf{Case of $K=+10_{33}$} by using the checkerboard coloring method as shown in Figure \ref{+10-33a} we can find the negative definite Goeritz matrix G such as :
\\
$$G = \begin{bmatrix}
-4 & 1 & 0 & 0 &0\\
1 & -3 & 1 & 0 & 0 \\
0 & 1 & -3 & 1 &0 \\
0 & 0 & 1 & -2 &1 \\
0 & 0 & 0 & 1 & -2\\
\end{bmatrix}.$$
\\
Figure \ref{+10-33b}. is the geometric representation of the Goertiz Matrix. The incidence graph of the knot $10_{33}$ is exactly the same with the incidence graph of the knot $-10_{33}$. Thus the embedding $\phi : (\mathbb{Z}^5,G\oplus [65])\hookrightarrow(\mathbb{Z}^5, -Id)$ does not exist and therefore that $\gamma_4(10_{33})\geq2$. As we have shown in Figure \ref{10-33}, there is a non-oriented band from $10_{33}$ to $9_{26}$, and since $\gamma_4(9_{26})=1$, we can conclude that $\gamma_4(10_{33})=2$.\ 
\\
\newpage
\begin{figure}[h]
\centering
\begin{subfigure}[b]{0.30\textwidth}
\includegraphics[width=\textwidth]{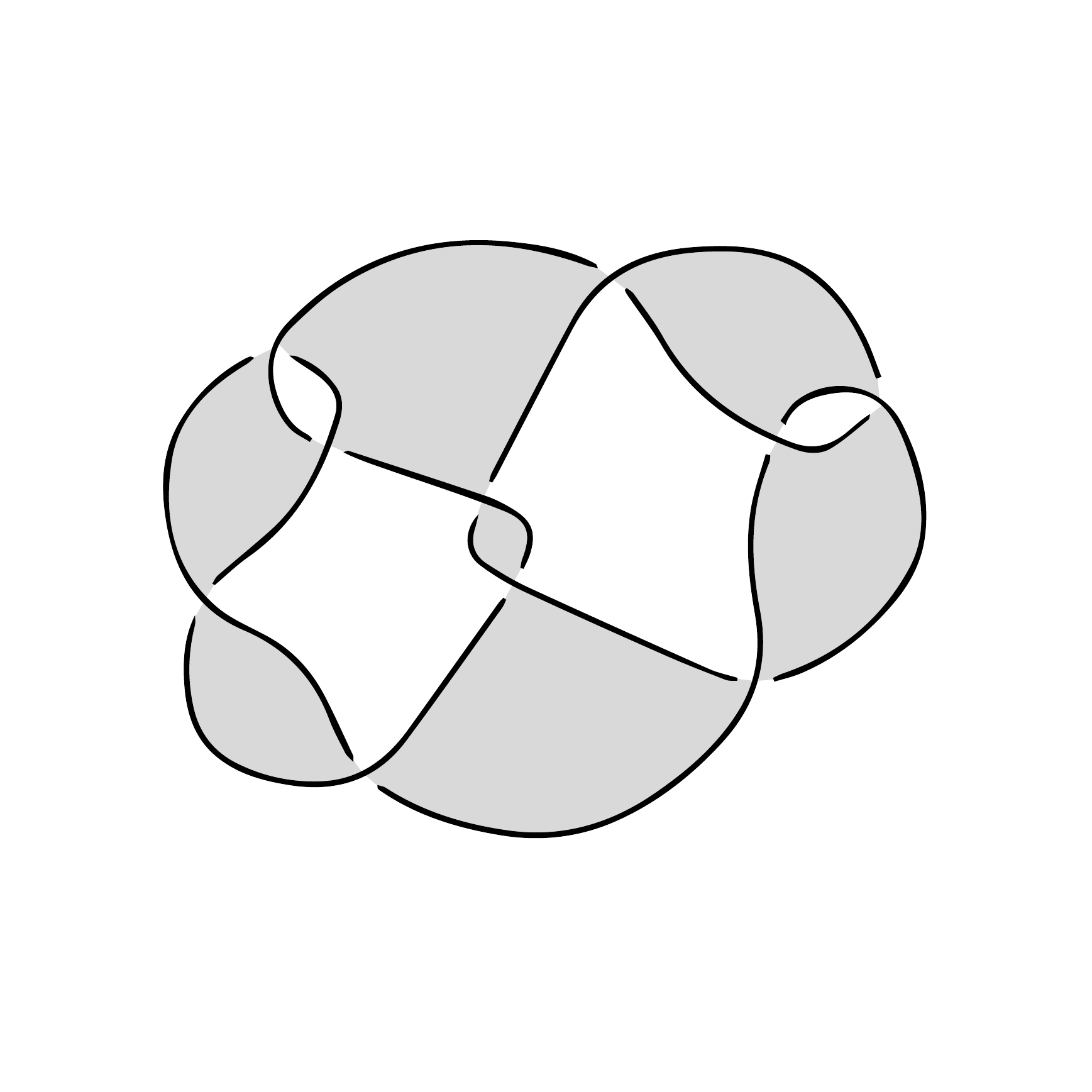}
\caption{Checkerboard diagram for $K=+10_{58}$} \label{+10-58a}
\end{subfigure}
\qquad \qquad \qquad \qquad
\begin{subfigure}[b]{0.30\textwidth}
\includegraphics[width=\textwidth]{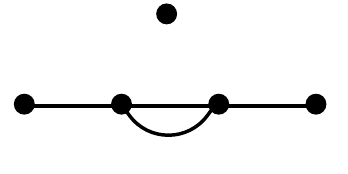}
\put(-130,10){$e_1$}
\put(-130,35){$-2$}
\put(-90,10){$e_2$}
\put(-95,35){$-5$}
\put(-50,10){$e_3$}
\put(-60,35){$-5$}
\put(-10,10){$e_4$}
\put(-20,35){$-2$}
\put(-70,50){$e_5$}
\put(-80,70){$-65$}
\caption{Incidence graph.}
\label{+10-58b}
\end{subfigure}
\caption{Case of $K=+10_{58}$}
\end{figure}\
\par
\textbf{Case of $K=+10_{58}$} by using the checkerboard coloring method as shown in Figure \ref{+10-58a} we can find the negative definite Goeritz matrix G such as :
\\
$$G = \begin{bmatrix}
-2 & 1 & 0 & 0 \\
1 & -5 & 2 & 0 \\
0 & 2 & -5 & 1 \\
0 & 0 & 1 & -2 \\
\end{bmatrix}.$$
\\
Figure \ref{+10-58b}. is the geometric representation of the Goertiz Matrix. Since  $\det (+10_{58}) = 65$   is square-free, we seek an embedding $\phi : (\mathbb{Z}^5,G\oplus [58])\hookrightarrow(\mathbb{Z}^5, -Id)$. If such a $\phi$ existed, there would be two cases for $\phi(e_2)$:
\begin{enumerate}
\item $\phi(e_2)=2f_1+f_2$ . In this case we have two possibilities for $\phi(e_1)$:
\begin{enumerate}
\item $\phi(e_1)=-f_1+f_2$, then let  $\phi(e_3)=\sum_{i=1}^{5}\lambda_if_i$, where $\lambda_i's$ are integers. Then $e_3.e_1=0$ and $e_3.e_2=2$ imply that
\begin{align*}
-2\lambda_1-\lambda_2&=2 \cr
\lambda_1-\lambda_2&=0 \cr
\end{align*}
By multiplying 2 to the second equation and add to the first equation we find that, $\lambda_2=-\frac{2}{3}$, which is impossible because $\lambda_2$ is an integer.
\item $\phi(e_1)=-f_2+f_4$ then $\phi(e_3)=2f_3-f_1$.
 Let $\phi(e_4)=\sum_{i=1}^{5}\lambda_if_i$, where $\lambda_i's$ are integers. Then for $i=1,2,3$ , $e_4.e_1=e_4.e_2=0$ and $e_4.e_3=1$ give us the following equations:
\begin{align*}
-2\lambda_1-\lambda_2& =0\cr
\lambda_2-\lambda_4&=0 \cr
-2\lambda_3+\lambda_1&=1 \cr
\end{align*}
By adding the first two equations, $\lambda_4=-2\lambda_1$ which is an even number and must be zero because $\lambda_1^2+...+\lambda_5^2=2$. Then $\lambda_1=\lambda_2=0$ and $\lambda_3=-\frac{1}{2}$ which is not possible.
\end{enumerate}
\item $\phi(e_2)= f_1+f_2+f_3+f_4+f_5$. This case  implies that $\lambda_1+\lambda_2+\lambda_3+\lambda_4+\lambda_5=0$ which is contradiction because:
\begin{align*}
65 & =\lambda_1^2+\lambda_2^2+\lambda_3^2+\lambda_4^2+\lambda_5^2 \cr
& \equiv(\lambda_1+\lambda_2+\lambda_3+\lambda_4+\lambda_5)^2 \pmod 2 \cr
&  \equiv 0 \pmod 2.
\end{align*}
\end{enumerate}
These two cases show that the embedding $\phi$ for $K=+10_{58}$ does not exists. Now we need to check $K=-10_{58}$ as follows:
\\
\begin{figure}[h]
\centering
\begin{subfigure}[b]{0.30\textwidth}
\includegraphics[width=\textwidth]{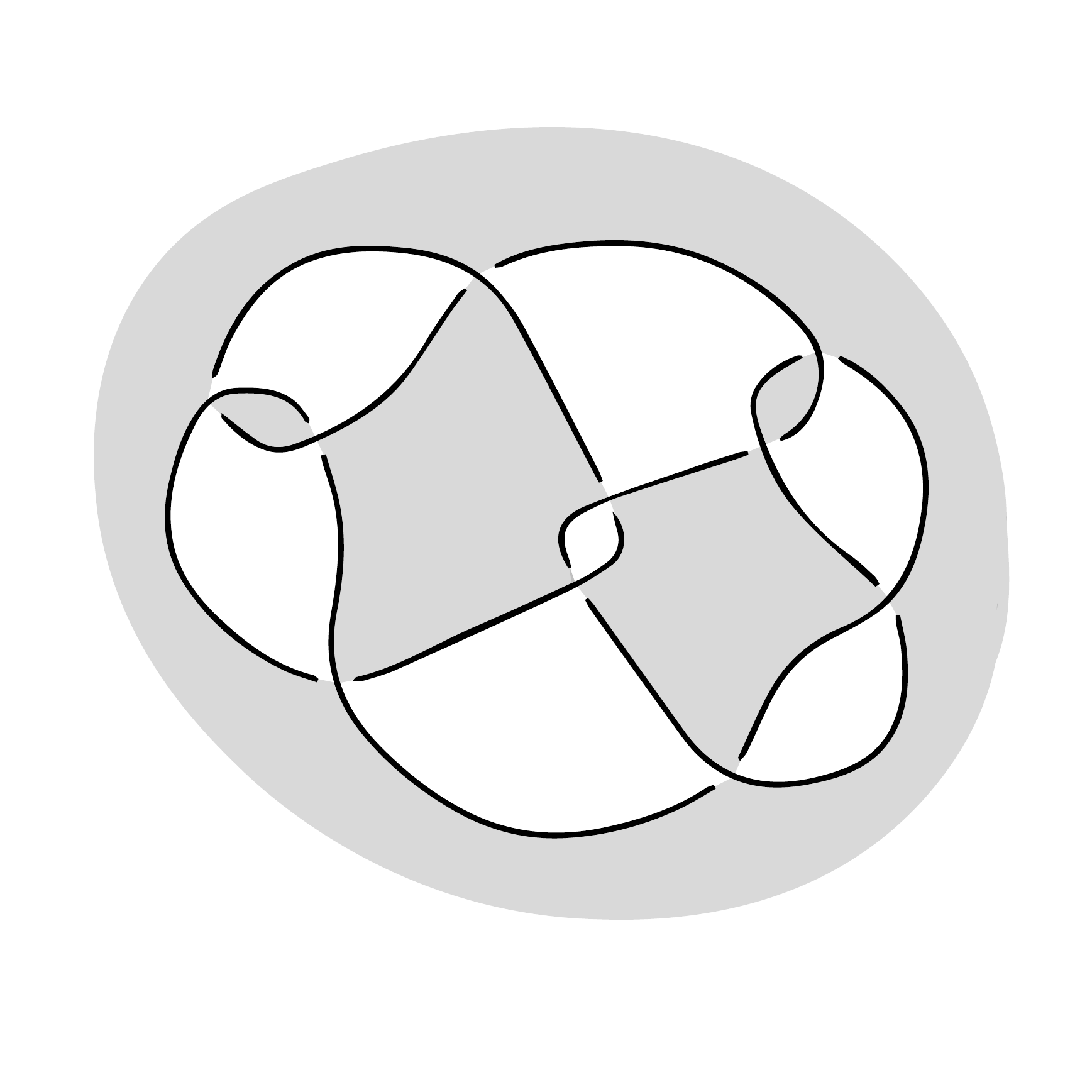}
\caption{Checkerboard diagram for $K=10_{58}$} \label{-10-58a}
\end{subfigure}
\qquad \qquad \qquad \qquad
\begin{subfigure}[b]{0.40\textwidth}
\includegraphics[width=\textwidth]{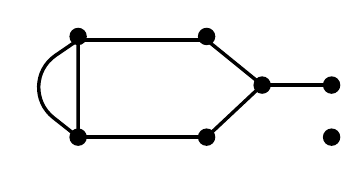}
\put(-150,10){$e_1$}
\put(-135,23){$-4$}
\put(-70,8){$e_5$}
\put(-87,24){$-2$}
\put(-20,8){$e_7$}
\put(-3,15){$-65$}
\put(-155,70){$e_2$}
\put(-135,58){$-3$}
\put(-70,75){$e_3$}
\put(-85,55){$-2$}
\put(-50,30){$e_4$}
\put(-47,50){$-3$}
\put(-20,35){$e_6$}
\put(-3,40){$-3$}
\caption{Incidence graph.}
\label{-10-58b}
\end{subfigure}
\caption{Case of $K=-10_{58}$}
\end{figure}\
\par
\begin{figure}[h]
\centering
\includegraphics[width=70mm]{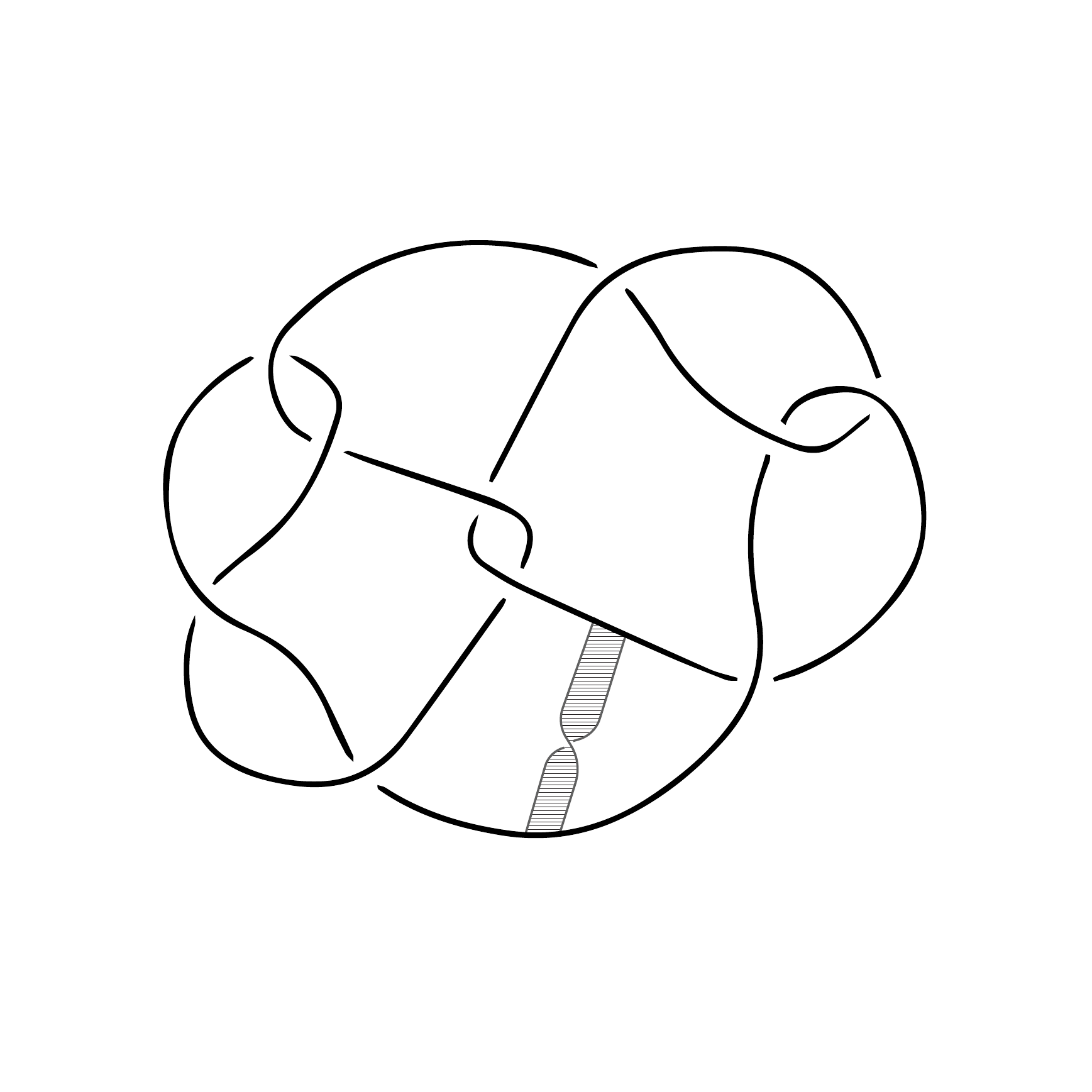}
\caption{$K=10_{58}\stackrel{-1}{\longrightarrow} 9_{25}$} \label{10-58}
\end{figure}

\textbf{Case of $K=-10_{58}$} by using the checkerboard coloring method as shown in Figure \ref{-10-58a} we can find the negative definite Goeritz matrix G such as :
\\
$$G = \begin{bmatrix}
-4 & 2 & 0 & 0 & 1 & 0 \\
2 & -3 & 1 & 0 & 0 & 0 \\
0 & 1 & -2 & 1 & 0 & 0 \\
0 & 0 & 1 & -3 & 1 & 1 \\
1 & 0 & 0 & 1 & -2 & 0 \\
0 & 0 & 0 & 1 & 0 & -3 \\
\end{bmatrix}.$$
\\
Figure \ref{-10-58b}. is the geometric representation of the Goertiz Matrix. Since  $\det (-10_{58}) = 65$   is square-free, we seek an embedding $\phi : (\mathbb{Z}^5,G\oplus [58])\hookrightarrow(\mathbb{Z}^5, -Id)$. If such a $\phi$ existed, suppose $\phi(e_2)=f_1+f_2+f_3$ then $\phi(e_3)$ must be $-f_3+f_4$. To find $\phi(e_1)$ we have two possibilities
\begin{enumerate}
\item $\phi(e_1)=-f_1-f_2+f_5+f_6$, then
\\
\begin{align*}
\phi(e_5)&=-f_5+f_7 \cr
\phi(e_4)&=-f_4+f_5-f_6 \cr
\end{align*}
Let $\phi(e_6)=\sum_{i=1}^{7}\lambda_if_i$, where $\lambda_i's$ are integers. Then for $i=1,2,3,5$ , $e_6.e_i=0$ and $e_6.e_4=1$ give us the following equations:
\begin{align*}
-\lambda_1-\lambda_2-\lambda_3&=0 \cr
\lambda_3-\lambda_4&=0 \cr
\lambda_1+\lambda_2-\lambda_5-\lambda_6&=0 \cr
\lambda_5-\lambda_7&=0 \cr
\lambda_4-\lambda_5+\lambda_6&=1 \cr
\end{align*}
Solving the equations, we obtain 
$$\lambda_3=\lambda_4, \quad \lambda_5=\lambda_7$$
Sum of the first, third and fifth equations lead us to $-\lambda_3+\lambda_4-2\lambda_5=1$, but $\lambda_3=\lambda_4$ and it implies that
 $$-2\lambda_5=1\quad \Longrightarrow \quad \lambda_5=-\frac{1}{2}$$
 But $\lambda_5$ is an integer, and this shows that case of $\phi(e_1)=-f_1-f_2+f_5+f_6$ cannot occur.
\item $\phi(e_1)=-f_1-f_3-f_4+f_5$ then
\begin{align*}
\phi(e_6)&=-f_1+f_2-f_5 \cr
\phi(e_5)&=f_6+f_7 \cr
\end{align*}
Let $\phi(e_4)=\sum_{i=1}^{7}\lambda_if_i$, where $\lambda_i's$ are integers. Then for $i=1,2,6$ , $e_4.e_i=0$, $e_4.e_5=1$ and $e_4.e_3=1$ give us the following equations:
\begin{align*}
-\lambda_1-\lambda_2-\lambda_3&=0 \cr
\lambda_3-\lambda_4&=1 \cr
\lambda_1+\lambda_3+\lambda_4-\lambda_5&=0 \cr
\lambda_1-\lambda_2+\lambda_5&=0 \cr
-\lambda_6-\lambda_7&=1 \cr
\end{align*}
By adding the first and third equations we obtain that $\lambda_5=0$, then based on the second equation there are two possibilities for $\lambda_3$:
\\
If $\lambda_3=0$ then $\lambda_1=\lambda_2=0$ and $\lambda_4=-1$. Based on the fifth equation, exactly one of $\lambda_6$ or $\lambda_7$ must be zero. Therefore only two of the $\lambda_i's$ can be nonzero, which is a contradiction because $\lambda_1^2+...+\lambda_7^2=3$.
\\
Now if $\lambda_3=1$, by adding the first and forth equations, we have $\lambda_2=-\frac{1}{2}$ which is a contradiction because $\lambda_2$ is an integer.
\end{enumerate}
These two cases show that the embedding $\phi$ does not exist and therefore that $\gamma_4(10_{58})\geq2$. As we have shown in Figure \ref{10-58}, there is a non-oriented band from $10_{58}$ to $9_{26}$, and since $\gamma_4(9_{26})=1$, we can conclude that $\gamma_4(10_{58})=2$.\ 

\newpage
\begin{figure}[h]
	\centering
	\begin{subfigure}[b]{0.3\textwidth}
		\includegraphics[width=\textwidth]{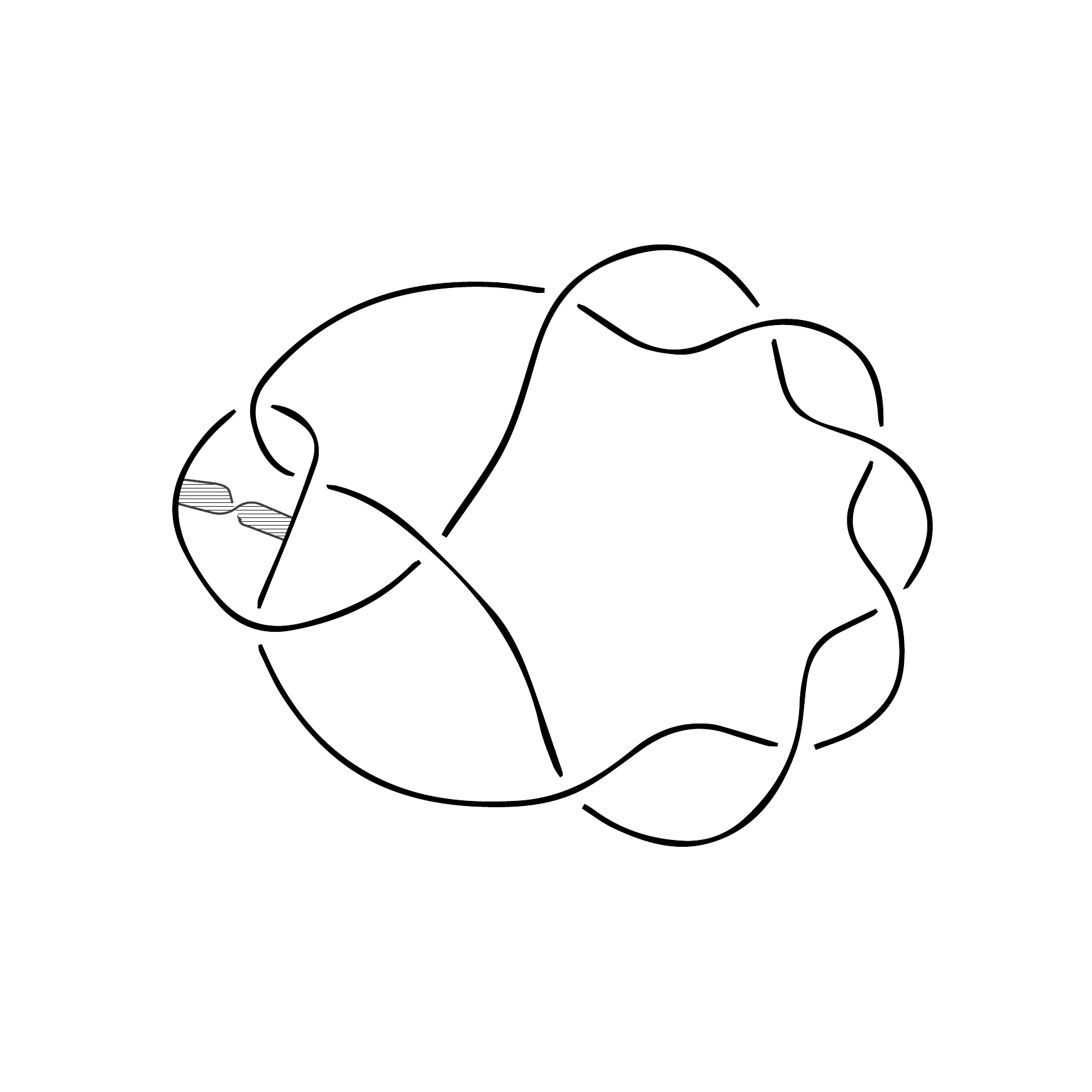}
		\caption{$10_{5}\stackrel{-1}{\longrightarrow} 9_{3}$}
		\label{FigureFor10-5}
	\end{subfigure}
~
	\begin{subfigure}[b]{0.3\textwidth}
		\includegraphics[width=\textwidth]{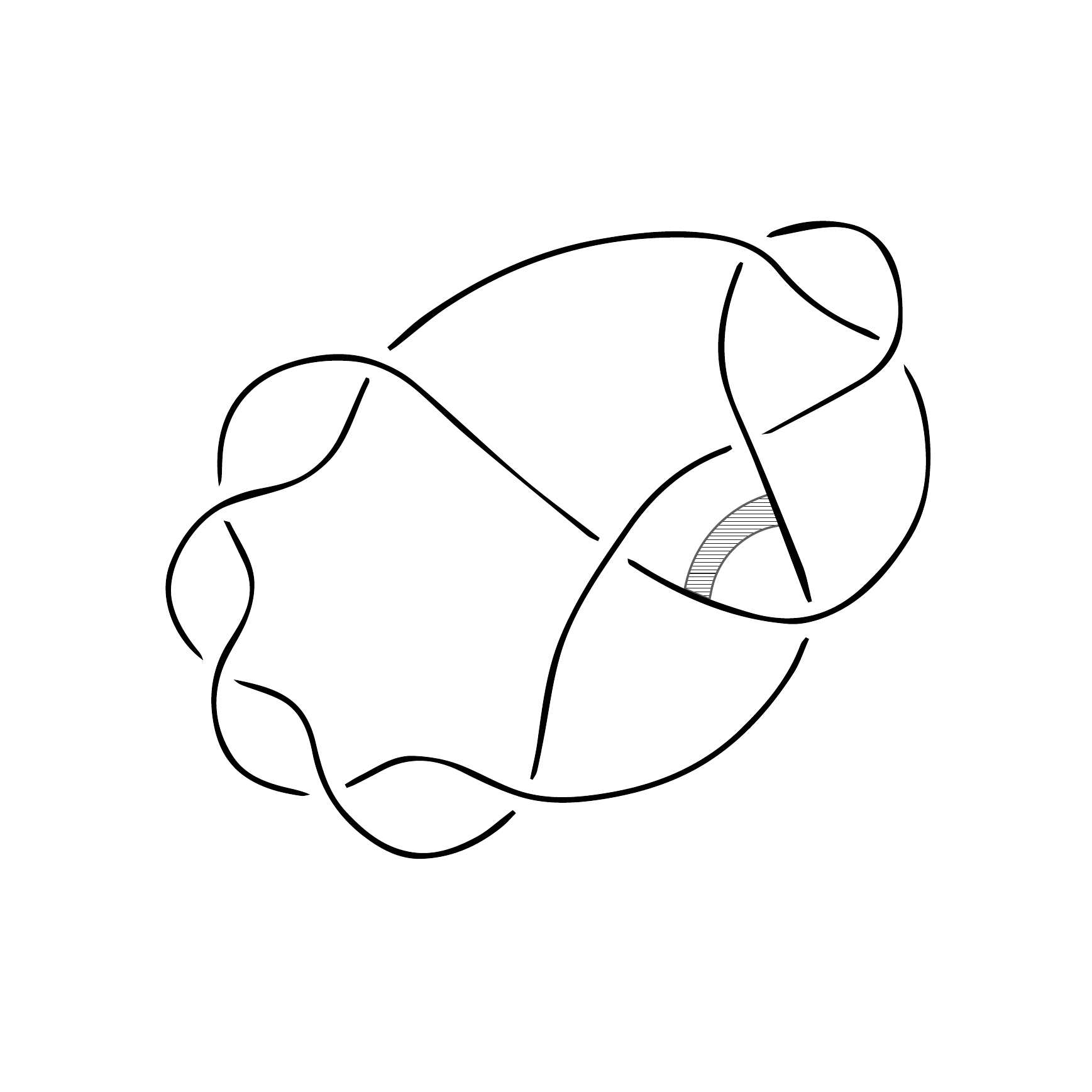}
		\caption{$10_{10}\stackrel{0}{\longrightarrow} 9_6$}
		\label{FigureFor10-10}
	\end{subfigure}
	~
	\begin{subfigure}[b]{0.3\textwidth}
		\includegraphics[width=\textwidth]{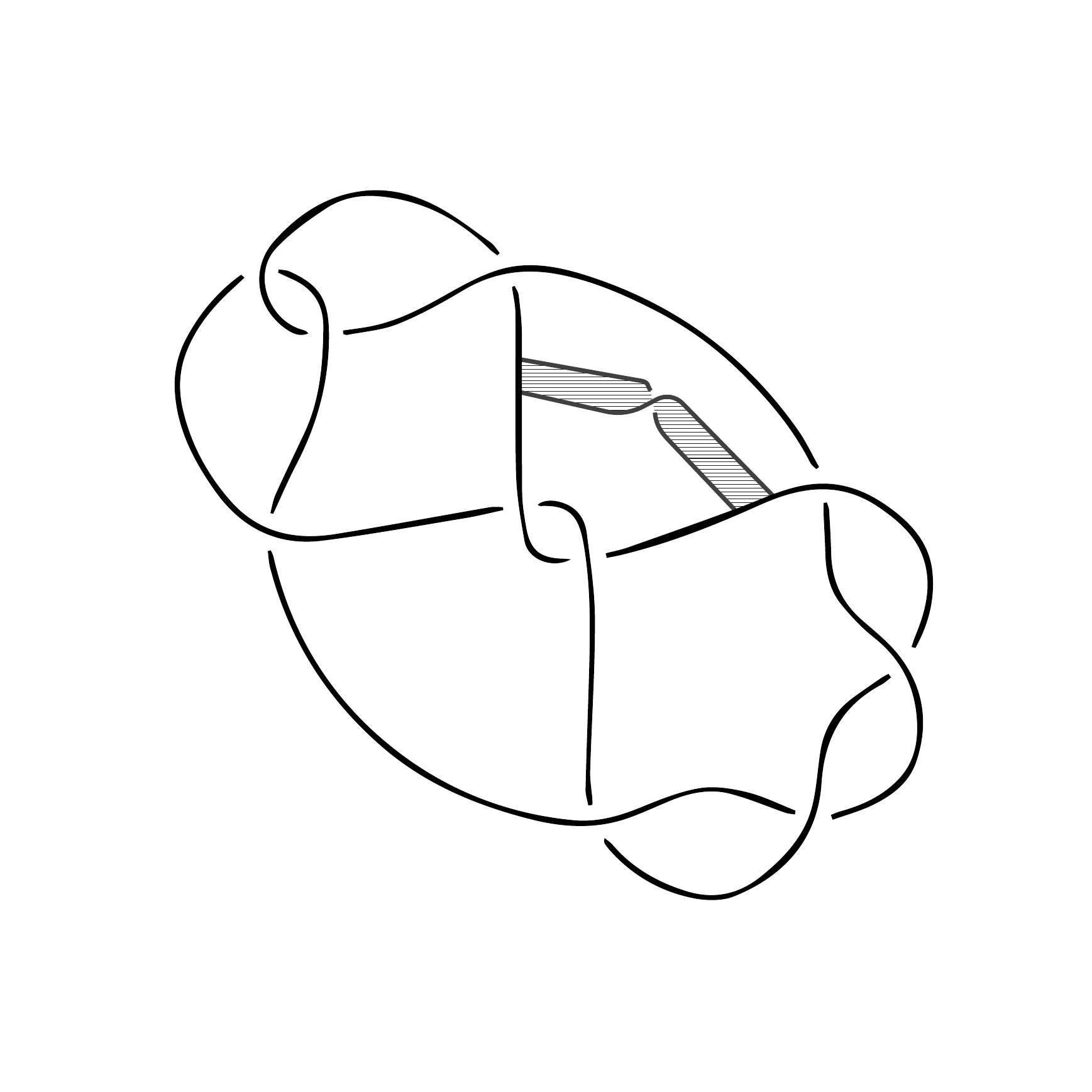}
		\caption{$10_{13}\stackrel{-1}{\longrightarrow} 8_{11}$}
		\label{FigureFor10-13}
	\end{subfigure}
	~
	\vskip3mm
	~
	\begin{subfigure}[b]{0.3\textwidth}
		\includegraphics[width=\textwidth]{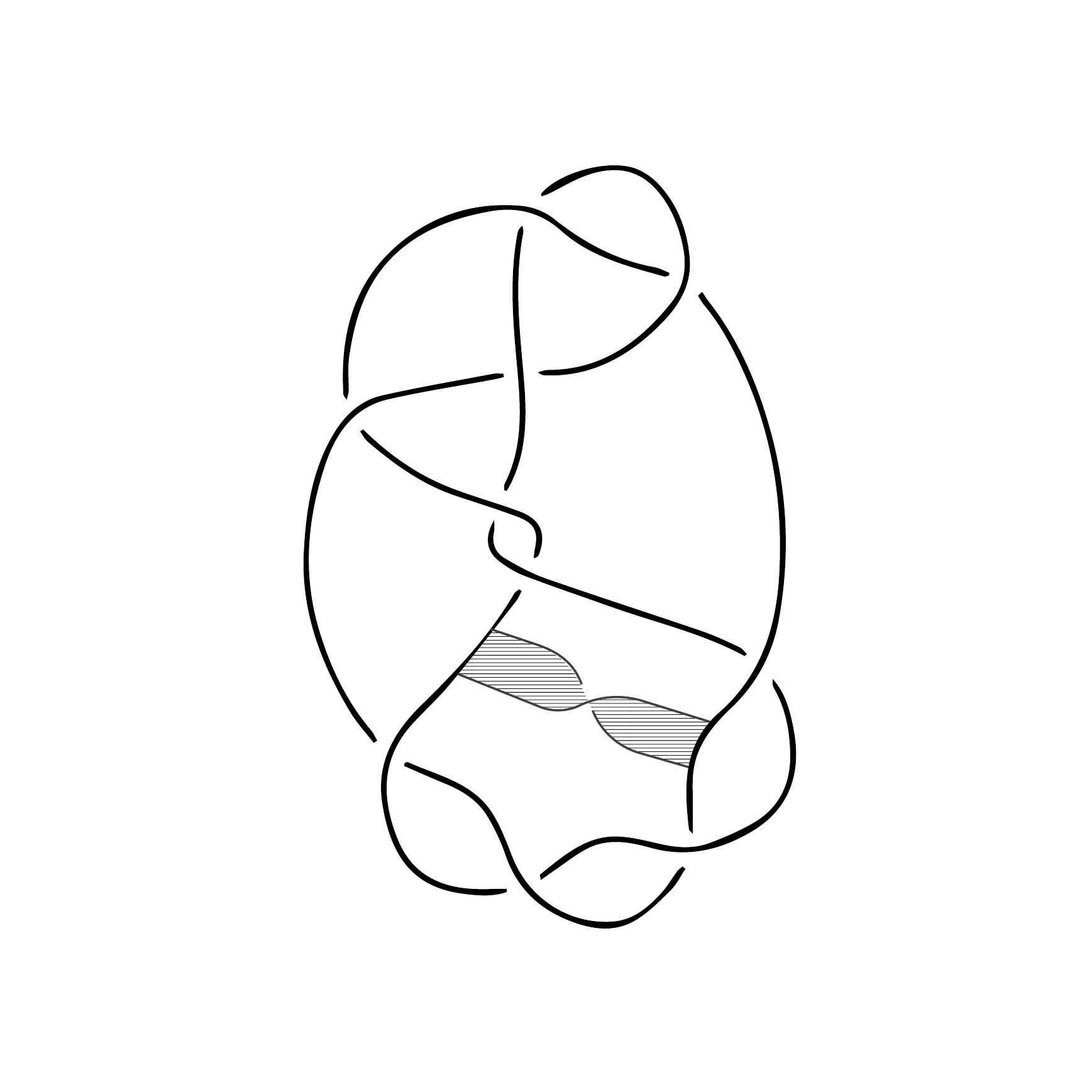}
		\caption{$10_{14}\stackrel{-1\phantom{i}}{\longrightarrow} 3_{1}$}
		\label{FigureFor10-14}
	\end{subfigure}
	~
	\begin{subfigure}[b]{0.3\textwidth}
		\includegraphics[width=\textwidth]{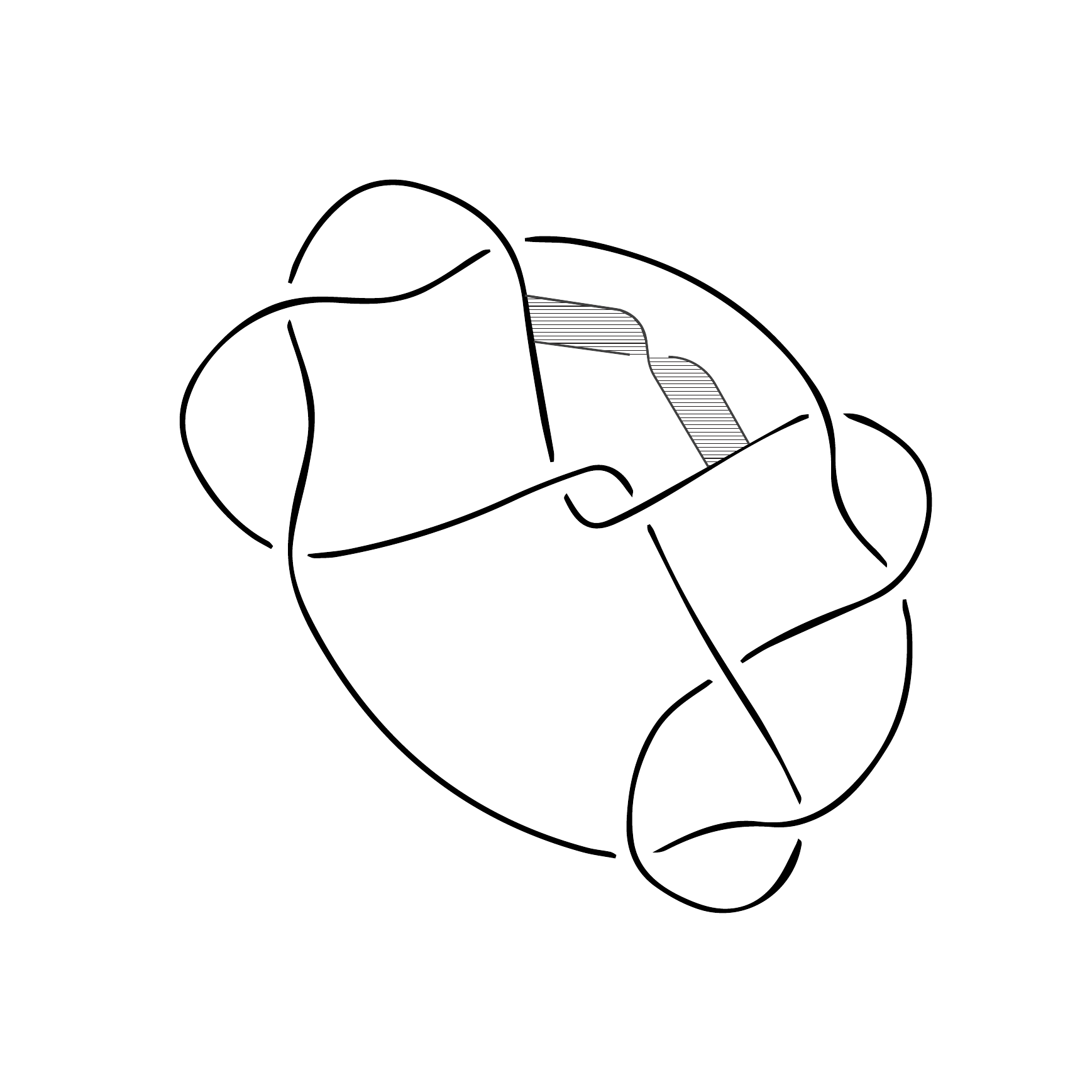}
		\caption{$10_{25}\stackrel{1}{\longrightarrow} 8_{14}$}
		\label{FigureFor10-25}
	\end{subfigure}
	~
	\begin{subfigure}[b]{0.3\textwidth}
		\includegraphics[width=\textwidth]{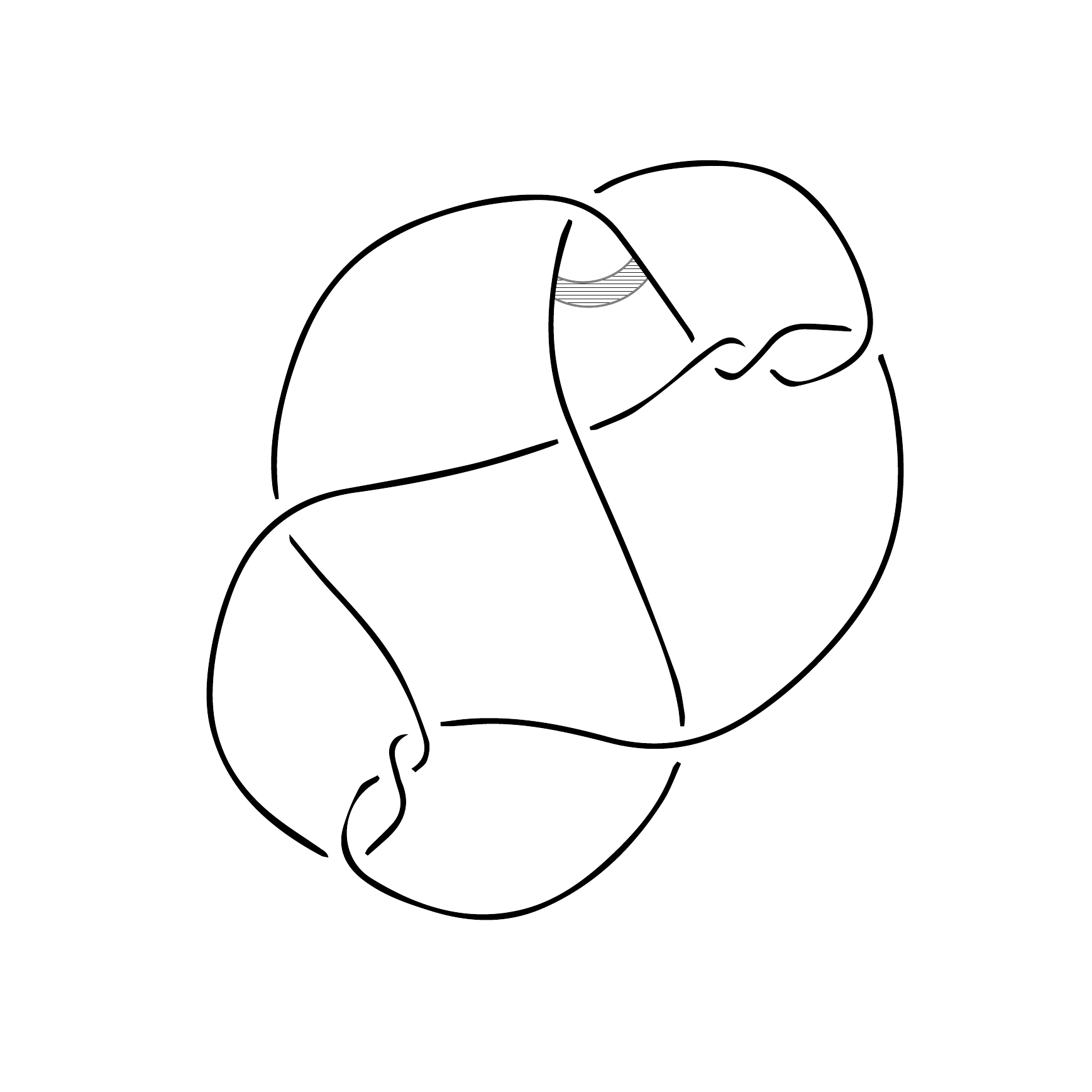}
		\caption{$10_{26}\stackrel{0}{\longrightarrow} 9_{9}$}
		\label{FigureFor10-26}
	\end{subfigure}
	~
	\vskip3mm
	~
	\begin{subfigure}[b]{0.3\textwidth}
		\includegraphics[width=\textwidth]{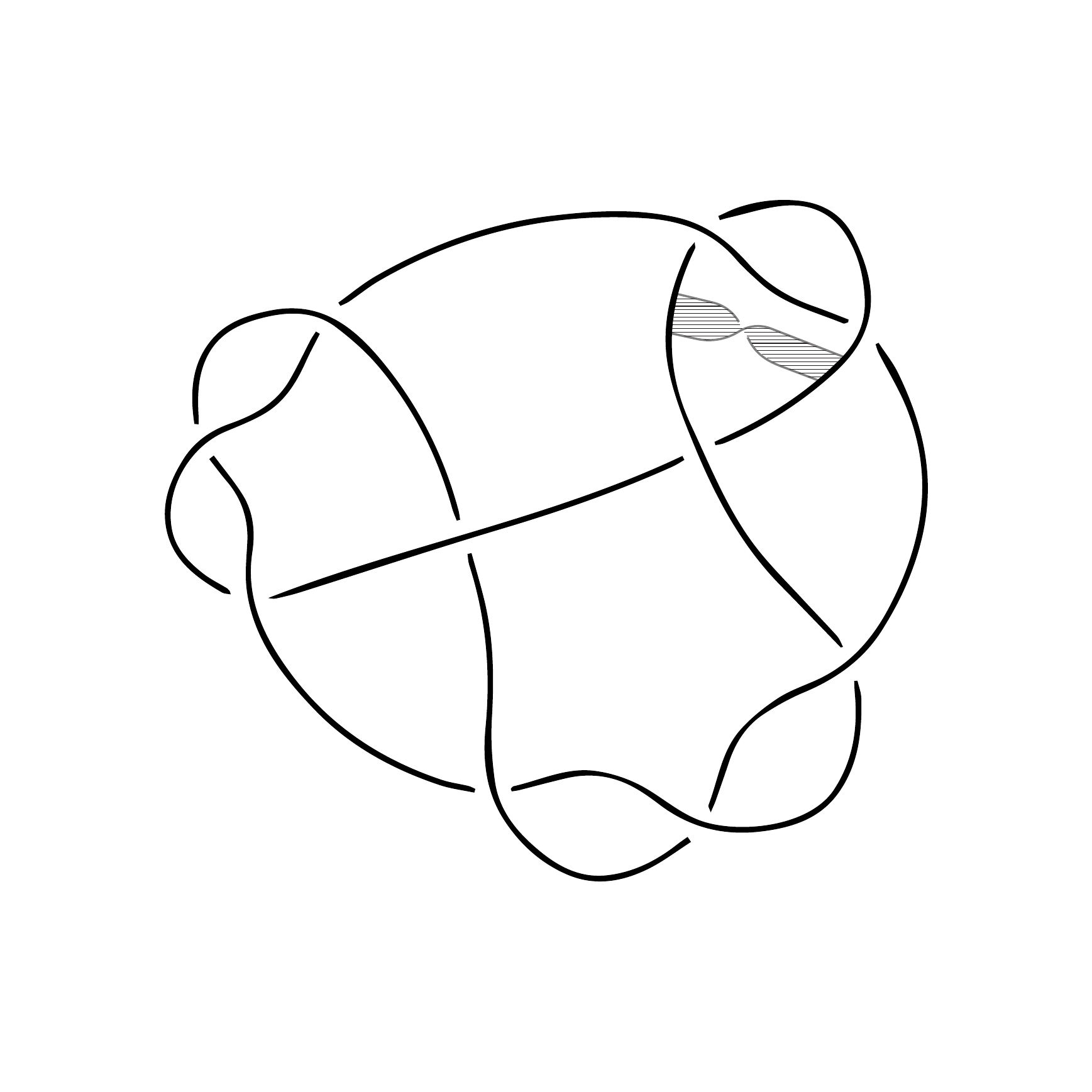}
		\caption{$10_{28}\stackrel{-1}{\longrightarrow} 9_{5}$}
		\label{FigureFor10-28}
	\end{subfigure}
~
	\begin{subfigure}[b]{0.3\textwidth}
		\includegraphics[width=\textwidth]{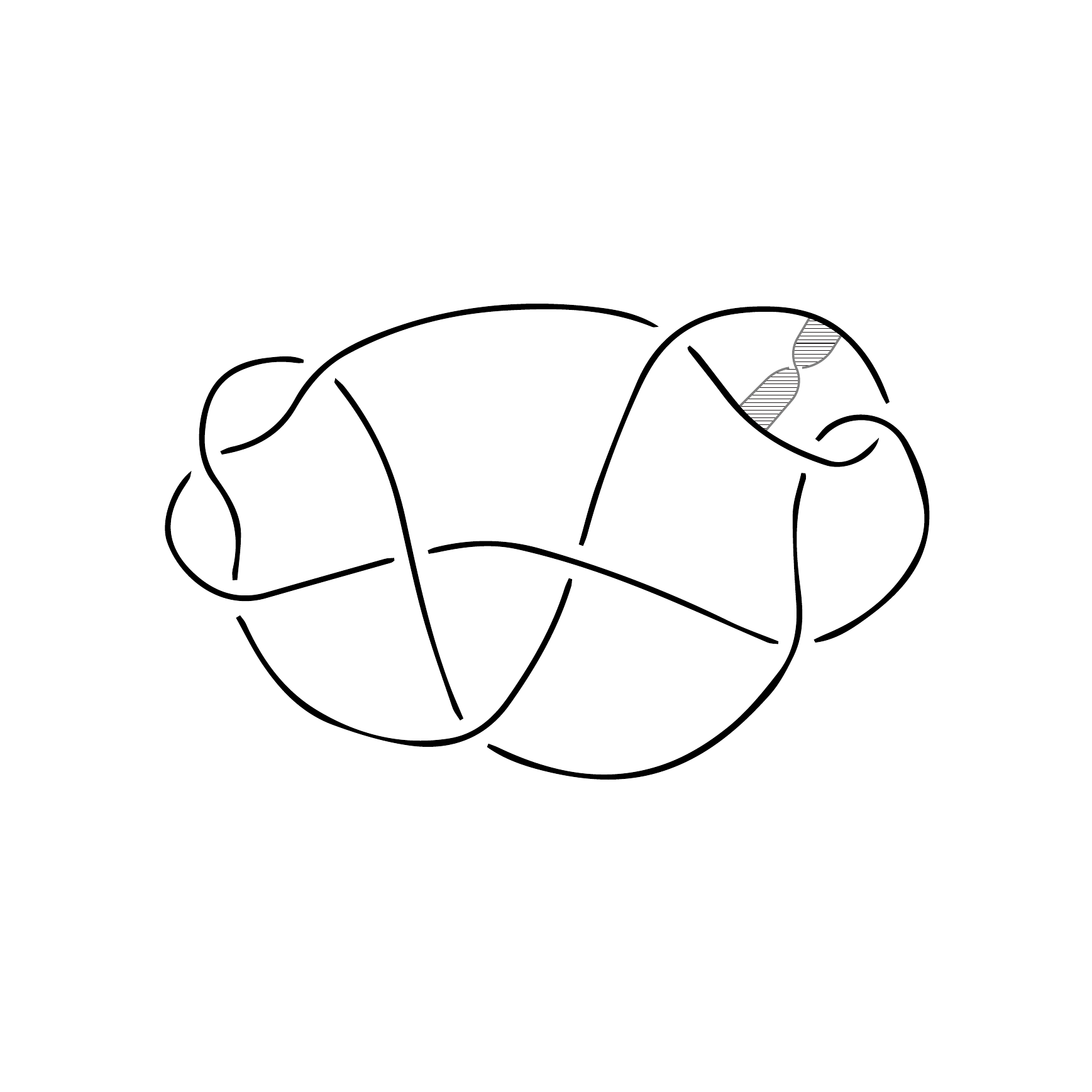}
		\caption{$10_{32}\stackrel{-1}{\longrightarrow} 9_{26}$}
		\label{FigureFor10-32}
	\end{subfigure}
	~
	\begin{subfigure}[b]{0.3\textwidth}
		\includegraphics[width=\textwidth]{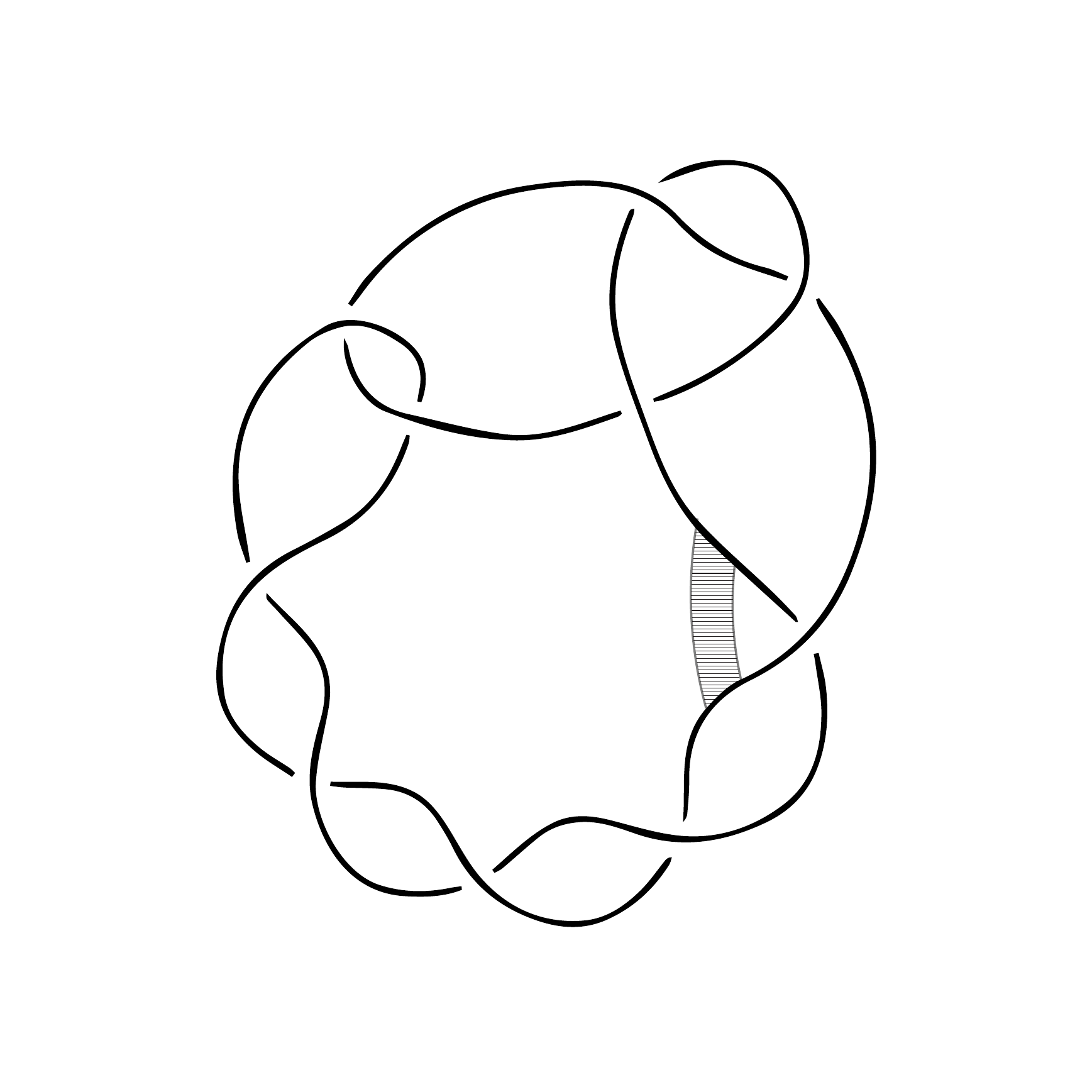}
		\caption{$10_{34}\stackrel{0}{\longrightarrow} 9_8$}
		\label{FigureFor10-34}
	\end{subfigure}
	~
	\vskip3mm
	~

	\caption{Non-oriented band moves from the knots $10_5$, $10_{10}$, $10_{13}$, $10_{14}$, $10_{25}$, $10_{26}$, $10_{28}$, $10_{32}$, $10_{34}$ to knots with $\gamma_4=1$.}\label{gamma4,1}
\end{figure}
\newpage
\begin{figure}[h]
	\centering
	\begin{subfigure}[b]{0.3\textwidth}
		\includegraphics[width=\textwidth]{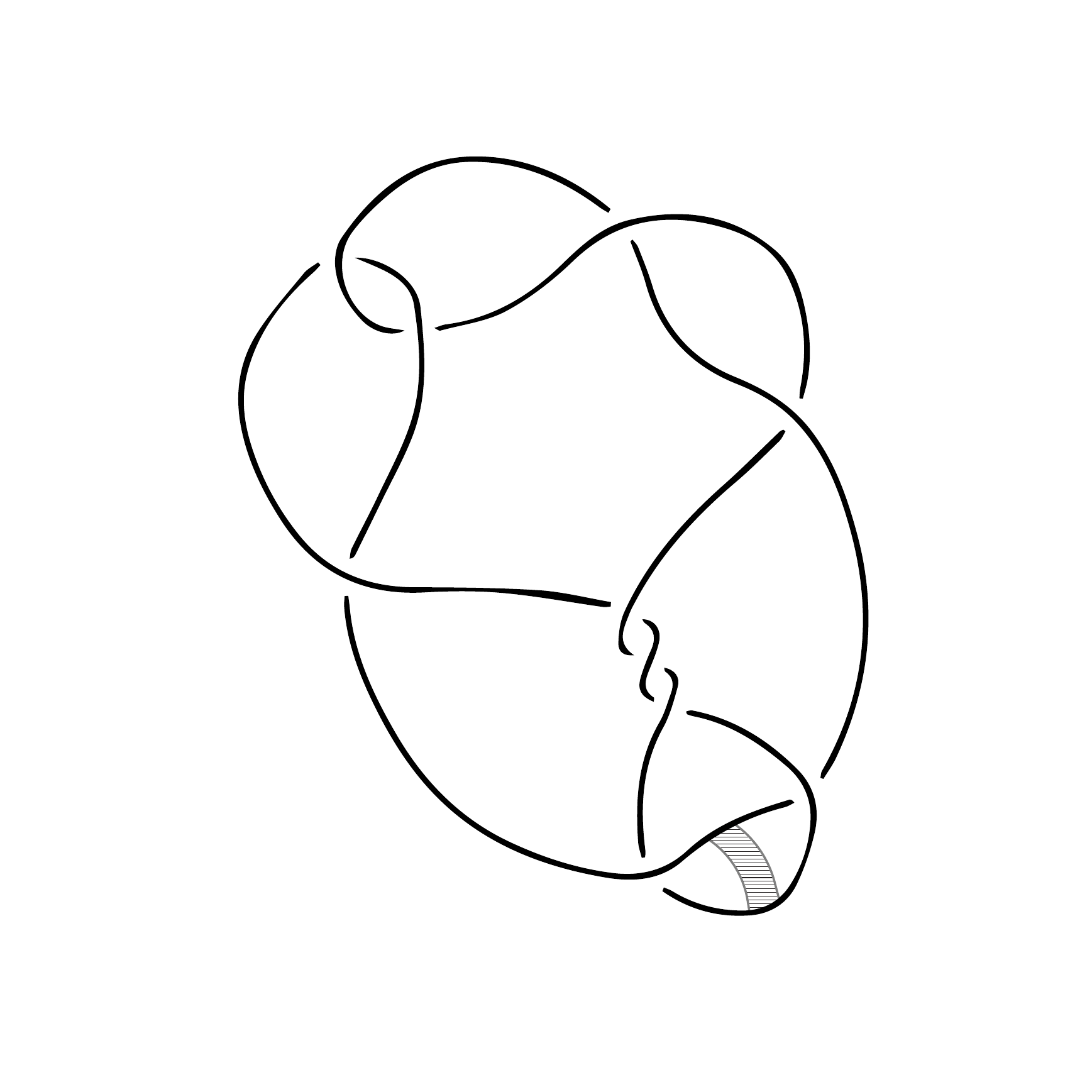}
		\caption{$10_{37}\stackrel{0}{\longrightarrow} 8_6$}
		\label{FigureFor10-37}
	\end{subfigure}
~
	\begin{subfigure}[b]{0.3\textwidth}
		\includegraphics[width=\textwidth]{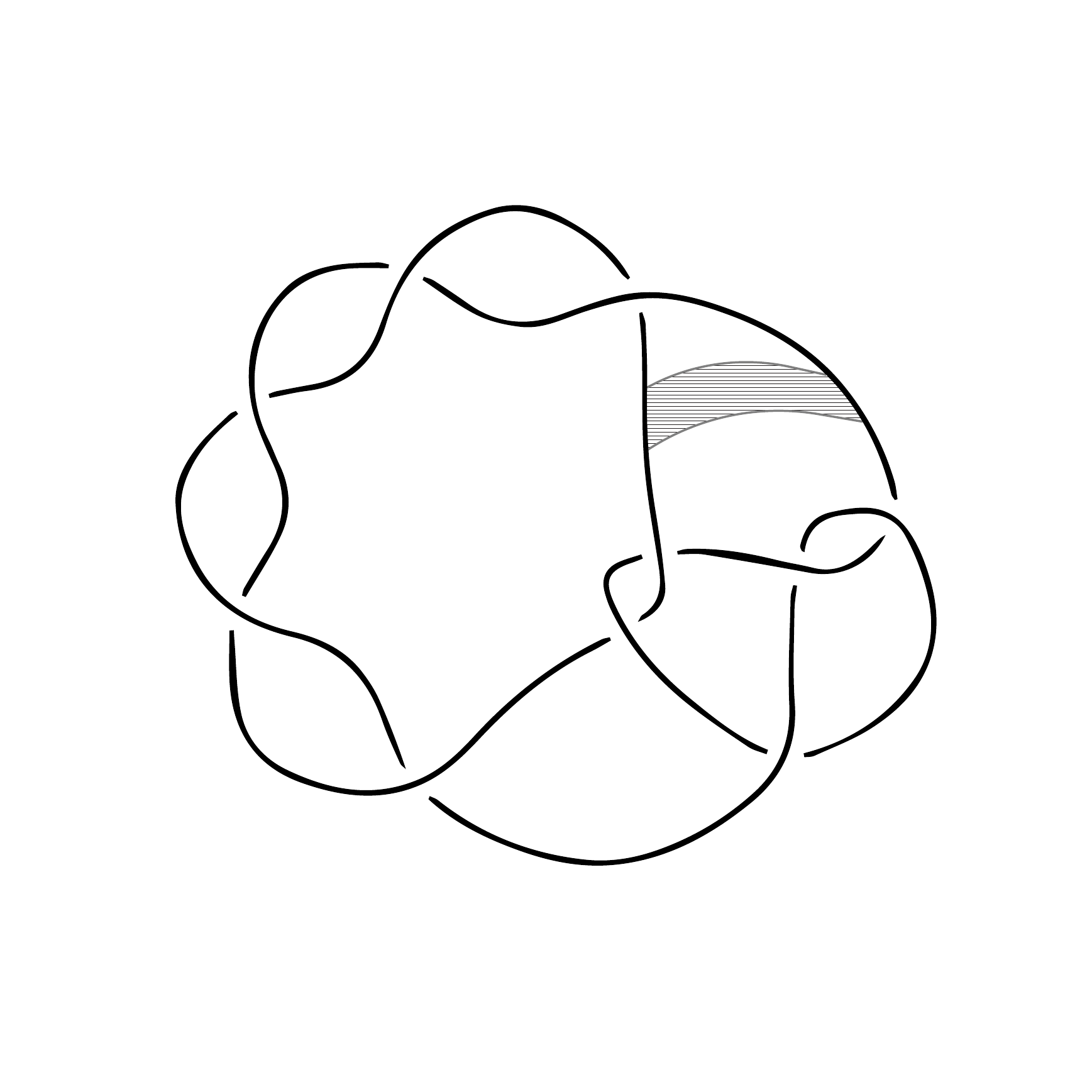}
		\caption{$10_{47}\stackrel{0}{\longrightarrow} 5_2$}
		\label{FigureFor10-47}
	\end{subfigure}
	~
	\begin{subfigure}[b]{0.3\textwidth}
		\includegraphics[width=\textwidth]{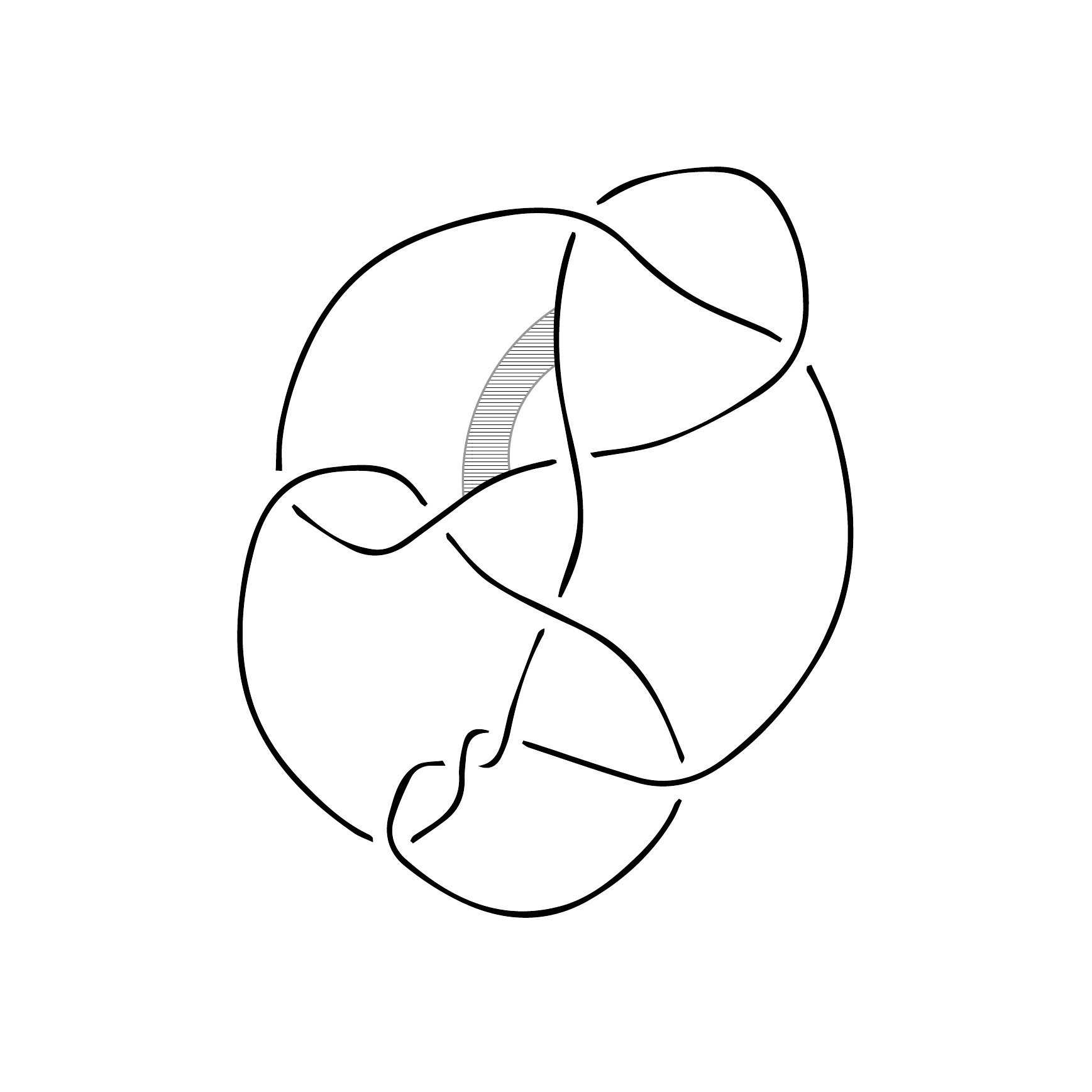}
		\caption{$10_{53}\stackrel{0\phantom{i}}{\longrightarrow} 9_{21}$}
		\label{FigureFor10-53}
	\end{subfigure}
	~
	\vskip3mm
	~
	\begin{subfigure}[b]{0.3\textwidth}
		\includegraphics[width=\textwidth]{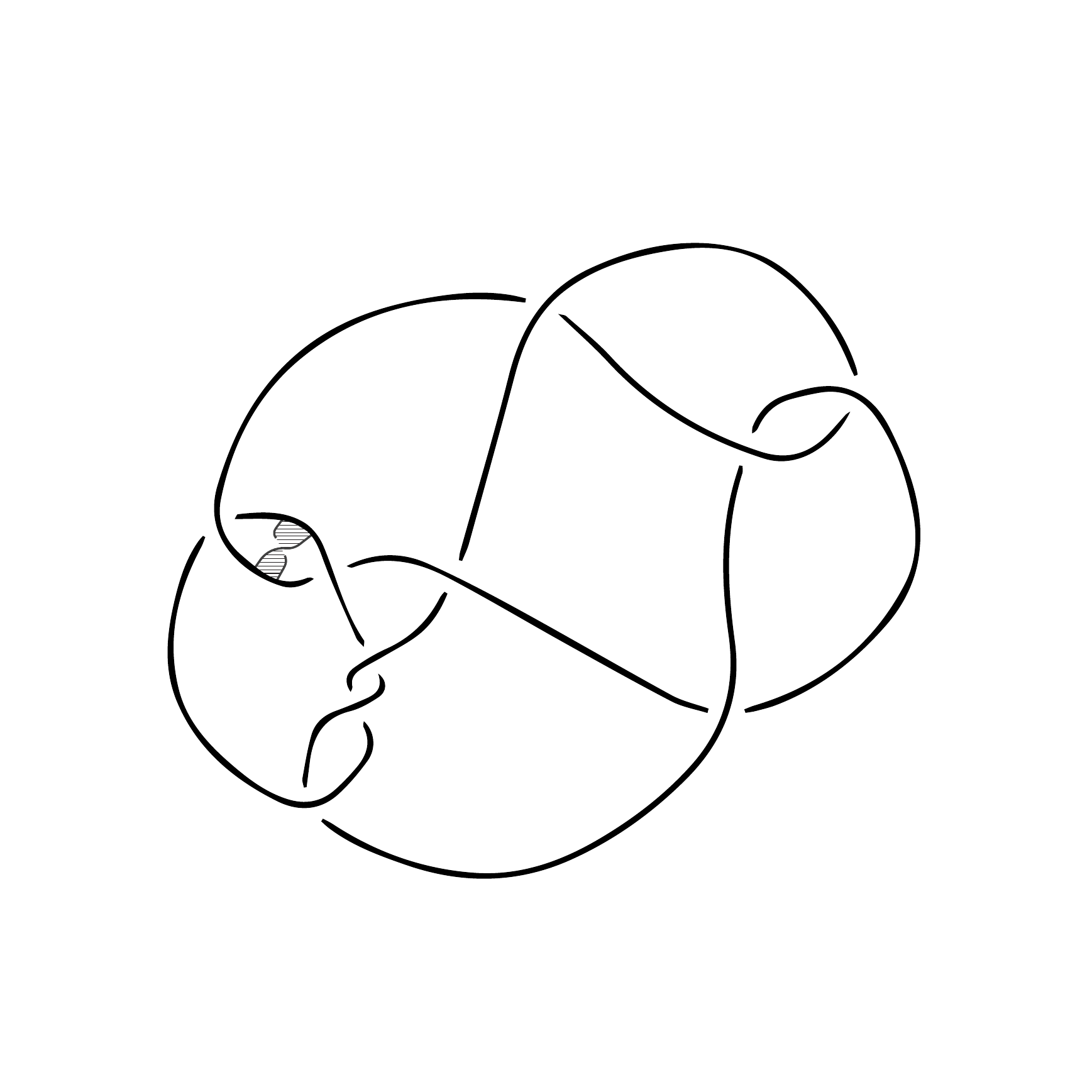}
		\caption{$10_{56}\stackrel{1}{\longrightarrow} 9_{21}$}
		\label{FigureFor10-56}
	\end{subfigure}
	~
	\begin{subfigure}[b]{0.3\textwidth}
		\includegraphics[width=\textwidth]{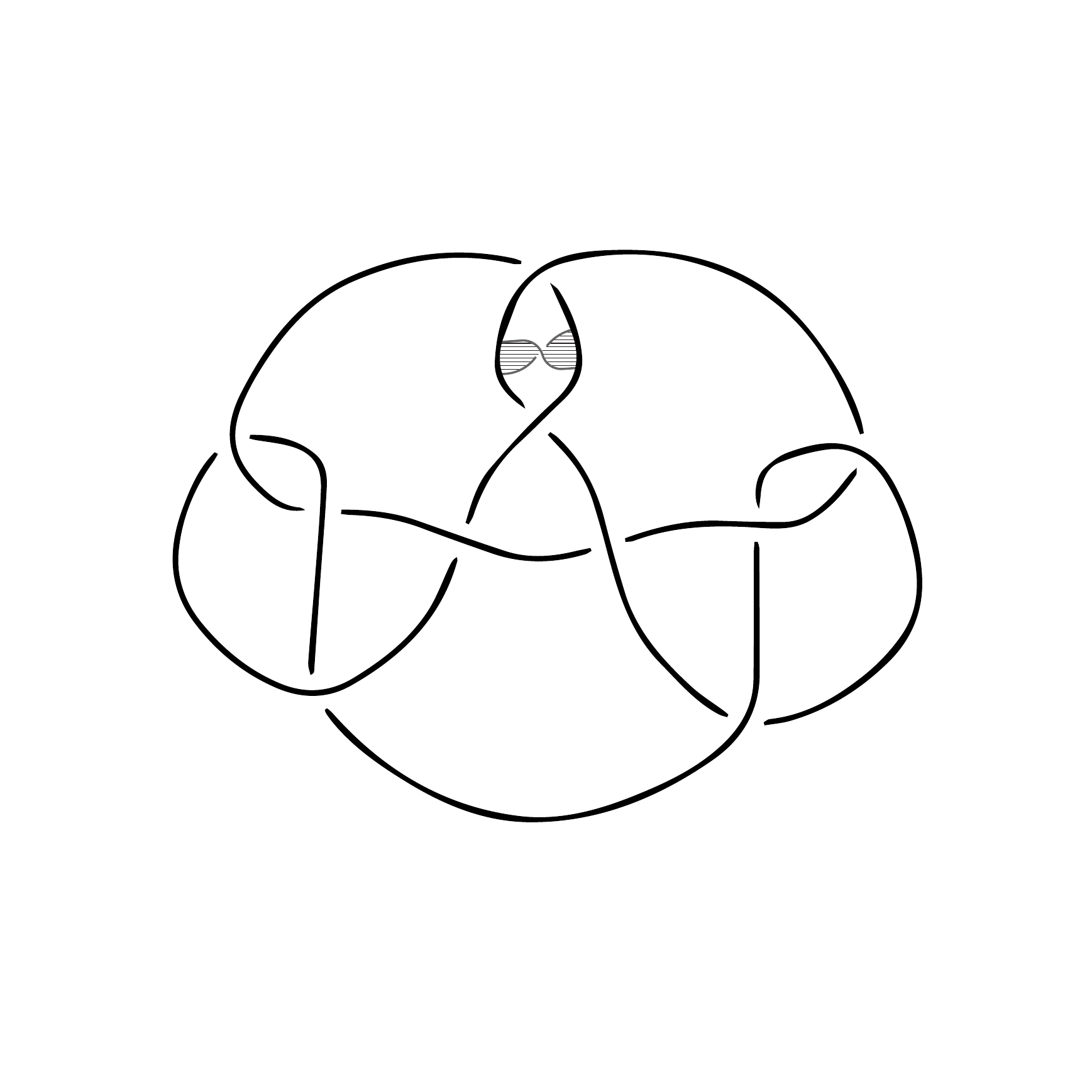}
		\caption{$10_{60}\stackrel{1\phantom{i}}{\longrightarrow} 9_{31}$}
		\label{FigureFor10-60}
	\end{subfigure}
	~
	\begin{subfigure}[b]{0.3\textwidth}
		\includegraphics[width=\textwidth]{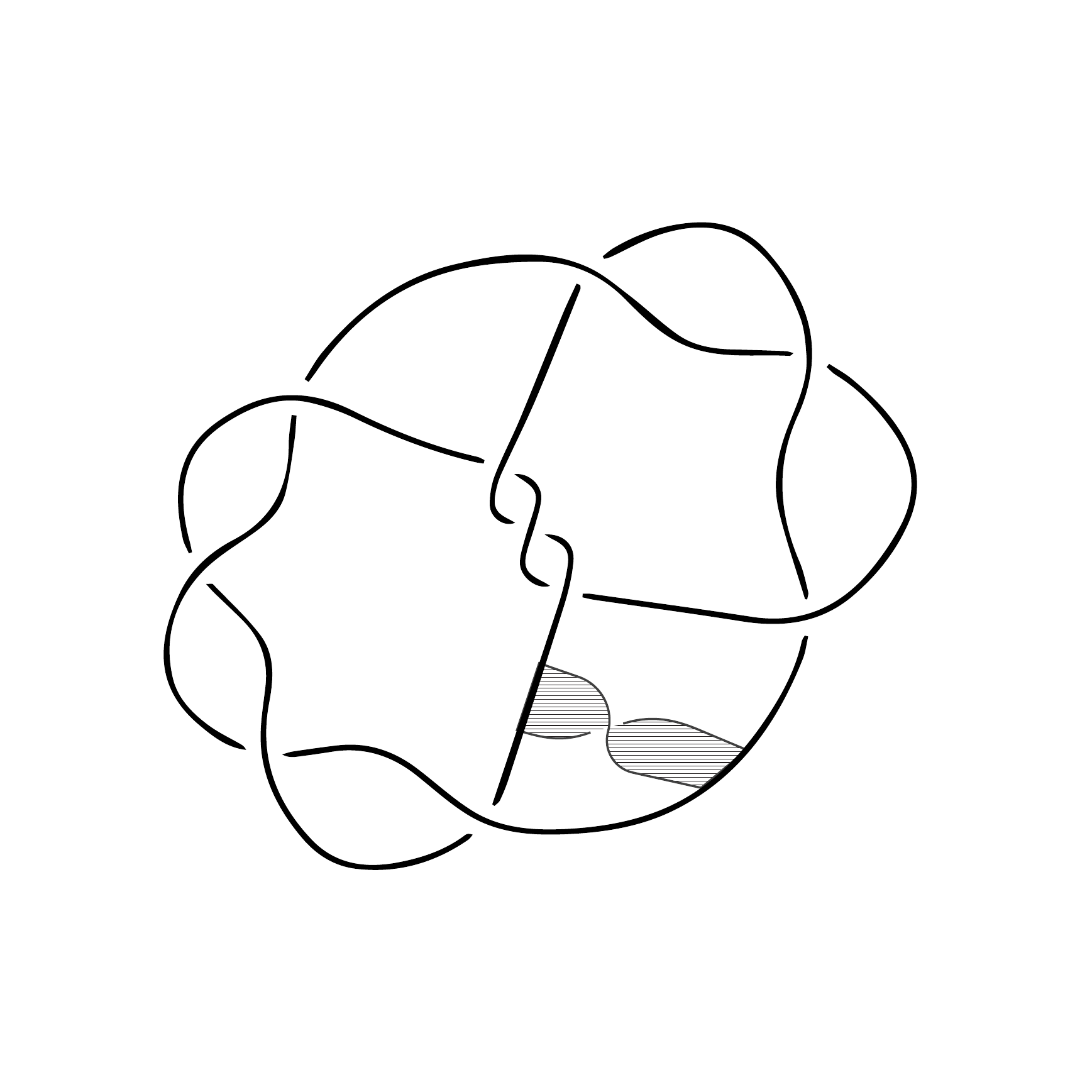}
		\caption{$10_{61}\stackrel{1}{\longrightarrow} 9_{35}$}
		\label{FigureFor10-61}
	\end{subfigure}
	~
	\vskip3mm
	~
	\begin{subfigure}[b]{0.3\textwidth}
		\includegraphics[width=\textwidth]{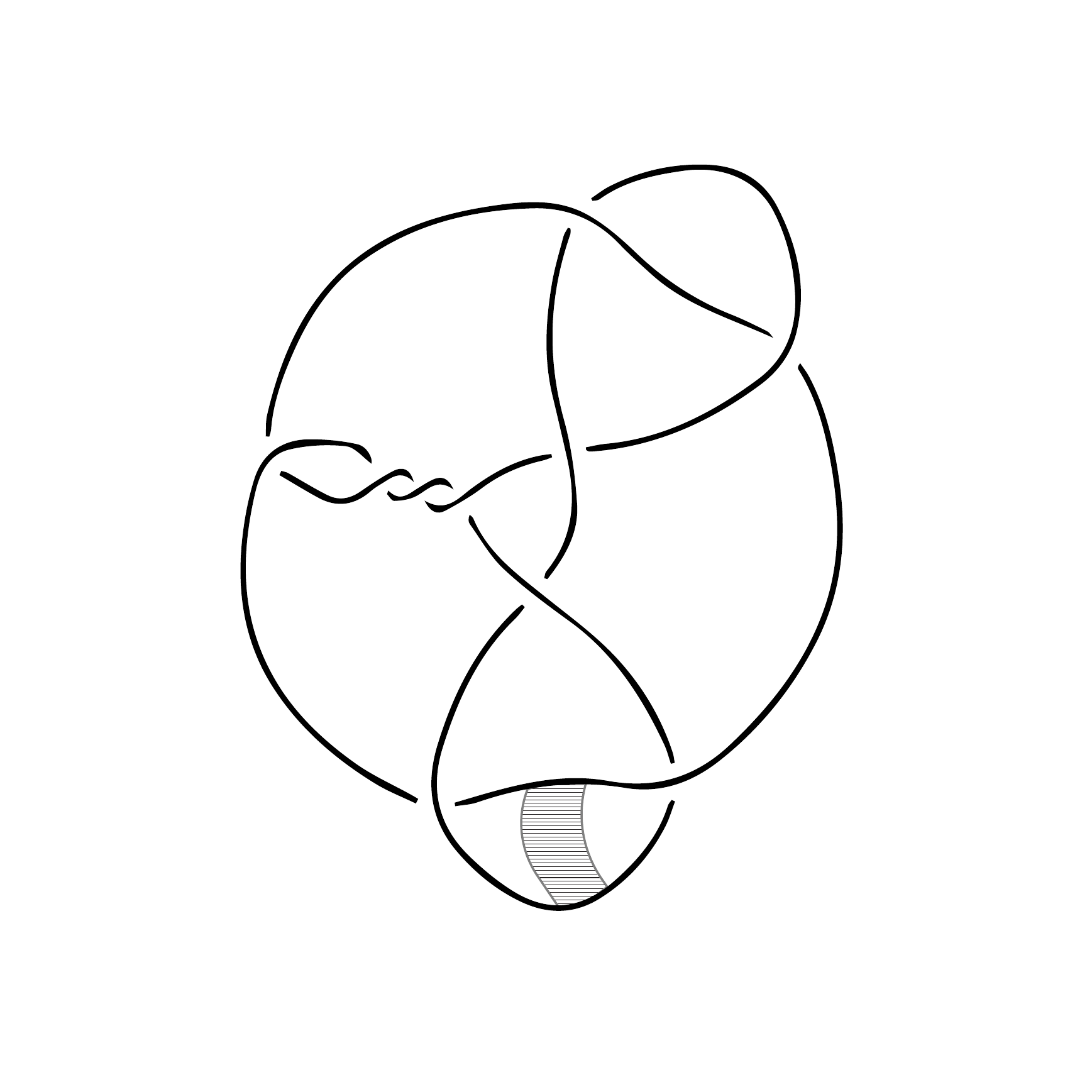}
		\caption{$10_{63}\stackrel{0\phantom{i}}{\longrightarrow} 8_7$}
		\label{FigureFor10-63}
	\end{subfigure}
	~
	\begin{subfigure}[b]{0.3\textwidth}
		\includegraphics[width=\textwidth]{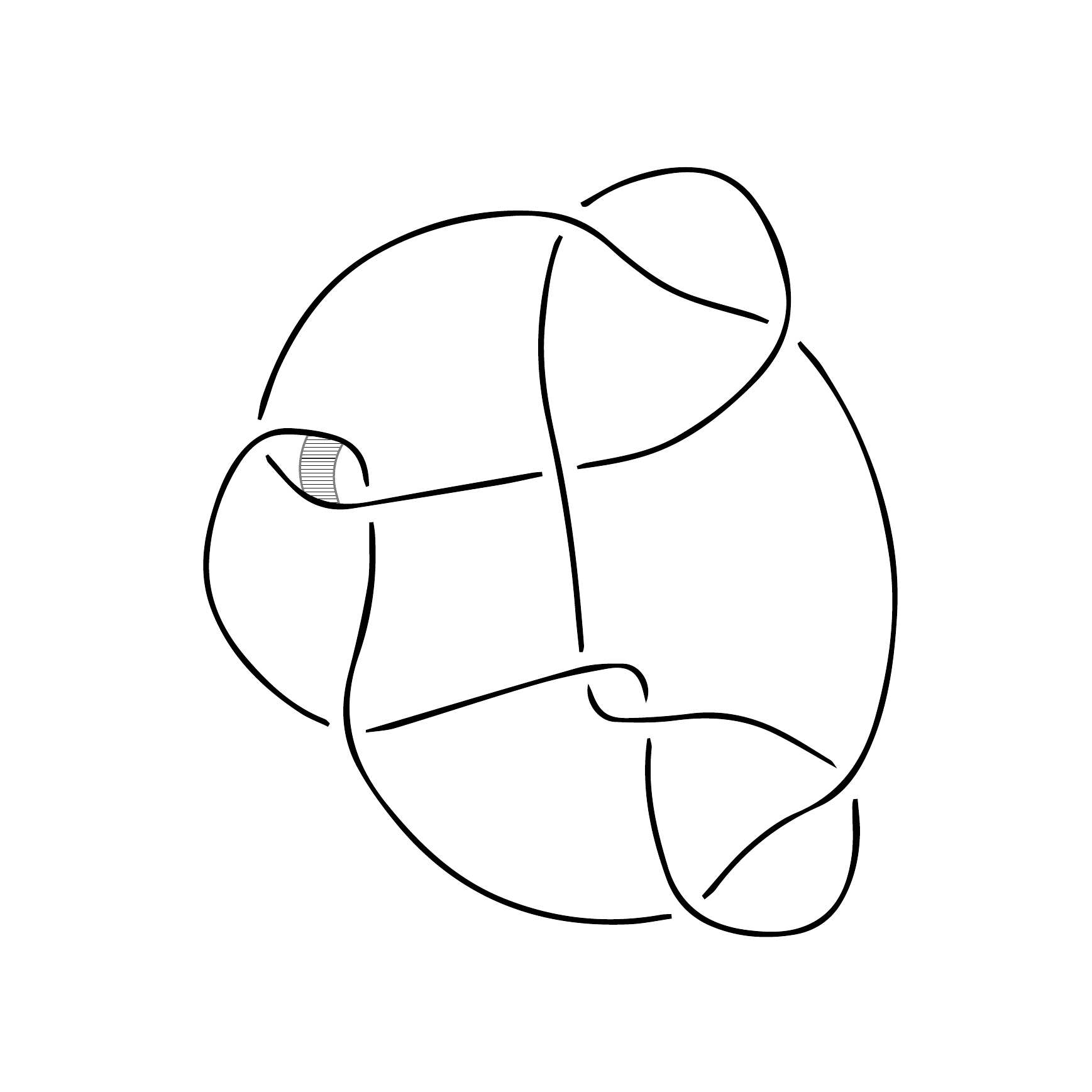}
		\caption{$10_{71}\stackrel{0}{\longrightarrow} 8_{14}$}
		\label{FigureFor10-71}
	\end{subfigure}
	~
	\begin{subfigure}[b]{0.3\textwidth}
		\includegraphics[width=\textwidth]{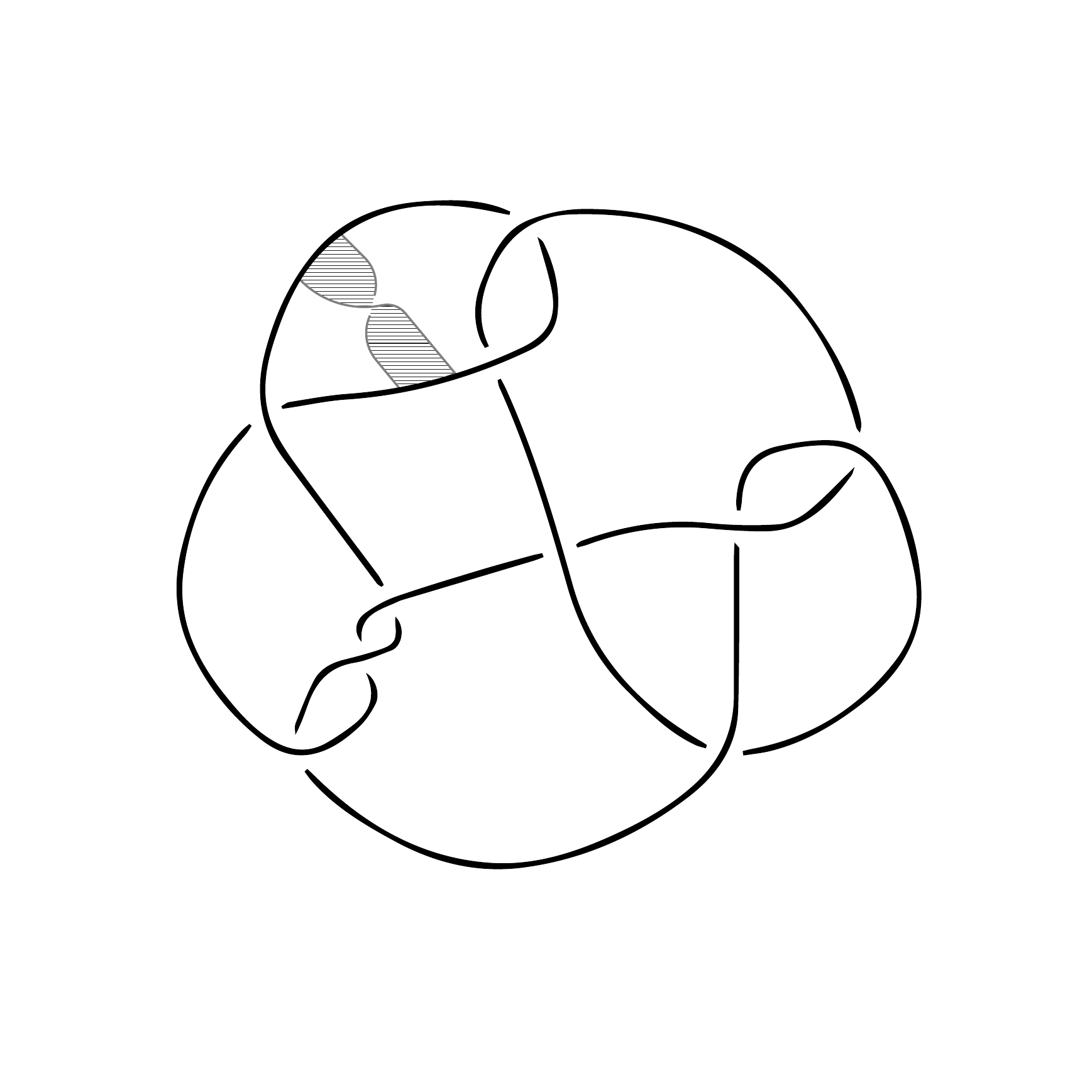}
		\caption{$10_{72}\stackrel{-1}{\longrightarrow} 9_{22}$}
		\label{FigureFor10-72}
	\end{subfigure}
	~
	\vskip3mm
	~
	\caption{Non-oriented band moves from the knots $10_{37}$, $10_{47}$, $10_{53}$, $10_{56}$, $10_{60}$, $10_{61}$, $10_{63}$, $10_{71}$, $10_{72}$ to knots with $\gamma_4=1$.}\label{gamma4,2}
\end{figure}
\newpage
\begin{figure}[h]
	\centering
	 \begin{subfigure}[b]{0.3\textwidth}
		\includegraphics[width=\textwidth]{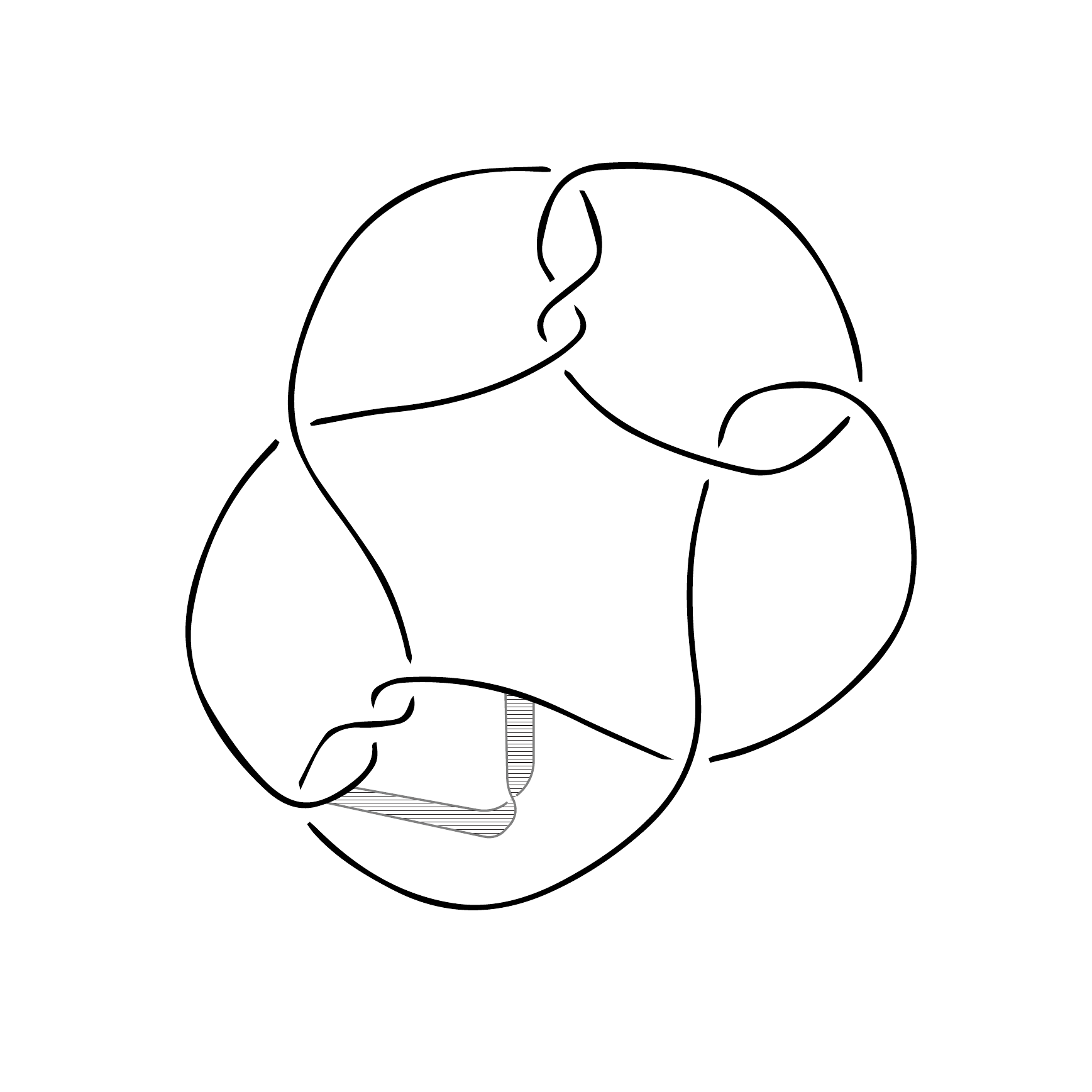}
		\caption{$10_{76}\stackrel{1}{\longrightarrow} 6_{2}$}
		\label{FigureFor10-76}
	\end{subfigure}
~
	\begin{subfigure}[b]{0.3\textwidth}
		\includegraphics[width=\textwidth]{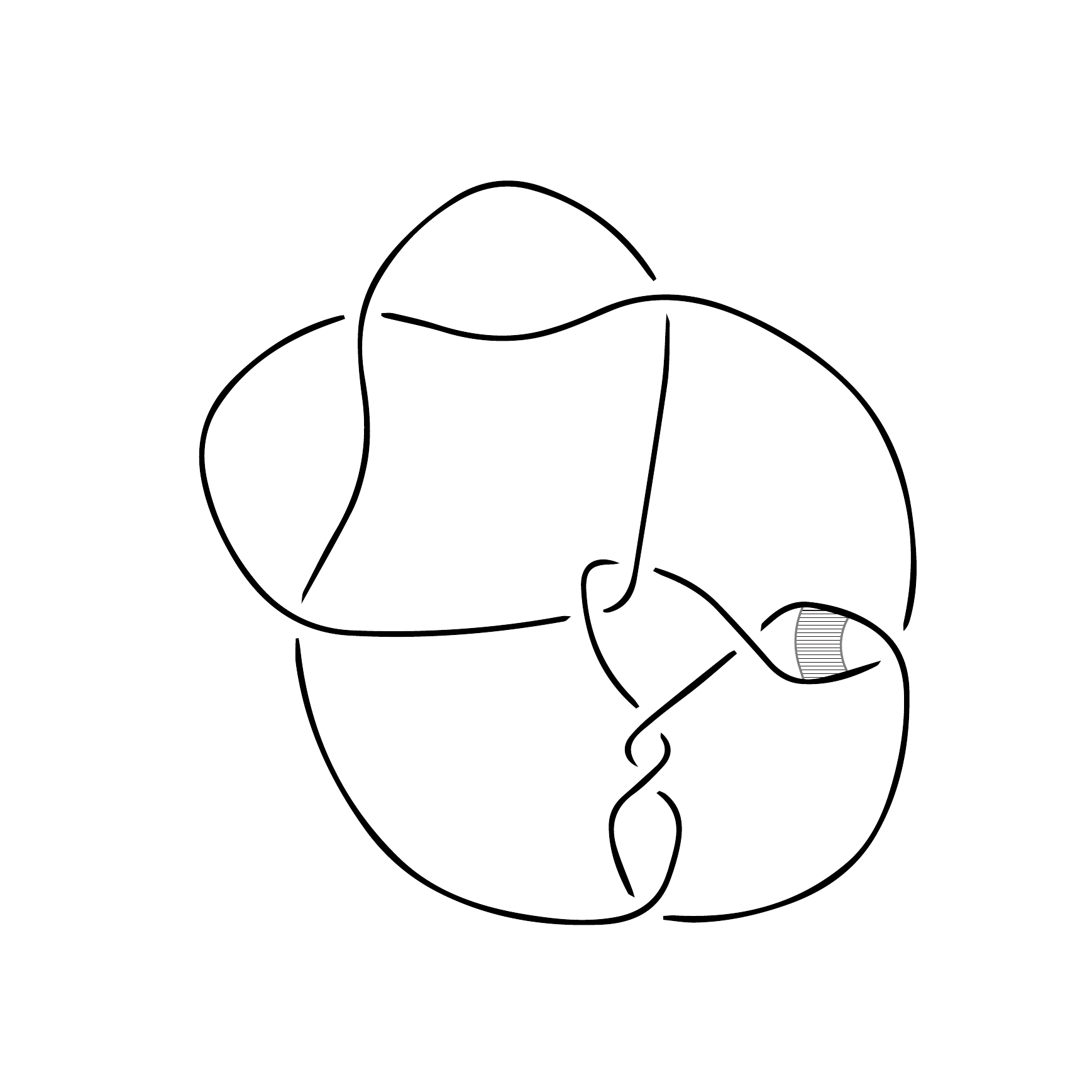}
		\caption{$10_{79}\stackrel{0}{\longrightarrow} 8_6$}
		\label{FigureFor10-79}
	\end{subfigure}
~
	\begin{subfigure}[b]{0.30\textwidth}
		\includegraphics[width=\textwidth]{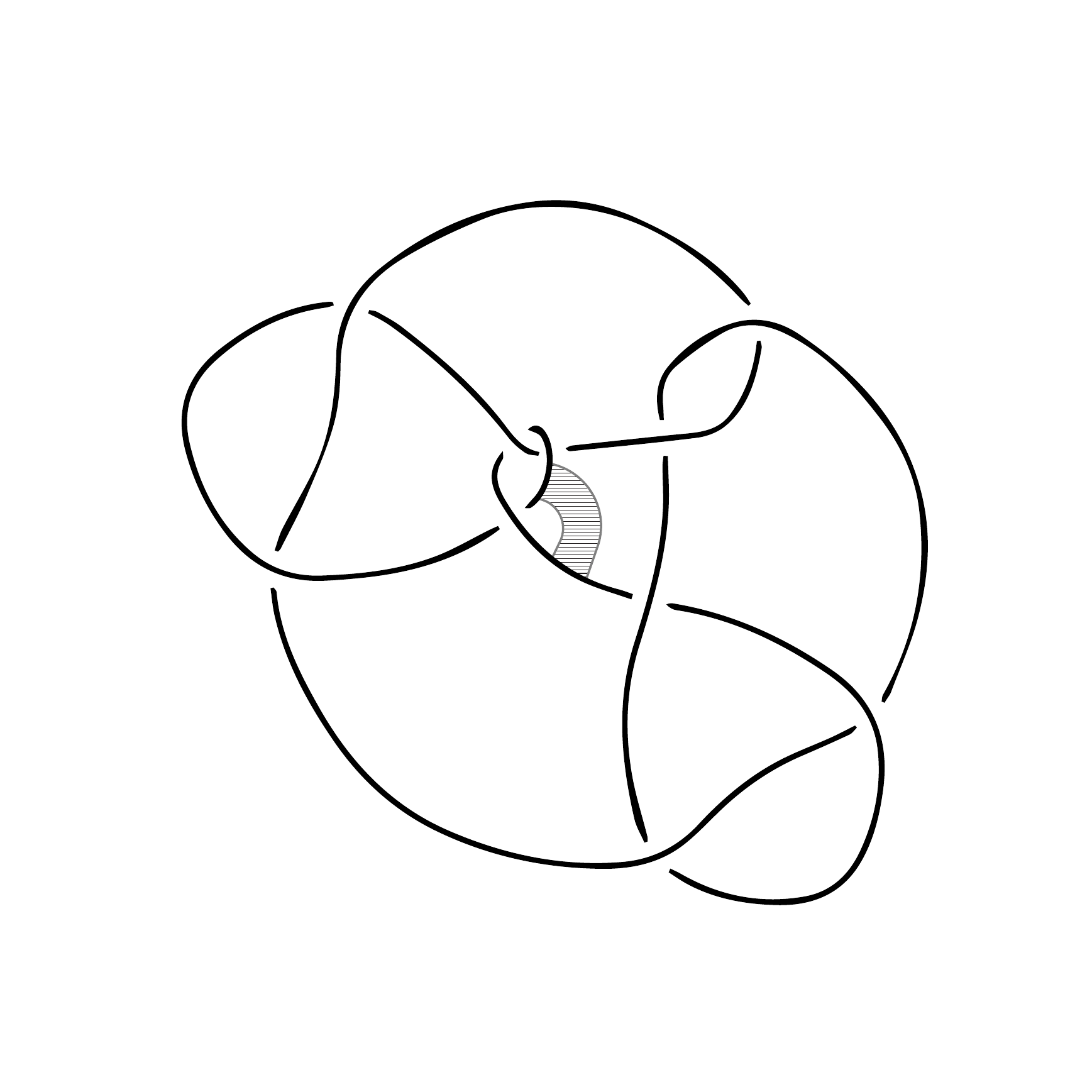}
		\caption{$10_{81}\stackrel{0}{\longrightarrow} 9_{25}$}
		\label{FigureFor10-81}
	\end{subfigure}
	~
	\vskip3mm
	~
	\begin{subfigure}[b]{0.30\textwidth}
		\includegraphics[width=\textwidth]{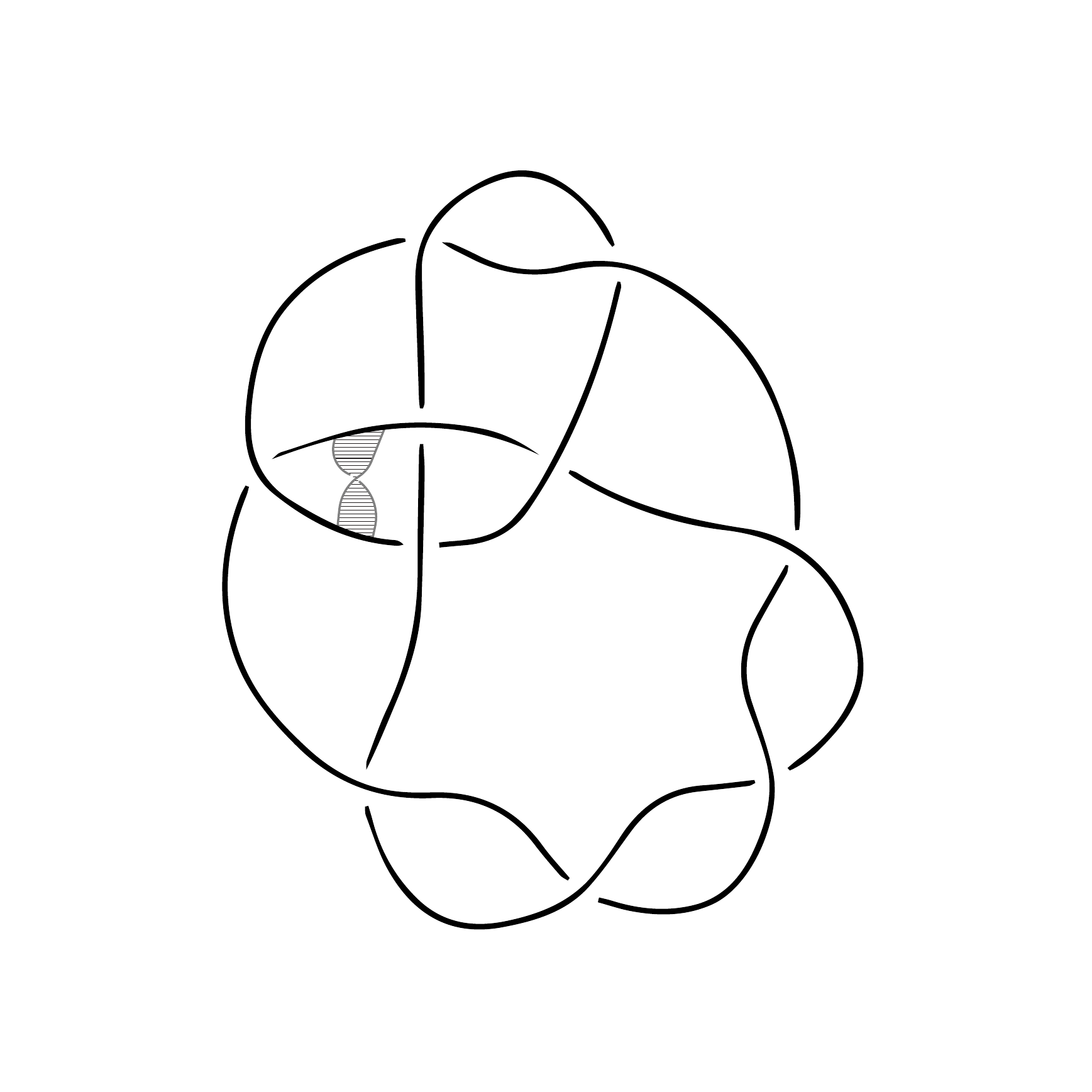}
		\caption{$10_{85}\stackrel{1}{\longrightarrow} 9_5$}
		\label{FigureFor10-85}
	\end{subfigure}
	~
	\begin{subfigure}[b]{0.3\textwidth}
		\includegraphics[width=\textwidth]{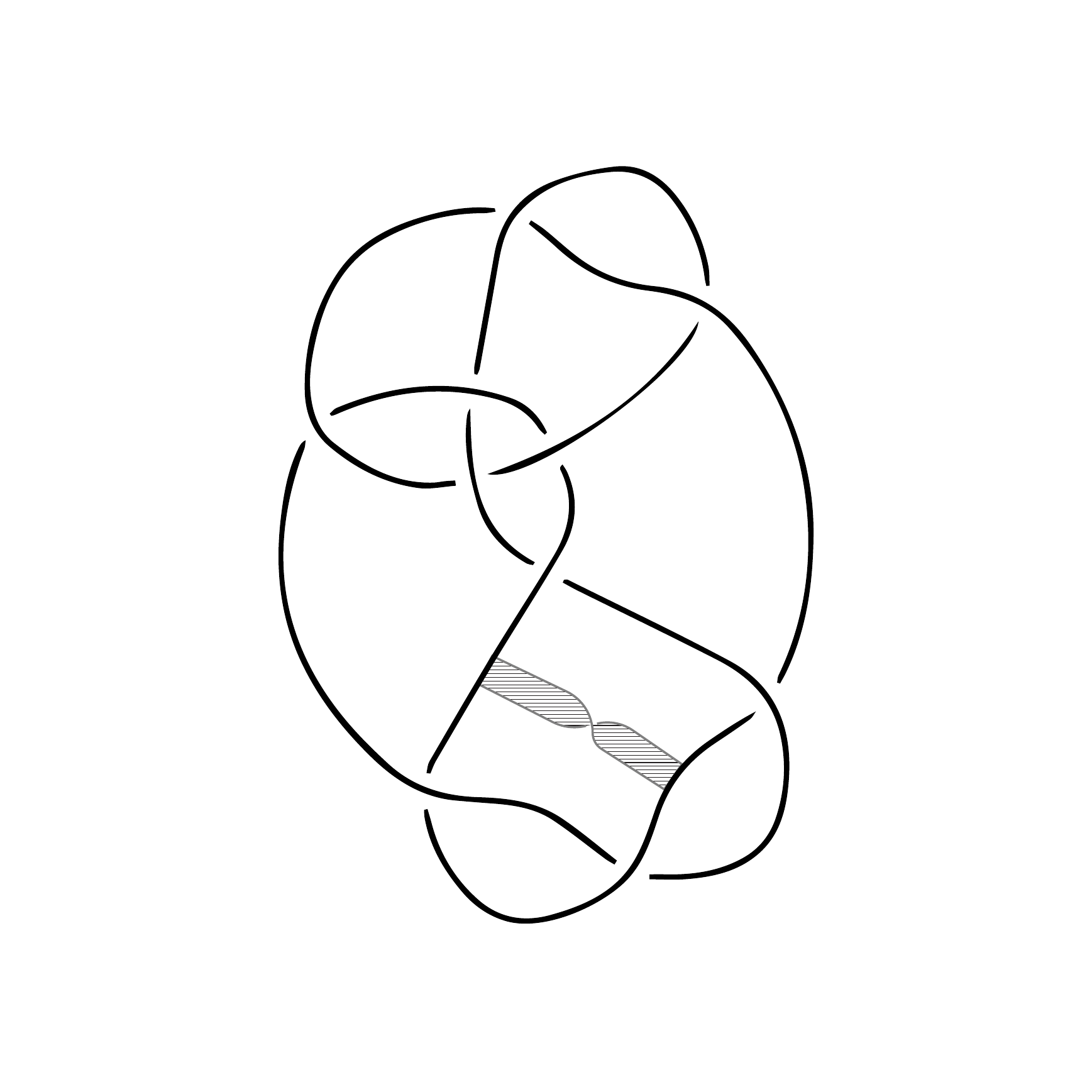}
		\caption{$10_{86}\stackrel{1}{\longrightarrow} 6_{2}$}
		\label{FigureFor10-86}
	\end{subfigure}%
~
	\begin{subfigure}[b]{0.30\textwidth}
		\includegraphics[width=\textwidth]{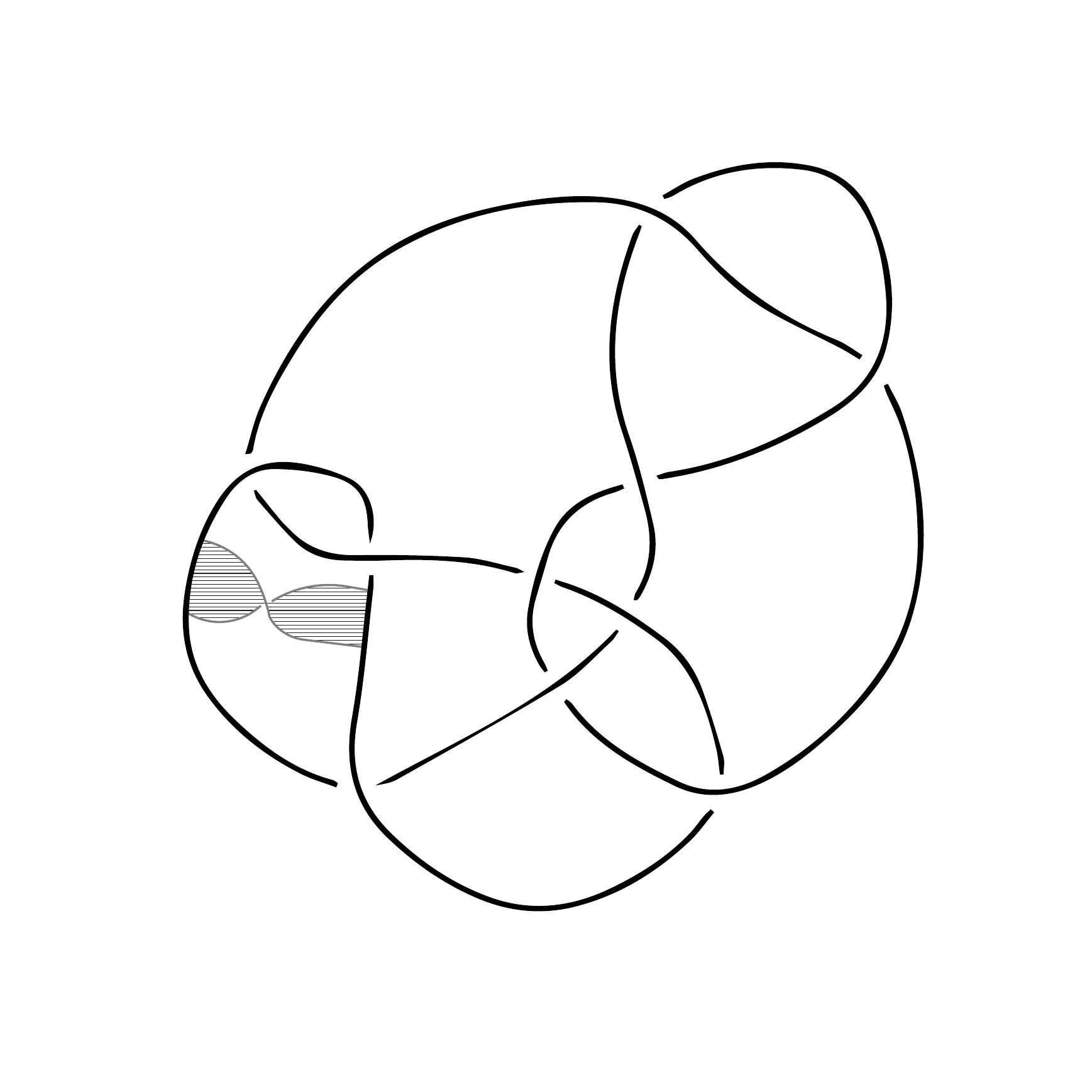}
		\caption{$10_{88}\stackrel{1}{\longrightarrow} 9_{32}$}
		\label{FigureFor10-88}
	\end{subfigure}
	\vskip3mm
	\begin{subfigure}[b]{0.30\textwidth}
		\includegraphics[width=\textwidth]{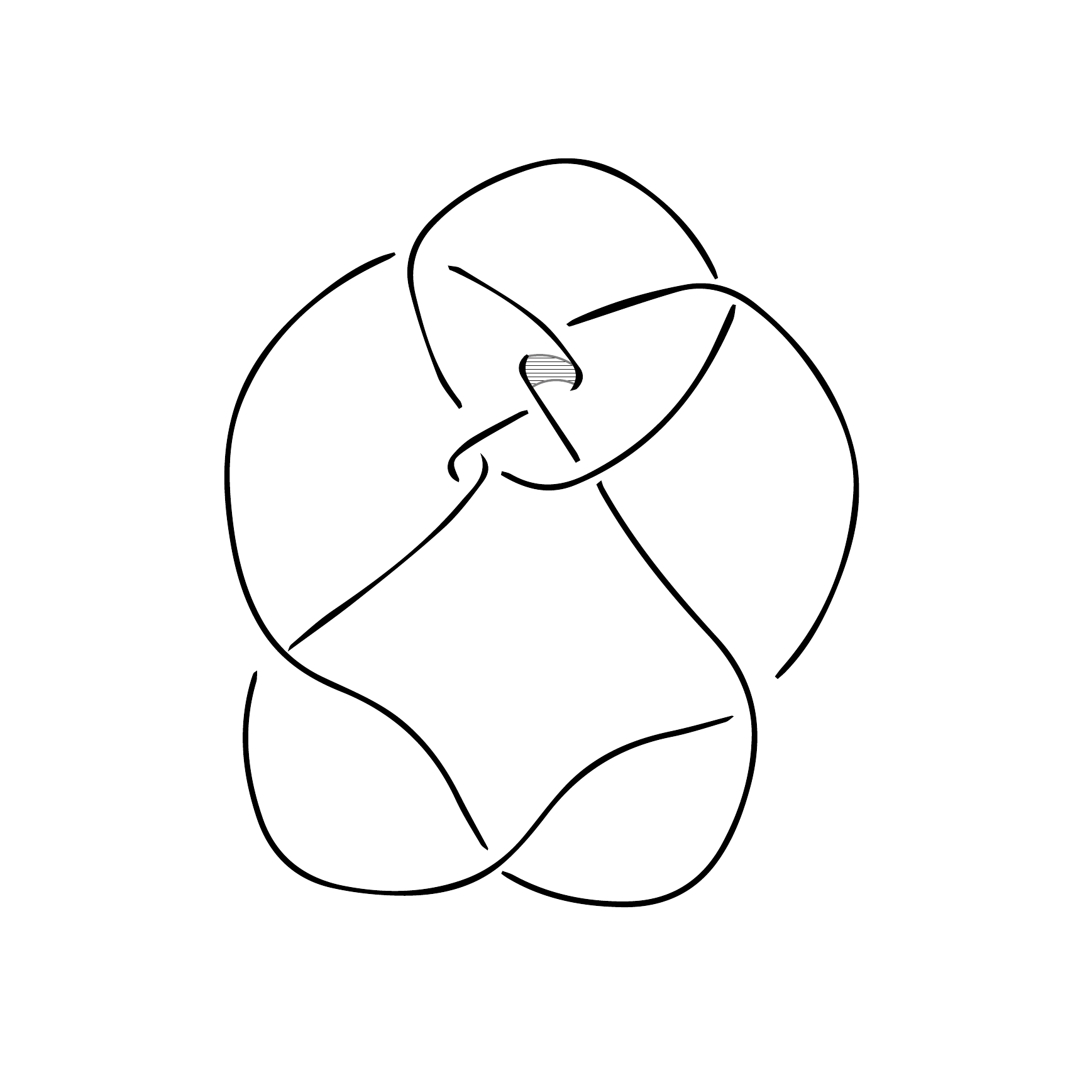}
		\caption{$10_{90}\stackrel{0\phantom{i}}{\longrightarrow} 8_{10}$}
		\label{FigureFor10-90}
	\end{subfigure}
	~
	\begin{subfigure}[b]{0.30\textwidth}
		\includegraphics[width=\textwidth]{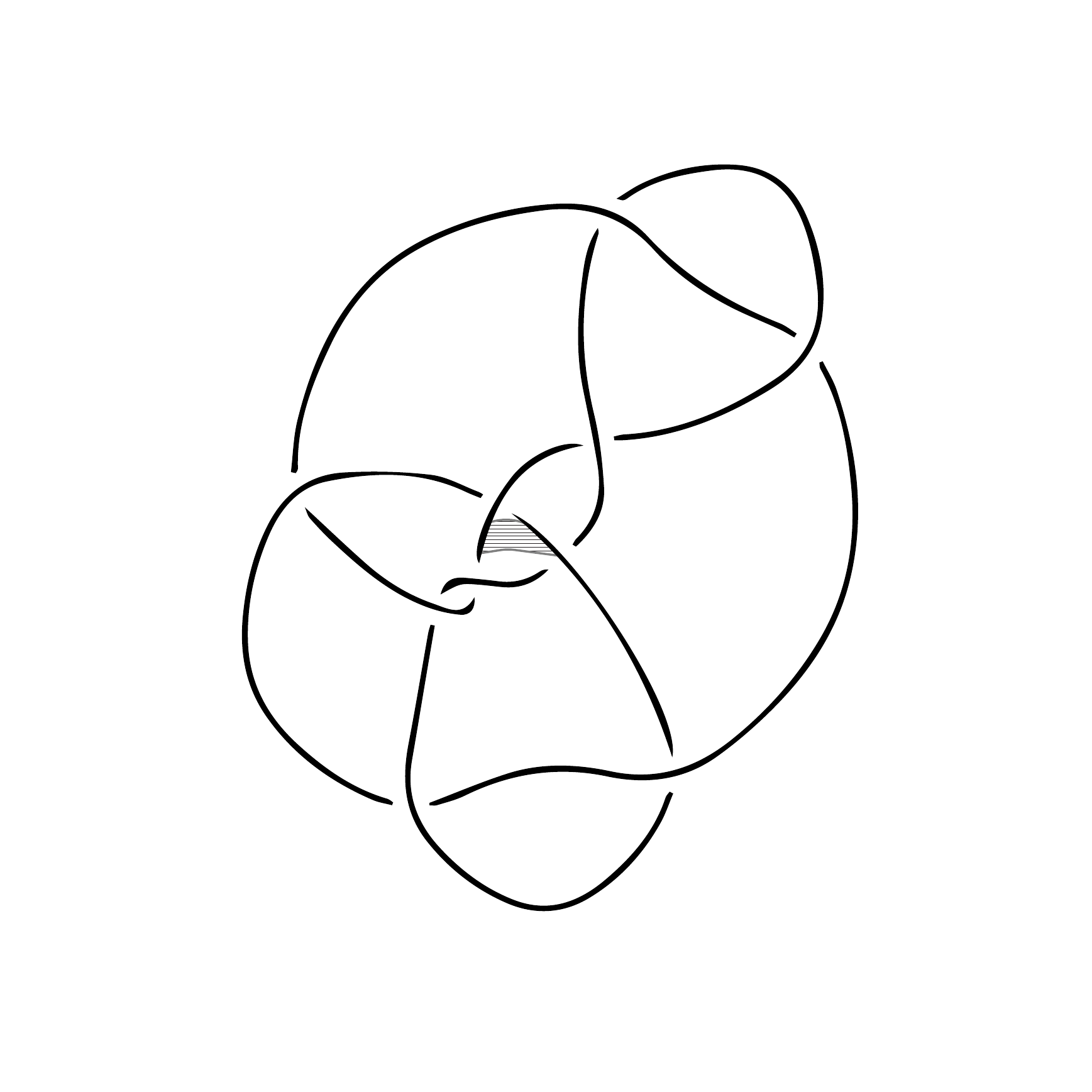}
		\caption{$10_{92}\stackrel{0}{\longrightarrow} 9_{22}$}
		\label{FigureFor10-92}
	\end{subfigure}
	~
	\begin{subfigure}[b]{0.3\textwidth}
		\includegraphics[width=\textwidth]{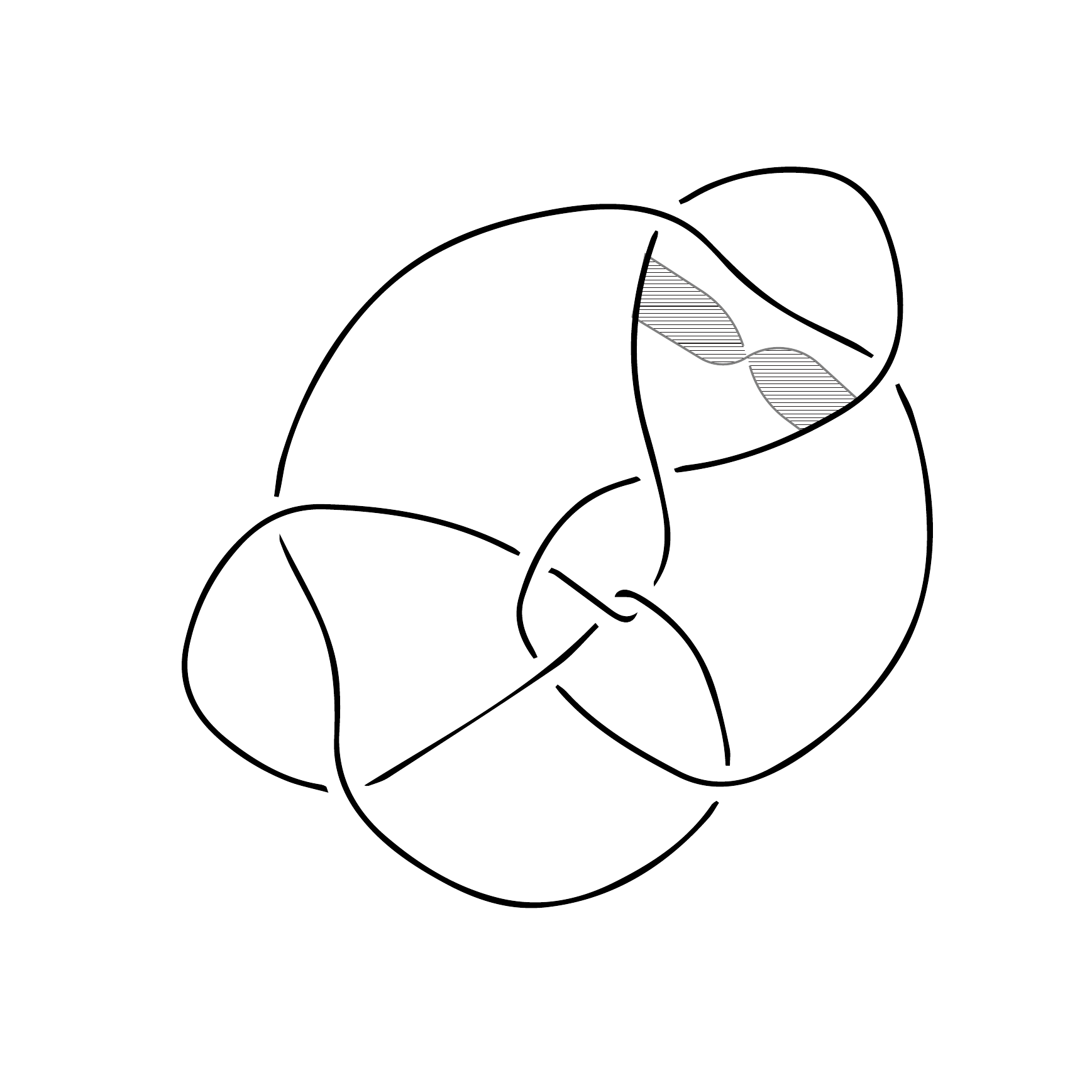}
		\caption{$10_{96}\stackrel{-1}{\longrightarrow} 9_{29}$}
		\label{FigureFor10-96}
	\end{subfigure}%
	   	~
		\vskip3mm
	\caption{Non-oriented band moves from the knots $10_{76}$, $10_{79}$, $10_{81}$, $10_{85}$, $10_{86}$, $10_{88}$, $10_{90}$, $10_{92}$, $10_{96}$ to knots with $\gamma_4=1$.}\label{gamma4,3}
\end{figure}
\newpage
\begin{figure}[h]
	\centering
	\begin{subfigure}[b]{0.30\textwidth}
		\includegraphics[width=\textwidth]{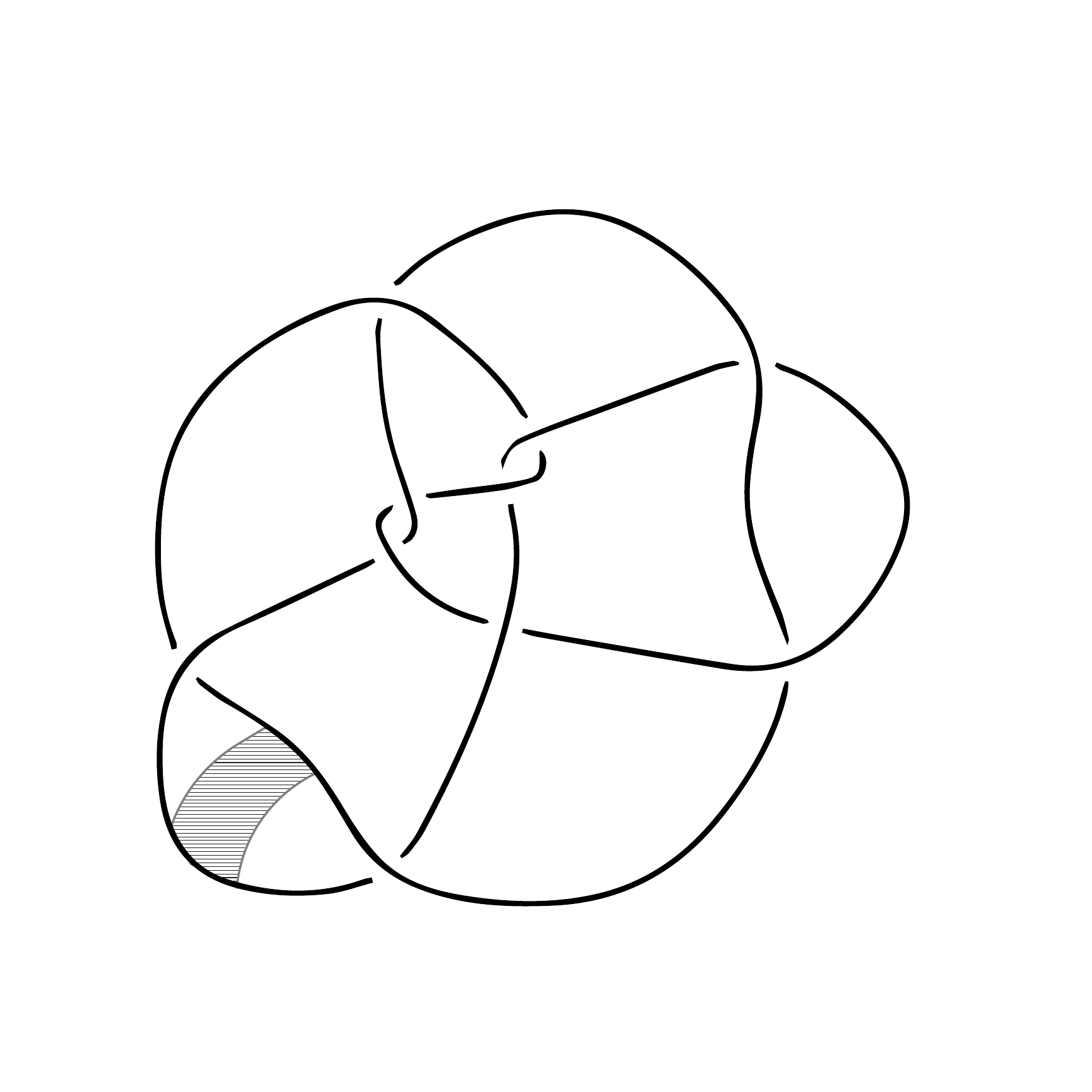}
		\caption{$10_{98}\stackrel{0}{\longrightarrow} 8_{10}$}
		\label{FigureFor10-98}
	\end{subfigure}
	~
	\begin{subfigure}[b]{0.30\textwidth}
		\includegraphics[width=\textwidth]{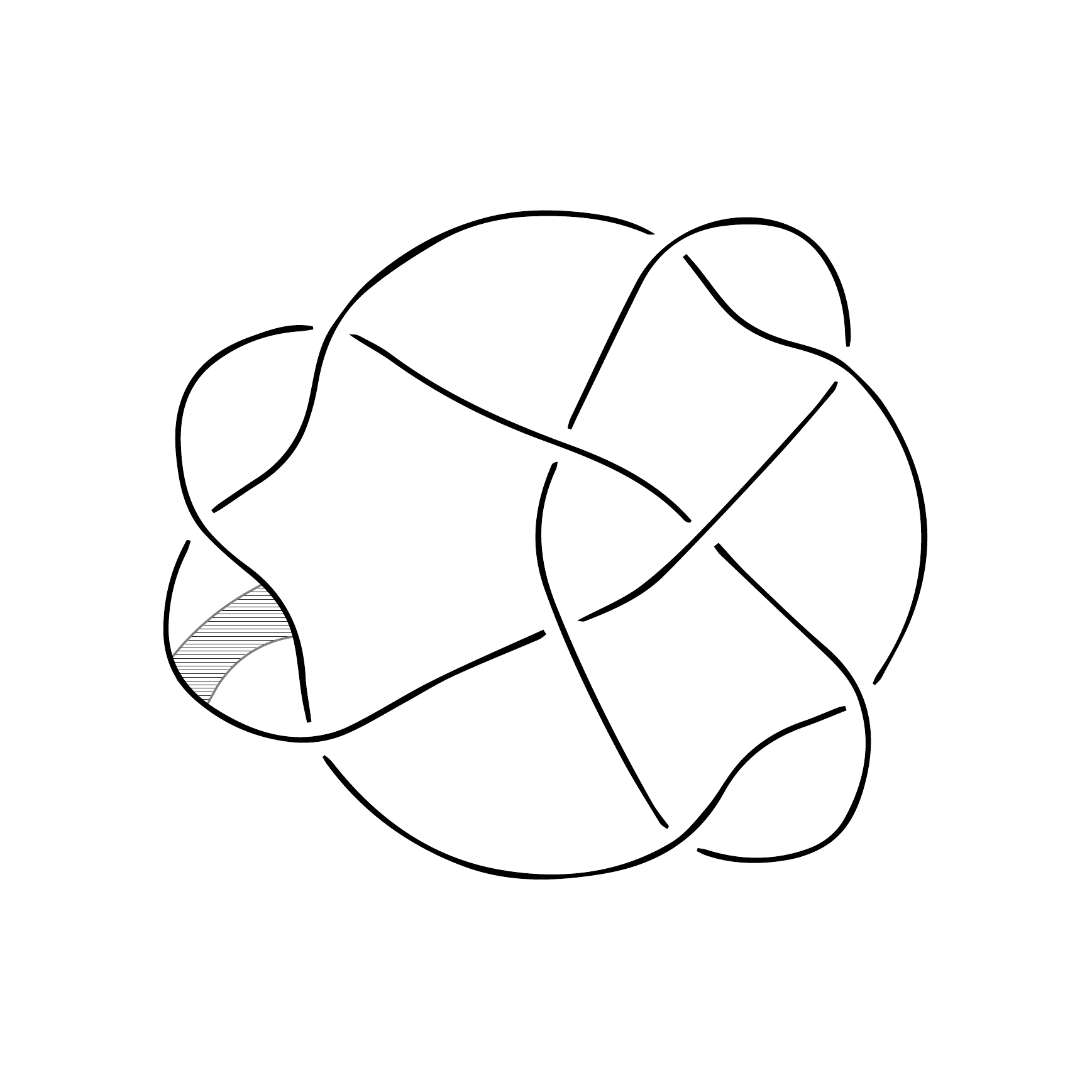}
		\caption{$10_{100}\stackrel{0\phantom{i}}{\longrightarrow} 7_4$}
		\label{FigureFor10-100}
	\end{subfigure}
	~
	\begin{subfigure}[b]{0.30\textwidth}
		\includegraphics[width=\textwidth]{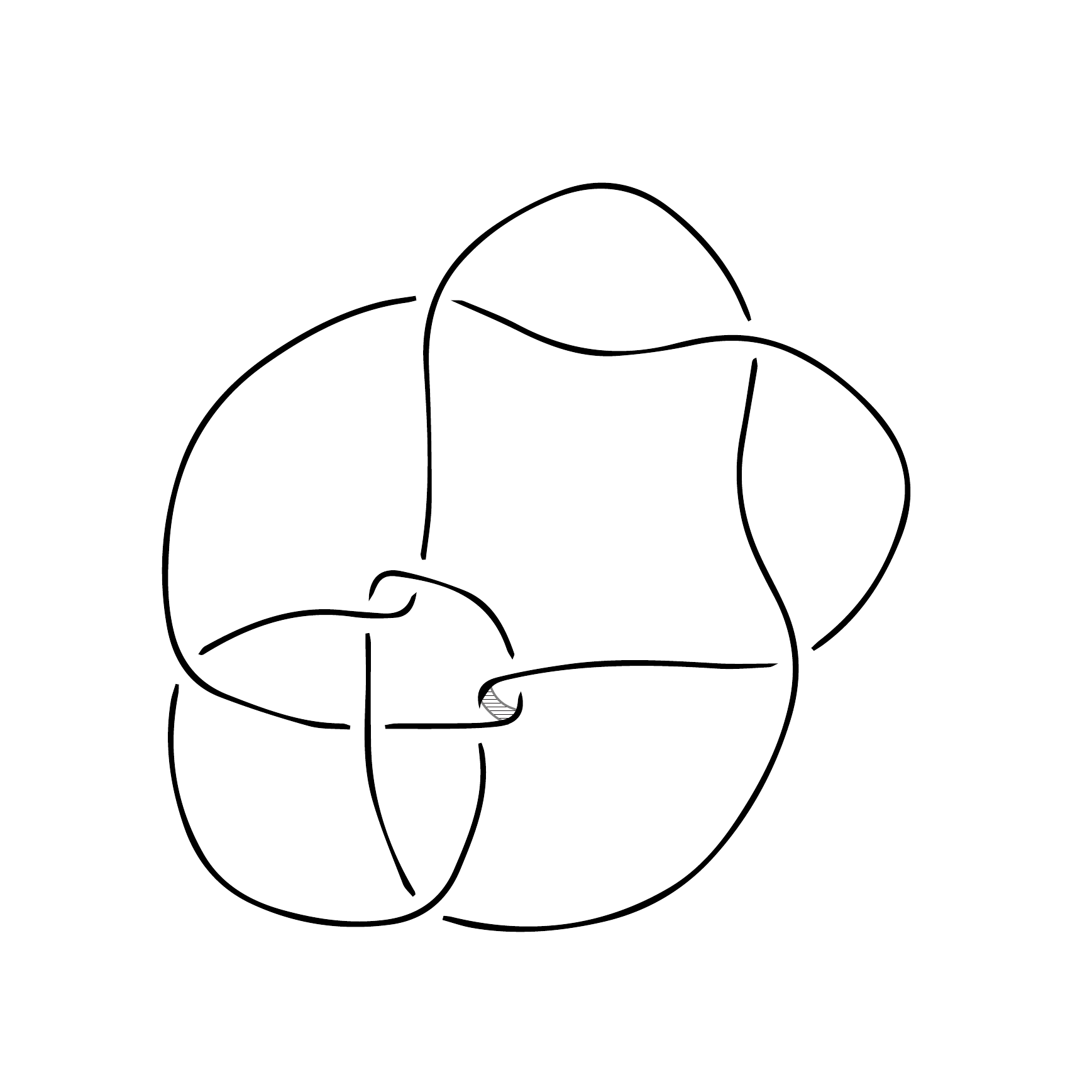}
		\caption{$10_{104}\stackrel{0}{\longrightarrow} 8_{11}$}
		\label{FigureFor10-104}
	\end{subfigure}
	~
	\vskip3mm
	~
	\begin{subfigure}[b]{0.30\textwidth}
		\includegraphics[width=\textwidth]{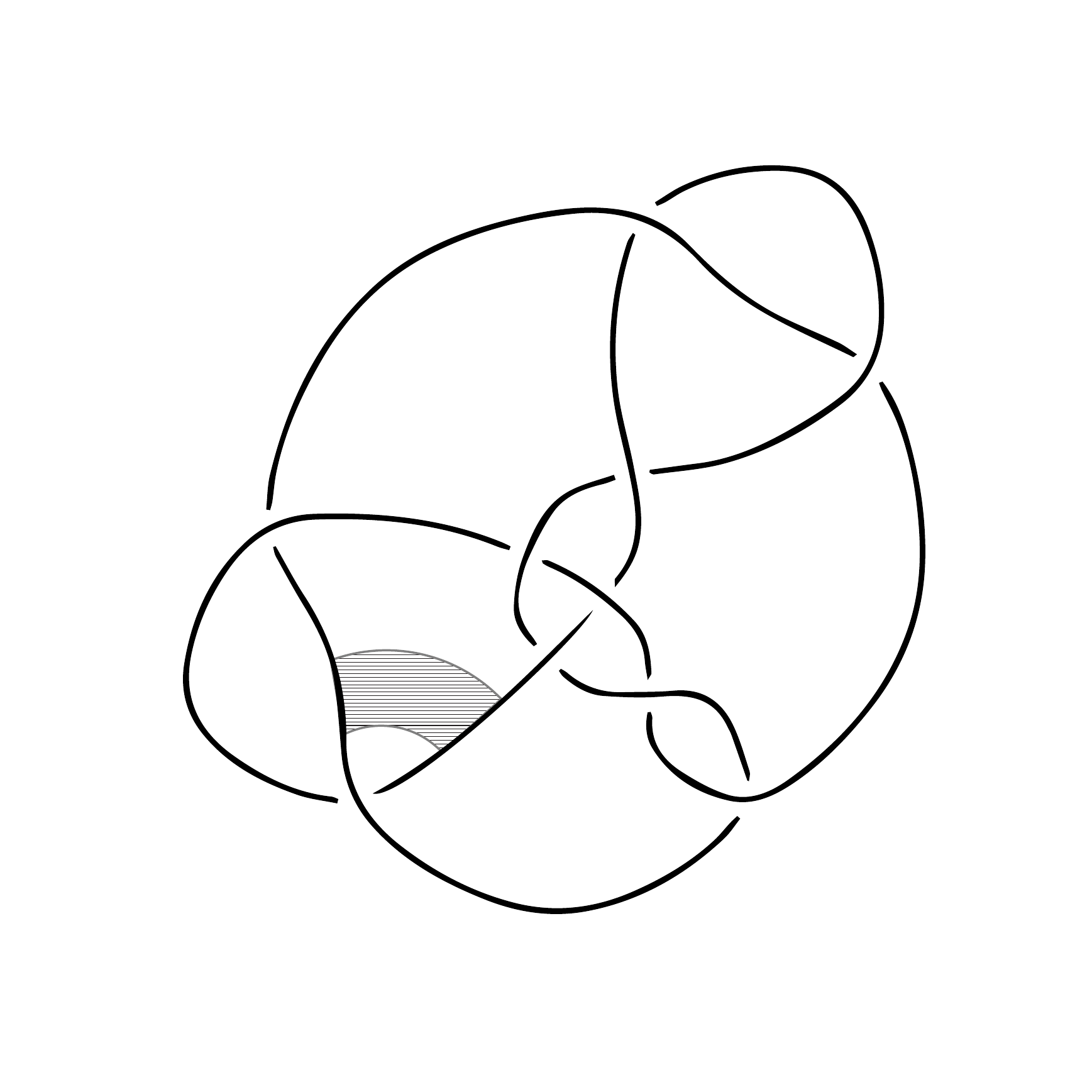}
		\caption{$10_{107}\stackrel{0}{\longrightarrow} 9_{32}$}
		\label{FigureFor10-107}
	\end{subfigure}
~
	\begin{subfigure}[b]{0.3\textwidth}
		\includegraphics[width=\textwidth]{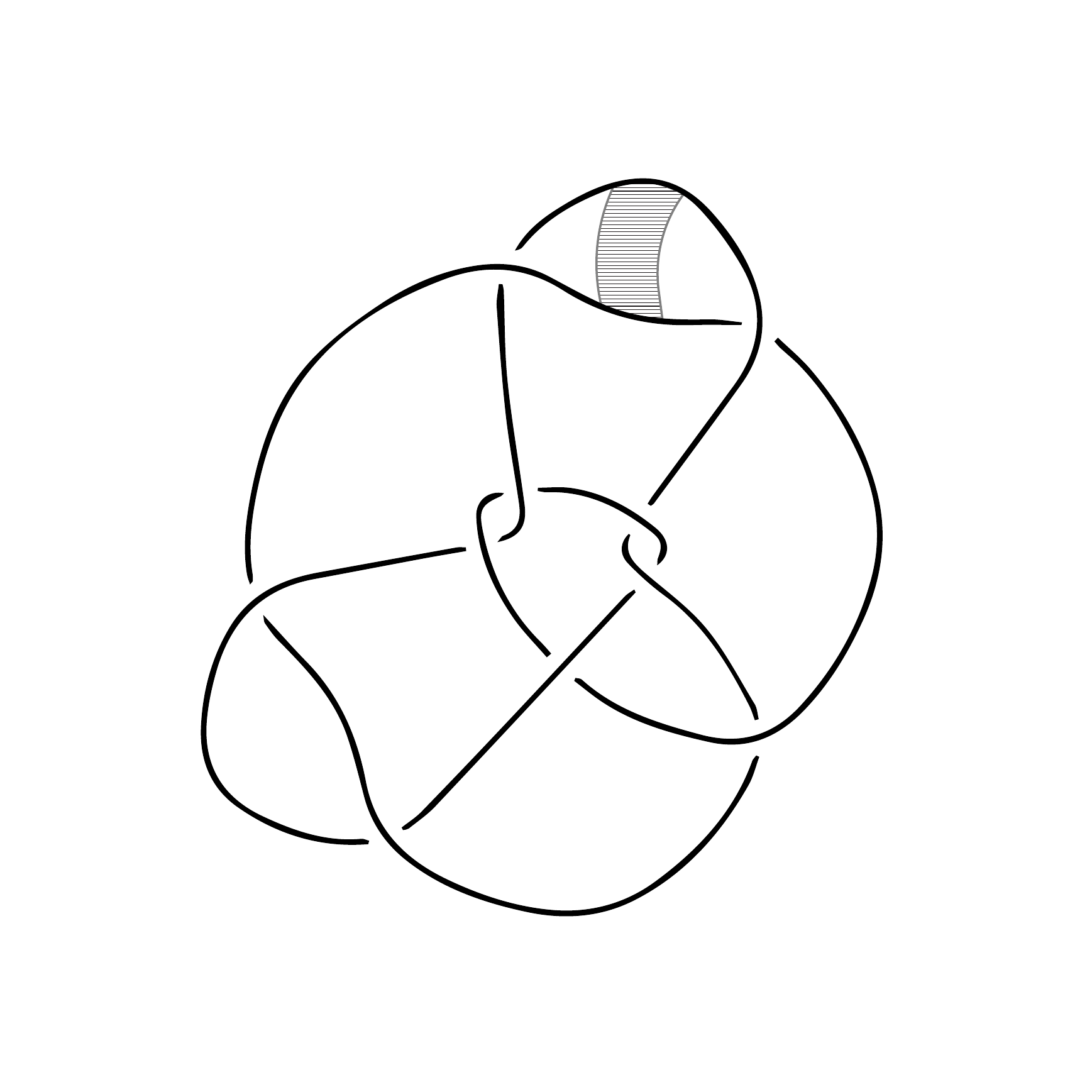}
		\caption{$10_{109}\stackrel{0}{\longrightarrow} 8_{14}$}
		\label{FigureFor10-109}
	\end{subfigure}
	~
	\begin{subfigure}[b]{0.3\textwidth}
		\includegraphics[width=\textwidth]{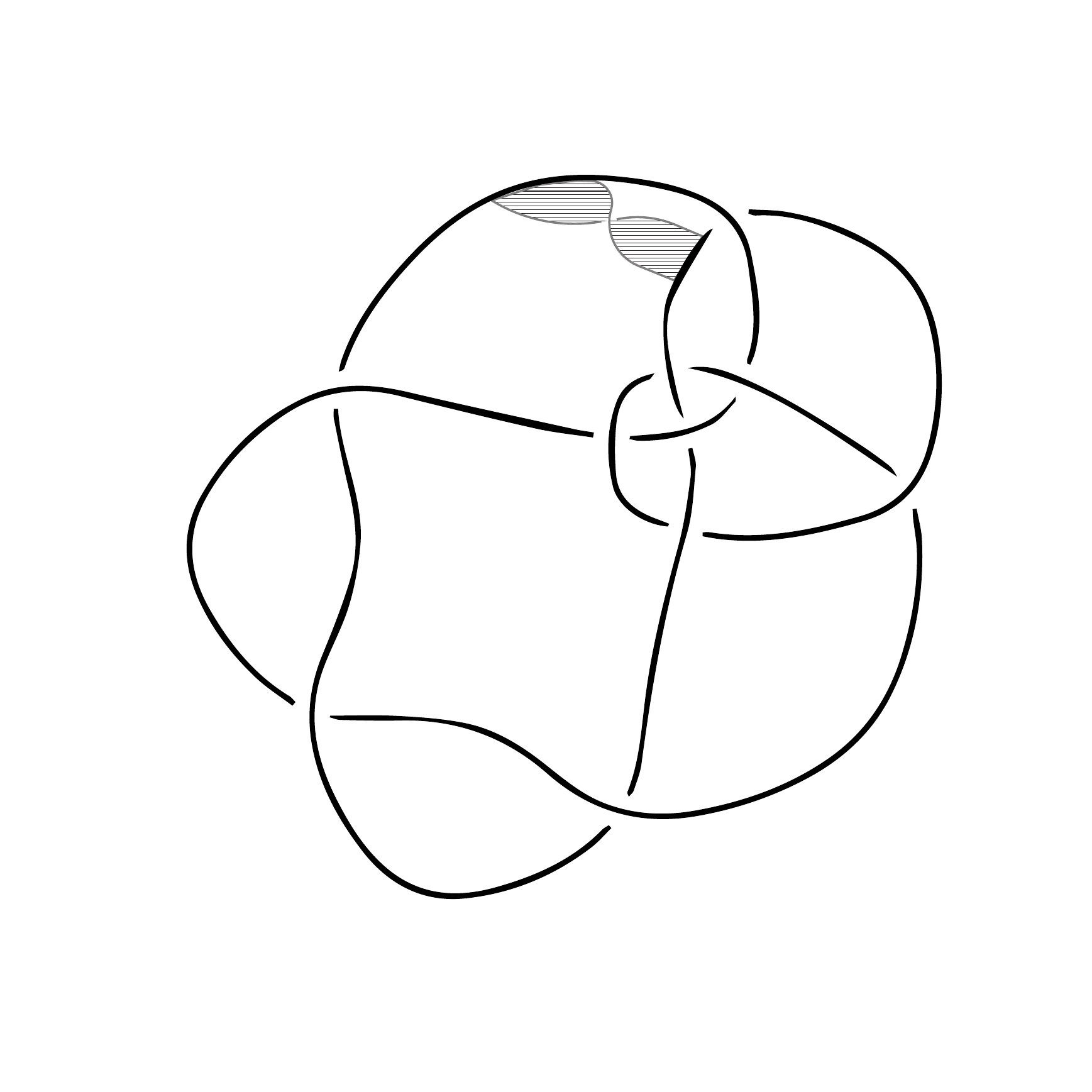}
		\caption{$10_{114}\stackrel{1}{\longrightarrow} 9_{26}$}
		\label{FigureFor10-114}
	\end{subfigure}
	~
	\vskip3mm
	~
	\begin{subfigure}[b]{0.3\textwidth}
		\includegraphics[width=\textwidth]{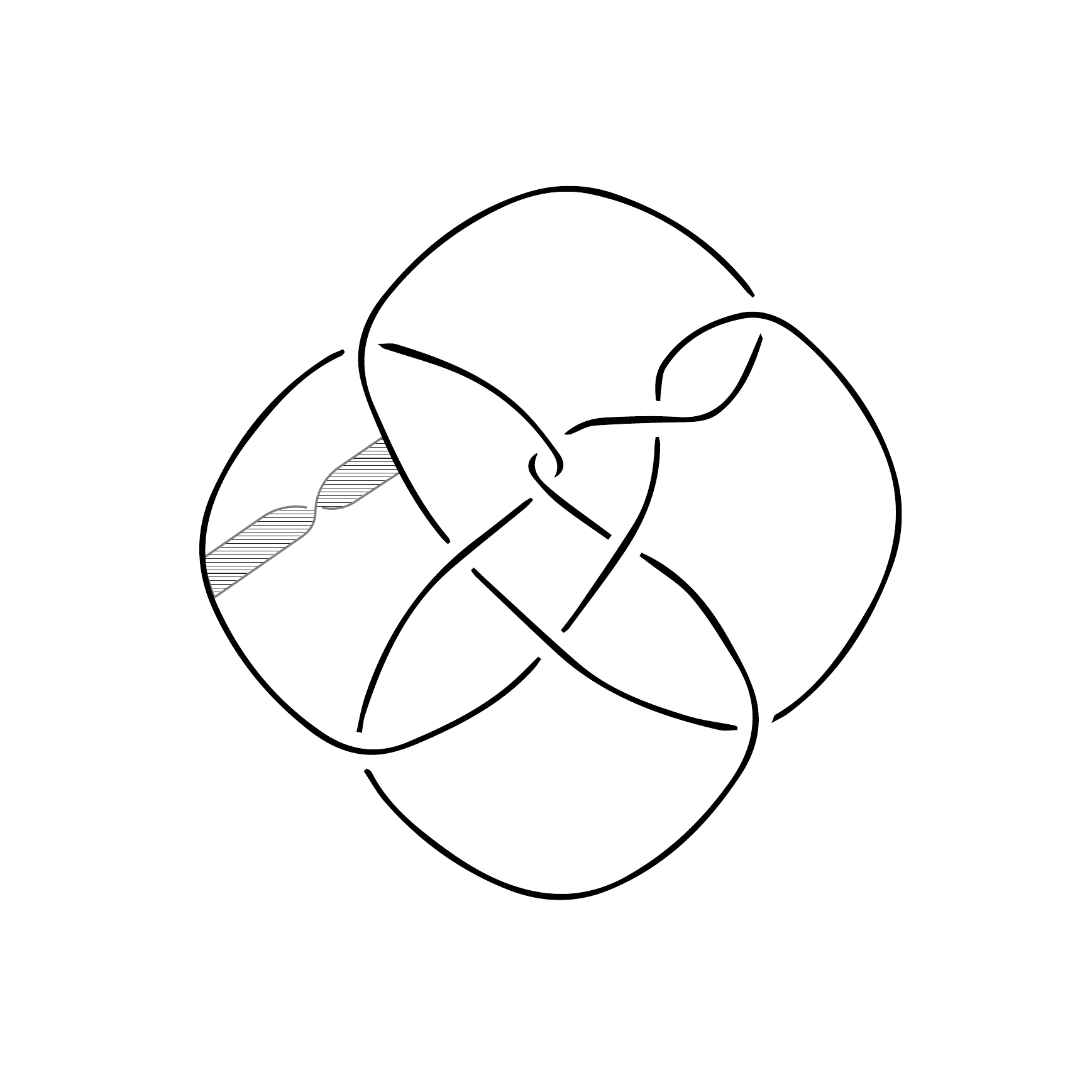}
		\caption{$10_{115}\stackrel{-1}{\longrightarrow} 9_{32}$}
		\label{FigureFor10-115}
	\end{subfigure}
	~
	\begin{subfigure}[b]{0.3\textwidth}
		\includegraphics[width=\textwidth]{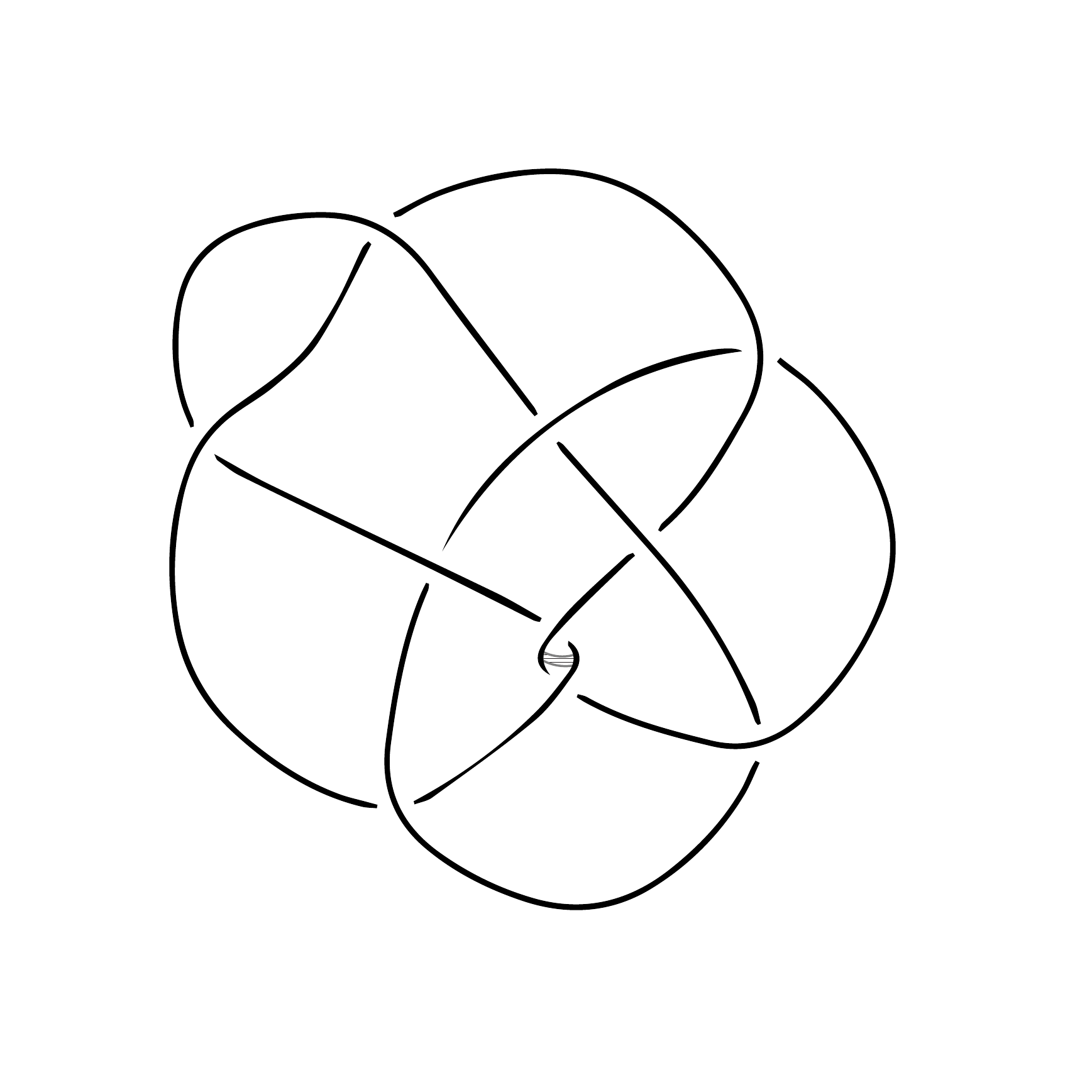}
		\caption{$10_{119}\stackrel{0}{\longrightarrow} 8_{16}$}
		\label{FigureFor10-119}
	\end{subfigure}
	~
	\begin{subfigure}[b]{0.3\textwidth}
		\includegraphics[width=\textwidth]{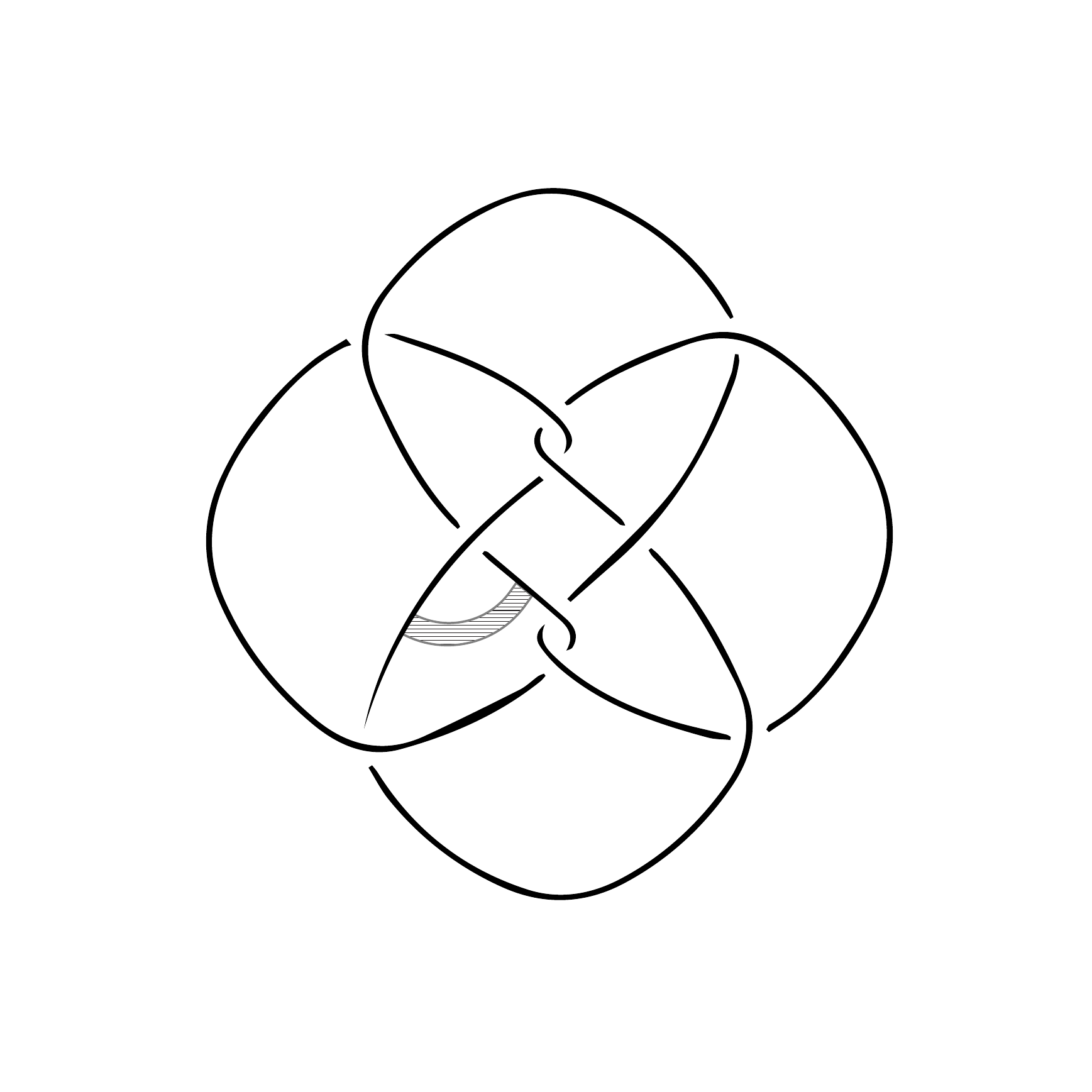}
		\caption{$10_{120}\stackrel{0}{\longrightarrow} 9_{31}$}
		\label{FigureFor10-120}
	\end{subfigure}
	~
		\vskip3mm
	\caption{Non-oriented band moves from the knots $10_{98}$, $10_{100}$, $10_{104}$, $10_{107}$, $10_{109}$, $10_{114}$, $10_{115}$, $10_{119}$, $10_{120}$ to knots with $\gamma_4=1$.}\label{gamma4,4}
\end{figure}
\newpage
\begin{figure}[h]
	\centering
	\begin{subfigure}[b]{0.3\textwidth}
		\includegraphics[width=\textwidth]{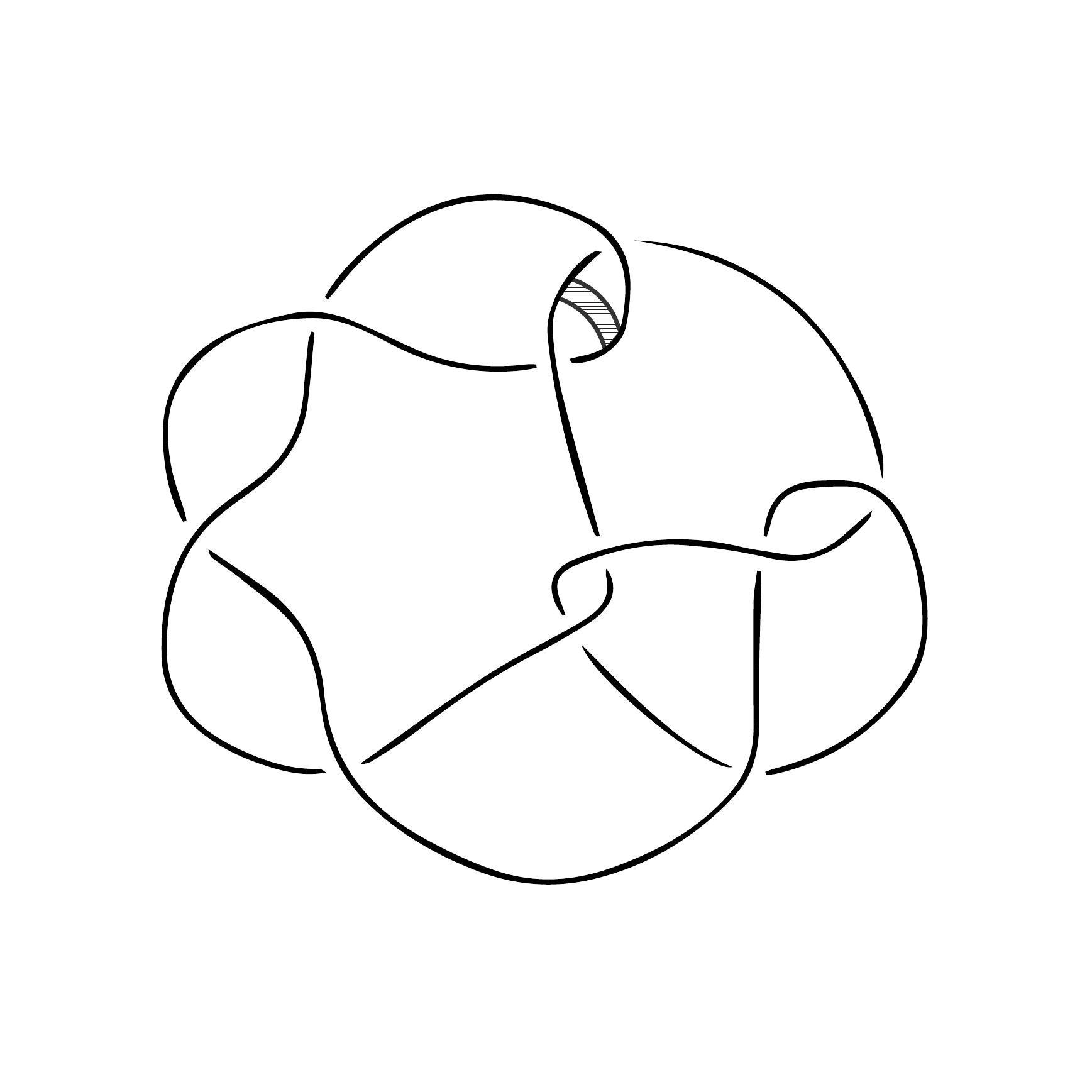}
		\caption{$10_{132}\stackrel{0}{\longrightarrow} 3_{1}$}
		\label{FigureFor10-132}
	\end{subfigure}
	~
	\begin{subfigure}[b]{0.3\textwidth}
		\includegraphics[width=\textwidth]{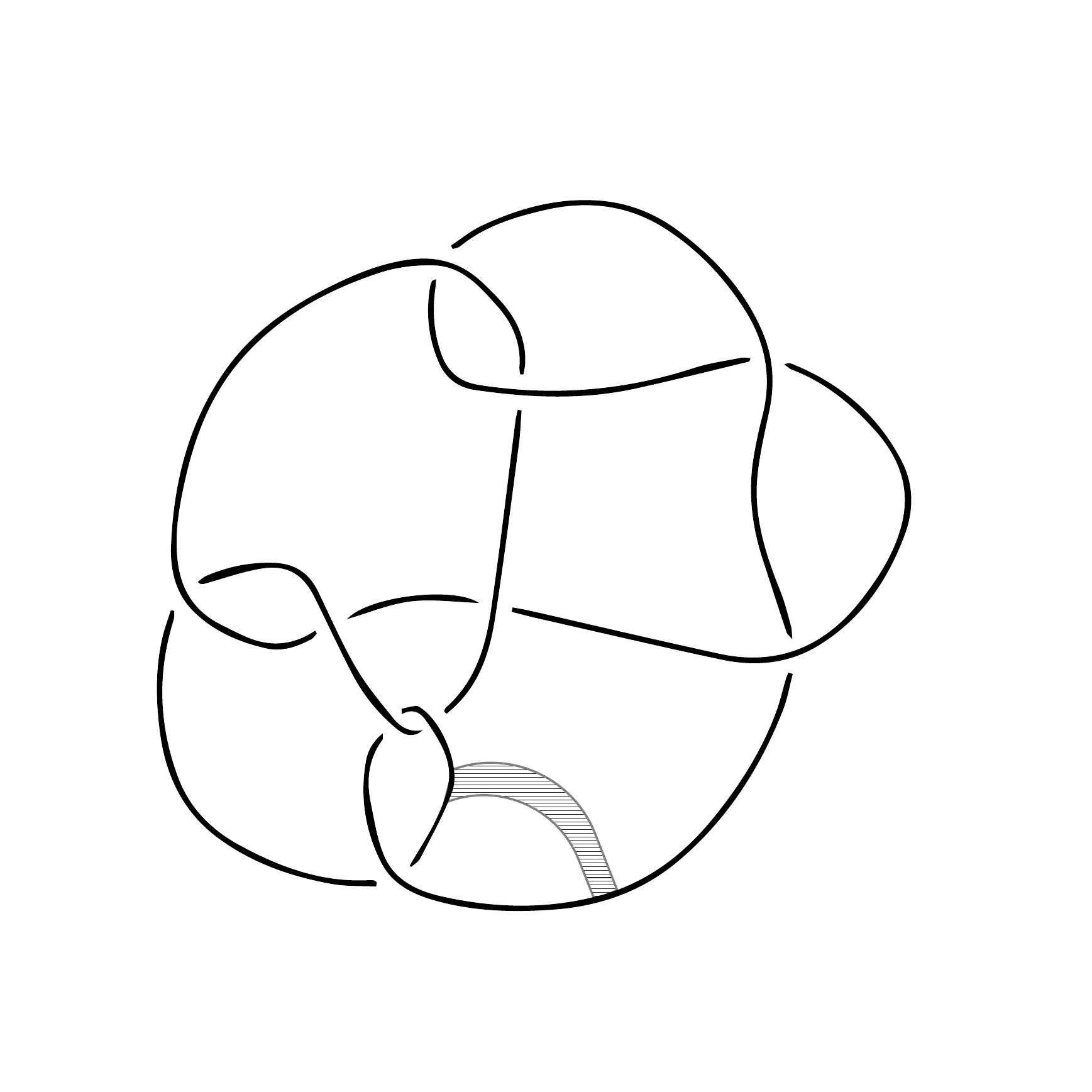}
		\caption{$10_{135}\stackrel{0}{\longrightarrow} 8_{14}$}
		\label{FigureFor10-135}
	\end{subfigure}
	~
	\begin{subfigure}[b]{0.3\textwidth}
		\includegraphics[width=\textwidth]{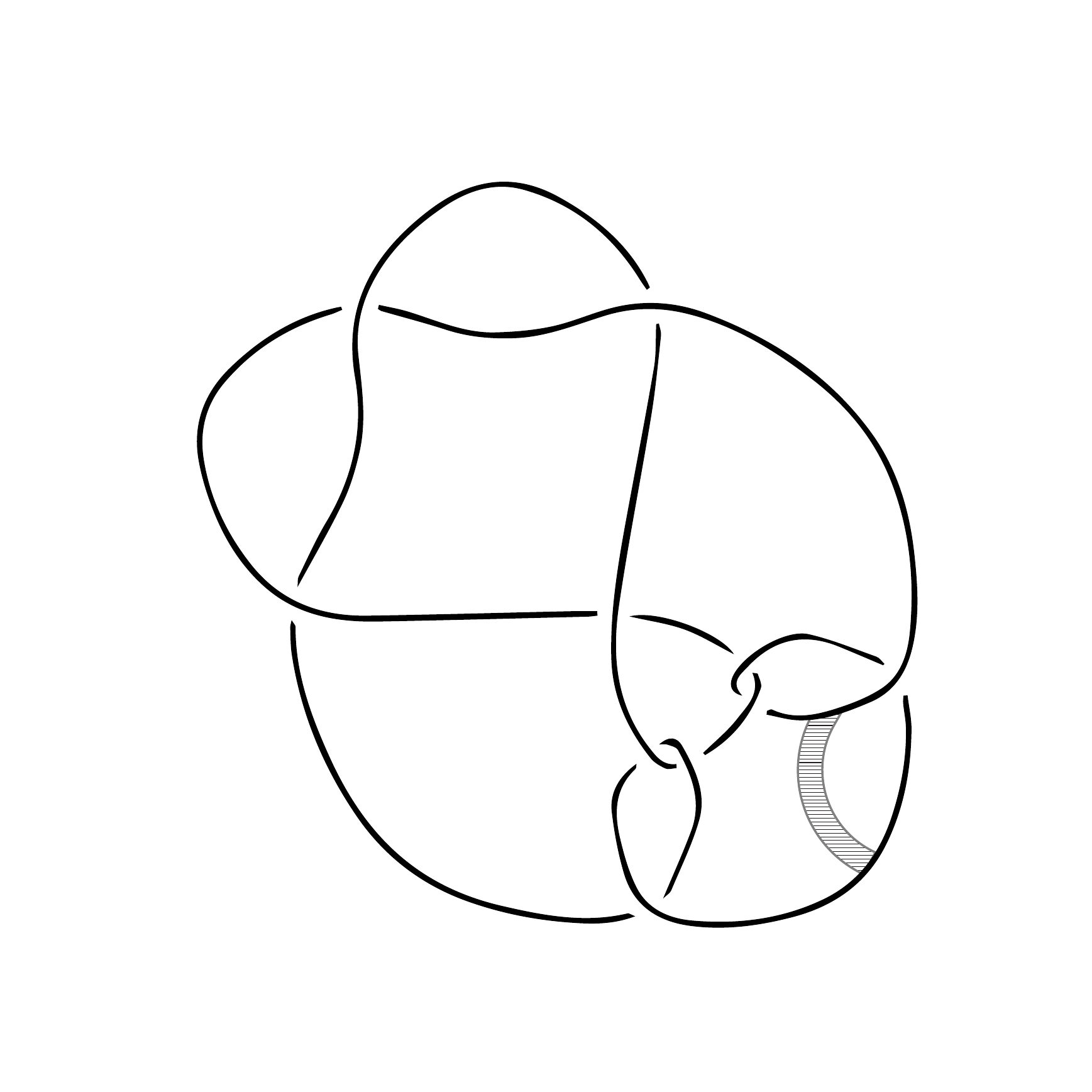}
		\caption{$10_{141}\stackrel{0\phantom{i}}{\longrightarrow} 8_{7}$}
		\label{FigureFor10-141}
	\end{subfigure}
	~
	\vskip3mm
	~
	\begin{subfigure}[b]{0.3\textwidth}
		\includegraphics[width=\textwidth]{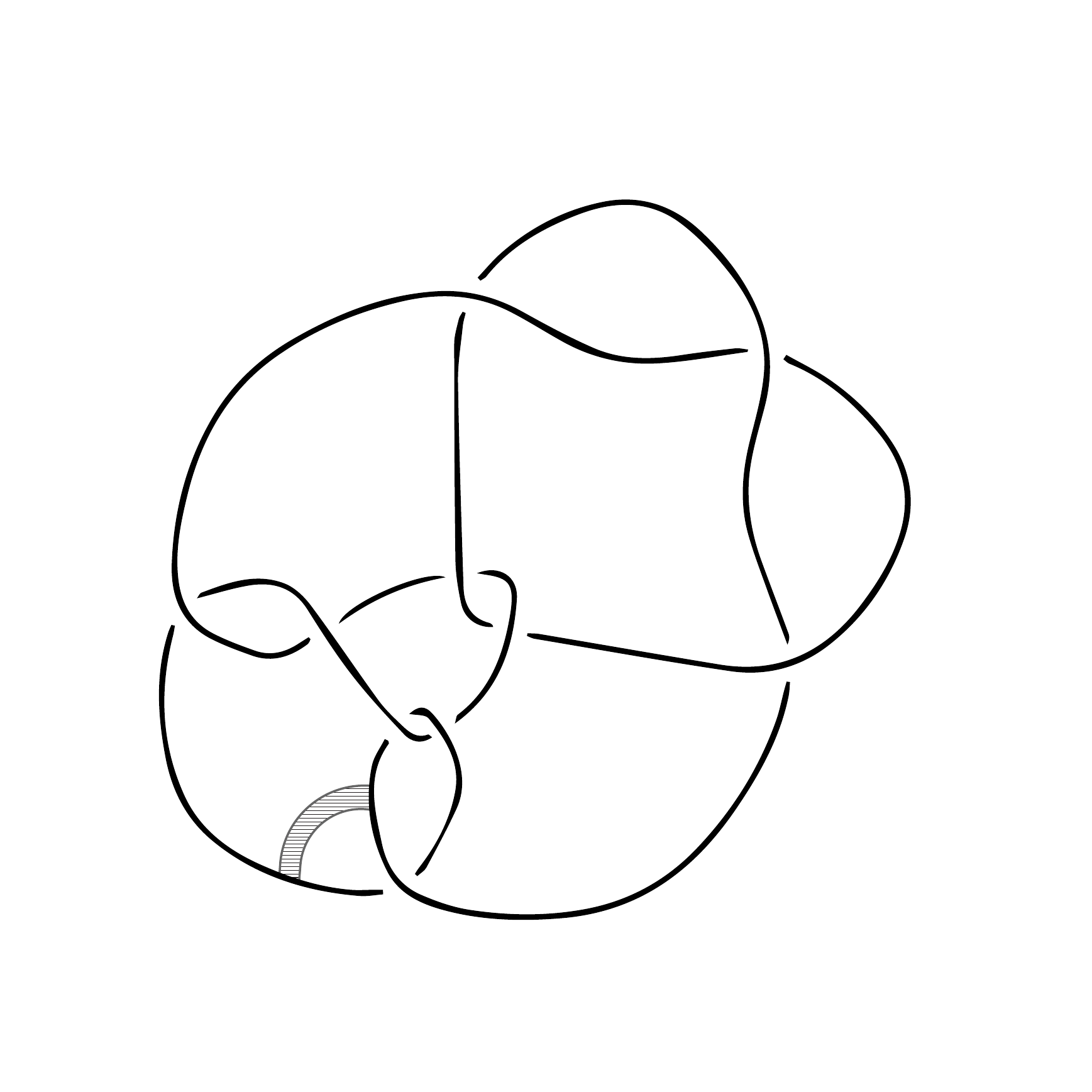}
		\caption{$10_{149}\stackrel{0}{\longrightarrow} 8_{10}$}
		\label{FigureFor10-149}
	\end{subfigure}
	~    
	\begin{subfigure}[b]{0.3\textwidth}
		\includegraphics[width=\textwidth]{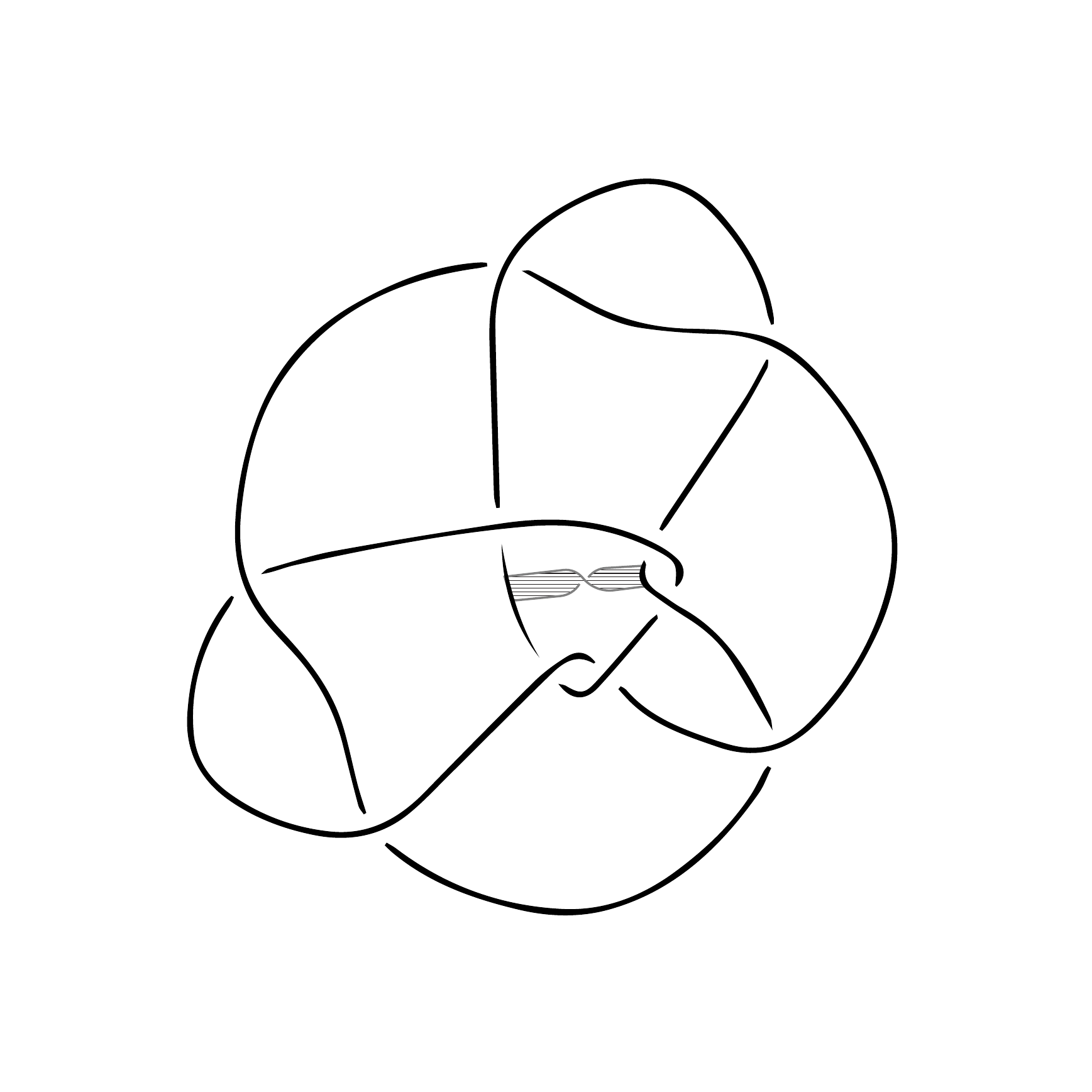}
		\caption{$10_{157}\stackrel{1}{\longrightarrow} 8_{16}$}
		\label{FigureFor10-157}
	\end{subfigure}
	~
	\begin{subfigure}[b]{0.3\textwidth}
		\includegraphics[width=\textwidth]{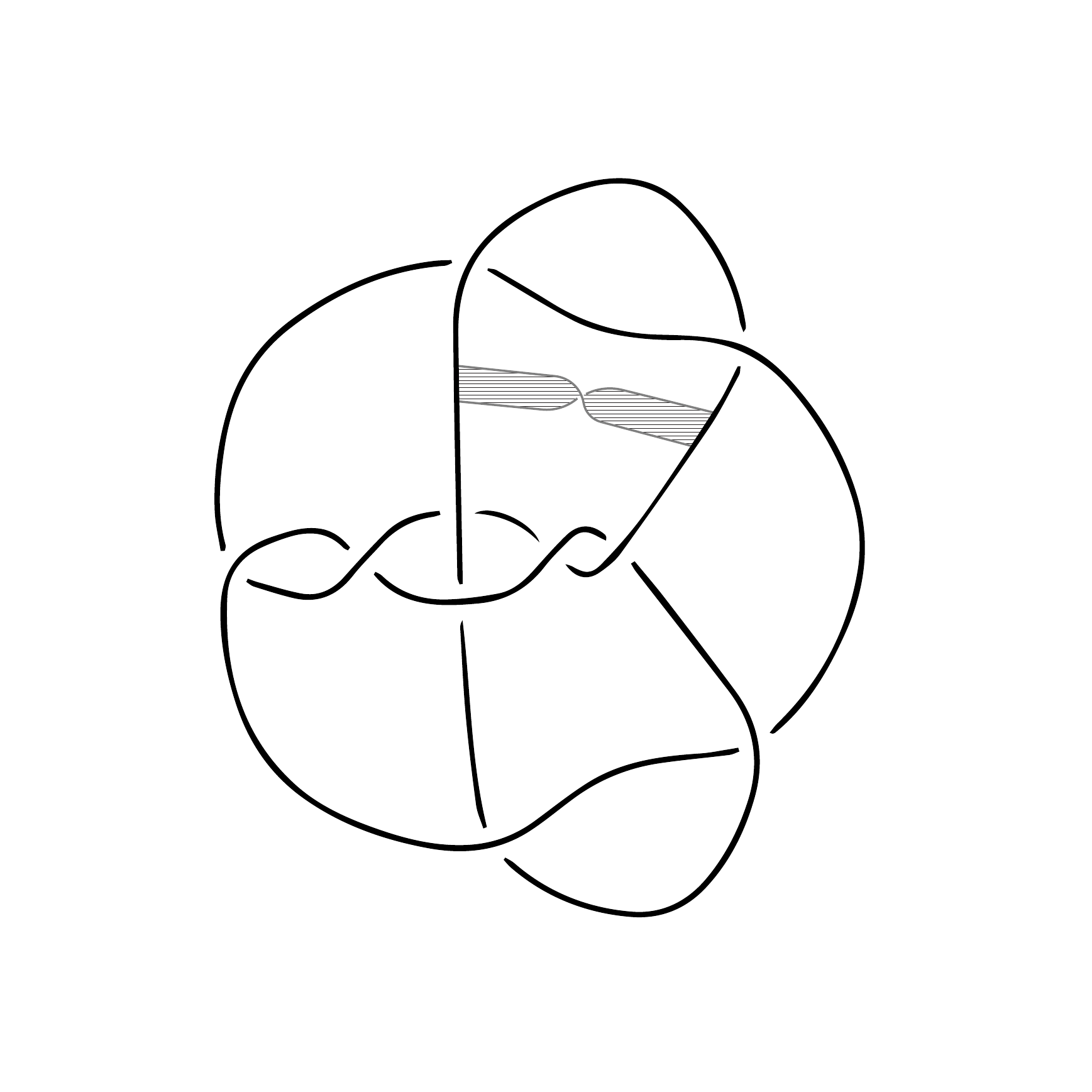}
		\caption{$10_{158}\stackrel{1}{\longrightarrow} 9_{45}$}
		\label{FigureFor10-158}
	\end{subfigure}
	~
	\vskip3mm
	~
	\begin{subfigure}[b]{0.3\textwidth}
		\includegraphics[width=\textwidth]{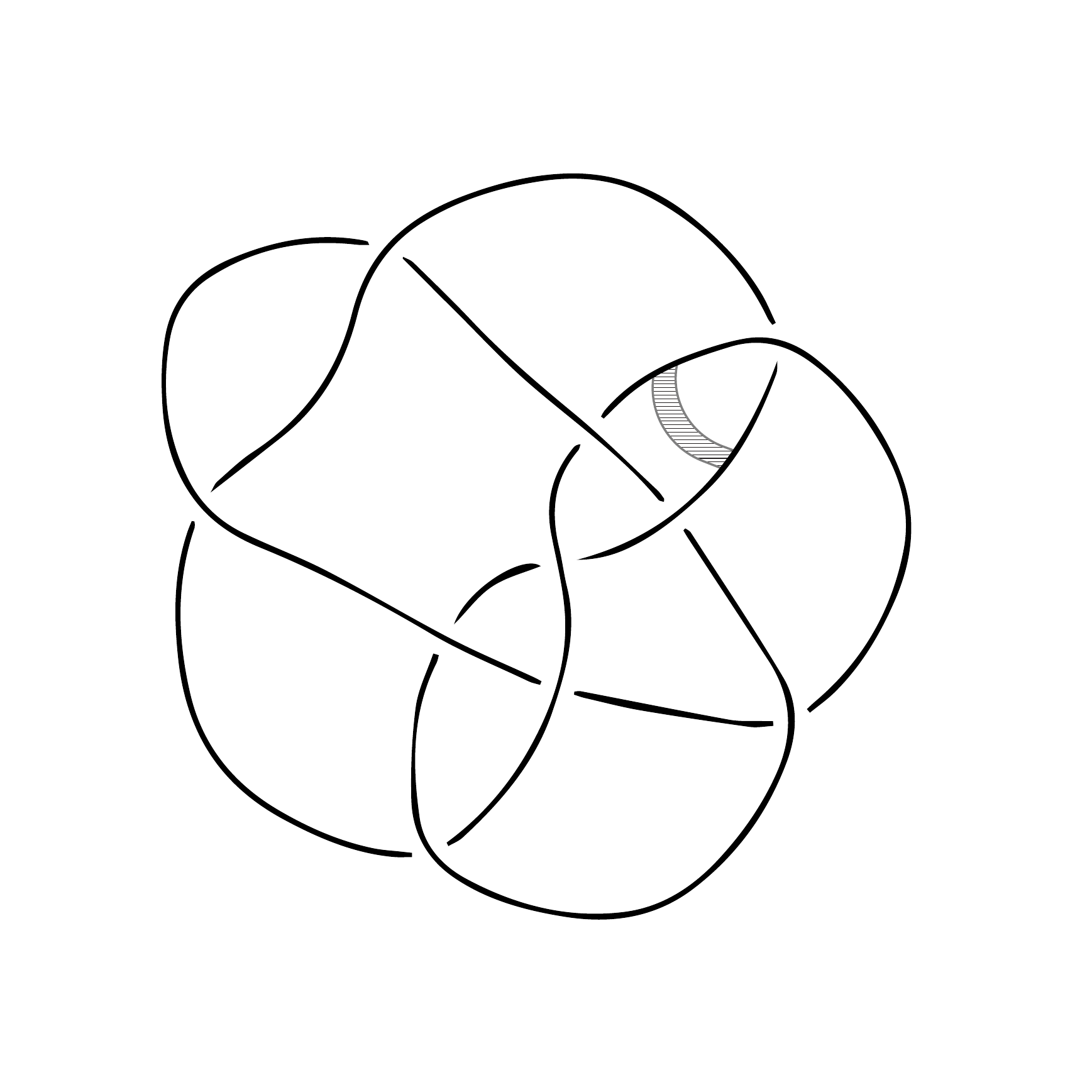}
		\caption{$10_{164}\stackrel{0}{\longrightarrow} 9_{45}$}
		\label{FigureFor10-164}
	\end{subfigure}
	\vskip3mm
		~

	\caption{Non-oriented band moves from the knots  $10_{132}$, $10_{135}$, $10_{141}$, $10_{149}$, $10_{157}$, $10_{158}$, $10_{164}$ to knots with $\gamma_4=1$.}\label{gamma4,5}
\end{figure}\
\newpage
\newpage
\newpage
\newpage
\newpage
\newpage
\newpage
	\newpage
\newpage
\newpage
\newpage
\end{document}